\documentclass[10pt,english]{article}
\usepackage{mathptmx}
\usepackage[T1]{fontenc}
\usepackage[utf8]{inputenc}
\usepackage[english]{babel}
\usepackage{amsmath}
\usepackage{amssymb}
\usepackage{amsthm}
\numberwithin{equation}{section}
\renewenvironment{cases}
 {\left\{\begin{array}{@{}l@{\qquad}l@{}}}
 {\end{array}\right.}
\usepackage{xurl}
\usepackage{setspace}
\usepackage{graphicx}
\usepackage[pass]{geometry}
\usepackage{afterpage}
\makeatletter
\newcommand{\restoretextgeometry}{%
  \restoregeometry
  \global\@colht=\textheight
  \global\@colroom=\textheight
  \global\vsize=\textheight
}
\makeatother
\AtBeginDocument{
  \setlength{\abovedisplayskip}{12pt}
  \setlength{\belowdisplayskip}{12pt}
  \setlength{\abovedisplayshortskip}{12pt}
  \setlength{\belowdisplayshortskip}{12pt}
}
\newtheorem{theorem}{Theorem}
\newtheorem{lemma}[theorem]{Lemma}
\newtheorem{mainresult}{Theorem}

\usepackage[pdfusetitle,
 bookmarks=true,bookmarksnumbered=false,bookmarksopen=false,
 breaklinks=true,pdfborder={0 0 1},backref=false,colorlinks=false]
 {hyperref}
\begin{document}
\title{Large deviations for generalized P\'olya urns\\
 with nonbinary increments}
\author{Simone Franchini{\normalsize\thanks{Correspondence: simone.franchini@yahoo.it}\thanks{CNR-ISTC, Via Gian Domenico Romagnosi 18, 00196 Rome, Italy}}}
\date{~}
\maketitle

\begin{abstract}
In this paper we show how to extend the large deviation principle
for the paths of the Hill, Lane and Sudderth urn model to the case in which
the urn increment is not a binary variable. In particular, we extend
Theorem~1 of [\emph{Stochastic Processes and their Applications} 127 (2017),
3372--3411] to urn processes with integer increments from zero to $K$. The increment probabilities depend on the mean
of the preceding increments and are continuous and strictly positive
on the entire state interval.
We derive the terminal entropy and uniform estimates for the probabilities of individual
attainable final values. We identify the zero-cost states and paths.
We prove the existence and strict convexity of the normalized limiting cumulant generating function
and provide conditions for duality with the entropy.
We characterize the realizable continuous concave entropies and construct
the corresponding urn probabilities.

\end{abstract}

\noindent\newpage{\large\tableofcontents{}}\newpage{}

\section{Introduction}
\label{sec:introduction}

The Hill, Lane and Sudderth (HLS) urn is a growth model with memory:
at each step a ball, black or white, is added, and the probability
of its color depends on the proportion of black balls already present
\cite{HLS,HLS_orig_02,HLS_orig_03}. The function assigning this probability is called the urn function.
The limit points, fluctuations and transitions of the model are studied
in \cite{HLS_orig_04,HLS_orig_05,HLS_orig_06}. The evolution of the proportion connects HLS urns to reinforced
processes and stochastic approximation \cite{HLS_orig_18}; for a continuous urn
function, the proportion converges almost surely to a fixed point \cite{HLS}.

In this paper we extend the theory to integer increments
$\sigma_n\in\Omega=\{0,\ldots,K\}$, with $K\geq1$ fixed.
We write $M_n$ for the sum of the first $n$ increments and set $\psi_n=M_n/n$:
the next increment $k$ has probability $\pi_k(\psi_n)$.
The history thus enters through the mean of the preceding increments.
The process is realized as a two-color urn by adding
$K$ balls at each step, of which $k$ are black when the increment is $k$.
With the initialization in Subsection~\ref{subsec:market_history},
after $n$ steps there are $M_n$ black balls out of $Kn$ and their proportion
is $\psi_n/K$. For $K=1$ we recover the HLS urn.
Throughout the paper the functions $\pi_k$ are continuous and strictly
positive on $[0,K]$, including the endpoints, and satisfy
$\sum_{k\in\Omega}\pi_k(\alpha)=1$ for every state.
Additional assumptions are stated in the individual results.

Large deviations quantify the probability of observing paths
and final states far from typical behavior. For the binary HLS urn,
Ref.~\cite{Franchini_URNS} proves the path LDP under the assumptions of its Eq.~(1.4)
(Theorem~1) and the logarithmic limit of point probabilities (Theorem~4).
The latter requires the combinatorial comparison in Subsection~3.2.1,
in addition to the terminal contraction that gives estimates for sets.
Here we extend both results to increments in $\Omega$.
Applications to increasing returns are developed in \cite{FB}.
The tools of the general theory \cite{Dembo_Zeitouni} complement the analytical approaches
for solvable urns \cite{HLS_orig_10,HLS_orig_11,HLS_orig_12,HLS_orig_14} and large deviations
in tree growth \cite{HLS_orig_13}. The theses~\cite{HLS_orig_15,HLS_orig_16} address, respectively, large
deviations of urns and polymers with self-intersections;
Ref.~\cite{HLS_orig_17} revisits the formulation of the HLS urn.

\subsection{Relation to other models}
\label{subsec:relation_other_models}

Generalized urns offer a common description of various reinforcement
mechanisms, but the relation to each model depends on the observed
variable and its update rule. In urns with a replacement
matrix, including the Bagchi--Pal model \cite{HLS_orig_19,HLS_orig_20}, the additions depend
on the color drawn. Some balanced two-color urns can be represented
by affine urn functions after a change of state (Ref.~\cite{Franchini_URNS}, Section~2.3).
To reduce an urn with more colors to the scalar process considered here,
the increment probabilities of the observed quantity must depend
on the composition only through $\psi_n$.

An explicit example is the \emph{elephant random walk} \cite{ERW=000020shcutz=000020trimper,ERW=000020UM=000020Baur=000020Berton}. After the
first step, a preceding time is chosen uniformly and its sign is repeated
with probability $p$, and reversed with probability $1-p$.
Encoding positive and negative signs by $1$ and $0$, the position is
$2M_n-n$ and the probability of a positive sign is
$p\psi_n+(1-p)(1-\psi_n)$. For $0<p<1$ both probabilities are
uniformly positive; our results therefore transfer
to the walk, with point probabilities on the attainable lattice
$\{-n,-n+2,\ldots,n\}$ and, for the initial law
\eqref{eq:initial_distribution}, a symmetric first step.
The relation in law to random-replacement urns is distinct from this
encoding (Ref.~\cite{ERW=000020UM=000020Baur=000020Berton}, Section~3, Eq.~(3)). Variants of memory, associated
observables and rules with multiple draws are studied in \cite{HLS_orig_23,HLS_orig_24,HLS_orig_25,HLS_orig_26,HLS_orig_28};
the preprint in Ref.~\cite{HLS_orig_28} was subsequently incorporated into the work
in Ref.~\cite{FB}. The work~\cite{HLS_orig_27} analyzes the mechanisms of rare fluctuations in systems with memory.

Nonlinear urn functions are central to Arthur's theory of increasing
returns: previous choices modify the attractiveness of
technologies and can select persistent outcomes \cite{HLS_orig_03,HLS_orig_29,HLS_orig_31,HLS_orig_32,HLS_orig_33,HLS_orig_34,HLS_orig_35,HLS_orig_36,HLS_orig_37}.
This framework includes developments on nonlinear feedback
\cite{HLS_orig_30,HLS_orig_41}, on the stabilization of leadership \cite{HLS_orig_38,HLS_orig_39}, on the distribution
of wages and returns on capital \cite{HLS_orig_40} and on social influence
\cite{HLS_orig_42,HLS_orig_43}. The large deviation theory in \cite{FB} also allows rare
adoption histories to be studied. In its model, described in
Sections~I and~II, a customer consults an odd-sized sample of previous
customers and follows the majority with probability $p$, and the opposite choice
with probability $1-p$. Consulting just one customer yields the
affine function $\pi_1(\alpha)=(1-p)+(2p-1)\alpha$; if $0<p<1$, both
probabilities are at least $\min\{p,1-p\}>0$ on $[0,1]$. A larger sample
modifies the reinforcement, but leaves the final choice binary.

Other connections come from physics and biology.
The cluster growth models studied by Klymko, Garrahan,
Whitelam and collaborators \cite{HLS_orig_46,HLS_orig_47}, and their analysis using large
deviations \cite{HLS_orig_44,HLS_orig_45}, relate reinforcement to symmetry
breaking and the coexistence of different growth behaviors.
In particular, the vanishing of the cost between stationary values of the
magnetization \cite{HLS_orig_44} has a counterpart in the phenomenon described here by
Theorem~\ref{thm:intro_zero_cost}.
The Khanin and Khanin model for neuronal polarity \cite{HLS_orig_48}
instead uses competitive growth of neurites, in which the
longer ones are more likely to grow further.
These examples motivate the study of nonlinear reinforcement. Applying
our theorems requires a scalar description with nonnegative increments
and the assumptions stated above; reversible growth dynamics
therefore require a separate treatment.

Urn theory also enters the study of the number of sites
visited by a random walk and of polymer
self-intersections \cite{HLS_orig_15,HLS_orig_16,HLS_orig_49,HLS_orig_50,HLS_orig_51,HLS_orig_52}. The microcanonical formulation in \cite{HLS_orig_49}
connects growth probabilities to large deviation properties
and also leads to a mean-field approximation based on the
Bagchi--Pal model. In this context one encounters the Stanley model,
the Rosenstock trapping problem \cite{HLS_orig_16,HLS_orig_51}, self-avoiding
walks and the corresponding continuous problem of the volume of the
\emph{Wiener sausage} \cite{HLS_orig_53}. These connections may concern
terminal distributions or an asymptotic limit: identification
with the law of all urn paths requires the probability
of the next increment to depend only on the scalar state. For the number
of visited sites, this property does not follow from the binary nature
of the increment alone: the probability of reaching a new site may depend
on the geometry of the path.

\subsection{Beyond binary increments}
\label{subsec:beyond_binary_increments}

Allowing increments from zero to $K$ makes it possible to describe additions
of different sizes while retaining a two-color urn and a scalar state.
This extension connects to an already developed literature on
multiple additions. Launay~\cite{Launay_simultaneous} studies simultaneous
draws with reinforcement updated after the entire block;
Morcrette and Mahmoud~\cite{HLS_orig_14} treat balanced urns with
random replacements, including binomial and uniform additions.
For urns with rules depending on a sample of several balls,
Kuba and Mahmoud~\cite{Kuba_Mahmoud_multidrawing} study the class
with an affine conditional mean; Lasmar, Mailler and
Selmi~\cite{Lasmar_Mailler_Selmi} also treat rules without this
restriction.
These models do not generally coincide with the class considered here:
the probabilities may depend on absolute counts or on the total
number of balls, and some increments may have zero probability.
Here every increment has a continuous and strictly positive probability,
determined by the mean $\psi_n$ alone.

One possible interpretation is the number of units purchased by a
customer, up to $K$, including no purchase. Models of
innovation diffusion~\cite{HLS_orig_54} illustrate the role of
feedback in successive adoptions; purchases of additional units are also studied
in diffusion modeling
\cite{Steffens_multiple_units}. In this interpretation $\psi_n$
is the mean quantity per customer, including nonbuyers, and dependence
on the mean alone is an assumption of our model. Describing the shares
of several products would instead require a state that distinguishes them.

For large deviations, too, there are precedents beyond the binary case.
Budhiraja and Waterbury~\cite{Budhiraja_Waterbury} prove an LDP
for the empirical frequencies of reinforced processes on a finite set,
when the distribution of the new symbol is obtained by applying
a strictly positive stochastic matrix to the preceding frequencies.
Our subclass in which all
$\pi_k$ are affine falls within that model: the LDP for the final state
is then obtained by contraction of the frequencies.
Budhiraja, Waterbury and Zoubouloglou~\cite{Budhiraja_Waterbury_Zoubouloglou}
also consider dependence on the last symbol and more general
communication conditions, still with affine dependence on the
frequencies for the full LDP. In the present paper we treat
continuous probabilities $\pi_k$, including nonaffine ones, and describe both
the paths and the probabilities of individual final values.

Coghi and Garrahan~\cite{Coghi_Garrahan_spectral} develop a spectral
approach to self-interacting jump processes and discrete-time chains,
coupling tilted operators to Hamilton--Jacobi characteristics through
a WKBJ expansion whose terms are assumed sufficiently smooth.
For the scalar urn model considered here, the path and count LDPs
and the uniform point-probability estimates are proved for continuous,
strictly positive urn probabilities.
Our results also include an inverse characterization of realizable
concave entropies.

The transition from binary to multiple increments requires distinguishing
the cumulative sum from the frequencies that produce it.
Over a horizon of $N$ steps, the path $\varphi$ interpolates
the cumulative sums $M_n/N$ at times $n/N$;
the counts $\gamma_k$ interpolate, at the same times, the number
of the first $n$ increments equal to $k$, divided by $N$.
For $\tau\in[0,1]$ we therefore have
\begin{equation}
\sum_{k\in\Omega}\gamma_k(\tau)=\tau,\qquad\sum_{k\in\Omega}k\gamma_k(\tau)=\varphi(\tau).
\label{eq:style_s1_count_constraints}
\end{equation}
For $K=1$, $\gamma_1=\varphi$ and $\gamma_0(\tau)=\tau-\varphi(\tau)$;
for $K>1$, the same limiting path may correspond to different
counts. Example~\eqref{eq:distinct_count_paths}--\eqref{eq:style_s3_distinct_second_action}
shows that histories with the same scalar limit can have different probabilities
already at exponential order in $N$.
At finite time the interpolated path reconstructs all increments:
the information loss occurs in the uniform limit.
The scalar rate is therefore obtained by minimizing the rate of the compatible
counts, as made explicit in Theorem~\ref{thm:intro_path_ldp}.

The mean increment
\begin{equation}
\bar\pi(\alpha)=\sum_{k\in\Omega}k\pi_k(\alpha)
\label{eq:style_s1_mean_increment}
\end{equation}
determines the equation for zero-cost paths, but deviations
require the entire increment distribution.
In the family~\eqref{eq:three_increment_double_root_probabilities}--\eqref{eq:style_s9_2_042},
with $K=2$ and $0<c<1/12$, we always have $\bar\pi(\alpha)=1$,
whereas the terminal entropy at endpoint $2$ is
$\phi(2)=\log(1/3+2c)$ and varies with $c$, as in
\eqref{eq:style_s9_2_050}.
Section~\ref{sec:main_results} collects the full statements
and their respective assumptions.
Subsection~\ref{subsec:main_binary_comparison} compares the scope
of the extension with the binary theory.

\section{Main results}
\label{sec:main_results}
\label{subsec:results_and_organization}

Theorems~\ref{thm:intro_path_ldp}--\ref{thm:intro_zero_cost}
describe the cost of paths, terminal probabilities and zero-cost
paths. Theorems~\ref{thm:intro_cumulant}
and~\ref{thm:intro_cumulant_duality} connect the entropy to the cumulant generating function;
Theorem~\ref{thm:intro_inverse_entropy} solves the inverse problem
in the concave class. Finally, Theorem~\ref{thm:intro_affine} collects
the formulas and regularity results for affine urns.
Examples and figures follow the results they illustrate.
Their full calculations are collected in Section~\ref{sec:examples};
Subsection~\ref{subsec:main_binary_comparison} states the scope and
limitations of the extension.

\subsection{Model, assumptions and notation}
\label{subsec:main_model}

Fix an integer $K\geq1$ and set $\Omega=\{0,\ldots,K\}$.
At step $n$, $K$ balls are added, of which $\sigma_n\in\Omega$
are black and the others white. For a history $\sigma=(\sigma_1,\ldots,\sigma_N)$ of $N$ steps,
with $N\geq1$ an integer, define
\begin{equation}
M_0=0,\qquad M_n=\sum_{s=1}^{n}\sigma_s,
\label{eq:display_001}
\end{equation}
\begin{equation}
\psi_n=M_n/n,\qquad 1\leq n\leq N.
\label{eq:style_s2_display_001_2}
\end{equation}
The functions $\pi_k:[0,K]\to(0,1)$ are continuous and satisfy
$\sum_{k\in\Omega}\pi_k(\alpha)=1$ for every $\alpha\in[0,K]$.
We denote the law of the process by $\mathbb P$ and expectation
by $\mathbb E$. The first increment is uniform on $\Omega$ and,
for $1\leq n<N$,
\begin{equation}
\mathbb P(\sigma_{n+1}=k\mid\sigma_1,\ldots,\sigma_n)
=\pi_k(\psi_n),\qquad k\in\Omega.
\label{eq:display_002}
\end{equation}

The cumulative path $\varphi(\sigma)$ and the counts
$\gamma(\sigma)=(\gamma_0(\sigma),\ldots,\gamma_K(\sigma))$
are the linear interpolations of the values
\begin{equation}
\varphi(n/N\mid\sigma)=M_n/N,\qquad 0\leq n\leq N,
\label{eq:display_003}
\end{equation}
\begin{equation}
\gamma_k(n/N\mid\sigma)=\frac1N\#\{s\in\{1,\ldots,n\}:\sigma_s=k\},\qquad 0\leq n\leq N.
\label{eq:style_s2_display_003_2}
\end{equation}
with zero counts for $n=0$.
Let $Q$ be the space of absolutely continuous functions
$\varphi:[0,1]\to\mathbb R$, vanishing at zero, with
$0\leq\partial_\tau\varphi\leq K$ almost everywhere, and set
$Q(\alpha)=\{\varphi\in Q:\varphi(1)=\alpha\}$ for $\alpha\in[0,K]$.
A vector $\gamma=(\gamma_0,\ldots,\gamma_K)$ is admissible if its
components are continuous, nondecreasing, vanish at zero and
$\sum_{k\in\Omega}\gamma_k(\tau)=\tau$ for every $\tau\in[0,1]$.
The components are therefore absolutely continuous: the derivatives
$\partial_\tau\gamma_k$ are the local increment frequencies,
nonnegative and summing to one almost everywhere. The associated path is
$\varphi=\sum_{k\in\Omega}k\gamma_k$.
We use the uniform distance for paths and the maximum of the uniform
distances between components for counts.
For every $\varphi\in Q$, the continuous mean is
$\psi(\tau)=\varphi(\tau)/\tau$ for $0<\tau\leq1$;
we do not impose a limit at zero.
The definitions and constraints are developed in
Section~\ref{sec:notation}.

\subsection{Large deviations and terminal entropy}
\label{subsec:main_ldp_entropy}

\begin{mainresult}[Large deviations of paths and counts]
\label{thm:intro_path_ldp}
The laws of $\varphi(\sigma)$ satisfy a large
deviation principle (LDP) with speed $N$ on $Q$.
The good rate function is finite on all of $Q$ and is
\begin{equation}
I(\varphi)=\int_0^1
\sup_{\lambda\in\mathbb{R}}
\left\{\lambda\,\partial_\tau\varphi(\tau)
-\log\!\left[\sum_{k\in\Omega}
\pi_k\!\left(\frac{\varphi(\tau)}{\tau}\right)
e^{k\lambda}\right]\right\}\,d\tau.
\label{eq:intro_path_rate}
\end{equation}
The laws of the counts $\gamma(\sigma)$ satisfy an LDP
on the space of admissible vectors, with speed $N$
and good rate function
\begin{equation}
I_{\mathrm c}(\gamma):=
\int_0^1\sum_{k\in\Omega}
\partial_\tau\gamma_k(\tau)
\log\frac{\partial_\tau\gamma_k(\tau)}
{\pi_k\!\left(\tau^{-1}\sum_{j\in\Omega}j\gamma_j(\tau)\right)}
\,d\tau.
\label{eq:intro_count_rate}
\end{equation}
For every $\varphi\in Q$, the two rates satisfy
\begin{equation}
I(\varphi)=
\min_{\substack{\gamma\ \mathrm{admissible}\\
\sum_{k\in\Omega}k\gamma_k=\varphi}}
I_{\mathrm c}(\gamma),
\label{eq:main_rate_contraction}
\end{equation}
with the minimum attained; the constraint is an equality of functions on $[0,1]$.
In the integrals, division by $\tau$ is used only for $\tau>0$;
a term with $\partial_\tau\gamma_k=0$ is zero.
\end{mainresult}

The rate $I_{\mathrm c}$ is the cost of the counts, whereas $I$ is the minimum
cost among all counts producing the same cumulative path.
Theorems~\ref{thm:count_ldp} and~\ref{thm:path_ldp} contain the
full estimates. The proof starts from the LDP for uniform counts,
derived from Ref.~\cite{Dembo_Zeitouni}. The change of measure uses continuity
of the action and the discrete approximation, which extend Lemmas~13 and~14
of Ref.~\cite{Franchini_URNS}. Contraction, with the local minimum in
Lemma~\ref{lem:local_frequency_minimum} and the optimal counts
\eqref{eq:optimal_count_construction}, gives \eqref{eq:intro_path_rate}.
For $K=1$ we recover Theorem~1 of Ref.~\cite{Franchini_URNS},
Eq.~(2.5)--(2.7), in the strictly positive class.

\begin{mainresult}[Terminal entropy and point probabilities]
\label{thm:intro_terminal_entropy}
The terminal entropy, defined by
\begin{equation}
\phi(\alpha)=-\min_{\varphi\in Q(\alpha)}I(\varphi),
\qquad 0\leq\alpha\leq K,
\label{eq:intro_terminal_variational}
\end{equation}
is continuous, nonpositive and has maximum zero.
The laws of $M_N/N$ satisfy a large
deviation principle with speed $N$ and good rate function
$-\phi$. Moreover,
\begin{equation}
\lim_{N\to\infty}\max_{r\in\{0,\ldots,KN\}}
\left|\frac1N\log\mathbb{P}(M_N=r)
-\phi\!\left(\frac rN\right)\right|=0,
\label{eq:intro_point_probabilities}
\end{equation}
and at the endpoints we have
\begin{equation}
\phi(0)=\log\pi_0(0),
\qquad \phi(K)=\log\pi_K(K).
\label{eq:intro_entropy_endpoints}
\end{equation}
In particular, for every sequence
$r_N\in\{0,\ldots,KN\}$ with $r_N/N\to\alpha$,
the limit of $N^{-1}\log\mathbb{P}(M_N=r_N)$
is $\phi(\alpha)$.
\end{mainresult}

The variational formula determines the logarithmic limit of point
probabilities, including states near the boundaries.
Theorem~\ref{thm:terminal_entropy} uses terminal contraction
and proves continuity by correcting the endpoint of the minimizing
path. Theorem~\ref{thm:point_probabilities} passes from
intervals to individual states by comparing histories with
nearby cumulative sums in Lemma~\ref{lem:point_mass_comparison}, which
controls weights and multiplicities. This extends Theorem~4 of
Ref.~\cite{Franchini_URNS}, with the argument in Subsection~3.2.1,
Eq.~(3.75)--(3.77). Theorem~\ref{thm:terminal_endpoints} derives
the endpoint values from the local costs at velocities $0$ and $K$.

The rate functions and limits in Theorems~\ref{thm:intro_path_ldp}
and~\ref{thm:intro_terminal_entropy} remain unchanged with an
initial history of fixed length $m$, including a random one with a law
independent of $N$, and with continuous,
strictly positive and normalized probabilities $\pi_{n,k}$ converging uniformly
to $\pi_k$, with no required rate.
Theorems~\ref{thm:finite_initial_history} and~\ref{thm:inhomogeneous_urn}
use exponential equivalence with the pushforward law that replaces
the prefix, the comparison of tails and the uniform comparison of logarithmic
densities with common initialization.
The maximum in the point probabilities remains restricted to attainable
states (Theorem~\ref{thm:point_probabilities}): for a prefix
of length $m$ and cumulative sum $j$, these are
$\{j,\ldots,j+K(N-m)\}$. These extensions correspond to
Corollaries~2 and~3 of Ref.~\cite{Franchini_URNS}, with the necessary
support conditions for pointwise limits.

\medskip
\noindent\emph{Scalar paths and increment counts.}
For $K=1$, $\gamma_1=\varphi$ and
$\gamma_0(\tau)=\tau-\varphi(\tau)$ determine the counts.
For $K>1$, a velocity in $(0,K)$ leaves $K-1$ frequency degrees
of freedom. Figure~\ref{fig:count_path_contraction} compares histories
of increments $1$ with histories alternating $0$ and $2$:
their scalar limits agree but their counts differ. Information is lost
only in the uniform limit; at finite $N$, the path determines every
increment. The minimizing counts in Theorem~\ref{thm:intro_path_ldp}
are nevertheless unique. Lemma~\ref{lem:local_frequency_minimum}
gives their exponential tilt, including its limits at velocities $0$ and $K$.

\afterpage{%
\newgeometry{left=20mm,right=20mm,top=18mm,bottom=18mm}
\begingroup
\singlespacing
\setlength{\abovecaptionskip}{10pt}
\setlength{\belowcaptionskip}{0pt}

\begin{figure}[p]
\centering
\includegraphics[width=\linewidth,height=0.72\textheight,keepaspectratio]{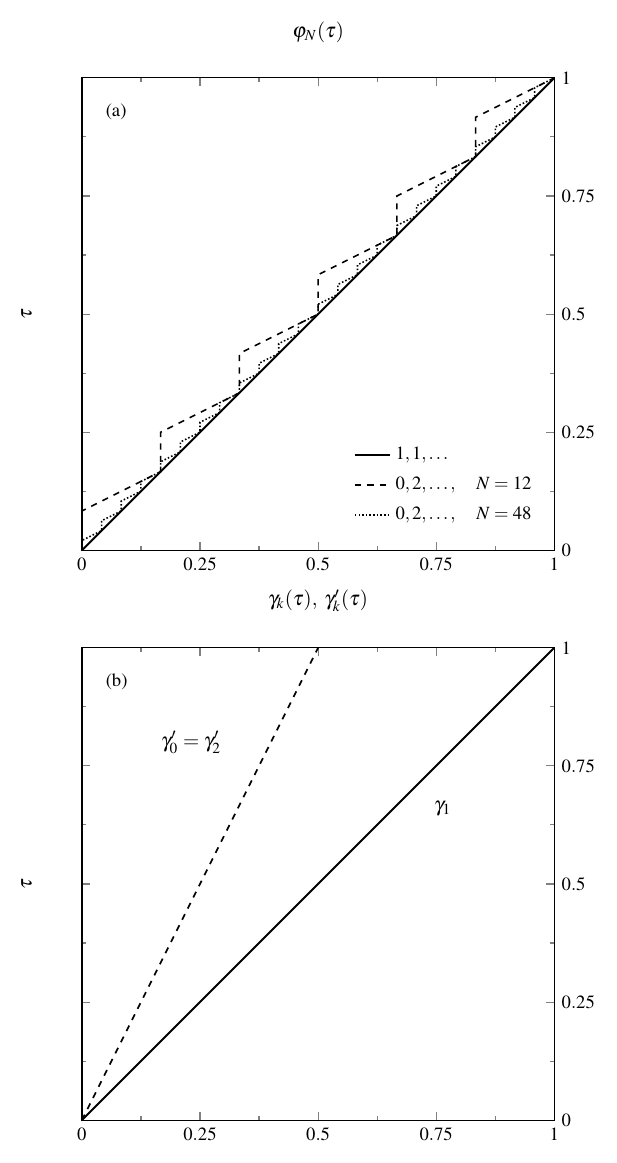}
\caption[Scalar paths and increment counts]{Different increment counts with the same limiting cumulative path, for $K=2$.
(a) Histories containing only increments $1$ give the line
$\varphi_N(\tau)=\tau$ (solid); histories alternating $0$ and $2$
give the dashed ($N=12$) and dotted ($N=48$) paths,
whose uniform distance from that line is $1/N$.
(b) Their limiting count vectors are $\gamma(\tau)=(0,\tau,0)$
and $\gamma'(\tau)=(\tau/2,0,\tau/2)$,
where the prime labels the second vector:
only the nonzero components are drawn, solid and dashed, respectively.
Both satisfy $\gamma_1+2\gamma_2=\gamma'_1+2\gamma'_2=\tau$;
time increases upward in both panels.
For the constant urn probabilities $\pi=(2/5,1/5,2/5)$,
the limits of $N^{-1}\log\mathbb P(\sigma)$ for these individual
histories are $\log(1/5)$ and $\log(2/5)$, respectively
\eqref{eq:distinct_count_actions}--\eqref{eq:style_s3_distinct_second_action}.
These are distinct from the count costs:
\eqref{eq:intro_count_rate} gives $I_{\mathrm c}(\gamma)=\log5$
and $I_{\mathrm c}(\gamma')=\log(5/4)$.
Minimization over all compatible counts in
Theorem~\ref{thm:intro_path_ldp} gives $I(\varphi)=0$
for their common scalar limit $\varphi(\tau)=\tau$.}
\label{fig:count_path_contraction}
\end{figure}
\clearpage
\endgroup
\restoretextgeometry
}

\subsection{Zero-cost paths}
\label{subsec:main_zero_cost}

For $\alpha\in[0,K]$, we call a path
$\varphi\in Q(\alpha)$ optimal if it attains the minimum in
\eqref{eq:intro_terminal_variational}.
We define the mean increment and its set
of fixed points, as in Section~\ref{sec:notation},
by
\begin{equation}
\bar\pi(\alpha)=\sum_{k\in\Omega}k\pi_k(\alpha),
\label{eq:intro_mean_fixed_points}
\end{equation}
\begin{equation}
C=\{\alpha\in[0,K]:\bar\pi(\alpha)=\alpha\}.
\label{eq:style_s2_intro_mean_fixed_points_2}
\end{equation}
Positivity implies that $C$ is a nonempty compact set
contained in $(0,K)$.

\begin{mainresult}[Zero-entropy interval and optimal paths]
\label{thm:intro_zero_cost}
For every $\alpha\in[0,K]$ we have
\begin{equation}
\phi(\alpha)=0
\qquad\Longleftrightarrow\qquad
\alpha\in[\inf C,\sup C].
\label{eq:intro_zero_entropy}
\end{equation}
For every final state in this interval,
the optimal paths are precisely those
$\varphi\in Q(\alpha)$ satisfying
\begin{equation}
\partial_\tau\varphi(\tau)
=\bar\pi\!\left(\frac{\varphi(\tau)}{\tau}\right)
\qquad\text{for almost every }\tau\in(0,1).
\label{eq:intro_zero_cost_equation}
\end{equation}
Their minimizing counts are uniquely
determined by
\begin{equation}
\gamma_k^*(\tau)
=\int_0^\tau\pi_k\!\left(\frac{\varphi(s)}s\right)ds,
\qquad k\in\Omega,\qquad 0\leq\tau\leq1.
\label{eq:intro_zero_cost_counts}
\end{equation}
The mean $\psi$ is nondecreasing or nonincreasing,
remains in $[\inf C,\sup C]$ and converges to a point of
$C$ as $\tau\downarrow0$.
Outside $[\inf C,\sup C]$ we have $\phi(\alpha)<0$.

If $\alpha\in C$, there exists an optimal path
in $Q(\alpha)$ other than $\varphi(\tau)=\alpha\tau$
if and only if at least one of the following two alternatives holds.
The first is that there exists $u\in[\inf C,\alpha)$ such that
\begin{equation}
\bar\pi(v)-v\geq0,\qquad u\leq v\leq\alpha,
\label{eq:intro_fixed_point_arrival_left}
\end{equation}
\begin{equation}
\int_u^\alpha\frac{dv}{\bar\pi(v)-v}<\infty.
\label{eq:style_s2_intro_fixed_point_arrival_left_2}
\end{equation}
The second is that there exists $u\in(\alpha,\sup C]$ such that
\begin{equation}
v-\bar\pi(v)\geq0,\qquad \alpha\leq v\leq u,
\label{eq:intro_fixed_point_arrival_right}
\end{equation}
\begin{equation}
\int_\alpha^u\frac{dv}{v-\bar\pi(v)}<\infty.
\label{eq:style_s2_intro_fixed_point_arrival_right_2}
\end{equation}
The integrals are understood as extended Lebesgue
integrals, with the reciprocal of zero equal to $+\infty$.
The straight line is therefore the unique optimal path
when neither condition is satisfied.
\end{mainresult}

The mean $\bar\pi$ determines the entire zero-entropy interval
and its optimal paths; the counts instead depend on all
components of $\pi$. Lemma~\ref{lem:zero_cost_paths} uses
equality in relative entropy, and
Theorem~\ref{thm:zero_entropy_interval} constructs the solutions
by separating variables in $\tau\partial_\tau\psi=\bar\pi(\psi)-\psi$
between fixed points. These extend Lemma~18 and Corollary~5
of Ref.~\cite{Franchini_URNS}.
Theorems~\ref{thm:component_optimal_paths}
and~\ref{thm:fixed_point_component_arrival} describe crossings
of components and waiting at fixed points, extending
Corollaries~6 and~7 of the same source.
Theorem~\ref{thm:fixed_point_global_uniqueness} and
Lemma~\ref{lem:fixed_point_integral_characterization} give the general
criterion without assuming that $C$ is finite. The integrals in
\eqref{eq:style_s2_intro_fixed_point_arrival_left_2} and
\eqref{eq:style_s2_intro_fixed_point_arrival_right_2} measure the arrival time
at the fixed point in the variable $\log\tau$; after arrival the mean
may remain constant.

\medskip
\noindent\emph{Smooth and non-Lipschitz examples.}
The zero-entropy interval, monotone means, and possible waiting,
departure and arrival at fixed points already occur in the binary
theory \cite[Corollaries~5--7 and Figure~2.2]{Franchini_URNS}.
In the smooth cubic example of Figure~\ref{fig:cubic_zero_cost},
every endpoint in $[1/2,3/2]$ admits a zero-cost path, but nonconstant
means neither wait at nor reach a fixed point at positive times.
The continuous, non-Lipschitz example in
Figure~\ref{fig:finite_arrival_paths} instead permits waiting at $1$,
motion to $1/2$ or $3/2$, and remaining there after arrival
(Subsection~\ref{subsec:nonlipschitz_example}).
Neither urn has explicitly time-dependent probabilities.

\medskip
\noindent\emph{Turning proportions and accumulated composition.}
The additional behavior for $K>1$ concerns the accumulated
composition. For a zero-cost path and its minimizing counts, set
$q_k(\tau)=\gamma_k^*(\tau)/\tau$ for $0<\tau\leq1$.
Equation~\eqref{eq:intro_zero_cost_counts} gives
\begin{equation}
\tau\partial_\tau q_k(\tau)=\pi_k(\psi(\tau))-q_k(\tau),\qquad k\in\Omega,
\label{eq:main_count_proportion_dynamics}
\end{equation}
\begin{equation}
\sum_{k\in\Omega}q_k(\tau)=1,\qquad\sum_{k\in\Omega}kq_k(\tau)=\psi(\tau).
\label{eq:style_s2_main_count_proportion_dynamics_2}
\end{equation}
A proportion changes direction when its current increment
probability crosses its accumulated value.
For $K=1$, $q_1=\psi$ and $q_0=1-\psi$, so monotonicity of the
zero-cost mean excludes such a turn. For $K>1$, individual
proportions can turn while their weighted mean remains monotone.
In the left-branch family
\eqref{eq:nonlipschitz_counts}--\eqref{eq:style_s9_nonlipschitz_counts_line3},
$q_0$ first increases and then decreases, while $q_1$ first decreases
and then increases. Figure~\ref{fig:component_memory} shows these turns, crossings
between family members, and distinct terminal proportions. The cumulative counts themselves
remain strictly increasing because
$\partial_\tau\gamma_k^*=\pi_k(\psi)>0$.
The binary exclusion is specific to zero-cost paths; arbitrary binary
paths can already have nonmonotone proportions.

Counts can distinguish zero-cost scalar paths after coalescence.
If two paths of the same urn coincide on $[\tau_0,1]$, with
$0<\tau_0<1$, their count derivatives agree there. For constants
$d_k$ and every $\tau\in[\tau_0,1]$, therefore,
\begin{equation}
\gamma_k^*(\tau)-\widetilde\gamma_k^*(\tau)=d_k,\qquad q_k(\tau)-\widetilde q_k(\tau)=d_k/\tau,
\label{eq:main_count_coalescence_offsets}
\end{equation}
\begin{equation}
\sum_{k\in\Omega}d_k=0,\qquad\sum_{k\in\Omega}kd_k=0.
\label{eq:style_s2_main_count_coalescence_offsets_2}
\end{equation}
For $K=1$ every offset vanishes; for $K>1$,
\eqref{eq:nonlipschitz_count_memory} gives nonzero offsets retaining
earlier composition. The discrete process remains Markovian in
$(n,M_n)$: equal states at equal times give the same subsequent law.
Finite-time coalescence requires failure of backward uniqueness;
the example uses a non-Lipschitz mean. Turning proportions do not
require coalescence.

\afterpage{%
\newgeometry{left=20mm,right=20mm,top=18mm,bottom=18mm}
\begingroup
\singlespacing
\setlength{\abovecaptionskip}{10pt}
\setlength{\belowcaptionskip}{0pt}

\begin{figure}[p]
\centering
\includegraphics[width=\linewidth,height=0.72\textheight,keepaspectratio]{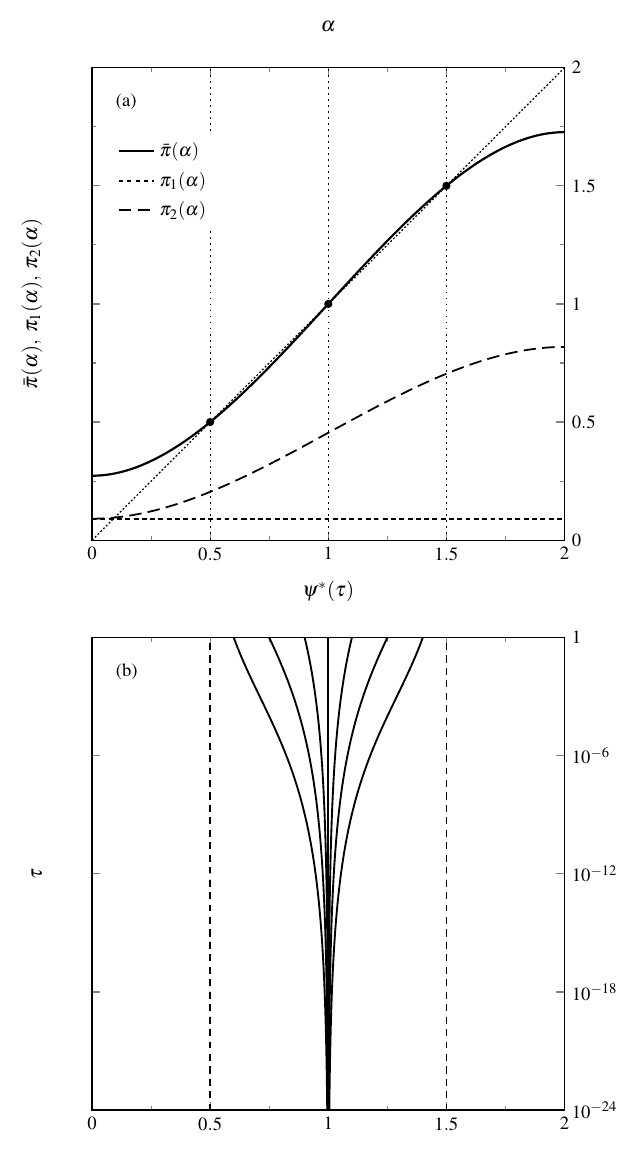}
\caption[Fixed points and zero-cost paths]{Zero-cost paths for the cubic urn with $K=2$ and probabilities
$\pi_0(\alpha)=(9-6\alpha^2+2\alpha^3)/11$,
$\pi_1(\alpha)=1/11$ and
$\pi_2(\alpha)=(1+6\alpha^2-2\alpha^3)/11$ on $[0,2]$.
(a) The mean increment
$\bar\pi(\alpha)=(3+12\alpha^2-4\alpha^3)/11$ (solid)
is shown together with $\pi_1$ (short dashes) and $\pi_2$ (long dashes).
The intersections of the mean with the dotted diagonal give the fixed
points $1/2$, $1$ and $3/2$.
(b) The optimal means $\psi^*(\tau)=\varphi^*(\tau)/\tau$
from \eqref{eq:cubic_optimal_average}, ending at
$\alpha=1/2,\,0.6,\,0.75,\,0.9,\,1,\,1.1,\,1.25,\,1.4,\,3/2$.
Time increases upward on a logarithmic scale from $10^{-24}$ to $1$;
the horizontal axes of the two panels are aligned.
The three fixed points give vertical lines, with the two outer lines dashed.
Every nonconstant mean tends to $1$ as $\tau\downarrow0$
and moves monotonically toward its terminal value:
it neither waits at nor reaches a fixed point at a positive time.
The logarithmic scale resolves the variation governed by $\tau^{1/11}$.
Thus every endpoint in $[1/2,3/2]$ admits a zero-cost path,
although only three endpoints are fixed points.
The exact formulas and uniqueness are proved in
Theorem~\ref{thm:cubic_zero_cost_paths}, illustrating
Theorem~\ref{thm:intro_zero_cost}.}
\label{fig:cubic_zero_cost}
\end{figure}
\clearpage

\begin{figure}[p]
\centering
\includegraphics[width=\linewidth,height=0.72\textheight,keepaspectratio]{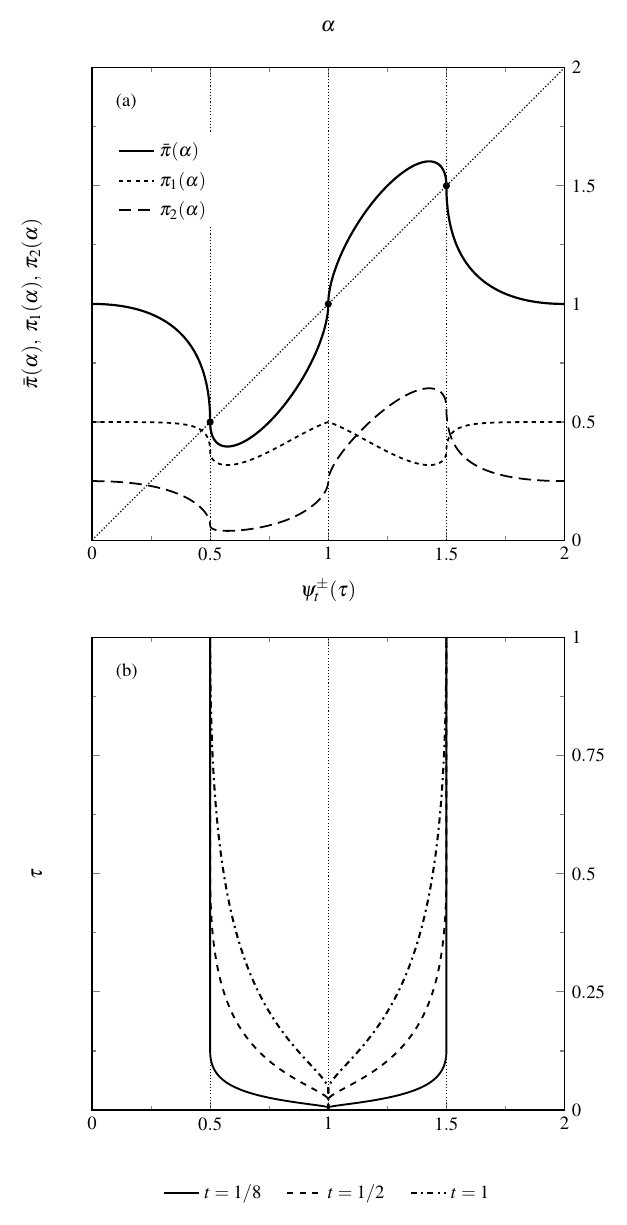}
\caption[Waiting and finite arrival for a fixed urn]{Waiting, motion and arrival at a fixed point for a fixed urn with $K=2$.
The probabilities are binomial:
$\pi_0=(1-\bar\pi/2)^2$, $\pi_1=\bar\pi(1-\bar\pi/2)$ and
$\pi_2=(\bar\pi/2)^2$, with the mean function
\eqref{eq:nonlipschitz_mean}.
This mean rescales the binary example in Figure~2.2 of
Ref.~\cite{Franchini_URNS}; all three probabilities are continuous
and bounded below by $((3-\sqrt2)/8)^2$.
(a) The mean $\bar\pi$ (solid) is shown together with
$\pi_1$ (short dashes) and $\pi_2$ (long dashes), as in Figure~\ref{fig:cubic_zero_cost}.
The intersections of the mean with the dotted diagonal give the
fixed points $1/2$, $1$ and $3/2$.
(b) The exact optimal means $\psi_t^\pm$ in
\eqref{eq:nonlipschitz_paths}--\eqref{eq:style_s9_nonlipschitz_paths_4}, with time increasing upwards.
The solid, dashed and dash-dotted curves correspond to arrival
times $t=1/8$, $1/2$ and $1$, respectively.
Each mean waits at $1$ until $s=e^{-\pi}t$, then moves to
$1/2$ on the left branch or $3/2$ on the right branch.
It reaches that point at $t$ and remains there afterwards.
Finite reciprocal-drift integrals permit both departure and arrival
at positive times, as in
Theorem~\ref{thm:fixed_point_component_arrival}.
The cumulative paths $\varphi_t^\pm(\tau)=\tau\psi_t^\pm(\tau)$
are exact zero-cost minimizers.
The urn probabilities themselves have no time dependence.
The two panels have the same horizontal scale, and dotted vertical
lines mark the fixed points.}
\label{fig:finite_arrival_paths}
\end{figure}
\clearpage

\begin{figure}[p]
\centering
\includegraphics[width=\linewidth,height=0.72\textheight,keepaspectratio]{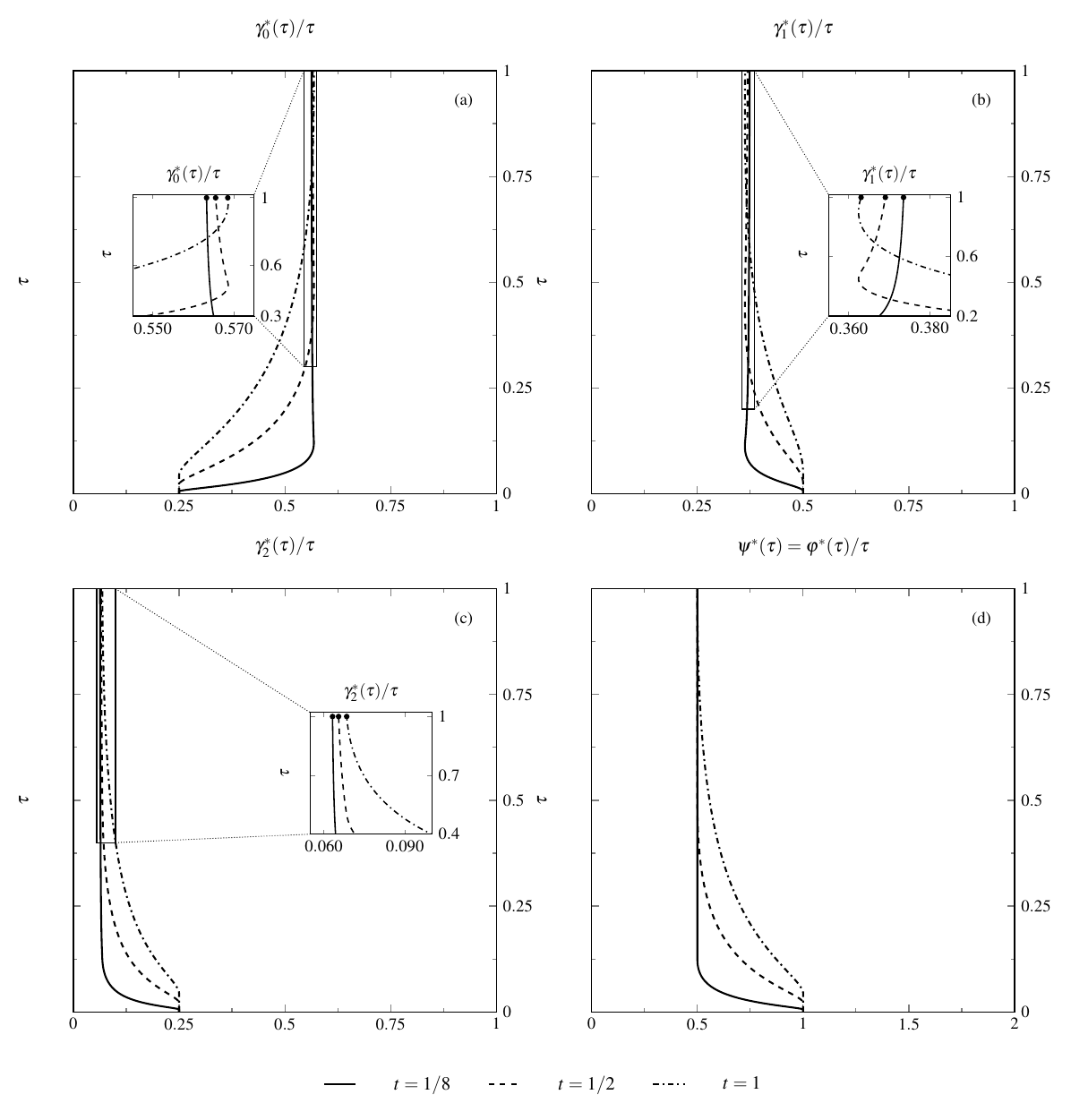}
\caption[Increment counts and memory of the earlier path]{Increment counts can distinguish optimal paths whose scalar means
have already merged.
For the fixed binomial urn in Figure~\ref{fig:finite_arrival_paths},
panels (a)--(c) show the cumulative count proportions
$\gamma_{k,t}^*(\tau)/\tau$, for $k=0,1,2$; panel (d) shows
$\psi_t^-(\tau)=\varphi_t^-(\tau)/\tau$.
The exact formulas are
\eqref{eq:nonlipschitz_paths}--\eqref{eq:style_s9_nonlipschitz_paths_4} and \eqref{eq:nonlipschitz_counts}--\eqref{eq:style_s9_nonlipschitz_counts_line3}.
The solid, dashed and dash-dotted curves use the same arrival times
$t=1/8$, $1/2$ and $1$ in every panel, with time increasing upwards.
All means leave $1$ at $s=e^{-\pi}t$ and reach $1/2$ at $t$.
The counts satisfy $\sum_k\gamma_{k,t}^*/\tau=1$ and
$\psi_t^-=\sum_k k\gamma_{k,t}^*/\tau$; they are cumulative
proportions, not the instantaneous probabilities
$\pi_k(\psi_t^-)=\partial_\tau\gamma_{k,t}^*$.
The square insets magnify the marked rectangles, revealing crossings
and changes of direction in (a)--(b), and curvature in (c).
The mean is monotone, whereas the individual proportions can turn and cross.
Their different time ranges include the relevant features of each
component; dots mark the distinct endpoints at $\tau=1$.
After arrival, $\psi_t^-=1/2$, while the normalized count vector is
$(9/16,3/8,1/16)+(t-s)(1,-2,1)/(160\tau)$.
At $\tau=1$, later arrivals give larger proportions for $k=0,2$
and a smaller proportion for $k=1$.
Thus the counts retain a contribution from the earlier motion,
even on the common terminal portion of the scalar path, as in
\eqref{eq:nonlipschitz_count_memory}.
The values shown at $\tau=0$ are the continuous limits
$(1/4,1/2,1/4)$ and $1$.}
\label{fig:component_memory}
\end{figure}
\clearpage
\endgroup
\restoretextgeometry
}

\subsection{Cumulant generating function and duality}
\label{subsec:main_cumulant_duality}

For every real parameter $\lambda$, define
$\zeta_N(\lambda)=N^{-1}\log\mathbb E[e^{\lambda M_N}]$,
where $\mathbb E$ is expectation under the law $\mathbb P$.
This is the cumulant generating function normalized
by the number of steps. The following result determines
its limit and connects it to the terminal entropy.

\begin{mainresult}[Cumulant generating function and limiting equation]
\label{thm:intro_cumulant}
Under the general assumptions on the urn functions, for every
$\lambda\in\mathbb{R}$ the following limit exists
\begin{equation}
\zeta(\lambda)
=\lim_{N\to\infty}\zeta_N(\lambda)
=\max_{\alpha\in[0,K]}
\{\lambda\alpha+\phi(\alpha)\}.
\label{eq:intro_cumulant_variational}
\end{equation}
Convergence is uniform on compact subsets of $\mathbb{R}$.
The function $\zeta$ is strictly convex, nondecreasing
and $K$-Lipschitz, with $\zeta(0)=0$.
We denote the left and right derivatives by $\zeta'_-$ and $\zeta'_+$,
respectively. For every $\lambda\in\mathbb R$ and for either choice of sign, we have
\begin{equation}
e^{\zeta(\lambda)}=\sum_{k\in\Omega}\pi_k\bigl(\zeta'_\pm(\lambda)\bigr)e^{k\lambda}.
\label{eq:intro_cumulant_equation}
\end{equation}
Moreover,
\begin{equation}
0<\zeta'_-(\lambda)\leq\zeta'_+(\lambda)<K,\qquad \lambda\in\mathbb R,
\label{eq:intro_cumulant_boundaries}
\end{equation}
\begin{equation}
\zeta'_-(0)=\inf C,\qquad\zeta'_+(0)=\sup C,
\label{eq:style_s2_intro_cumulant_boundaries_2}
\end{equation}
\begin{equation}
\lim_{\lambda\to-\infty}\zeta(\lambda)=\log\pi_0(0),
\label{eq:style_s2_intro_cumulant_boundaries_3}
\end{equation}
\begin{equation}
\lim_{\lambda\to+\infty}\{\zeta(\lambda)-K\lambda\}=\log\pi_K(K).
\label{eq:style_s2_intro_cumulant_boundaries_4}
\end{equation}
In particular, $\zeta$ is differentiable at zero if and
only if $C$ contains a single point.
\end{mainresult}

The variational formula derives the cumulant generating function from the entropy;
the slopes at zero identify the endpoints of the zero-entropy segment.
The proofs are in Theorems~\ref{thm:cumulant_variational}
and~\ref{thm:cumulant_limit_equation}, in
Lemma~\ref{lem:cumulant_boundary_values} and in
\eqref{eq:cumulant_derivatives_at_zero}--\eqref{eq:style_s7_004}.
The proof of the equation uses Lemma~\ref{lem:cumulant_changed_ldp}:
weighting histories by $e^{\lambda M_N}$ and normalizing,
the terminal state concentrates at $\zeta'(\lambda)$ at points
of differentiability. Summing the recursion
\eqref{eq:cumulant_exact_recursion} yields
\eqref{eq:cumulant_recursion_telescope}, developing
Ref.~\cite{Franchini_URNS}, Subsections~3.3 and~3.3.1 and Lemma~20.
The limit along points of differentiability gives the equation
for the one-sided derivatives; \eqref{eq:cumulant_fixed_state_strict_convexity}
rules out affine segments and proves strict convexity.

To determine the cumulant generating function from its equation and derive
the entropy from the cumulant generating function, we add the following injectivity
assumptions. They require that, for a fixed parameter $\lambda$,
distinct states give distinct values of the exponential sum,
on the interval specified for each sign of $\lambda$.

\begin{mainresult}[Uniqueness of the cumulant generating function and duality]
\label{thm:intro_cumulant_duality}
Under the general assumptions, suppose that the map
\begin{equation}
\alpha\longmapsto
\sum_{k\in\Omega}\pi_k(\alpha)e^{k\lambda}
\label{eq:intro_cumulant_injectivity_map}
\end{equation}
is injective on $[0,\inf C]$ for every $\lambda<0$
and on $[\sup C,K]$ for every $\lambda>0$.
Then $\zeta\in C^1(\mathbb{R}\setminus\{0\})$ and is
the unique function
$f\in C(\mathbb{R})\cap C^1(\mathbb{R}\setminus\{0\})$
satisfying
\begin{equation}
e^{f(\lambda)}=\sum_{k\in\Omega}\pi_k\bigl(f'(\lambda)\bigr)e^{k\lambda},\qquad\lambda\neq0,
\label{eq:intro_cumulant_boundary_problem}
\end{equation}
\begin{equation}
0\leq f'(\lambda)\leq\inf C,\qquad\lambda<0,
\label{eq:style_s2_intro_cumulant_boundary_problem_2}
\end{equation}
\begin{equation}
\sup C\leq f'(\lambda)\leq K,\qquad\lambda>0.
\label{eq:style_s2_intro_cumulant_boundary_problem_3}
\end{equation}
The equation and the global constraints imply $f(0)=0$ and
\begin{equation}
\lim_{\lambda\to-\infty}f(\lambda)=\log\pi_0(0),
\label{eq:style_s2_candidate_left_limit}
\end{equation}
\begin{equation}
\lim_{\lambda\to+\infty}\{f(\lambda)-K\lambda\}=\log\pi_K(K),
\label{eq:style_s2_candidate_right_limit}
\end{equation}
without assuming convexity of $f$.
The derivative $\zeta'$ is continuous and strictly increasing
on each half-line; its values range over $(0,\inf C)$
for $\lambda<0$ and $(\sup C,K)$ for $\lambda>0$. For every $\alpha\in[0,K]$,
\begin{equation}
\phi(\alpha)
=\inf_{\lambda\in\mathbb{R}}
\{\zeta(\lambda)-\lambda\alpha\}.
\label{eq:intro_entropy_duality}
\end{equation}
The function $\phi$ is concave and absolutely continuous
on $[0,K]$, of class $C^1$ on $(0,K)$ and strictly
concave on $[0,\inf C]$ and $[\sup C,K]$.
Concavity is strict on all of $[0,K]$ if and only
if $C$ contains a single point. Moreover,
\begin{equation}
\phi'\bigl(\zeta'(\lambda)\bigr)=-\lambda,
\qquad \lambda\neq0.
\label{eq:intro_entropy_parametric_derivative}
\end{equation}
For $0<\alpha<K$ the dual infimum is a unique minimum;
for $\alpha=0$ and $\alpha=K$ it is obtained only in the
limits $\lambda\to-\infty$ and $\lambda\to+\infty$.
\end{mainresult}

Theorems~\ref{thm:cumulant_boundary_uniqueness} and~\ref{thm:entropy_shape}
also prove the two separate results: injectivity for $\lambda<0$
suffices for uniqueness on $(-\infty,0)$ and duality on $[0,\inf C]$;
injectivity for $\lambda>0$ suffices for uniqueness on $(0,\infty)$
and duality on $[\sup C,K]$.
The slope constraints, imposed on the entire half-line,
also make the limits at the origin and at infinity follow from the equation.
Injectivity identifies the two one-sided derivatives of the cumulant generating function.
For uniqueness, two solutions are compared at a point where
their difference has a nonzero maximum or minimum,
as in \eqref{eq:cumulant_solution_comparison}.
The ranges of $\zeta'$ and
Lemma~\ref{lem:entropy_dual_equality} give the inverse formula
and hence concavity of the entropy.
Lemma~\ref{lem:cumulant_probability_sums} provides a sufficient
criterion for injectivity through sums of probabilities.
For $K=1$ we recover Theorem~9 and Corollary~11
of Ref.~\cite{Franchini_URNS}.

\subsection{The inverse problem}
\label{subsec:main_inverse}

The inverse problem prescribes the terminal entropy and seeks the urn
probabilities that realize it. The following characterization and
construction require no injectivity assumption.

\begin{mainresult}[Characterization and construction of concave entropies]
\label{thm:intro_inverse_entropy}
Let $f:[0,K]\to\mathbb{R}$ be continuous and concave.
There exist continuous urn functions, strictly
positive on $[0,K]$ and normalized on $\Omega$, whose
terminal entropy is $f$, if and only if
\begin{equation}
f\in C^1((0,K)),\qquad \max_{\alpha\in[0,K]}f(\alpha)=0,
\label{eq:intro_inverse_conditions}
\end{equation}
\begin{equation}
\lim_{\alpha\downarrow0}f'(\alpha)=+\infty,\qquad\lim_{\alpha\uparrow K}f'(\alpha)=-\infty.
\label{eq:style_s2_intro_inverse_conditions_2}
\end{equation}
When these conditions hold, we define the scalar parameter
$\pi_f$ by
\begin{equation}
\pi_f(\alpha)=\frac{e^{[f(\alpha)-\alpha f'(\alpha)]/K}-1}{e^{-f'(\alpha)}-1},\qquad 0<\alpha<K,\qquad f'(\alpha)\neq0,
\label{eq:intro_inverse_parameter}
\end{equation}
\begin{equation}
\pi_f(\alpha)=\alpha/K,\qquad 0<\alpha<K,\qquad f'(\alpha)=0,
\label{eq:style_s2_intro_inverse_parameter_2}
\end{equation}
\begin{equation}
\pi_f(0)=1-e^{f(0)/K},\qquad\pi_f(K)=e^{f(K)/K}.
\label{eq:style_s2_intro_inverse_parameter_3}
\end{equation}
An explicit realization is given by the binomial probabilities
\begin{equation}
\pi_k(\alpha)
=\binom{K}{k}\pi_f(\alpha)^k
\bigl(1-\pi_f(\alpha)\bigr)^{K-k},
\qquad k\in\Omega.
\label{eq:intro_inverse_binomial}
\end{equation}
The function $\pi_f$ is continuous and takes values in $(0,1)$.
The fixed-point set of the urn thus constructed is
$C=\{\alpha\in[0,K]:f(\alpha)=0\}$.
More generally, for a function $f$ satisfying
\eqref{eq:intro_inverse_conditions}--\eqref{eq:style_s2_intro_inverse_conditions_2}, every family
of urn functions under the general assumptions satisfying
\begin{equation}
e^{f(\alpha)-\alpha f'(\alpha)}
=\sum_{k\in\Omega}\pi_k(\alpha)e^{-k f'(\alpha)},
\qquad 0<\alpha<K,
\label{eq:intro_inverse_equation}
\end{equation}
has terminal entropy equal to $f$.
\end{mainresult}

Strict convexity of the cumulant generating function excludes interior
corners and finite endpoint slopes of a realizable $f$
(Theorem~\ref{thm:inverse_binomial_construction}).
Theorem~\ref{thm:inverse_verification} proves
$I(\varphi)\geq-f(\varphi(1))$ in
\eqref{eq:inverse_action_lower_bound} and constructs an attaining path
through \eqref{eq:inverse_optimal_mean_equation}--\eqref{eq:style_s7_034}.
For $K=1$, the construction recovers
\cite[Corollary~10, Eq.~(2.34)]{Franchini_URNS}.

At each fixed state where $f'\neq0$, the compatible probability vectors
have dimension $K-1$ (Lemma~\ref{lem:inverse_probability_freedom}).
For $K>1$, perturbations
\eqref{eq:inverse_global_perturbation}--\eqref{eq:style_s7_048}
give distinct continuous, uniformly positive urns on $[0,K]$
with the same entropy.
For $K=1$, the inverse equation fixes the urn where $f'\neq0$
and reduces to normalization on an interval of maxima.
Continuity forces $\bar\pi(c)=c$ at any zero $c$ approached by states
with $f'\neq0$. For strictly concave entropy, each interior state is
the unique maximizer at its dual parameter; thus
\eqref{eq:intro_cumulant_equation} makes the inverse equation necessary
and determines the binary urn by continuity.

\medskip
\noindent\emph{The same entropy with different path costs.}
The ternary family
\eqref{eq:three_increment_analytic_family}--\eqref{eq:style_s9_2_053},
$0<c<1/4$, satisfies the inverse equation for a common strictly
concave entropy
(Theorem~\ref{thm:three_increment_analytic_family}).
Figure~\ref{fig:shared_entropy_path_cost} shows that the cost of
\eqref{eq:three_increment_comparison_path} nevertheless increases
strictly with $c$. This test path ends at $1$ but has positive cost,
so it is not optimal. The optimal paths are discussed after
Theorem~\ref{thm:intro_affine}.

\afterpage{%
\newgeometry{left=20mm,right=20mm,top=18mm,bottom=18mm}
\begingroup
\singlespacing
\setlength{\abovecaptionskip}{10pt}
\setlength{\belowcaptionskip}{0pt}

\begin{figure}[p]
\centering
\includegraphics[width=\linewidth,height=0.72\textheight,keepaspectratio]{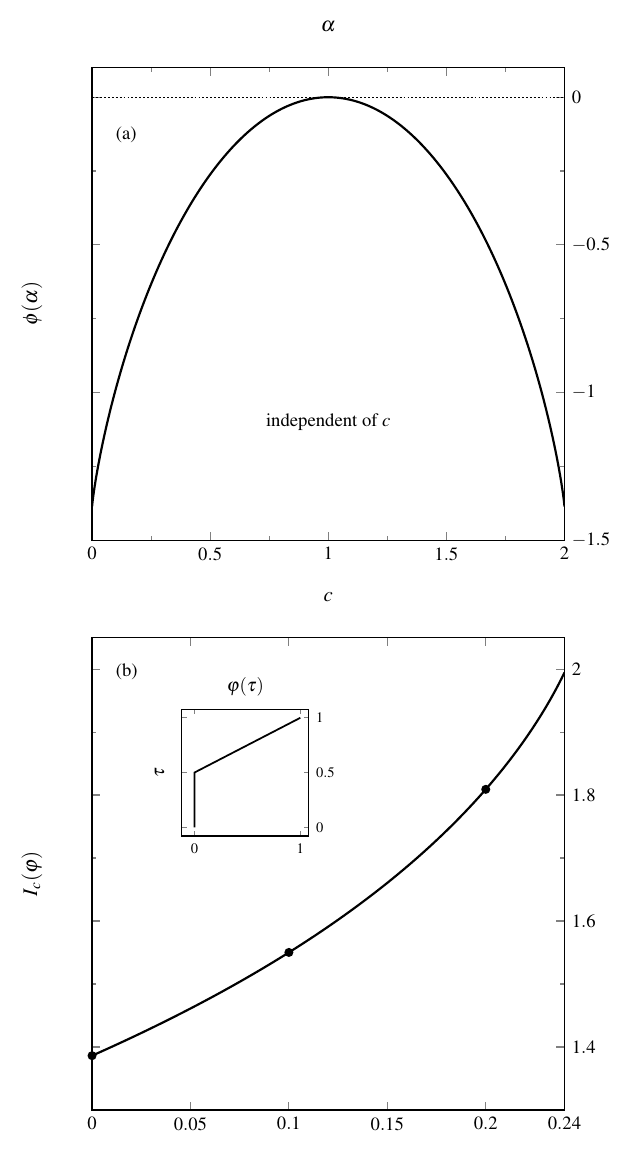}
\caption[Common terminal entropy and different path costs]{Identical terminal entropy with different path costs for the $K=2$ urns
$\pi_0(\alpha)=1/4-c\alpha/2$, $\pi_1(\alpha)=1/2+c$ and
$\pi_2(\alpha)=1/4-c+c\alpha/2$, where $0\leq\alpha\leq2$
and $0\leq c<1/4$.
(a) The terminal entropy
$\phi(\alpha)=-\alpha\log\alpha-(2-\alpha)\log(2-\alpha)$,
with $0\log0=0$, is independent of $c$
(Theorem~\ref{thm:three_increment_analytic_family}).
(b) The cost $I_c(\varphi)$ of the fixed path
$\varphi(\tau)=0$ for $0\leq\tau\leq1/2$ and
$\varphi(\tau)=2\tau-1$ for $1/2\leq\tau\leq1$
increases strictly with $c$, as given by
\eqref{eq:three_increment_comparison_path_cost}.
The inset shows this path with time increasing upward;
the marked costs correspond to $c=0,0.1,0.2$.
This path ends at $1$ but has positive cost, whereas
$\min_{\varphi\in Q(1)}I_c(\varphi)=-\phi(1)=0$:
it is a test path, not an optimal path.
The panel plots $0\leq c\leq0.24$; at $c=0$ the urn has
constant probabilities $(1/4,1/2,1/4)$.
Equality of the terminal entropies therefore does not determine
the cost of an individual path.}
\label{fig:shared_entropy_path_cost}
\end{figure}
\clearpage
\endgroup
\restoretextgeometry
}

\subsection{Affine urns}
\label{subsec:main_affine_examples}

Consider affine urn functions with real coefficients
$a_k,b_k$, for $k\in\Omega$, of the form
\begin{equation}
\pi_k(\alpha)=a_k+b_k\alpha,\qquad k\in\Omega,\qquad 0\leq\alpha\leq K,
\label{eq:intro_affine_probabilities}
\end{equation}
\begin{equation}
\sum_{k\in\Omega}a_k=1,\qquad\sum_{k\in\Omega}b_k=0,
\label{eq:style_s2_intro_affine_probabilities_2}
\end{equation}
\begin{equation}
a_k>0,\qquad a_k+Kb_k>0,\qquad k\in\Omega.
\label{eq:style_s2_intro_affine_probabilities_3}
\end{equation}
The inequalities at the two boundaries are equivalent to strict positivity
of the probabilities on all of $[0,K]$. For $\lambda\in\mathbb R$ define
the exponential sums
\begin{equation}
A(\lambda)=\sum_{k\in\Omega}a_ke^{k\lambda},
\qquad
B(\lambda)=\sum_{k\in\Omega}b_ke^{k\lambda}.
\label{eq:intro_affine_sums}
\end{equation}
These are the functions in \eqref{eq:affine_exponential_sums}--\eqref{eq:style_s8_exponential_slope_sum}.
The following result combines
Theorems~\ref{thm:affine_cumulant_equation},
\ref{thm:affine_integral_representation},
\ref{thm:affine_entropy}
and~\ref{thm:affine_second_derivatives}.

\begin{mainresult}[Cumulant generating function and entropy for affine urns]
\label{thm:intro_affine}
Under \eqref{eq:intro_affine_probabilities}--\eqref{eq:style_s2_intro_affine_probabilities_3},
the cumulant generating function is the unique solution
$\zeta\in C^1(\mathbb{R})$ of
\begin{equation}
e^{\zeta(\lambda)}
=A(\lambda)+B(\lambda)\zeta'(\lambda),
\qquad\lambda\in\mathbb R,
\label{eq:intro_affine_equation}
\end{equation}
with the conditions
\begin{equation}
\lim_{\lambda\to-\infty}\zeta(\lambda)=\log a_0,
\label{eq:intro_affine_boundaries}
\end{equation}
\begin{equation}
\lim_{\lambda\to+\infty}\{\zeta(\lambda)-K\lambda\}=\log(a_K+Kb_K).
\label{eq:style_s2_intro_affine_boundaries_2}
\end{equation}
It is strictly convex and
$\{\zeta'(\lambda):\lambda\in\mathbb{R}\}=(0,K)$.
If $B\not\equiv0$, let $J=(\lambda_-,\lambda_+)$
be a maximal interval containing no zeros of $B$, with
$-\infty\leq\lambda_-<\lambda_+\leq+\infty$.
For $\lambda\in J$, when $B>0$ on $J$, we have
\begin{equation}
e^{-\zeta(\lambda)}
=\int_\lambda^{\lambda_+}
\frac{1}{B(t)}
\exp\!\left\{-\int_\lambda^t
\frac{A(u)}{B(u)}\,du\right\}\,dt,
\label{eq:intro_affine_quadrature_positive}
\end{equation}
whereas, if $B<0$ on $J$,
\begin{equation}
e^{-\zeta(\lambda)}
=\int_{\lambda_-}^{\lambda}
\frac{1}{-B(t)}
\exp\!\left\{-\int_\lambda^t
\frac{A(u)}{B(u)}\,du\right\}\,dt.
\label{eq:intro_affine_quadrature_negative}
\end{equation}
The integrals converge absolutely,
even when the integration endpoint is infinite
or is a multiple zero of $B$.
At every zero $\lambda_0$ of $B$ the matching data are
\begin{equation}
\zeta(\lambda_0)=\log A(\lambda_0),
\qquad
\zeta'(\lambda_0)
=\frac{A'(\lambda_0)}
{A(\lambda_0)-B'(\lambda_0)},
\label{eq:intro_affine_matching}
\end{equation}
with $A(\lambda_0)-B'(\lambda_0)>0$.
Under the affine assumptions alone, duality
\eqref{eq:intro_entropy_duality} holds for every $\alpha\in[0,K]$
and identity \eqref{eq:intro_entropy_parametric_derivative}
holds for every $\lambda\in\mathbb R$, including zero.
The entropy is strictly concave and absolutely continuous
on $[0,K]$ and belongs to $C^2((0,K))$.
For $0<\alpha<K$ the minimum is unique and is characterized
by $\zeta'(\lambda)=\alpha$.
Away from the zeros of $B$, the cumulant generating function is real analytic,
$\zeta''>0$ and
\begin{equation}
\phi''\bigl(\zeta'(\lambda)\bigr)
=-\frac{1}{\zeta''(\lambda)}<0.
\label{eq:intro_affine_dual_derivatives}
\end{equation}
Assume that $B\not\equiv0$ and that $\lambda_0$ is a zero of $B$.
Then the curvature of $\zeta$ is finite and positive at
$\lambda_0$ if $2B'(\lambda_0)<A(\lambda_0)$;
it instead diverges to $+\infty$ when approaching
$\lambda_0$ if $2B'(\lambda_0)\geq A(\lambda_0)$.
The divergence is logarithmic in the equality case
and has exponent $A(\lambda_0)/B'(\lambda_0)-2$
in the strict inequality case.
In particular,
\begin{equation}
\phi''\bigl(\zeta'(\lambda_0)\bigr)=0
\qquad\Longleftrightarrow\qquad
2B'(\lambda_0)\geq A(\lambda_0).
\label{eq:intro_affine_zero_curvature}
\end{equation}
If $B\equiv0$, the probabilities are constant and
$\zeta=\log A$ is real analytic on $\mathbb{R}$
with $\zeta''>0$; the entropy
is real analytic on $(0,K)$ with $\phi''<0$.
\end{mainresult}

Theorem~\ref{thm:affine_cumulant_equation} proves uniqueness without
the injectivity assumptions of Theorem~\ref{thm:intro_cumulant_duality}.
The substitution $e^{-\zeta}$ linearizes the equation, and the sign
of $B$ selects the integration endpoint in
Theorem~\ref{thm:affine_integral_representation}.
Lemma~\ref{lem:affine_coefficient_zeros} supplies the positive
matching denominator; Theorems~\ref{thm:affine_entropy}
and~\ref{thm:affine_second_derivatives} give duality and curvature.
For $K=1$, these formulas recover Corollaries~11--12 of
Ref.~\cite{Franchini_URNS}.

\medskip
\noindent\emph{Affine examples and optimal paths.}
Figures~\ref{fig:constant_urn_paths} and~\ref{fig:binary_affine_paths}
compare constant probabilities with the binary affine urn.
In the constant case, every optimal mean is constant and the cumulative
path is a straight line. In the binary affine example, the optimal
means for atypical endpoints are curved and tend to $0$ or $1$ as
$\tau\downarrow0$. Only the path through the unique fixed point has
zero cost. Lemma~\ref{lem:affine_plot_paths} constructs the optimal
paths for all the affine examples, including those with positive cost.

For a nonconstant binary affine urn, $B(\lambda)=b_1(e^\lambda-1)$
has only the zero at the origin. With $K>1$, additional or multiple
zeros require the matching in Theorem~\ref{thm:intro_affine}.
In Figure~\ref{fig:extra_root_affine_paths}, the additional zero
$\lambda=\log2$ gives the constant optimal mean $40/27$ with
positive cost, distinct from the zero-cost mean $24/25$.
The cumulant generating function is $C^{27}$ but not $C^{28}$ there.
This loss of regularity does not involve divergent curvature:
Lemma~\ref{lem:three_increment_nonzero_roots} excludes such divergence
at nonzero zeros of $B$ for positive affine urns with $K=2$.

Figure~\ref{fig:double_root_affine_paths} shows the complementary
double-zero example. Its mean increment is constant, but its
probabilities depend on the state and its optimal means for atypical
endpoints are nonconstant. A constant mean increment therefore does
not imply the straight optimal paths of the constant-probability urn.
The quadratures and matching for both ternary examples are given in
Subsection~\ref{subsec:ternary_affine_examples}.

For $K=1$, $\bar\pi=\pi_1$ determines the urn. For $K=2$, the family
\eqref{eq:three_increment_double_root_probabilities}--\eqref{eq:style_s9_2_042}
has $\bar\pi\equiv1$ for every $0<c<1/12$, but its endpoint entropy
$\phi(2)=\log(1/3+2c)$ varies with $c$.
The same scalar zero-cost paths can therefore coexist with different
deviation costs.

Figure~\ref{fig:analytic_affine_paths} returns to the family with
common entropy in Figure~\ref{fig:shared_entropy_path_cost}.
It pairs the urn functions with the optimal means for $c=1/8$.
Formula~\eqref{eq:analytic_affine_optimal_mean} shows that these paths
also depend on $c$, although their minimum costs are unchanged.

\afterpage{%
\newgeometry{left=20mm,right=20mm,top=18mm,bottom=18mm}
\begingroup
\singlespacing
\setlength{\abovecaptionskip}{10pt}
\setlength{\belowcaptionskip}{0pt}

\begin{figure}[p]
\centering
\includegraphics[width=\linewidth,height=0.72\textheight,keepaspectratio]{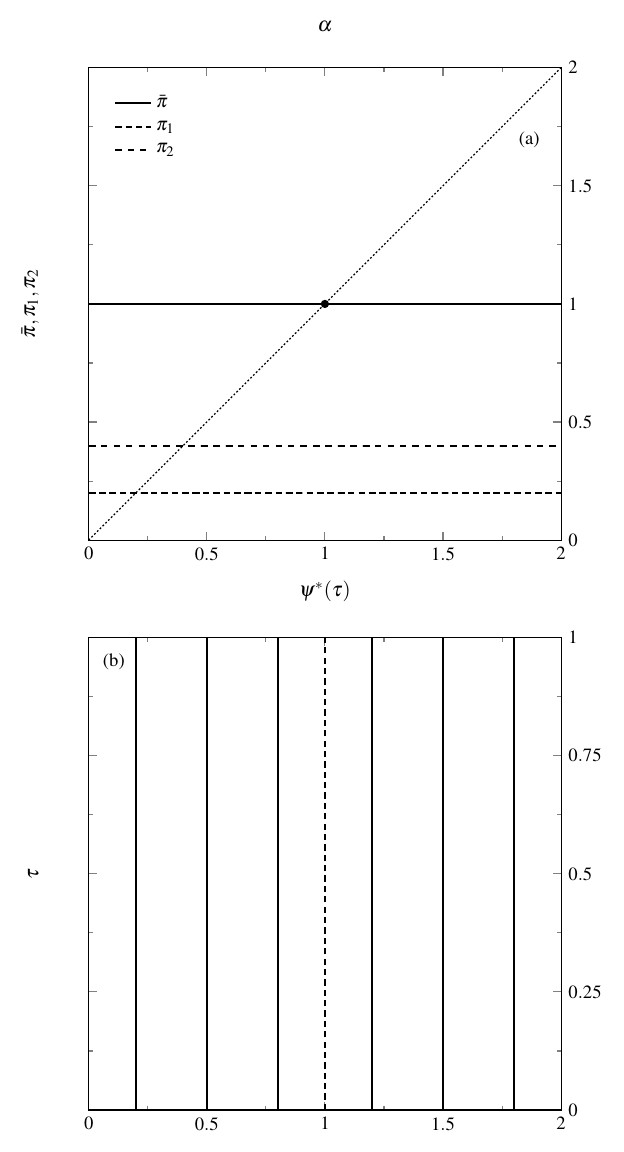}
\caption[Constant urn probabilities and optimal paths]{The constant urn
$\pi=(2/5,1/5,2/5)$ used in Figure~\ref{fig:count_path_contraction}.
(a) The probabilities $\pi_1=1/5$ and $\pi_2=2/5$, together with the
mean increment $\bar\pi=1$; the dotted diagonal identifies the unique
fixed point.
(b) Optimal means $\psi^*(\tau)=\alpha$ for terminal states
$\alpha=0.2,0.5,0.8,1,1.2,1.5,1.8$, with time increasing upward.
The corresponding cumulative paths are $\varphi^*(\tau)=\alpha\tau$.
Only the dashed line at $\alpha=1$ has zero cost; the other straight
paths minimize a positive cost under their respective terminal constraints.
Here $B\equiv0$, so the parameter in
Lemma~\ref{lem:affine_plot_paths} is constant along each path.
In contrast to Figure~\ref{fig:count_path_contraction}, which compares
two prescribed count histories with the same scalar limit, this panel
shows scalar minimizers for different final means.}
\label{fig:constant_urn_paths}
\end{figure}
\clearpage

\begin{figure}[p]
\centering
\includegraphics[width=\linewidth,height=0.72\textheight,keepaspectratio]{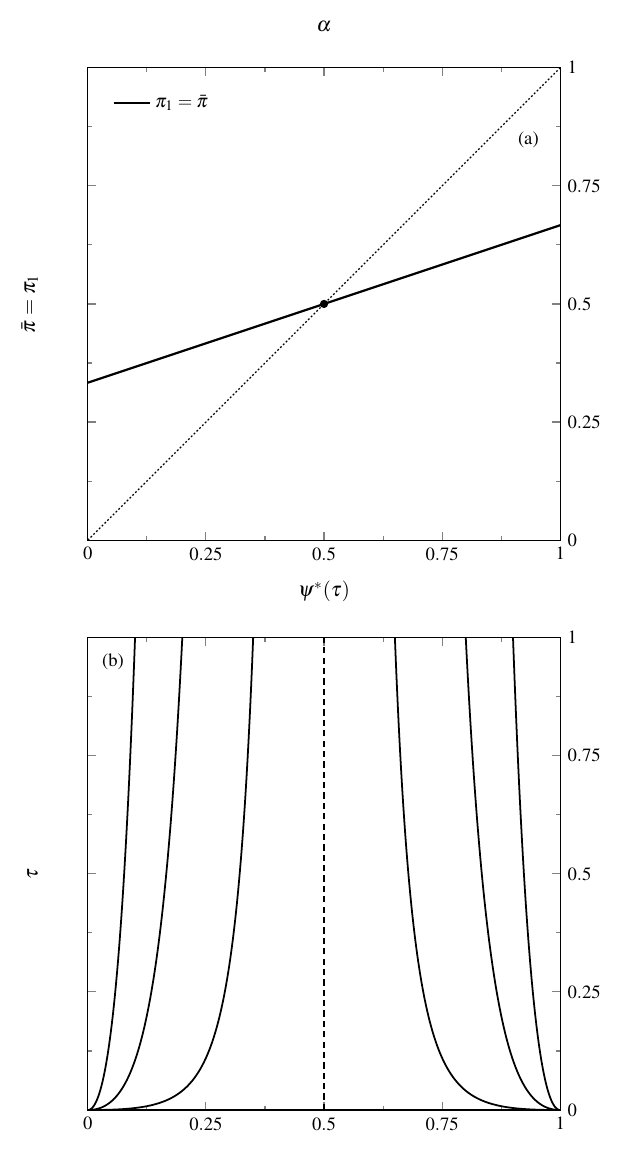}
\caption[Binary affine urn and optimal paths]{The binary affine example
$a=b=1/3$ in Subsection~\ref{subsec:binary_affine_examples}.
(a) The urn probability $\pi_1(\alpha)=(1+\alpha)/3$, which is also
the mean increment.
Its intersection with the dotted diagonal is the fixed point $1/2$.
(b) Optimal means ending at $\alpha=0.1,0.2,0.35,0.5,0.65,0.8,0.9$,
from \eqref{eq:binary_affine_optimal_mean} and its reflection.
Time increases upward, and the two panels share the horizontal scale.
The dashed mean $1/2$ corresponds to the unique zero-cost path.
For a terminal value below $1/2$, the optimal mean increases from $0$;
for one above $1/2$, it decreases from $1$.
Each solid curve minimizes a strictly positive cost for its own endpoint.
The values at $\tau=0$ are limits of the means; all cumulative paths
$\varphi^*(\tau)=\tau\psi^*(\tau)$ start at zero.}
\label{fig:binary_affine_paths}
\end{figure}
\clearpage

\begin{figure}[p]
\centering
\includegraphics[width=\linewidth,height=0.72\textheight,keepaspectratio]{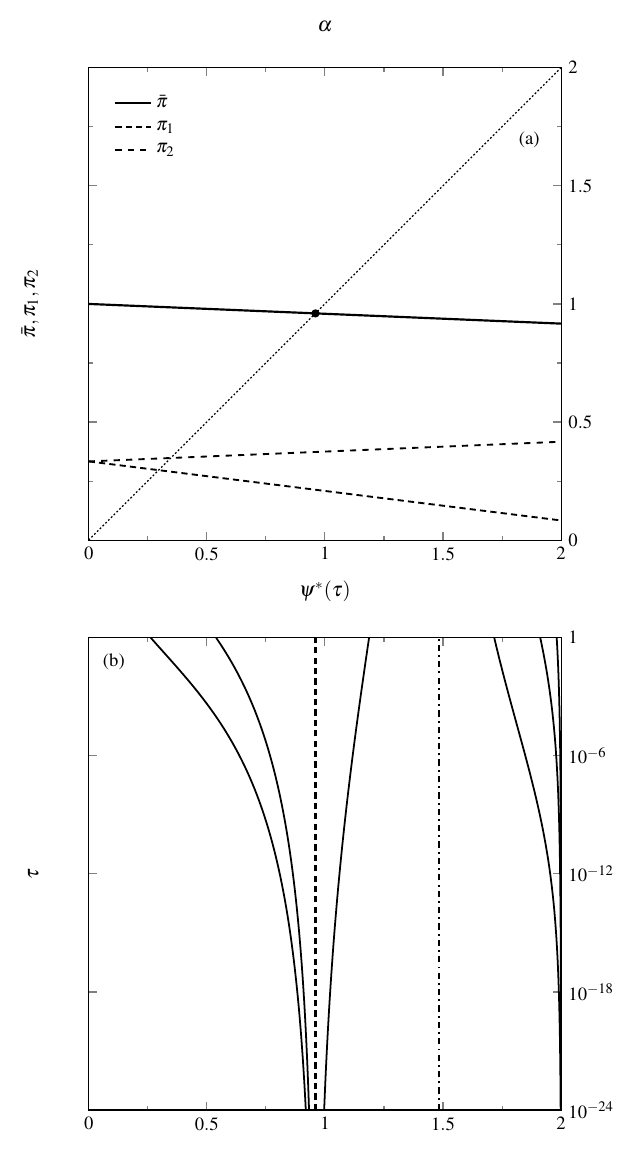}
\caption[Affine urn with an additional zero of the exponential slope sum]{
The ternary affine urn
\eqref{eq:three_increment_extra_root_probabilities}--\eqref{eq:style_s9_2_025}.
(a) The probabilities $\pi_1,\pi_2$ and the decreasing mean
$\bar\pi(\alpha)=1-\alpha/24$; the unique fixed point is $24/25$.
(b) Optimal means constructed by Lemma~\ref{lem:affine_plot_paths}
from the quadratures
\eqref{eq:three_increment_extra_root_lower_integral}--\eqref{eq:three_increment_extra_root_upper_integral}.
From left to right, the terminal means are approximately
$0.262,0.539,0.960,1.187,1.481,1.717,1.912,1.981$.
Time increases upward on a logarithmic scale from $10^{-24}$ to $1$.
The dashed line at $24/25$ has zero cost.
The dash-dotted line at $40/27$ is also an optimal constant mean, because
$B(\log2)=0$, but its cost is
$(40/27)\log2-\log(7/3)>0$.
Means ending below $40/27$ tend to $24/25$ as $\tau\downarrow0$;
those ending above it tend to $2$.
The slowly changing curves on the left are resolved by the logarithmic
time scale; the diagram does not identify zeros of $B$ with zero-cost states.}
\label{fig:extra_root_affine_paths}
\end{figure}
\clearpage

\begin{figure}[p]
\centering
\includegraphics[width=\linewidth,height=0.72\textheight,keepaspectratio]{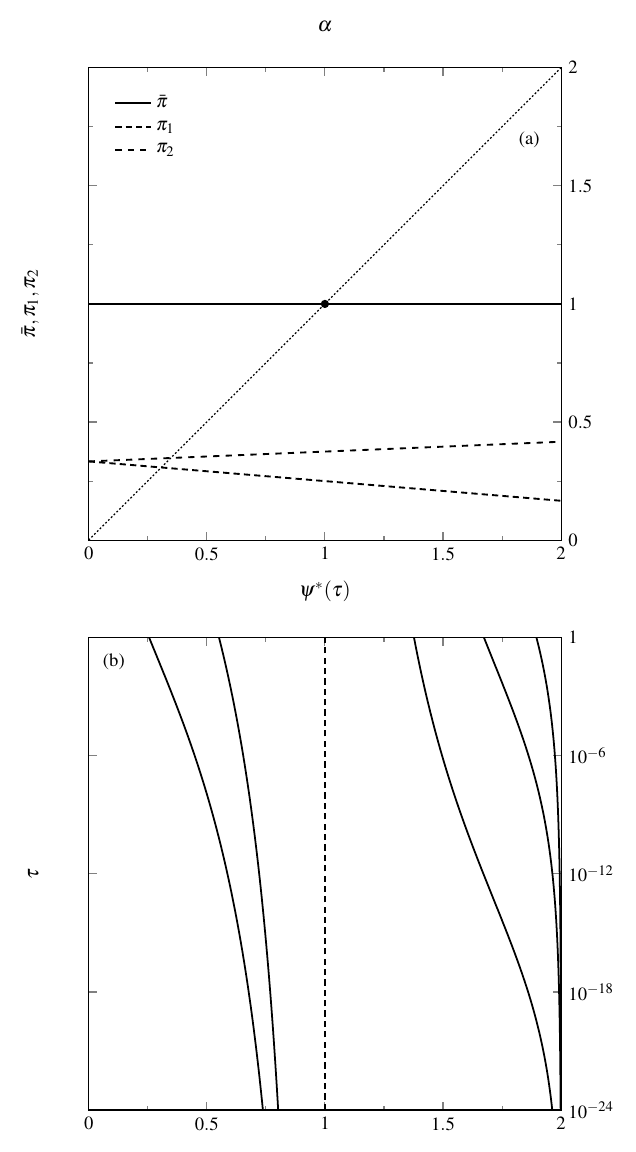}
\caption[Constant mean with state-dependent probabilities]{The ternary
double-zero family
\eqref{eq:three_increment_double_root_probabilities}--\eqref{eq:style_s9_2_042},
with $c=1/24$.
(a) The probabilities $\pi_1=1/3-\alpha/12$ and
$\pi_2=1/3+\alpha/24$ vary with the state, although $\bar\pi=1$.
(b) Optimal means obtained from Lemma~\ref{lem:affine_plot_paths}
and the quadratures
\eqref{eq:three_increment_double_root_lower_integral}--\eqref{eq:three_increment_double_root_upper_integral}.
From left to right, the terminal means are approximately
$0.256,0.552,1,1.377,1.673,1.896$.
Time increases upward on a logarithmic scale from $10^{-24}$ to $1$.
Only the dashed mean $1$ has zero cost.
For endpoints below $1$, the mean decreases from the limiting value $1$;
for endpoints above $1$, it decreases from the limiting value $2$.
These limits refer to $\tau\downarrow0$ and need not be attained
at the lower edge of the plot.
The positive-cost curves are asymmetric even though the mean increment
is constant and the probabilities of $0$ and $2$ are equal at each state:
the state dependence of their common value breaks reflection symmetry
about $\alpha=1$.}
\label{fig:double_root_affine_paths}
\end{figure}
\clearpage

\begin{figure}[p]
\centering
\includegraphics[width=\linewidth,height=0.72\textheight,keepaspectratio]{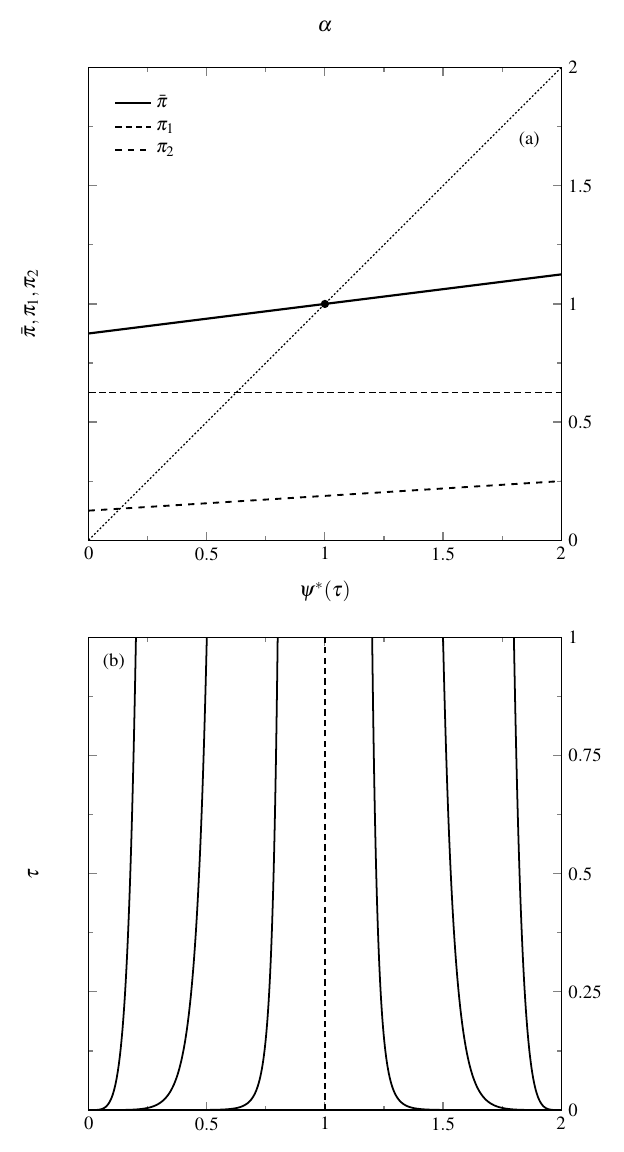}
\caption[Optimal paths in the family with common terminal entropy]{The
affine family
\eqref{eq:three_increment_analytic_family}--\eqref{eq:style_s9_2_053},
with $c=1/8$.
(a) The probabilities $\pi_1,\pi_2$ and the mean
$\bar\pi(\alpha)=7/8+\alpha/8$, whose unique fixed point is $1$.
(b) Optimal means from \eqref{eq:analytic_affine_optimal_mean},
ending at $\alpha=0.2,0.5,0.8,1,1.2,1.5,1.8$.
Time increases upward, with the same horizontal scale as in panel (a).
The dashed mean $1$ has zero cost; every solid curve has positive cost.
Means ending below $1$ increase from $0$, and those ending above $1$
decrease from $2$. The values at zero are limits of the means,
while every cumulative path starts at zero.
For any fixed positive time, the formula tends to the constant mean
$\alpha$ as $c\downarrow0$, the limit of constant urn probabilities.
Changing $c$ therefore changes the optimal paths without changing
their minimum costs, as well as changing the cost of the fixed test
path in Figure~\ref{fig:shared_entropy_path_cost}.}
\label{fig:analytic_affine_paths}
\end{figure}
\clearpage
\endgroup
\restoretextgeometry
}

\subsection{Scope of the extension}
\label{subsec:main_binary_comparison}

Beyond the count description, Theorem~\ref{thm:intro_terminal_entropy}
controls individual attainable point probabilities uniformly,
including sequences approaching the endpoints; a terminal LDP
for sets alone does not imply this.

The extension has a limitation. Ref.~\cite{Franchini_URNS} allows
binary probabilities $0$ and $1$: Subsections~3.1.3--3.1.4 use path
surgery for forbidden increments and singular logarithmic costs,
and Corollary~2 treats initial-condition-dependent reachability.
We do not extend those results. Here strict positivity and
compactness give a common positive lower bound, controlling
logarithms and permitting history modifications in the point
probability comparisons. With vanishing probabilities, allowed
increments may depend on the state and modified histories may
be impossible, with infinite cost. Simply allowing zeros or
removing positive regularization therefore does not extend our
arguments. These zeros differ from the fixed points of $\bar\pi$
already covered by Theorem~\ref{thm:intro_zero_cost}.

Nonbinary large deviation results also precede the present work,
as discussed in Subsection~\ref{subsec:beyond_binary_increments}.
The affine subclass belongs to the finite-state framework of
Budhiraja and Waterbury~\cite{Budhiraja_Waterbury}; broader
communication and dependence on the last symbol are treated in
\cite{Budhiraja_Waterbury_Zoubouloglou}.
In continuous time, Budhiraja and Coghi~\cite{Budhiraja_Coghi}
obtain a joint empirical-measure and flux LDP under affine-rate
and communication assumptions. Coghi and
Garrahan~\cite{Coghi_Garrahan_SIJP} develop tilting constructions
for more general rate maps under additional mixing and
time-scale separation assumptions.
Here the continuous positive probabilities may be nonlinear;
the conclusions combine path and count costs, uniform point
asymptotics, the inverse entropy problem and the global affine
cumulant analysis.

Section~\ref{sec:notation} develops the construction and path spaces;
Sections~\ref{sec:large_deviations}--\ref{sec:optimal_trajectories}
prove the results on paths, terminal probabilities and zero cost.
Section~\ref{sec:cumulants} treats the cumulant generating function,
duality and the inverse problem, followed by the affine analysis in
Section~\ref{sec:affine_urns}.
Section~\ref{sec:examples} supplies the specializations to $K=1,2$
and the explicit calculations for the examples.

\section{Notation and foundations}
\label{sec:notation}

We define the urn process and its path spaces,
following the notation of Ref.~\cite{Franchini_URNS} and
Ref.~\cite{FB}.

\subsection{Market history and initial condition}
\label{subsec:market_history}

Let $K\geq1$ and $N\geq1$ be integers, with $K$ fixed and $N$
the number of observed steps. Set $S:=\{1,\ldots,N\}$,
$\hat{\Omega}:=\{1,\ldots,K\}$ and
$\Omega:=\hat{\Omega}\cup\{0\}$. The history is the indexed
family
\begin{equation}
\sigma:=\{\sigma_n\in\Omega:\,n\in S\}\in\Omega^S,
\label{eq:market_history}
\end{equation}
\begin{equation}
\Omega^S:=\prod_{n\in S}\Omega^{\{n\}},
\label{eq:style_s3_history_product}
\end{equation}
which preserves the order of the increments. In the market
interpretation, $\sigma_n$ is the number of units purchased
by the $n$th customer, with zero for no purchase.
Set $M_0:=0$ and define the cumulative sum and mean by
\begin{equation}
M_n:=\sum_{s=1}^{n}\sigma_s,\qquad n\in S,
\label{eq:total_and_average}
\end{equation}
\begin{equation}
\psi_n:=\frac{M_n}{n},\qquad n\in S.
\label{eq:style_s3_discrete_average}
\end{equation}
We have $M_n\in\{0,\ldots,Kn\}$ and $\psi_n\in[0,K]$;
the mean is defined only for $n\geq1$.

We denote the law of the process by $\mathbb P$ and
expectation by $\mathbb E$. The first increment is uniform:
\begin{equation}
\mathbb P(\sigma_1=k)=\frac{1}{K+1},
\qquad k\in\Omega.
\label{eq:initial_distribution}
\end{equation}
In particular, $M_1=\psi_1=\sigma_1$.

\subsection{Urn functions and transition probabilities}
\label{subsec:urn_functions}

The components of the vector $\pi:=\{\pi_k:\,k\in\Omega\}$
are continuous and strictly positive on the entire interval:
\begin{equation}
\pi_k\in C([0,K]),\qquad k\in\Omega,
\label{eq:urn_positivity}
\end{equation}
\begin{equation}
\pi_k(\alpha)>0,\qquad k\in\Omega,\qquad\alpha\in[0,K].
\label{eq:style_s3_urn_strict_positivity}
\end{equation}
We also assume
\begin{equation}
\sum_{k\in\Omega}\pi_k(\alpha)=1,\qquad\alpha\in[0,K],
\label{eq:urn_normalization}
\end{equation}
\begin{equation}
\pi_0(\alpha)=1-\sum_{k\in\hat\Omega}\pi_k(\alpha),\qquad\alpha\in[0,K].
\label{eq:style_s3_urn_zero_component}
\end{equation}
These assumptions hold throughout the paper. By compactness
of $[0,K]$ and finiteness of $\Omega$,
\begin{equation}
\min_{k\in\Omega}\min_{\alpha\in[0,K]}\pi_k(\alpha)>0;
\label{eq:uniform_positivity}
\end{equation}
hence the logarithms $\log\pi_k$ are bounded and uniformly
continuous.

For $1\leq n<N$ the transition is
\begin{equation}
\mathbb P\bigl(\sigma_{n+1}=k\mid 
\sigma_1,\ldots,\sigma_n\bigr)
=\pi_k(\psi_n),
\qquad k\in\Omega.
\label{eq:transition_probability}
\end{equation}
The transition depends on the past through $\psi_n$
and, together with \eqref{eq:initial_distribution},
determines the law of the process. Every finite history
has positive probability.

The mean increment function is
\begin{equation}
\bar\pi(\alpha):=\sum_{k\in\hat{\Omega}}k\pi_k(\alpha),
\qquad \alpha\in[0,K].
\label{eq:mean_increment}
\end{equation}
It is continuous and satisfies $0<\bar\pi<K$, since
increments $0$ and $K$ have positive probability.
The relation
\begin{equation}
\psi_{n+1}=\frac{n\psi_n+\sigma_{n+1}}{n+1},\qquad 1\leq n<N,
\label{eq:style_s3_discrete_average_recursion}
\end{equation}
and \eqref{eq:transition_probability} give, for $1\leq n<N$,
\begin{equation}
\mathbb E\bigl(\psi_{n+1}-\psi_n\mid 
\sigma_1,\ldots,\sigma_n\bigr)
=\frac{\bar\pi(\psi_n)-\psi_n}{n+1}.
\label{eq:conditional_average_increment}
\end{equation}

The fixed-point set is
\begin{equation}
C:=\{\alpha\in[0,K]:\,\bar\pi(\alpha)=\alpha\}.
\label{eq:fixed_points}
\end{equation}
Continuity and the opposite signs of $\bar\pi(\alpha)-\alpha$
at $0$ and $K$ imply that $C$ is nonempty, compact and
\begin{equation}
0<\inf C\leq\sup C<K.
\label{eq:fixed_point_bounds}
\end{equation}
We do not require $C$ to be finite or to have only isolated points.
For $K=1$, $\bar\pi=\pi_1$ and $\pi_0=1-\bar\pi$;
for $K>1$, the mean and normalization do not generally determine
the entire vector $\pi$.

\subsection{History probabilities}
\label{subsec:history_probabilities}

For $k,k'\in\Omega$, the Kronecker delta $\delta_k(k')$ is $1$ if $k'=k$ and $0$ otherwise.
The probability of a history is
\begin{equation}
W(\sigma):=\mathbb P(\{\sigma\})
=\frac{1}{K+1}
\prod_{n=1}^{N-1}\pi_{\sigma_{n+1}}(\psi_n).
\label{eq:history_weight}
\end{equation}
An empty product is $1$ and an empty sum is $0$; in particular,
for $N=1$ we have $W(\sigma)=1/(K+1)$. Summing successively over the
increments, from last to first, \eqref{eq:urn_normalization} and
\eqref{eq:initial_distribution} give
\begin{equation}
\sum_{\sigma\in\Omega^S}W(\sigma)=1.
\label{eq:style_s3_history_total_mass}
\end{equation}
For every $E\subseteq\Omega^S$ we therefore have
\begin{equation}
\mathbb P(E)=\sum_{\sigma\in E}W(\sigma).
\label{eq:style_s3_history_event_mass}
\end{equation}
The action of a history is the logarithm of its probability:
\begin{equation}
\mathcal A(\sigma):=\log W(\sigma)
=-\log(K+1)
+\sum_{n=1}^{N-1}\log\pi_{\sigma_{n+1}}(\psi_n).
\label{eq:discrete_action}
\end{equation}
It is finite and nonpositive and retains the contribution
of the initialization, including for $N=1$.

\subsection{Continuous paths}
\label{subsec:continuous_trajectories}

To compare histories of different lengths, we linearly
interpolate the values $M_n/N$ in time $\tau\in[0,1]$.
For $n=0,\ldots,N-1$ and $n/N\leq\tau\leq(n+1)/N$, set
\begin{equation}
N\varphi(\tau\mid\sigma)
:=M_n+(N\tau-n)\sigma_{n+1}.
\label{eq:path_interpolation}
\end{equation}
The formulas agree at the shared interval endpoints.
We denote the resulting function by $\varphi(\sigma)$;
at the nodes and in the interiors of the respective intervals we have
\begin{equation}
\varphi(n/N\mid\sigma)=\frac{M_n}{N},\qquad 0\leq n\leq N,
\label{eq:interpolation_nodes}
\end{equation}
\begin{equation}
\partial_\tau\varphi(\tau\mid\sigma)=\sigma_{n+1},\qquad \frac nN<\tau<\frac{n+1}{N}.
\label{eq:style_s3_path_derivative}
\end{equation}
In particular, $\varphi(0\mid\sigma)=0$ and
$\varphi(1\mid\sigma)=\psi_N$.

We introduce the space
\begin{equation}
Q:=\bigl\{\varphi\in C([0,1]):\ \varphi(0)=0,\ 0\leq\varphi(\tau')-\varphi(\tau)\leq K(\tau'-\tau)\bigr\},
\label{eq:path_space}
\end{equation}
where the constraint holds for every $0\leq\tau\leq\tau'\leq1$.
Equivalently, the elements of $Q$ are the absolutely
continuous functions vanishing at $0$, with
$0\leq\partial_\tau\varphi\leq K$ almost everywhere:
the increment constraint implies the Lipschitz property,
while the converse follows by integrating the derivative.
We have $0\leq\varphi(\tau)\leq K\tau$.

We equip $Q$ with the uniform distance
\begin{equation}
\|\varphi-\varphi'\|_\infty
:=\sup_{\tau\in[0,1]}|\varphi(\tau)-\varphi'(\tau)|.
\label{eq:uniform_distance}
\end{equation}
The functions in $Q$ are uniformly bounded and equicontinuous by their common
Lipschitz bound. Since uniform convergence preserves the value
at zero and the inequalities in \eqref{eq:path_space}, $Q$ is closed;
the Arzelà--Ascoli theorem therefore gives its compactness.

The continuous mean is
\begin{equation}
\psi(\tau):=\frac{\varphi(\tau)}{\tau},
\qquad 0<\tau\leq1.
\label{eq:continuous_average}
\end{equation}
It is continuous on $(0,1]$, takes values in $[0,K]$ and, for
$\varphi=\varphi(\sigma)$, satisfies $\psi(n/N)=\psi_n$
for $1\leq n\leq N$.
We do not impose a limit as $\tau\downarrow0$.
For every $\varepsilon\in(0,1)$ and $\varphi,\varphi'\in Q$,
\begin{equation}
\sup_{\tau\in[\varepsilon,1]}
\left|\frac{\varphi(\tau)}{\tau}
-\frac{\varphi'(\tau)}{\tau}\right|
\leq\frac{\|\varphi-\varphi'\|_\infty}{\varepsilon}.
\label{eq:average_uniform_bound}
\end{equation}
Uniform convergence of paths therefore controls
the means on every $[\varepsilon,1]$.

For $\alpha\in[0,K]$, set
\begin{equation}
Q(\alpha):=\{\varphi\in Q:\,\varphi(1)=\alpha\}.
\label{eq:terminal_path_space}
\end{equation}
This set is nonempty, because it contains
$\varphi(\tau)=\alpha\tau$, and compact, because
evaluation at $1$ is continuous.
A history of length $N$ can have final mean
$\alpha$ only if $N\alpha$ is an integer.

\subsection{Counts of individual increments}
\label{subsec:increment_counts}

For $k\in\Omega$, $n=0,\ldots,N-1$ and
$n/N\leq\tau\leq(n+1)/N$, define the interpolated count
\begin{equation}
N\gamma_k(\tau\mid\sigma)
:=\sum_{s=1}^{n}\delta_k(\sigma_s)
+(N\tau-n)\delta_k(\sigma_{n+1}).
\label{eq:count_interpolation}
\end{equation}
The formulas agree at the shared endpoints.
We write $\gamma(\sigma)
=(\gamma_0(\sigma),\ldots,\gamma_K(\sigma))$.
At the nodes and in the interiors of the respective intervals,
\begin{equation}
N\gamma_k(n/N\mid\sigma)=\sum_{s=1}^{n}\delta_k(\sigma_s),\qquad 0\leq n\leq N,
\label{eq:count_derivative}
\end{equation}
\begin{equation}
\partial_\tau\gamma_k(\tau\mid\sigma)=\delta_k(\sigma_{n+1}),\qquad \frac nN<\tau<\frac{n+1}{N}.
\label{eq:style_s3_count_slope}
\end{equation}
For every $\tau\in[0,1]$,
\begin{equation}
\sum_{k\in\Omega}\gamma_k(\tau\mid\sigma)=\tau,
\label{eq:count_identities}
\end{equation}
\begin{equation}
\sum_{k\in\Omega}k\gamma_k(\tau\mid\sigma)=\varphi(\tau\mid\sigma).
\label{eq:style_s3_count_scalar_identity}
\end{equation}

In the continuous setting, a vector
$\gamma=(\gamma_0,\ldots,\gamma_K)$ is
admissible if its components are continuous, nondecreasing and
\begin{equation}
\gamma_k(0)=0,\qquad k\in\Omega,
\label{eq:count_constraints}
\end{equation}
\begin{equation}
\sum_{k\in\Omega}\gamma_k(\tau)=\tau,\qquad 0\leq\tau\leq1.
\label{eq:style_s3_count_total_constraint}
\end{equation}
Monotonicity and the sum constraint give
\begin{equation}
0\leq\gamma_k(\tau')-\gamma_k(\tau)\leq\tau'-\tau,
\qquad 0\leq\tau\leq\tau'\leq1.
\label{eq:display_004}
\end{equation}
The components are therefore Lipschitz with constant $1$.
Equivalently, they are absolutely continuous, vanish
at zero and satisfy
\begin{equation}
\partial_\tau\gamma_k(\tau)\geq0,\qquad\text{almost everywhere},
\label{eq:local_frequencies}
\end{equation}
\begin{equation}
\sum_{k\in\Omega}\partial_\tau\gamma_k(\tau)=1,\qquad\text{almost everywhere}.
\label{eq:style_s3_local_frequency_total}
\end{equation}
These derivatives are the local increment frequencies
and may vanish.

We use the uniform distance
\begin{equation}
\|\gamma-\gamma'\|_\infty
:=\max_{k\in\Omega}\|\gamma_k-\gamma_k'\|_\infty.
\label{eq:count_uniform_distance}
\end{equation}
Since $K\geq1$, the set of admissible vectors is
the closed subset of the compact set $Q^{K+1}$ defined
by the constraint $\sum_k\gamma_k(\tau)=\tau$ for every $\tau$;
it is therefore compact.

The associated cumulative path and mean are
\begin{equation}
\varphi(\tau):=\sum_{k\in\Omega}k\gamma_k(\tau),
\label{eq:counts_to_path}
\end{equation}
\begin{equation}
\psi(\tau):=\frac{\varphi(\tau)}{\tau},\qquad\tau>0.
\label{eq:style_s3_counts_to_mean}
\end{equation}
Relations \eqref{eq:local_frequencies} and \eqref{eq:style_s3_local_frequency_total} give
$0\leq\partial_\tau\varphi\leq K$ almost everywhere,
hence $\varphi\in Q$. The projection
$\gamma\mapsto\varphi$ is continuous, since
\begin{equation}
\|\varphi-\varphi'\|_\infty
\leq\sum_{k\in\Omega}k\|\gamma_k-\gamma_k'\|_\infty
\leq\frac{K(K+1)}2
\|\gamma-\gamma'\|_\infty.
\label{eq:count_projection_bound}
\end{equation}
It is also surjective: for $\varphi\in Q$, the vector
\begin{equation}
\gamma_0(\tau)=\tau-\frac{\varphi(\tau)}K,
\label{eq:count_representation}
\end{equation}
\begin{equation}
\gamma_K(\tau)=\frac{\varphi(\tau)}K,
\label{eq:style_s3_count_representation_last}
\end{equation}
with zero intermediate components is admissible.
Monotonicity of $\gamma_0$ follows from the upper bound
in \eqref{eq:path_space}; the other constraints follow
directly from the definitions. For $K=1$ the representation
is unique: $\gamma_1=\varphi$ and
$\gamma_0(\tau)=\tau-\varphi(\tau)$.

The action of the counts is
\begin{equation}
\Phi(\gamma):=\sum_{k\in\Omega}\int_0^1\log\pi_k\bigl(\psi(\tau)\bigr)\,d\gamma_k(\tau),
\label{eq:count_action}
\end{equation}
\begin{equation}
\Phi(\gamma)=\int_0^1\sum_{k\in\Omega}\partial_\tau\gamma_k(\tau)\log\pi_k\bigl(\psi(\tau)\bigr)\,d\tau.
\label{eq:style_s3_count_action_density}
\end{equation}
The equality follows from absolute continuity of the counts.
The integrands are defined for $\tau>0$; the value
at the single point $0$ does not affect the integrals.
Boundedness of the logarithms and the identity
\begin{equation}
\sum_{k\in\Omega}\int_0^1d\gamma_k=1
\label{eq:style_s3_count_mass}
\end{equation}
give
\begin{equation}
\min_{k\in\Omega}\min_{\beta\in[0,K]}\log\pi_k(\beta)
\leq\Phi(\gamma)\leq0.
\label{eq:count_action_bounds}
\end{equation}

For fixed $N$, the interpolated path determines all
increments through \eqref{eq:style_s3_path_derivative}
and hence all counts. In the uniform limit, however,
a scalar path may correspond to different counts.
For $K=2$, the admissible vectors
\begin{equation}
\gamma(\tau)=(0,\tau,0),
\label{eq:distinct_count_paths}
\end{equation}
\begin{equation}
\gamma'(\tau)=(\tau/2,0,\tau/2)
\label{eq:style_s3_distinct_count_second}
\end{equation}
both have cumulative path $\varphi(\tau)=\tau$.
The first comes from histories with only increments $1$;
the second is the uniform limit of the counts of
histories alternating $0$ and $2$. The scalar paths
of these histories lie at distance from the line $\varphi(\tau)=\tau$
of at most $1/N$,
even when $N$ is odd. The actions are
\begin{equation}
\Phi(\gamma)=\log\pi_1(1),
\label{eq:distinct_count_actions}
\end{equation}
\begin{equation}
\Phi(\gamma')=\frac12\log\pi_0(1)+\frac12\log\pi_2(1).
\label{eq:style_s3_distinct_second_action}
\end{equation}
For the constant vector $\pi=(2/5,1/5,2/5)$, the
normalized actions $\mathcal A(\sigma)/N$ of the
two families also converge, respectively, to $\log(1/5)$
and $\log(2/5)$, while the scalar paths have the same
limit. In general, therefore, no continuous functional
of the scalar path alone can approximate
$\mathcal A/N$, as $N$ increases, uniformly over all
histories: the action $\Phi$ requires the entire count vector.
Figure~\ref{fig:count_path_contraction} compares these two families
and their limiting counts.

\section{Large deviations and scaling limit}
\label{sec:large_deviations}

In this section we prove the large deviation principle
for the paths of the process. The method follows that of
Ref.~\cite{Franchini_URNS}, Subsections~3.1.1--3.1.2: the
change of measure is applied to the increment counts,
while the result for the cumulative path is obtained by
contraction. Dependence on the entire vector of urn functions
remains explicit in both representations.

\subsection{Large deviation principle for paths}
\label{subsec:path_ldp}

For $\varphi\in Q$ define the functional
\begin{equation}
I(\varphi):=\int_0^1
\sup_{\lambda\in\mathbb{R}}
\left\{\lambda\,\partial_\tau\varphi(\tau)
-\log\left[\sum_{k\in\Omega}
\pi_k\!\left(\frac{\varphi(\tau)}{\tau}\right)e^{k\lambda}\right]
\right\}\,d\tau.
\label{eq:path_rate}
\end{equation}
The supremum is computed separately for almost every $\tau$. The
ratio $\varphi(\tau)/\tau$ is used for $\tau>0$, and the value
of the integrand at the single point $0$ does not affect the integral.

\begin{theorem}[Large deviations for paths]
\label{thm:path_ldp}
Let $K\geq1$ be fixed and let the functions $\pi_k$, $k\in\Omega$,
be continuous and strictly positive on $[0,K]$, with
\begin{equation}
\sum_{k\in\Omega}\pi_k(\alpha)=1,\qquad\alpha\in[0,K].
\label{eq:style_s4_path_ldp_normalization}
\end{equation}
For the process with initial condition
\eqref{eq:initial_distribution} and transitions
\eqref{eq:transition_probability}, the laws of the paths
$\varphi(\sigma)$ satisfy, as $N$ increases, a large
deviation principle with speed $N$ on the space $Q$ equipped
with the uniform distance \eqref{eq:uniform_distance}.

The rate function is $I$ defined in \eqref{eq:path_rate}.
It is nonnegative, finite on all of $Q$, lower semicontinuous
and has compact sublevel sets. For every Borel set $B\subseteq Q$
fixed as $N$ varies, we have
\begin{equation}
\liminf_{N\to\infty}\frac{1}{N}
\log\mathbb{P}\bigl(\varphi(\sigma)\in B\bigr)
\geq-\inf_{\varphi\in B^\circ}I(\varphi)
\label{eq:path_ldp_lower}
\end{equation}
and
\begin{equation}
\limsup_{N\to\infty}\frac{1}{N}
\log\mathbb{P}\bigl(\varphi(\sigma)\in B\bigr)
\leq-\inf_{\varphi\in\overline{B}}I(\varphi),
\label{eq:path_ldp_upper}
\end{equation}
where $B^\circ$ and $\overline{B}$ are the interior and closure relative
to $Q$. The conventions $\inf\varnothing=+\infty$ and
$\log0=-\infty$ are understood.

For every $\varphi\in Q$, the rate function also admits the
representation
\begin{equation}
I(\varphi)=
\min_{\substack{\gamma\ \mathrm{admissible}\\
\sum_{k\in\Omega}k\gamma_k=\varphi}}
\left\{\int_0^1\sum_{k\in\Omega}
\partial_\tau\gamma_k(\tau)
\log\partial_\tau\gamma_k(\tau)\,d\tau
-\Phi(\gamma)\right\},
\label{eq:path_rate_counts}
\end{equation}
or, equivalently,
\begin{equation}
I(\varphi)=
\min_{\substack{\gamma\ \mathrm{admissible}\\
\sum_{k\in\Omega}k\gamma_k=\varphi}}
\int_0^1\sum_{k\in\Omega}
\partial_\tau\gamma_k(\tau)
\log\frac{\partial_\tau\gamma_k(\tau)}
{\pi_k(\varphi(\tau)/\tau)}\,d\tau.
\label{eq:path_rate_relative_entropy}
\end{equation}
The admissible vectors are those of subsection
\ref{subsec:increment_counts}; the constraint on their weighted sum
is an equality of functions on $[0,1]$. The minimum is attained.
The convention $0\log0=0$ is adopted; in particular, a term
of the integrand in \eqref{eq:path_rate_relative_entropy} vanishes
when $\partial_\tau\gamma_k=0$.
\end{theorem}

Nonnegativity already follows from \eqref{eq:path_rate}, choosing
$\lambda=0$ in the supremum. Moreover, the representation
\eqref{eq:path_rate_counts}, the inequality
\begin{equation}
\sum_{k\in\Omega}\alpha_k\log\alpha_k\leq0
\label{eq:style_s4_frequency_entropy_bound}
\end{equation}
for nonnegative frequencies summing to one, and the bound
\eqref{eq:count_action_bounds} give
\begin{equation}
0\leq I(\varphi)
\leq-\min_{k\in\Omega}\min_{\beta\in[0,K]}\log\pi_k(\beta)
<\infty,
\qquad \varphi\in Q.
\label{eq:path_rate_bounds}
\end{equation}
The count constraint is always feasible, as shown by the
construction \eqref{eq:count_representation}--\eqref{eq:style_s3_count_representation_last}. Setting $I=+\infty$
on $C([0,1])\setminus Q$, the same principle also holds on the
space $C([0,1])$ with the uniform topology, since all laws
are supported on the compact set $Q$.

For a Borel set $B\subseteq Q$ fixed as $N$ varies,
if
\begin{equation}
\inf_{B^\circ}I=\inf_{\overline B}I,
\label{eq:style_s4_ldp_set_matching}
\end{equation}
the two estimates give
\begin{equation}
\lim_{N\to\infty}\frac1N
\log\mathbb P\bigl(\varphi(\sigma)\in B\bigr)
=-\inf_{\varphi\in B}I(\varphi).
\label{eq:display_005}
\end{equation}

For $K=1$, the counts compatible with $\varphi$ are uniquely
$\gamma_1=\varphi$ and $\gamma_0(\tau)=\tau-\varphi(\tau)$.
The representation \eqref{eq:path_rate_counts} therefore recovers
the functional of Theorem~1 of Ref.~\cite{Franchini_URNS},
Eq.~(2.5)--(2.7), in the strictly positive class considered here.

The proof proceeds through the uniform reference process,
continuity and approximation of the action, and the LDP for
counts. It concludes in subsection
\ref{subsec:scalar_contraction} through contraction and
minimization over local frequencies, constructing a
measurable minimizer before identifying
\eqref{eq:path_rate}.

\subsection{Reference process and change of measure}
\label{subsec:reference_process}

Let $\mathbb P_0$ be the law on $\Omega^S$ under which the
increments are independent and uniform on $\Omega$, and let
$\mathbb E_0$ be the corresponding expectation. For every history,
\begin{equation}
\mathbb P_0(\{\sigma\})=(K+1)^{-N}.
\label{eq:reference_history_probability}
\end{equation}
The law of the first increment agrees with
\eqref{eq:initial_distribution}. From
\eqref{eq:history_weight}--\eqref{eq:discrete_action} we obtain
the exact density
\begin{equation}
\frac{d\mathbb P}{d\mathbb P_0}(\sigma)
=(K+1)^N e^{\mathcal A(\sigma)}
=\prod_{n=1}^{N-1}
\left[(K+1)\pi_{\sigma_{n+1}}(\psi_n)\right].
\label{eq:reference_density}
\end{equation}
The product contains $N-1$ factors because the two laws
have the same initialization; for $N=1$ it equals $1$.
Denoting the indicator of an event by $\mathbb I$,
for every $E\subseteq\Omega^S$ we obtain
\begin{equation}
\mathbb P(E)=(K+1)^N\mathbb E_0
\left[e^{\mathcal A(\sigma)}\mathbb I(\sigma\in E)\right].
\label{eq:reference_change_event}
\end{equation}
This identity holds for every $N$, including for events depending
on $N$. We will apply it to the counts, after uniformly comparing
$\mathcal A/N$ with the continuous action $\Phi$.

\subsection{Continuity and approximation of the action}
\label{subsec:action_continuity}

The following two lemmas extend to the count space
Lemmas~13 and~14 of Ref.~\cite{Franchini_URNS}.
Set
\begin{equation}
\|\log\pi\|_\infty:=
\max_{k\in\Omega}\sup_{\alpha\in[0,K]}|\log\pi_k(\alpha)|
\label{eq:log_urn_norm}
\end{equation}
and use the common modulus of continuity
\begin{equation}
\omega(r):=\max_{k\in\Omega}
\sup_{\substack{\alpha,\beta\in[0,K]\\|\alpha-\beta|\leq r}}
|\log\pi_k(\alpha)-\log\pi_k(\beta)|,
\qquad r\geq0.
\label{eq:log_urn_modulus}
\end{equation}
The assumptions on the urn functions imply
$\|\log\pi\|_\infty<\infty$, while $\omega$ is nondecreasing,
with $\omega(0)=0$ and $\omega(r)\to0$ as $r\downarrow0$.
The definition also includes $r>K$.

\begin{lemma}[Continuity of the action on counts]
\label{lem:count_action_continuity}
The functional $\Phi$ defined in \eqref{eq:count_action} is
uniformly continuous on the set of admissible vectors,
equipped with the distance \eqref{eq:count_uniform_distance}.
\end{lemma}

\begin{proof}
Let $\gamma^{(j)}\to\gamma$ uniformly
in the count space, and let $\psi_j$ and $\psi$ be the corresponding
means. For every admissible vector and $0<\varepsilon<1$,
positivity of the measures $d\gamma_k$ and the identity
\eqref{eq:style_s3_count_total_constraint} give
\begin{equation}
\left|\sum_{k\in\Omega}\int_0^\varepsilon
\log\pi_k(\psi(\tau))\,d\gamma_k(\tau)\right|
\leq\varepsilon\|\log\pi\|_\infty.
\label{eq:action_initial_interval}
\end{equation}
For fixed $\varepsilon$, on $[\varepsilon,1]$ we have
$\psi_j\to\psi$ uniformly, and hence also
$g_{j,k}:=\log\pi_k(\psi_j)\to g_k:=\log\pi_k(\psi)$
uniformly for every $k$.
The Stieltjes measures $d\gamma_k^{(j)}$ converge weakly
to $d\gamma_k$ on this interval: for $C^1$ test functions it suffices
to integrate by parts and use uniform convergence of the primitives;
for continuous test functions we conclude by uniform approximation,
since all masses are at most one.
The counts are Lipschitz, so these measures have no atoms,
including at the endpoints of the restricted interval. Therefore
\begin{equation}
\left|\int_\varepsilon^1g_{j,k}\,d\gamma_k^{(j)}-\int_\varepsilon^1g_k\,d\gamma_k\right|
\leq\|g_{j,k}-g_k\|_\infty 
+\left|\int_\varepsilon^1g_k\,d\bigl(\gamma_k^{(j)}-\gamma_k\bigr)\right|
\longrightarrow0.
\label{eq:display_006}
\end{equation}
Summing over the finitely many components and using
\eqref{eq:action_initial_interval} for both vectors,
the contribution to $\Phi(\gamma^{(j)})-
\Phi(\gamma)$ from the segment $[0,\varepsilon]$
has absolute value at most $2\varepsilon\|\log\pi\|_\infty$.
First $j\to\infty$ and then $\varepsilon\downarrow0$ prove
continuity of $\Phi$. Compactness of the count space,
proved in Subsection~\ref{subsec:increment_counts},
implies its uniform continuity.
\end{proof}

\begin{lemma}[Uniform approximation of the discrete action]
\label{lem:discrete_action_approximation}
For every $N\geq1$ we have
\begin{equation}
\sup_{\sigma\in\Omega^S}
\left|\frac{\mathcal A(\sigma)}N-\Phi\bigl(\gamma(\sigma)\bigr)\right| 
\leq\frac{\log(K+1)+\|\log\pi\|_\infty}N
+\frac1N\sum_{n=1}^{N-1}\omega\!\left(\frac K{n+1}\right),
\label{eq:discrete_action_uniform_bound}
\end{equation}
\begin{equation}
\frac{\log(K+1)+\|\log\pi\|_\infty}N
+\frac1N\sum_{n=1}^{N-1}\omega\!\left(\frac K{n+1}\right)
\longrightarrow0,\qquad N\to\infty.
\label{eq:style_s4_discrete_action_bound_limit}
\end{equation}
\end{lemma}

\begin{proof}
For a history $\sigma$, on the first interval
$0<\tau\leq1/N$ the mean equals $\sigma_1$. From
\eqref{eq:style_s3_count_slope} we obtain the exact identity
\begin{equation}
\Phi(\gamma(\sigma))
=\frac1N\log\pi_{\sigma_1}(\sigma_1)
+\sum_{n=1}^{N-1}\int_{n/N}^{(n+1)/N}
\log\pi_{\sigma_{n+1}}(\psi(\tau))\,d\tau.
\label{eq:continuous_action_on_history}
\end{equation}
For $1\leq n<N$ and $n/N\leq\tau\leq(n+1)/N$,
\eqref{eq:path_interpolation} gives
\begin{equation}
\psi(\tau)-\psi_n=\frac{N\tau-n}{N\tau}(\sigma_{n+1}-\psi_n),
\label{eq:discrete_continuous_average_bound}
\end{equation}
\begin{equation}
|\psi(\tau)-\psi_n|\leq\frac K{n+1}.
\label{eq:style_s4_discrete_continuous_average_estimate}
\end{equation}
We use $|\sigma_{n+1}-\psi_n|\leq K$ and monotonicity of
$u/(n+u)$ for $0\leq u\leq1$.
In the comparison with \eqref{eq:discrete_action} divided by $N$,
the first interval contributes exactly
$-[\log(K+1)+\log\pi_{\sigma_1}(\sigma_1)]/N$;
each subsequent interval contributes in absolute value
at most $N^{-1}\omega(K/(n+1))$. Summing gives
\eqref{eq:discrete_action_uniform_bound}, with an upper bound
independent of the history. It tends to zero because
$\omega(K/(n+1))\to0$ and its arithmetic mean has the same limit.
For $N=1$ the subsequent sum is empty.
\end{proof}

\subsection{Large deviations for increment counts}
\label{subsec:count_ldp}

We derive the LDP for counts under $\mathbb P_0$ and
apply the change of measure. The procedure extends
to counts Lemma~15 and the argument in Subsection~3.1.2
of Ref.~\cite{Franchini_URNS}.

\begin{lemma}[Large deviations of uniform counts]
\label{lem:reference_count_ldp}
Under $\mathbb P_0$, the paths
$\gamma(\sigma)$ satisfy an LDP with speed
$N$ on the set of admissible vectors of subsection
\ref{subsec:increment_counts}, equipped with the distance
\eqref{eq:count_uniform_distance}. The rate function is
\begin{equation}
\log(K+1)+\int_0^1\sum_{k\in\Omega}
\partial_\tau\gamma_k(\tau)
\log\partial_\tau\gamma_k(\tau)\,d\tau.
\label{eq:reference_count_rate}
\end{equation}
It is lower semicontinuous, has compact sublevel sets
and takes values in $[0,\log(K+1)]$, with $0\log0=0$.
\end{lemma}

\begin{proof}
The vectors $(\delta_k(\sigma_n))_{k\in\Omega}$ are independent
and uniform on the $K+1$ canonical basis vectors. For
$(\lambda_k)_{k\in\Omega}\in\mathbb R^{K+1}$,
\begin{equation}
\log\mathbb E_0\exp\biggl(\sum_{k\in\Omega}\lambda_k\delta_k(\sigma_1)\biggr)
=\log\biggl[\frac1{K+1}\sum_{k\in\Omega}e^{\lambda_k}\biggr].
\label{eq:reference_count_cumulant}
\end{equation}
The function is finite for every parameter vector. For every
$(\alpha_k)_{k\in\Omega}\in\mathbb R^{K+1}$, its
Legendre transform is
\begin{equation}
\sup_{(\lambda_k)\in\mathbb R^{K+1}}
\biggl\{\sum_{k\in\Omega}\alpha_k\lambda_k
-\log\biggl[\frac1{K+1}\sum_{k\in\Omega}e^{\lambda_k}\biggr]\biggr\}.
\label{eq:reference_count_local_rate}
\end{equation}
For nonnegative frequencies summing to one, the value is
\begin{equation}
\log(K+1)+\sum_{k\in\Omega}\alpha_k\log\alpha_k.
\label{eq:style_global_reference_count_local_value}
\end{equation}
For every other vector, the value is $+\infty$.

Indeed, on the frequency simplex, concavity of the
logarithm gives
\begin{equation}
\sum_{\alpha_k>0}\alpha_k
\log\frac{e^{\lambda_k}}{\alpha_k}
\leq\log\sum_{\alpha_k>0}e^{\lambda_k}
\leq\log\sum_{k\in\Omega}e^{\lambda_k}.
\label{eq:display_007}
\end{equation}
This is the required upper bound; choosing
$\lambda_k=\log(\alpha_k+\varepsilon)$ and sending
$\varepsilon\downarrow0$ yields the stated value,
even if some frequencies vanish.
If the frequencies do not sum to one, the choice $\lambda_k=t$
for every $k$ makes the expression equal to
\begin{equation}
t\biggl(\sum_{k\in\Omega}\alpha_k-1\biggr),
\label{eq:style_s4_reference_transform_nonunit_mass}
\end{equation}
which is unbounded above. If instead the frequencies sum to one
and $\alpha_j<0$, choosing $\lambda_j=-t$ and setting the other
parameters to zero yields
\begin{equation}
-t\alpha_j-\log\frac{K+e^{-t}}{K+1}
\longrightarrow+\infty,\qquad t\to+\infty.
\label{eq:display_008}
\end{equation}

Theorem~5.1.2 of Ref.~\cite{Dembo_Zeitouni}, in dimension
$K+1$, gives the LDP for stepwise counts on the space of
bounded functions from $[0,1]$ to $\mathbb R^{K+1}$, with the
uniform distance. In the same space, for every history,
\begin{equation}
\max_{k\in\Omega}\sup_{\tau\in[0,1]}
\left|\gamma_k(\tau\mid\sigma)
-\frac1N\sum_{s=1}^{\lfloor N\tau\rfloor}
\delta_k(\sigma_s)\right|\leq\frac1N.
\label{eq:reference_count_interpolation}
\end{equation}
Thus, for every fixed positive distance, the probability
of exceeding it is zero for sufficiently large $N$.
Exponential equivalence and Theorem~4.2.13 of
Ref.~\cite{Dembo_Zeitouni} transfer the LDP to the interpolated
counts. The maximum norm over components and the
Euclidean norm induce the same uniform topology.

The integral cost is finite precisely on absolutely
continuous paths, vanishing at the origin, with nonnegative
derivatives summing to one almost everywhere. These are exactly
the admissible vectors of
\eqref{eq:count_constraints}--\eqref{eq:style_s3_local_frequency_total},
which form a closed set and contain all interpolated counts.
The LDP therefore restricts to this set; integrating
\eqref{eq:style_global_reference_count_local_value} yields
\eqref{eq:reference_count_rate}, with lower semicontinuity
and compact sublevel sets. The cost is nonnegative, by choosing
all parameters zero in the transform, and does not exceed
$\log(K+1)$ by \eqref{eq:style_s4_frequency_entropy_bound}.
\end{proof}

\begin{theorem}[Large deviations of urn counts]
\label{thm:count_ldp}
Under the assumptions of Theorem~\ref{thm:path_ldp}, the laws of
$\gamma(\sigma)$ satisfy an LDP with speed
$N$ on the set of admissible vectors, with the distance
\eqref{eq:count_uniform_distance} and rate function
\begin{equation}
I_{\mathrm c}(\gamma):=\int_0^1\sum_{k\in\Omega}
\partial_\tau\gamma_k(\tau)\log\partial_\tau\gamma_k(\tau)\,d\tau-\Phi(\gamma),
\label{eq:count_rate}
\end{equation}
\begin{equation}
I_{\mathrm c}(\gamma)=\int_0^1\sum_{k\in\Omega}\partial_\tau\gamma_k(\tau)
\log\frac{\partial_\tau\gamma_k(\tau)}{\pi_k\bigl(\psi(\tau)\bigr)}\,d\tau,
\label{eq:style_s4_count_rate_entropy}
\end{equation}
where $\psi$ is the mean of the path
\eqref{eq:counts_to_path}--\eqref{eq:style_s3_counts_to_mean}. The cost is finite, nonnegative,
lower semicontinuous, has compact sublevel sets and minimum zero.

For every Borel set $B$ fixed as $N$ varies, the estimates
\eqref{eq:path_ldp_lower}--\eqref{eq:path_ldp_upper} hold
for $\gamma(\sigma)$, replacing $Q$ with the
space of admissible counts and $I$ with the rate
$I_{\mathrm c}$ in \eqref{eq:count_rate}. Interior and closure are relative
to this space. We use the conventions
$\inf\varnothing=+\infty$, $\log0=-\infty$ and $0\log0=0$;
a term of the integrand vanishes if
$\partial_\tau\gamma_k=0$, and the value at the single point
$\tau=0$ does not affect the integral.
\end{theorem}

\begin{proof}
For every $\eta>0$, Lemma~\ref{lem:discrete_action_approximation}
and \eqref{eq:reference_change_event} give, for sufficiently large
$N$ and every Borel set $B$ of counts,
\begin{equation}
e^{N[\log(K+1)-\eta]}\,
\mathbb E_0\Bigl[e^{N\Phi(\gamma(\sigma))}\mathbb I\bigl(\gamma(\sigma)\in B\bigr)\Bigr]
\leq\mathbb P\bigl(\gamma(\sigma)\in B\bigr),
\label{eq:count_change_measure_bounds}
\end{equation}
\begin{equation}
\mathbb P\bigl(\gamma(\sigma)\in B\bigr)
\leq e^{N[\log(K+1)+\eta]}\,
\mathbb E_0\Bigl[e^{N\Phi(\gamma(\sigma))}\mathbb I\bigl(\gamma(\sigma)\in B\bigr)\Bigr].
\label{eq:style_s4_count_change_measure_upper}
\end{equation}
The threshold in $N$ does not depend on $B$.
By Lemma~\ref{lem:reference_count_ldp}, the reference process
satisfies an LDP with a good rate function, denoted
here by $R_0$, given by \eqref{eq:reference_count_rate}.
Lemma~\ref{lem:count_action_continuity} and
\eqref{eq:count_action_bounds} ensure that $\Phi$ is
continuous and bounded. The integral estimates
on open and closed sets from Ref.~\cite{Dembo_Zeitouni},
Exercise~4.3.11, therefore apply; the positive-tail condition is
automatic because $\Phi\leq0$.
For every open set $G$ and closed set $F$, fixed as $N$ varies,
they give
\begin{equation}
\liminf_{N\to\infty}\frac1N\log\mathbb E_0
\Bigl[e^{N\Phi(\gamma(\sigma))}\mathbb I\bigl(\gamma(\sigma)\in G\bigr)\Bigr]
\geq\sup_{\gamma\in G}(\Phi-R_0)(\gamma),
\label{eq:display_009}
\end{equation}
\begin{equation}
\limsup_{N\to\infty}\frac1N\log\mathbb E_0
\Bigl[e^{N\Phi(\gamma(\sigma))}\mathbb I\bigl(\gamma(\sigma)\in F\bigr)\Bigr]
\leq\sup_{\gamma\in F}(\Phi-R_0)(\gamma).
\label{eq:style_s4_count_laplace_upper}
\end{equation}
The comparisons \eqref{eq:count_change_measure_bounds} and
\eqref{eq:style_s4_count_change_measure_upper}, followed by the limit
$\eta\downarrow0$, give the two estimates with cost
$R_0-\Phi-\log(K+1)$, which is exactly \eqref{eq:count_rate}.
For a general Borel set we use
$B^\circ\subseteq B\subseteq\overline B$.
Empty sets are covered by the conventions.

The representation \eqref{eq:style_s4_count_rate_entropy} follows from
\eqref{eq:style_s3_count_action_density}. The relative entropy between the frequencies
and the probabilities $\pi_k(\psi)$ is nonnegative. Moreover,
\begin{equation}
\sum_{k\in\Omega}\partial_\tau\gamma_k\log\partial_\tau\gamma_k\leq0,
\label{eq:style_s4_count_frequency_entropy}
\end{equation}
and $\Phi\geq-\|\log\pi\|_\infty$, hence
\begin{equation}
0\leq
\int_0^1\sum_{k\in\Omega}
\partial_\tau\gamma_k(\tau)
\log\partial_\tau\gamma_k(\tau)\,d\tau
-\Phi(\gamma)
\leq\|\log\pi\|_\infty.
\label{eq:count_rate_bounds}
\end{equation}
Subtracting the continuous function $\log(K+1)+\Phi$ from the reference
cost preserves lower semicontinuity;
compactness of the space gives compact sublevel sets.
The upper bound applied to the entire space, of probability
one, implies that the infimum of the cost is at most zero.
Nonnegativity and compactness show that the minimum
is attained and equals zero.
\end{proof}

In the next subsection we will contract the LDP through
\eqref{eq:counts_to_path} and solve the minimization over
local frequencies, completing Theorem~\ref{thm:path_ldp}.

\subsection{Contraction to the cumulative path}
\label{subsec:scalar_contraction}

The projection
\begin{equation}
\gamma\longmapsto\sum_{k\in\Omega}k\gamma_k
\label{eq:style_s4_count_projection}
\end{equation}
is continuous and surjective onto $Q$, by
\eqref{eq:count_projection_bound} and
\eqref{eq:count_representation}--\eqref{eq:style_s3_count_representation_last}.
The contraction principle, Ref.~\cite{Dembo_Zeitouni},
Theorem~4.2.1, Eq.~(4.2.2), applied to
Theorem~\ref{thm:count_ldp}, therefore gives a scalar LDP with
a good rate function equal to the infimum of the cost
\eqref{eq:count_rate} over compatible counts.
The following lemma allows us to identify this infimum with
$I$ and construct the count vector attaining it.

\begin{lemma}[Minimization over local frequencies]
\label{lem:local_frequency_minimum}
For every $\alpha,\beta\in[0,K]$ we have
\begin{equation}
\min_{\substack{\alpha_k\geq0,\ \sum_{k\in\Omega}\alpha_k=1\\
\sum_{k\in\Omega}k\alpha_k=\alpha}}
\sum_{k\in\Omega}\alpha_k\log\frac{\alpha_k}{\pi_k(\beta)} 
=\sup_{\lambda\in\mathbb R}
\biggl\{\lambda\alpha-\log\biggl[\sum_{k\in\Omega}\pi_k(\beta)e^{k\lambda}\biggr]\biggr\}.
\label{eq:local_frequency_duality}
\end{equation}
The vector attaining the minimum is unique. We denote its
components by $\alpha_k^*(\alpha,\beta)$,
$k\in\Omega$. For $0<\alpha<K$, let
$\lambda(\alpha,\beta)$ be the unique real solution of
\begin{equation}
\alpha=
\frac{\sum_{k\in\Omega}k\pi_k(\beta)
 e^{k\lambda(\alpha,\beta)}}
{\sum_{k\in\Omega}\pi_k(\beta)
 e^{k\lambda(\alpha,\beta)}}.
\label{eq:local_frequency_parameter}
\end{equation}
Then
\begin{equation}
\alpha_k^*(\alpha,\beta)=
\frac{\pi_k(\beta)e^{k\lambda(\alpha,\beta)}}
{\sum_{j\in\Omega}\pi_j(\beta)e^{j\lambda(\alpha,\beta)}},
\qquad 0<\alpha<K,
\label{eq:local_optimal_frequencies}
\end{equation}
while at the endpoints
\begin{equation}
\alpha_k^*(0,\beta)=\delta_k(0),
\label{eq:local_optimal_endpoint_frequencies}
\end{equation}
\begin{equation}
\alpha_k^*(K,\beta)=\delta_k(K).
\label{eq:style_s4_local_optimal_last_frequencies}
\end{equation}
The corresponding minimum values are
$-\log\pi_0(\beta)$ and $-\log\pi_K(\beta)$.
For $0<\alpha<K$, the supremum in
\eqref{eq:local_frequency_duality} is attained uniquely
at $\lambda(\alpha,\beta)$; for $\alpha=0$ and $\alpha=K$
it is obtained, respectively, in the limits
$\lambda\to-\infty$ and $\lambda\to+\infty$ and is not
attained by any real parameter.

The functions $\alpha_k^*$ are continuous, and hence
uniformly continuous, on the square $[0,K]^2$.
The minimum value in
\eqref{eq:local_frequency_duality} is also continuous in
$(\alpha,\beta)$ on this square.
\end{lemma}

\begin{proof}
The derivative with respect to $\lambda$ of the ratio in \eqref{eq:local_frequency_parameter}
equals the variance of the distribution
with probabilities proportional to $\pi_k(\beta)e^{k\lambda}$.
This variance is positive, since all probabilities are
positive and $K\geq1$. The ratio therefore increases strictly
from $0$ to $K$: for $0<\alpha<K$ there exists a unique parameter
$\lambda(\alpha,\beta)$, the unique maximizer of the dual expression,
because the derivative of the expression to be maximized is
$\alpha$ minus that ratio.

The vector \eqref{eq:local_optimal_frequencies} sums to one
and has mean $\alpha$. Set $\lambda=\lambda(\alpha,\beta)$.
For every admissible competing vector,
\begin{equation}
\sum_{k\in\Omega}\alpha_k\log\frac{\alpha_k}{\pi_k(\beta)}
=\sum_{k\in\Omega}\alpha_k\log\frac{\alpha_k}{\alpha_k^*(\alpha,\beta)} 
+\lambda\alpha
-\log\biggl[\sum_{k\in\Omega}\pi_k(\beta)e^{k\lambda}\biggr].
\label{eq:local_entropy_decomposition}
\end{equation}
The relative entropy on the right-hand side is nonnegative and
vanishes exactly when $\alpha_k=\alpha_k^*(\alpha,\beta)$
for every $k$. The identity, also valid for zero competing
frequencies with the convention $0\log0=0$, proves duality and
uniqueness of the interior minimum.

For $\alpha=0$ and $\alpha=K$, the constraints force, respectively,
the vectors $\delta_k(0)$ and $\delta_k(K)$, with costs
$-\log\pi_0(\beta)$ and $-\log\pi_K(\beta)$.
In the dual expression these values are obtained in the limits
$\lambda\to-\infty$ and $\lambda\to+\infty$;
for $\alpha=K$ it suffices to factor out $e^{K\lambda}$ from the sum.
Strict positivity of the other terms makes both
values unattainable for a finite parameter.

For interior continuity, fix $(\alpha,\beta)$ with
$0<\alpha<K$ and $\eta>0$. The ratio in
\eqref{eq:local_frequency_parameter}, evaluated at the two fixed
parameters $\lambda(\alpha,\beta)-\eta$ and $\lambda(\alpha,\beta)+\eta$,
is respectively less and greater than $\alpha$.
By continuity of the $\pi_k$, these
inequalities still hold for $(\alpha',\beta')$ close
to $(\alpha,\beta)$. Monotonicity then gives
$|\lambda(\alpha',\beta')-\lambda(\alpha,\beta)|<\eta$.
This yields joint continuity of the parameter and, by
\eqref{eq:local_optimal_frequencies}, of the frequencies
on $(0,K)\times[0,K]$.

The sum and mean constraints also give,
uniformly in $\beta$,
\begin{equation}
\max_{k\in\Omega}|\alpha_k^*(\alpha,\beta)-\delta_k(0)|\leq\alpha,
\label{eq:local_frequency_boundary_continuity}
\end{equation}
\begin{equation}
\max_{k\in\Omega}|\alpha_k^*(\alpha,\beta)-\delta_k(K)|\leq K-\alpha.
\label{eq:style_s4_local_frequency_last_continuity}
\end{equation}
The frequencies are thus continuous on the entire compact
square, and hence uniformly continuous.
Substituting them into the cost, continuity of $r\log r$,
including at $r=0$, and positivity of the $\pi_k$ also give
continuity of the minimum.
\end{proof}

\begin{proof}[Conclusion of the proof of Theorem~\ref{thm:path_ldp}]
Fix $\varphi\in Q$. The derivatives of any compatible
count vector are, almost everywhere, nonnegative frequencies
summing to one with mean $\partial_\tau\varphi(\tau)$.
Lemma~\ref{lem:local_frequency_minimum}, applied with
$\alpha=\partial_\tau\varphi(\tau)$ and
$\beta=\varphi(\tau)/\tau$, shows that the corresponding cost
is at least $I(\varphi)$. The local minimum is continuous in its two
arguments and lies between $0$ and $\|\log\pi\|_\infty$:
the upper bound follows from
\eqref{eq:style_s4_frequency_entropy_bound}.
The integrand of $I$ is therefore measurable and bounded.

To attain the lower bound, choose a Borel
version of $\partial_\tau\varphi$ taking values in $[0,K]$,
setting it equal to zero where the derivative does not exist.
Use $\psi(\tau)=\varphi(\tau)/\tau$ for $\tau>0$ and
set $\psi(0)=0$, without requiring a limit at the origin.
Define
\begin{equation}
\gamma_k^*(\tau):=
\int_0^\tau
\alpha_k^*\bigl(\partial_s\varphi(s),\psi(s)\bigr)\,ds,
\qquad k\in\Omega,\qquad 0\leq\tau\leq1.
\label{eq:optimal_count_construction}
\end{equation}
The integrands are Borel and lie between zero and one,
by continuity of the $\alpha_k^*$. The counts are therefore
absolutely continuous, nondecreasing and vanish at zero.
Integrating the frequency constraints and using absolute
continuity of $\varphi$, with $\varphi(0)=0$, yields
\begin{equation}
\sum_{k\in\Omega}\gamma_k^*(\tau)=\tau,\qquad 0\leq\tau\leq1,
\label{eq:optimal_count_constraints}
\end{equation}
\begin{equation}
\sum_{k\in\Omega}k\gamma_k^*(\tau)=\varphi(\tau),\qquad 0\leq\tau\leq1.
\label{eq:style_s4_optimal_count_projection}
\end{equation}
The constructed vector is therefore admissible and compatible
with $\varphi$. Local duality, integrated over time, gives
\begin{equation}
\int_0^1\sum_{k\in\Omega}\partial_\tau\gamma_k^*(\tau)
\log\frac{\partial_\tau\gamma_k^*(\tau)}{\pi_k\bigl(\varphi(\tau)/\tau\bigr)}\,d\tau
=I(\varphi).
\label{eq:optimal_count_value}
\end{equation}
The rate obtained by contraction therefore agrees with
\eqref{eq:path_rate}, and both representations
\eqref{eq:path_rate_counts} and \eqref{eq:path_rate_relative_entropy}
hold, with an attained minimum. Nonnegativity and the bound
\eqref{eq:path_rate_bounds} follow from the bounds on the local cost;
lower semicontinuity and compactness of the sublevel sets
are given by the contraction principle. The estimates for
Borel sets follow by applying those for open and closed sets
to the interior and closure, completing the theorem.

For the same path $\varphi$, the minimizing count vector
is unique. Indeed, for another minimizer the difference
between the local cost and the local minimum is nonnegative and has
zero integral, hence vanishes almost everywhere.
Local uniqueness forces the same derivatives as in
\eqref{eq:optimal_count_construction}; absolute continuity
and the zero initial value give equality of the counts
on all of $[0,1]$.
\end{proof}

\subsection{Initial histories and nonhomogeneous urn functions}
\label{subsec:initial_and_inhomogeneous}

We extend to the case $K\geq1$ Corollaries~2 and~3 of
Ref.~\cite{Franchini_URNS}, in the strictly positive class;
for initial conditions the reference is Lemma~16
of Subsection~3.1.5.

Fix, independently of $N$, an integer $m\geq1$
and an initial history
\begin{equation}
\sigma^*:=(\sigma_1^*,\ldots,\sigma_m^*)
\in\Omega^{\{1,\ldots,m\}}.
\label{eq:prescribed_initial_history}
\end{equation}
For $N\geq m$, let $\mathbb P^*$ be the law prescribing these
first $m$ increments and using \eqref{eq:transition_probability}
for $m\leq n<N$. It agrees with $\mathbb P$ conditioned
on the initial history, which has positive probability.
Paths and counts also include the prefix and
are interpolated as in \eqref{eq:path_interpolation}
and \eqref{eq:count_interpolation}.

\begin{theorem}[Finite initial histories]
\label{thm:finite_initial_history}
Under the assumptions of Theorem~\ref{thm:path_ldp},
the laws of the counts under $\mathbb{P}^*$
satisfy the LDP of Theorem~\ref{thm:count_ldp},
with speed $N$ and the same rate function
\eqref{eq:count_rate}. The cumulative paths
satisfy the LDP of Theorem~\ref{thm:path_ldp},
with the same functional $I$.

The conclusions also hold if the first $m$
increments have any fixed distribution,
independent of $N$, and subsequent transitions
are given by \eqref{eq:transition_probability}.
In particular, they hold for the original law
conditioned on $M_m=j$, for every
$j\in\{0,\ldots,Km\}$.
\end{theorem}

\begin{proof}
The map $T$ replaces the first $m$ increments with $\sigma^*$
and leaves the tail unchanged. Set $Q_N:=\mathbb P\circ T^{-1}$.
Coupling $\sigma$ with law $\mathbb P$ to $\sigma':=T\sigma$,
for every history we obtain
\begin{equation}
\|\gamma(\sigma')-\gamma(\sigma)\|_\infty\leq\frac mN,
\label{eq:initial_history_path_bounds}
\end{equation}
\begin{equation}
\|\varphi(\sigma')-\varphi(\sigma)\|_\infty\leq\frac{Km}N.
\label{eq:style_s4_initial_scalar_path_bound}
\end{equation}
The bounds at the nodes extend by linear interpolation.
The counts are therefore exponentially equivalent and
Theorem~4.2.13 of Ref.~\cite{Dembo_Zeitouni} transfers to $Q_N$
the LDP of Theorem~\ref{thm:count_ldp}, with cost
\eqref{eq:count_rate}.

The pushforward law $Q_N$ is not generally $\mathbb P^*$,
which on histories with prefix $\sigma^*$ equals
\begin{equation}
\mathbb{P}^*(\{\sigma\})=
\prod_{n=m}^{N-1}\pi_{\sigma_{n+1}}(\psi_n).
\label{eq:initial_history_weights}
\end{equation}
For a tail $v=(v_{m+1},\ldots,v_N)$ and a prefix
$u\in\Omega^{\{1,\ldots,m\}}$, the cumulative sums of $uv$ and $\sigma^*v$
differ by at most $Km$, and the means at time $n\geq m$
by at most $Km/n$. The modulus \eqref{eq:log_urn_modulus} gives,
uniformly in $u,v$,
\begin{equation}
e^{-b_N}\leq\frac{\mathbb P(v\mid u)}{\mathbb P(v\mid\sigma^*)}\leq e^{b_N},
\label{eq:display_010}
\end{equation}
\begin{equation}
b_N:=\sum_{n=m}^{N-1}\omega(Km/n)=o(N).
\label{eq:style_s4_initial_history_log_error}
\end{equation}
Here the tail probabilities are conditioned on the prefixes;
the limit follows by Cesàro averaging for fixed $m$.
For $N=m$ the product is one and the sum is empty. Moreover,
\begin{equation}
Q_N(\{\sigma^*v\})=\sum_u\mathbb P(u)\mathbb P(v\mid u),
\label{eq:display_011}
\end{equation}
\begin{equation}
\mathbb P^*(\{\sigma^*v\})=\mathbb P(v\mid\sigma^*).
\label{eq:style_s4_initial_conditioned_tail}
\end{equation}
The weights $\mathbb P(u)$ sum to one: the comparison also holds for $Q_N(\{\sigma^*v\})/\mathbb P^*(\{\sigma^*v\})$.
Both laws vanish outside this support;
summing, for every event $E\subseteq\Omega^S$ we obtain
\begin{equation}
e^{-b_N}\mathbb P^*(E)\leq Q_N(E)\leq e^{b_N}\mathbb P^*(E).
\label{eq:display_012}
\end{equation}
Since $b_N=o(N)$, the count LDP transfers from
$Q_N$ to $\mathbb P^*$. The contraction in subsection
\ref{subsec:scalar_contraction} gives the scalar LDP.

A fixed law on the first $m$ increments is a finite
mixture of these laws. For the lower bound, one
component of positive weight suffices; for the upper bound,
the components are summed. The logarithms of the positive weights,
divided by $N$, tend to zero and all components have
the same rate: both LDPs therefore hold.
Conditioning on $M_m=j$ is such a mixture over
initial histories with sum $j$, with weights independent of $N$;
for every $j\in\{0,\ldots,Km\}$ at least one exists
and all have positive probability.
\end{proof}

The uniform distribution of the first increment can therefore
be replaced by any distribution.
For an initial cumulative sum $M_m=j$, every integer $r$ satisfies
\begin{equation}
\mathbb P(M_N=r\mid M_m=j)>0
\qquad\Longleftrightarrow\qquad
r\in\{j,j+1,\ldots,j+K(N-m)\}.
\label{eq:initial_history_reachable_totals}
\end{equation}
Indeed, every integer between zero and $K(N-m)$ is a sum of $N-m$
increments in $\Omega$, and every sequence has positive
probability. The two endpoints, divided by $N$, tend to $0$ and $K$,
but point probabilities remain constrained to this support:
$M_N=0$ requires $j=0$, while $M_N=KN$ requires $j=Km$.

We now turn to urn functions depending on the step.
For $n\geq1$ and $k\in\Omega$, let $\pi_{n,k}$ be continuous
and strictly positive on $[0,K]$, with
\begin{equation}
\sum_{k\in\Omega}\pi_{n,k}(\alpha)=1,
\qquad n\geq1,\qquad\alpha\in[0,K].
\label{eq:inhomogeneous_normalization}
\end{equation}
Let $\widetilde{\mathbb P}$ be the law with a uniform first
increment and transitions
\begin{equation}
\widetilde{\mathbb P}\bigl(\sigma_{n+1}=k\mid\sigma_1,\ldots,\sigma_n\bigr)
=\pi_{n,k}(\psi_n),\qquad 1\leq n<N.
\label{eq:inhomogeneous_transitions}
\end{equation}
The index $n$ thus identifies the transition from step
$n$ to step $n+1$.

\begin{theorem}[Nonhomogeneous urn functions]
\label{thm:inhomogeneous_urn}
Suppose that the functions $\pi_{n,k}$ satisfy
the assumptions just stated and that
\begin{equation}
\max_{k\in\Omega}\sup_{\alpha\in[0,K]}
|\pi_{n,k}(\alpha)-\pi_k(\alpha)|
\longrightarrow0
\qquad\text{as }n\to\infty,
\label{eq:inhomogeneous_uniform_convergence}
\end{equation}
where the limiting functions $\pi_k$ satisfy the
assumptions of Theorem~\ref{thm:path_ldp}.
Then, under $\widetilde{\mathbb{P}}$, the counts
and cumulative paths satisfy, respectively,
the LDPs of Theorems~\ref{thm:count_ldp} and
\ref{thm:path_ldp}, with speed $N$ and the
rate functions constructed from $\pi$.

The same conclusions hold with a prescribed finite
initial history or a fixed
initial distribution on the first $m$
increments, as in
Theorem~\ref{thm:finite_initial_history}.
No rate of convergence is required
in \eqref{eq:inhomogeneous_uniform_convergence}.
\end{theorem}

\begin{proof}
Set
\begin{equation}
p_*:=\min_{k\in\Omega}\min_{\alpha\in[0,K]}\pi_k(\alpha)>0.
\label{eq:style_s4_limit_probability_minimum}
\end{equation}
By uniform convergence
\eqref{eq:inhomogeneous_uniform_convergence}, the $\pi_{n,k}$
are at least $p_*/2$ for sufficiently large $n$.
The mean value theorem applied to the logarithm then gives
\begin{equation}
\Delta_n:=\max_{k\in\Omega}\sup_{\alpha\in[0,K]}|\log\pi_{n,k}(\alpha)-\log\pi_k(\alpha)|,
\label{eq:inhomogeneous_log_convergence}
\end{equation}
\begin{equation}
\Delta_n\leq\frac2{p_*}\max_{k\in\Omega}\sup_{\alpha\in[0,K]}
|\pi_{n,k}(\alpha)-\pi_k(\alpha)|\longrightarrow0.
\label{eq:style_s4_inhomogeneous_log_estimate}
\end{equation}
The finitely many initial terms excluded from the bound are nevertheless
finite, by continuity and positivity on the compact interval.

Fix an integer $m\geq1$ and a distribution $\mu$ on the first
$m$ increments, independent of $N$ and common to both processes.
Denote by $\mathbb P^\mu$ and $\widetilde{\mathbb P}^\mu$
the laws with this initialization and homogeneous and nonhomogeneous
transitions, respectively. On the common support, consisting of the
histories with $\mu(\sigma_1,\ldots,\sigma_m)>0$, the initial factor
cancels and
\begin{equation}
\frac{d\widetilde{\mathbb{P}}^\mu}{d\mathbb{P}^\mu}(\sigma)
=\prod_{n=m}^{N-1}
\frac{\pi_{n,\sigma_{n+1}}(\psi_n)}
{\pi_{\sigma_{n+1}}(\psi_n)}.
\label{eq:inhomogeneous_initial_density}
\end{equation}
Outside this support both laws vanish.
Convergence of the logarithms and Cesàro averaging give
\begin{equation}
\sup_{\sigma\in\operatorname{supp}\mathbb P^\mu}
\left|\frac1N\log\frac{d\widetilde{\mathbb P}^\mu}{d\mathbb P^\mu}(\sigma)\right|
\leq\frac1N\sum_{n=m}^{N-1}\Delta_n\longrightarrow0.
\label{eq:inhomogeneous_density_bound}
\end{equation}
For $N=m$ the product is one and the sum is empty.
For every $\eta>0$ and sufficiently large $N$, it therefore follows,
simultaneously for all events $E\subseteq\Omega^S$,
including those depending on $N$, that
\begin{equation}
e^{-N\eta}\mathbb{P}^\mu(E)
\leq\widetilde{\mathbb{P}}^\mu(E)
\leq e^{N\eta}\mathbb{P}^\mu(E).
\label{eq:inhomogeneous_event_bounds}
\end{equation}
For fixed $m$ the threshold is independent of both $E$ and $\mu$.
Theorem~\ref{thm:finite_initial_history} supplies the homogeneous
LDPs: applying the comparison to open and closed sets of counts
or paths and sending $\eta\downarrow0$ yields
the same LDPs for the nonhomogeneous laws.
This includes the original law, with $m=1$ and uniform $\mu$,
the prescribed prefix, with point-mass $\mu$, and any fixed initial
mixture, without requiring a rate of convergence.
\end{proof}

\section{Entropy of the final state}
\label{sec:terminal_entropy}

The rescaled final state is the mean
$\psi_N=M_N/N=\varphi(1\mid\sigma)$.
The path LDP therefore allows us to
study the probabilities that $\psi_N$ belongs
to an interval through evaluation at the final
time. In this section we derive the corresponding
variational representation, pass to point
probabilities and determine the values
at the endpoints of the domain.

\subsection{Variational representation}
\label{subsec:terminal_variational}

For every $\alpha\in[0,K]$, the set $Q(\alpha)$
defined in \eqref{eq:terminal_path_space} is
nonempty and compact. Since $I$ is finite and
lower semicontinuous, it attains its
minimum there. Define
\begin{equation}
\phi(\alpha):=
-\min_{\varphi\in Q(\alpha)}I(\varphi),
\qquad 0\leq\alpha\leq K.
\label{eq:terminal_entropy_variational}
\end{equation}
We adopt the sign convention of
Ref.~\cite{Franchini_URNS}, Theorem~4,
Eq.~(2.12)--(2.13): $\phi$ is a nonpositive
entropy, while the rate function of the
final state is $-\phi$. We will first prove
the interpretation of $\phi$ through interval
probabilities; identification with the limit
of point probabilities will be proved
in Subsection~\ref{subsec:point_probabilities}.

\begin{theorem}[Entropy of the final state]
\label{thm:terminal_entropy}
Under the assumptions of Theorem~\ref{thm:path_ldp},
the function $\phi$ defined in
\eqref{eq:terminal_entropy_variational} is
continuous on $[0,K]$ and satisfies
\begin{equation}
-\|\log\pi\|_\infty\leq\phi(\alpha)\leq0,\qquad\alpha\in[0,K].
\label{eq:terminal_entropy_bounds}
\end{equation}
\begin{equation}
\max_{\alpha\in[0,K]}\phi(\alpha)=0.
\label{eq:style_s5_entropy_maximum}
\end{equation}
The laws of $\psi_N$ satisfy an LDP on
$[0,K]$ with speed $N$ and good rate
function $-\phi$: for every Borel set
$B\subseteq[0,K]$ independent of $N$,
\begin{equation}
\liminf_{N\to\infty}\frac{1}{N}
\log\mathbb{P}(\psi_N\in B)
\geq\sup_{\alpha\in B^\circ}\phi(\alpha),
\label{eq:terminal_ldp_lower}
\end{equation}
\begin{equation}
\limsup_{N\to\infty}\frac{1}{N}
\log\mathbb{P}(\psi_N\in B)
\leq\sup_{\alpha\in\overline B}\phi(\alpha).
\label{eq:terminal_ldp_upper}
\end{equation}
Interior and closure are relative to $[0,K]$;
we use the conventions $\sup\varnothing=-\infty$
and $\log0=-\infty$.

The representation in terms of counts is
\begin{equation}
\phi(\alpha)=
\max_{\substack{\gamma\ \mathrm{admissible}\\
\sum_{k\in\Omega}k\gamma_k(1)=\alpha}}
\left\{\Phi(\gamma)
-\int_0^1\sum_{k\in\Omega}
\partial_\tau\gamma_k(\tau)
\log\partial_\tau\gamma_k(\tau)\,d\tau\right\}.
\label{eq:terminal_entropy_counts}
\end{equation}
The maximum is attained and the convention
$0\log0=0$ is retained. In this formula
only the final value of the weighted sum
of counts is fixed.
\end{theorem}

\begin{proof}
The map $\varphi\mapsto\varphi(1)$ is continuous
in the uniform distance and its image is
$[0,K]$. The contraction principle,
Theorem~4.2.1 of Ref.~\cite{Dembo_Zeitouni},
applied to Theorem~\ref{thm:path_ldp},
gives the LDP for the final state with cost
$\min_{\varphi\in Q(\alpha)}I(\varphi)=-\phi(\alpha)$.
This yields
\eqref{eq:terminal_ldp_lower}--\eqref{eq:terminal_ldp_upper}
and goodness of the rate function.
In particular, $-\phi$ is lower
semicontinuous, that is, $\phi$ is upper
semicontinuous. The pointwise inequalities in
\eqref{eq:terminal_entropy_bounds} follow from
\eqref{eq:path_rate_bounds}.

To prove the other semicontinuity, fix
$\alpha\in[0,K]$ and a minimizing path
$\varphi\in Q(\alpha)$, so that
$I(\varphi)=-\phi(\alpha)$.
For every $\beta\in[0,K]$ define a
path $\varphi'$ by
\begin{equation}
\varphi'(\tau):= \frac{\beta}{\alpha}\varphi(\tau),\qquad 0\leq\beta<\alpha,
\label{eq:terminal_path_adjustment}
\end{equation}
\begin{equation}
\varphi'(\tau):= \varphi(\tau),\qquad \beta=\alpha,
\label{eq:style_global_terminal_path_adjustment_2}
\end{equation}
\begin{equation}
\varphi'(\tau):= \varphi(\tau)+
\frac{\beta-\alpha}{K-\alpha}
\bigl(K\tau-\varphi(\tau)\bigr),\qquad \alpha<\beta\leq K.
\label{eq:style_global_terminal_path_adjustment_3}
\end{equation}
In the first case $\alpha>0$, while in the third
$\alpha<K$, so all denominators are
positive. The two constructions are convex
combinations of $\varphi$ with the zero
path and with $K\tau$, respectively. They therefore preserve
the constraints defining $Q$ and give
$\varphi'(1)=\beta$; hence
$\varphi'\in Q(\beta)$.

As $\beta\to\alpha$, \eqref{eq:terminal_path_adjustment}--\eqref{eq:style_global_terminal_path_adjustment_3} gives
convergence of the derivatives almost everywhere and of the means for every
$\tau>0$. At the endpoints, $\varphi$ is respectively zero or equal
to $K\tau$ and the constructed path is $\varphi'(\tau)=\beta\tau$.
The local cost in \eqref{eq:path_rate} is continuous and lies between
$0$ and $\|\log\pi\|_\infty$ by
Lemma~\ref{lem:local_frequency_minimum} and the bounds in
Subsection~\ref{subsec:scalar_contraction}. By dominated convergence,
$I(\varphi')\to I(\varphi)$; since $\varphi'\in Q(\beta)$,
\begin{equation}
\liminf_{\beta\to\alpha}\phi(\beta)
\geq-I(\varphi)=\phi(\alpha).
\label{eq:display_013}
\end{equation}
The upper semicontinuity already proved establishes continuity
of $\phi$ on $[0,K]$, uniform by compactness.

For \eqref{eq:terminal_entropy_counts}, we minimize the cost
\eqref{eq:count_rate} first over counts with fixed scalar projection
$\varphi$ and then over $\varphi\in Q(\alpha)$.
The inner minimum is $I(\varphi)$ and is attained by
\eqref{eq:path_rate_relative_entropy}; the outer one is
$-\phi(\alpha)$ and is attained by compactness and semicontinuity.
This gives equality and attainment in the representation
in terms of counts.

Finally, the upper bound of the LDP applied
to the entire interval $[0,K]$, which has
probability one, gives
$0\leq\sup_{\alpha\in[0,K]}\phi(\alpha)$.
Since $\phi\leq0$ and is continuous, its
maximum is zero. This proves
\eqref{eq:style_s5_entropy_maximum}.
\end{proof}

Continuity of the entropy allows us to derive
exact limits for nondegenerate intervals.
This is the step corresponding to
Lemma~17 of Ref.~\cite{Franchini_URNS},
Eq.~(3.68)--(3.70).
For $0\leq a<b\leq K$, the interior of $[a,b]$
relative to $[0,K]$ is dense in $[a,b]$.
Continuity therefore gives
\begin{equation}
\sup_{\alpha\in[a,b]^\circ}\phi(\alpha)
=\max_{\alpha\in[a,b]}\phi(\alpha),
\qquad 0\leq a<b\leq K.
\label{eq:terminal_interval_continuity}
\end{equation}
The two LDP bounds agree and we obtain
\begin{equation}
\lim_{N\to\infty}\frac{1}{N}
\log\mathbb{P}(a\leq\psi_N\leq b)
=\max_{\alpha\in[a,b]}\phi(\alpha).
\label{eq:terminal_interval_limit}
\end{equation}

Similarly, for fixed $\alpha\in[0,K]$ and
$\delta>0$, the window
$\{\beta\in[0,K]:|\beta-\alpha|<\delta\}$
is open in the relative topology and its
closure is
$[\max\{0,\alpha-\delta\},\min\{K,\alpha+\delta\}]$.
The supremum of $\phi$ on the open window
agrees with the maximum on its closure.
Therefore, with $\delta$ fixed independently
of $N$,
\begin{equation}
\lim_{N\to\infty}\frac{1}{N}
\log\mathbb{P}(|\psi_N-\alpha|<\delta)
=\max_{\substack{\beta\in[0,K]\\|\beta-\alpha|\leq\delta}}
\phi(\beta).
\label{eq:terminal_window_limit}
\end{equation}
Then letting $\delta$ tend to zero,
continuity of $\phi$ gives
\begin{equation}
\phi(\alpha)=
\lim_{\delta\downarrow0}\,
\lim_{N\to\infty}\frac{1}{N}
\log\mathbb{P}(|\psi_N-\alpha|<\delta),
\qquad \alpha\in[0,K].
\label{eq:terminal_local_entropy}
\end{equation}
For $\alpha=0$ and $\alpha=K$ as well, the windows
are neighborhoods relative to the domain, and the two
formulas remain valid.

Theorems~\ref{thm:finite_initial_history} and
\ref{thm:inhomogeneous_urn} give the same
path LDP for the initializations
and nonhomogeneous urn functions considered
in Subsection~\ref{subsec:initial_and_inhomogeneous}.
Contraction at the final time and the arguments
just given show that Theorem~\ref{thm:terminal_entropy}
and the limits
\eqref{eq:terminal_interval_limit}--\eqref{eq:terminal_local_entropy}
also hold for these laws,
with the same function $\phi$ constructed
from the limiting urn functions.

The limit for the point event
$M_N=\lfloor N\alpha\rfloor$ requires a
further comparison between interval
probabilities and probabilities of individual states.
We will prove it in Subsection~\ref{subsec:point_probabilities},
following the discrete method of
Subsection~3.2.1 of Ref.~\cite{Franchini_URNS}.
For prescribed initial histories, this
comparison will preserve the attainability
constraint
\eqref{eq:initial_history_reachable_totals}.

\subsection{Point probabilities}
\label{subsec:point_probabilities}

To pass from intervals to individual states,
we adapt the discrete comparison in
Ref.~\cite{Franchini_URNS}, Subsection~3.2.1,
in particular Eq.~(3.75)--(3.77), to extend
Theorem~4, Eq.~(2.12)--(2.13), to increments in $\Omega$.
We modify the increments of a history to
give it a different final cumulative sum and
control both the change in its weight
and the number of histories that can produce
the same modified history. The second
bound accounts for all values
in $\Omega=\{0,\ldots,K\}$.

\begin{lemma}[Comparison of point probabilities]
\label{lem:point_mass_comparison}
Under the assumptions of Theorem~\ref{thm:path_ldp},
let $r,s\in\{0,\ldots,KN\}$ and
$h:=|r-s|$. With the norm and modulus defined
in \eqref{eq:log_urn_norm} and
\eqref{eq:log_urn_modulus}, we have
\begin{equation}
\left|\log
\frac{\mathbb{P}(M_N=r)}{\mathbb{P}(M_N=s)}\right|
\leq2h\|\log\pi\|_\infty
+\sum_{n=1}^{N-1}\omega\!\left(\frac hn\right)
+\log\left[\sum_{\ell=0}^{\min\{h,N\}}\binom N\ell K^\ell\right].
\label{eq:point_mass_finite_comparison}
\end{equation}
For the combinatorial contribution, set
\begin{equation}
\eta_K(\delta):=-\delta\log\delta-(1-\delta)\log(1-\delta)+\delta\log K,
\qquad0<\delta<1.
\label{eq:style_global_combinatorial_entropy}
\end{equation}
In particular, if $0<\delta\leq1/2$ and $|r-s|\leq\delta N$, then
\begin{equation}
\frac1N\left|\log\frac{\mathbb P(M_N=r)}{\mathbb P(M_N=s)}\right|
\leq2\|\log\pi\|_\infty(\delta+\sqrt\delta)+\omega(\sqrt\delta)+\eta_K(\delta).
\label{eq:point_mass_uniform_comparison}
\end{equation}
The right-hand side tends to zero as
$\delta\downarrow0$ and is independent of
$N,r,s$.

The same estimates hold for the law
$\mathbb{P}^*$ with a prescribed initial history
of length $m$ and cumulative sum $j$, provided that
$r,s\in\{j,\ldots,j+K(N-m)\}$.
\end{lemma}

\begin{proof}
The case $r=s$ is immediate. First assume
$s>r$ and consider a history
$\sigma$ with $M_N=r$. We construct $\sigma'$
by scanning the indices in increasing order and
adding to each increment the amount still
needed, up to the allowed maximum $K$.
Explicitly, for $n=1,\ldots,N$, set
\begin{equation}
\sigma_n':=\sigma_n+
\min\left\{K-\sigma_n,\,
h-\sum_{i=1}^{n-1}(\sigma_i'-\sigma_i)\right\}.
\label{eq:point_edit_increase}
\end{equation}
The amount still to be added is always
nonnegative, and the total available capacity
is $KN-r\geq s-r=h$. The construction therefore
ends with
\begin{equation}
\sum_{n=1}^{N}\sigma_n'=s.
\label{eq:style_s5_edited_total}
\end{equation}
If $s<r$, we traverse the same indices and subtract,
at each position, the minimum of $\sigma_n$ and the amount
still to be removed, namely
\begin{equation}
h-\sum_{i<n}(\sigma_i-\sigma_i').
\label{eq:style_s5_remaining_decrease}
\end{equation}
The total available amount is $r\geq h$, so the procedure
ends with cumulative sum $s$ and all increments remain in $\Omega$.

In both constructions all changes
have the same sign and
\begin{equation}
\sum_{n=1}^{N}|\sigma_n'-\sigma_n|=h,
\label{eq:point_edit_size}
\end{equation}
\begin{equation}
|\{n\in S:\sigma_n'\neq\sigma_n\}|\leq\min\{h,N\}.
\label{eq:style_s5_edited_positions}
\end{equation}
The second inequality uses the fact that
each modified increment changes by at least
one unit. Denoting the mean of the modified
history by $\psi_n'$, for every $n\geq1$
we have
\begin{equation}
|\psi_n'-\psi_n|
=\frac1n\left|
\sum_{i=1}^{n}(\sigma_i'-\sigma_i)\right|
\leq\frac hn.
\label{eq:point_edit_average_bound}
\end{equation}

To compare the transition probabilities, add and subtract the logarithm of the probability of the modified increment evaluated at the original mean. Summing the resulting bounds using
\eqref{eq:point_edit_size} and \eqref{eq:style_s5_edited_positions}, we obtain
\begin{equation}
\left|\log\pi_{\sigma_{n+1}'}(\psi_n')
-\log\pi_{\sigma_{n+1}}(\psi_n)\right|
\leq \omega\!\left(\frac hn\right)
+2\|\log\pi\|_\infty\,\mathbb I(\sigma_{n+1}'\neq\sigma_{n+1}),
\label{eq:point_edit_transition_bound}
\end{equation}
\begin{equation}
|\mathcal{A}(\sigma')-\mathcal{A}(\sigma)|
\leq2h\|\log\pi\|_\infty
+\sum_{n=1}^{N-1}\omega\!\left(\frac hn\right).
\label{eq:point_edit_action_bound}
\end{equation}
In \eqref{eq:point_edit_transition_bound}, the first term controls the change in the mean; the second controls the change in the increment. The uniform initial
factor is identical for the two histories.
The modulus $\omega$ is defined for every nonnegative
argument, so the estimate also covers
the first steps, for which $h/n$
may exceed $K$.

Fix a resulting history $\sigma'$. Every history
that produces it differs from it
at most $\min\{h,N\}$ positions.
The histories that differ at exactly
$\ell$ positions number $\binom N\ell K^\ell$:
we choose the positions and, at each one, one
of the $K$ values different from the assigned value.
The multiplicity of the construction is therefore
at most the sum in
\eqref{eq:point_mass_finite_comparison}.
In particular,
\begin{equation}
\sum_{\substack{\sigma\in\Omega^S\\M_N=r}}
W(\sigma')
\leq
\left[\sum_{\ell=0}^{\min\{h,N\}}
\binom N\ell K^\ell\right]
\mathbb{P}(M_N=s).
\label{eq:point_edit_weight_multiplicity}
\end{equation}
The comparison \eqref{eq:point_edit_action_bound},
summed over histories with cumulative sum $r$ and
combined with
\eqref{eq:point_edit_weight_multiplicity},
gives the upper bound for
$\log\mathbb{P}(M_N=r)-\log\mathbb{P}(M_N=s)$
in
\eqref{eq:point_mass_finite_comparison}.
Repeating the construction from $s$ to $r$
gives the reverse bound and hence the
estimate with the absolute value.

Now assume $h\leq\delta N$, with
$0<\delta\leq1/2$. In the sum we separate
the integers $n<\sqrt\delta\,N$ from the
remaining ones. The first group contains at
most $\sqrt\delta\,N$ terms, each
bounded by $2\|\log\pi\|_\infty$;
in the second group we have $h/n\leq\sqrt\delta$.
Therefore
\begin{equation}
\frac1N\sum_{n=1}^{N-1}\omega\!\left(\frac hn\right)
\leq
2\|\log\pi\|_\infty\sqrt\delta+\omega(\sqrt\delta).
\label{eq:point_edit_modulus_sum}
\end{equation}
For the combinatorial term, since
$K\geq1$ and $0<\delta/(1-\delta)\leq1$,
the binomial theorem gives
\begin{equation}
\sum_{\ell=0}^{h}\binom N\ell K^\ell
\leq K^{\delta N}\left(\frac{1-\delta}{\delta}\right)^{\delta N}
\sum_{\ell=0}^{N}\binom N\ell\left(\frac{\delta}{1-\delta}\right)^\ell,
\label{eq:point_edit_multiplicity_bound}
\end{equation}
\begin{equation}
K^{\delta N}\Bigl(\frac{1-\delta}{\delta}\Bigr)^{\delta N}
\sum_{\ell=0}^{N}\binom N\ell\Bigl(\frac{\delta}{1-\delta}\Bigr)^\ell
=e^{N\eta_K(\delta)}.
\label{eq:style_s5_multiplicity_evaluation}
\end{equation}
The inequality uses
$K^\ell\leq K^{\delta N}$ and
$(\delta/(1-\delta))^\ell
\geq(\delta/(1-\delta))^{\delta N}$
for $\ell\leq h\leq\delta N$.
Substituting the estimates
\eqref{eq:point_edit_modulus_sum} and
\eqref{eq:point_edit_multiplicity_bound}--\eqref{eq:style_s5_multiplicity_evaluation}
into
\eqref{eq:point_mass_finite_comparison}
proves \eqref{eq:point_mass_uniform_comparison}.
Its right-hand side tends to zero because
$\omega(u)\to0$ as $u\downarrow0$ and
$\delta\log\delta\to0$.

With a prescribed initial history, we leave
the first $m$ increments unchanged and apply
the constructions only to positions
$m+1,\ldots,N$. If the initial cumulative sum is
$j$, the available capacities for increasing
and decreasing the total are, respectively,
$K(N-m)-(r-j)$ and $r-j$.
They suffice when both totals belong to
the stated attainable set.
The ratio of the conditional weights agrees
with the ratio of the weights of the two
original histories, because the factor corresponding
to the initial history is common.
The sum in the action estimate can start
at $n=m$, and the multiplicity is at most
\begin{equation}
\sum_{\ell=0}^{\min\{h,N-m\}}
\binom{N-m}{\ell}K^\ell
\leq
\sum_{\ell=0}^{\min\{h,N\}}
\binom N\ell K^\ell.
\label{eq:point_edit_initial_multiplicity}
\end{equation}
The same upper bounds, written in terms of $N$,
therefore also hold for $\mathbb{P}^*$.
\end{proof}
\begin{theorem}[Point probabilities and variational principle]
\label{thm:point_probabilities}
Under the assumptions of Theorem~\ref{thm:path_ldp},
let $r_N\in\{0,\ldots,KN\}$ be a sequence
with $r_N/N\to\alpha\in[0,K]$.
Then
\begin{equation}
\lim_{N\to\infty}\frac1N
\log\mathbb{P}(M_N=r_N)
=\phi(\alpha)
=-\min_{\varphi\in Q(\alpha)}I(\varphi).
\label{eq:point_probability_limit}
\end{equation}
In particular, for every $\alpha\in[0,K]$,
\begin{equation}
\lim_{N\to\infty}\frac1N
\log\mathbb{P}(M_N=\lfloor N\alpha\rfloor)
=\phi(\alpha).
\label{eq:point_probability_floor_limit}
\end{equation}
Moreover, convergence is uniform:
\begin{equation}
\lim_{N\to\infty}
\max_{r\in\{0,\ldots,KN\}}
\left|
\frac1N\log\mathbb{P}(M_N=r)
-\phi\!\left(\frac rN\right)
\right|=0.
\label{eq:point_probability_uniform_limit}
\end{equation}

The conclusions also hold for the
initializations and nonhomogeneous
urn functions in subsection
\ref{subsec:initial_and_inhomogeneous},
replacing $\mathbb{P}$ by the law
under consideration and restricting the states to
those with positive probability.
For a prescribed initial history of
length $m$ and cumulative sum $j$, this
restriction is precisely
\begin{equation}
r_N\in\{j,\ldots,j+K(N-m)\}.
\label{eq:point_probability_initial_support}
\end{equation}
The maximum in
\eqref{eq:point_probability_uniform_limit}
is also taken only over attainable states.
Formula
\eqref{eq:point_probability_floor_limit}
holds for these initializations when
$\lfloor N\alpha\rfloor$ is attainable
for all sufficiently large $N$.
\end{theorem}

\begin{proof}
Let $\mathbb Q$ be the original law or a law with a fixed
prefix. On both supports, the window limit
\eqref{eq:terminal_window_limit} holds by
Subsection~\ref{subsec:terminal_variational}, as does the comparison in
Lemma~\ref{lem:point_mass_comparison}. Fix attainable states
$r_N$ with $r_N/N\to\alpha$; in the case of a prefix the support is
\eqref{eq:point_probability_initial_support}.
For every $\delta>0$, when $N$ is
sufficiently large, the event
$\{M_N=r_N\}$ is contained in the window
$\{|\psi_N-\alpha|<\delta\}$.
Equation \eqref{eq:terminal_window_limit} gives
\begin{equation}
\limsup_{N\to\infty}\frac1N
\log\mathbb Q(M_N=r_N)
\leq
\max_{\substack{\beta\in[0,K]\\|\beta-\alpha|\leq\delta}}
\phi(\beta).
\label{eq:point_probability_upper_bound}
\end{equation}
Letting $\delta$ tend to zero,
continuity of $\phi$ gives the required
upper bound.

For the lower bound, fix
$\eta>0$ and choose $0<\delta\leq1/2$
so that the right-hand side of
\eqref{eq:point_mass_uniform_comparison}
is at most $\eta$.
For all sufficiently large $N$, we have
$|r_N/N-\alpha|<\delta/2$.
The window $|\psi_N-\alpha|<\delta/2$
contains at most $KN+1$ terminal states
and contains an attainable state, $r_N$, so it has positive mass.
Choose a state $s_N$ in the window with
maximal probability. Since there are at most $KN+1$ states,
its mass is at least that of the window divided by
$KN+1$; moreover, $|s_N-r_N|<\delta N$.
Lemma~\ref{lem:point_mass_comparison} and, with the
window fixed as $N\to\infty$,
\eqref{eq:terminal_window_limit} therefore give
\begin{equation}
\frac1N\log\mathbb Q(M_N=r_N)
\geq\frac1N\log\mathbb Q(|\psi_N-\alpha|<\delta/2)
-\frac{\log(KN+1)}N-\eta,
\label{eq:point_probability_lower_comparison}
\end{equation}
\begin{equation}
\liminf_{N\to\infty}\frac1N\log\mathbb Q(M_N=r_N)
\geq\max_{\substack{\beta\in[0,K]\\|\beta-\alpha|\leq\delta/2}}
\phi(\beta)-\eta\geq\phi(\alpha)-\eta.
\label{eq:point_probability_lower_bound}
\end{equation}
Since $\eta$ is arbitrary, this proves
\eqref{eq:point_probability_limit}.
The windows are relative to $[0,K]$, so
the proof also covers
$\alpha=0$ and $\alpha=K$.
Taking $r_N=\lfloor N\alpha\rfloor$, when it is attainable
for all sufficiently large $N$, gives
\eqref{eq:point_probability_floor_limit}.

To prove
\eqref{eq:point_probability_uniform_limit},
suppose that along a subsequence
the maximum stays above a positive
constant. For each index choose an
attainable state realizing that maximum.
By compactness of $[0,K]$, along
a further subsequence the ratios
$r_N/N$ converge to a value $\alpha$.
The preceding argument also applies along
this subsequence and gives
$N^{-1}\log\mathbb Q(M_N=r_N)\to\phi(\alpha)$.
On the other hand, continuity implies
$\phi(r_N/N)\to\phi(\alpha)$.
The difference therefore tends to zero,
contradicting the choice of the
states. This proves uniform
convergence.

Now consider a fixed distribution
on the first $m$ increments. It gives a
finite mixture of the laws with prescribed
initial histories, and we consider
only the components with positive weights.
For every $\eta>0$, the uniform convergence
just proved ensures that, for all sufficiently large $N$,
each component law $\mathbb Q$ satisfies
\begin{equation}
e^{N[\phi(r/N)-\eta]}
\leq\mathbb Q(M_N=r)
\leq e^{N[\phi(r/N)+\eta]}
\label{eq:style_s5_component_mass_bounds}
\end{equation}
on all of its attainable states, simultaneously for all components.
Writing $\mathbb Q_{\mathrm{mix}}$ for the mixture and $w_{\min}$
for its smallest positive weight, we obtain
\begin{equation}
w_{\min}e^{N[\phi(r/N)-\eta]}
\leq\mathbb Q_{\mathrm{mix}}(M_N=r)
\leq e^{N[\phi(r/N)+\eta]}
\label{eq:style_s5_mixture_mass_bounds}
\end{equation}
for every state attainable by the mixture.
The weight $w_{\min}$ is a strictly positive
constant independent of $N$.
Its logarithm, divided by $N$, tends
to zero; we therefore obtain the same
uniform convergence on all states
attainable by the mixture. This
argument includes conditioning
on the initial cumulative sum $M_m=j$ alone.

Finally, for nonhomogeneous urn
functions, the comparison
\eqref{eq:inhomogeneous_event_bounds}
also holds for each event
$\{M_N=r\}$, simultaneously in $r$.
The logarithms of the corresponding masses,
divided by $N$, therefore differ by
a quantity tending uniformly to
zero. The comparison already covers every
common initialization, by
\eqref{eq:inhomogeneous_initial_density}.
Thus all pointwise limits and
uniform convergence transfer, on the
common supports, to the
nonhomogeneous law.
\end{proof}

Formula
\eqref{eq:point_probability_floor_limit}
therefore identifies the entropy defined
in \eqref{eq:terminal_entropy_variational}
with the limit for the probabilities of
individual states, uniformly on their respective supports.

\subsection{Endpoint values}
\label{subsec:terminal_endpoints}

At the endpoints of the domain, the terminal constraint
determines the entire path. This
allows us to evaluate explicitly both
the variational principle and the probabilities
of the corresponding histories. The direct
endpoint calculation in the binary case
is recalled in Ref.~\cite{Franchini_URNS},
Subsection~3.2.1, after Eq.~(3.81).
The exact probabilities below use
the discrete initialization
\eqref{eq:initial_distribution}.

\begin{theorem}[Endpoint values]
\label{thm:terminal_endpoints}
Under the assumptions of Theorem~\ref{thm:path_ldp},
the terminal entropy satisfies
\begin{equation}
\phi(0)=\log\pi_0(0)<0,
\label{eq:terminal_endpoint_values}
\end{equation}
\begin{equation}
\phi(K)=\log\pi_K(K)<0.
\label{eq:style_s5_endpoint_value_K}
\end{equation}
The only path in $Q(0)$ is
$\varphi(\tau)=0$; the only path
in $Q(K)$ is $\varphi(\tau)=K\tau$.
They therefore attain the respective
minima in
\eqref{eq:terminal_entropy_variational}.

With the first increment uniform on
$\Omega$, for every $N\geq1$ the following
exact identities hold:
\begin{equation}
\mathbb P(M_N=0)=\frac1{K+1}\,\pi_0(0)^{N-1},
\label{eq:terminal_endpoint_exact_masses}
\end{equation}
\begin{equation}
\mathbb P(M_N=KN)=\frac1{K+1}\,\pi_K(K)^{N-1}.
\label{eq:style_s5_endpoint_exact_mass_K}
\end{equation}
\end{theorem}

\begin{proof}
If $\varphi\in Q(0)$, monotonicity and
the conditions $\varphi(0)=\varphi(1)=0$
force $\varphi(\tau)=0$ on $[0,1]$.
If $\varphi\in Q(K)$, the function
$K\tau-\varphi(\tau)$ is nondecreasing,
since its derivative is
$K-\partial_\tau\varphi\geq0$ almost
everywhere. It is zero at both $0$ and
$1$, and therefore
$\varphi(\tau)=K\tau$ on $[0,1]$.

At velocities $0$ and $K$ and for every
mean $\beta\in[0,K]$, Lemma~\ref{lem:local_frequency_minimum}
gives the local costs
$-\log\pi_0(\beta)$ and $-\log\pi_K(\beta)$.
Along the two paths, the means are respectively
$0$ and $K$: integration in \eqref{eq:path_rate}
and \eqref{eq:terminal_entropy_variational}
give \eqref{eq:terminal_endpoint_values}--\eqref{eq:style_s5_endpoint_value_K}.
The extreme frequencies \eqref{eq:local_optimal_endpoint_frequencies}--\eqref{eq:style_s4_local_optimal_last_frequencies},
imposed by the local constraints, also give the only
compatible counts, respectively
$\gamma_k(\tau)=\tau\delta_k(0)$ and
$\gamma_k(\tau)=\tau\delta_k(K)$ for $k\in\Omega$
and $0\leq\tau\leq1$.
Since $K\geq1$ and all probabilities are
positive and sum to one, $0<\pi_0(0),\pi_K(K)<1$:
both entropy values are finite and strictly negative.

At the discrete level, $M_N=0$ forces
$\sigma_n=0$ for every $n$, whereas
$M_N=KN$ forces $\sigma_n=K$ for
every $n$. The means along these histories
are constantly zero and $K$,
respectively. Multiplying the
transition probabilities and
including the initial factor
$1/(K+1)$ gives
\eqref{eq:terminal_endpoint_exact_masses}--\eqref{eq:style_s5_endpoint_exact_mass_K}.
For $N=1$, the powers have exponent
zero and recover the initial distribution.
\end{proof}

Now consider the prescribed
initial history
\eqref{eq:prescribed_initial_history},
of length $m\geq1$ and cumulative sum $j\in\{0,\ldots,Km\}$,
with $m$ and $j$ fixed as
$N$ varies. Under the law $\mathbb{P}^*$,
the endpoints of the support in
\eqref{eq:initial_history_reachable_totals}
are $j$ and $j+K(N-m)$. To reach them,
all increments after time $m$
must be zero or $K$,
respectively. From
\eqref{eq:initial_history_weights}
we obtain, for $N\geq m$,
\begin{equation}
\mathbb P^*(M_N=j)=\prod_{n=m}^{N-1}\pi_0\!\left(\frac jn\right),
\label{eq:terminal_initial_endpoint_masses}
\end{equation}
\begin{equation}
\mathbb P^*(M_N=j+K(N-m))=\prod_{n=m}^{N-1}\pi_K\!\left(K-\frac{Km-j}{n}\right).
\label{eq:style_s5_initial_endpoint_mass_K}
\end{equation}
For $N=m$, both products are
empty and equal to one. The means appearing
in the factors belong to $[0,K]$:
in the first case the cumulative sum remains $j$;
in the second, at time $n$, it equals
$j+K(n-m)=Kn-(Km-j)$.
The formulas depend on the initial
history only through $m$ and
$j$, and therefore also hold for the
original law conditioned only
on the event $M_m=j$.

The events $M_N=0$ and $M_N=KN$
are attainable under
$\mathbb{P}^*$ only if,
respectively, $j=0$ and $j=Km$.

For the nonhomogeneous urn functions
of Theorem~\ref{thm:inhomogeneous_urn},
the same identification of the extreme
histories gives, for $N\geq1$,
\begin{equation}
\widetilde{\mathbb P}(M_N=0)=\frac1{K+1}\prod_{n=1}^{N-1}\pi_{n,0}(0),
\label{eq:terminal_inhomogeneous_endpoint_masses}
\end{equation}
\begin{equation}
\widetilde{\mathbb P}(M_N=KN)=\frac1{K+1}\prod_{n=1}^{N-1}\pi_{n,K}(K).
\label{eq:style_s5_inhomogeneous_endpoint_mass_K}
\end{equation}
With the prescribed initial history,
we instead have, for $N\geq m$,
\begin{equation}
\widetilde{\mathbb P}^{*}(M_N=j)=\prod_{n=m}^{N-1}\pi_{n,0}\!\left(\frac jn\right),
\label{eq:terminal_inhomogeneous_initial_masses}
\end{equation}
\begin{equation}
\widetilde{\mathbb P}^{*}(M_N=j+K(N-m))=\prod_{n=m}^{N-1}\pi_{n,K}\!\left(K-\frac{Km-j}{n}\right).
\label{eq:style_s5_inhomogeneous_initial_mass_K}
\end{equation}
Here too, empty products equal one.
Strict positivity of the transitions preserves
the supports and attainability conditions
of the homogeneous case.

Finally, for an arbitrary fixed initial
distribution on histories of
length $m$, the probability of
$M_N=0$ is the initial weight of the
all-zero history multiplied by
$\pi_0(0)^{N-m}$; the probability of
$M_N=KN$ is the initial weight of the
all-$K$ history multiplied by
$\pi_K(K)^{N-m}$. In the nonhomogeneous
case the powers are replaced,
respectively, by the products
\begin{equation}
\prod_{n=m}^{N-1}\pi_{n,0}(0)
\label{eq:style_s5_zero_endpoint_product}
\end{equation}
and
\begin{equation}
\prod_{n=m}^{N-1}\pi_{n,K}(K).
\label{eq:style_s5_K_endpoint_product}
\end{equation}
The corresponding initial weight is positive
if and only if the endpoint is attainable;
if it is zero, the event is impossible for every
$N\geq m$ and the logarithm of its mass is $-\infty$.

For attainable states, all exponential limits
of the preceding masses follow from
Theorem~\ref{thm:point_probabilities}.
The limit \eqref{eq:point_probability_limit},
uniform convergence \eqref{eq:point_probability_uniform_limit}
on the respective supports, and the values
\eqref{eq:terminal_endpoint_values}--\eqref{eq:style_s5_endpoint_value_K} give,
for every sequence of attainable states $r_N$,
\begin{equation}
\lim_{N\to\infty}\frac1N\log\mathbb P(M_N=r_N)=\log\pi_0(0),\qquad\text{if }\frac{r_N}{N}\longrightarrow0,
\label{eq:terminal_endpoint_sequence_limits}
\end{equation}
\begin{equation}
\lim_{N\to\infty}\frac1N\log\mathbb P(M_N=r_N)=\log\pi_K(K),\qquad\text{if }\frac{KN-r_N}{N}\longrightarrow0.
\label{eq:style_s5_endpoint_sequence_limit_K}
\end{equation}
The conclusions hold, on their respective supports, for all laws
in Subsection~\ref{subsec:initial_and_inhomogeneous}.
In particular, the attainable endpoints with a prefix satisfy
$j/N\to0$ and $[j+K(N-m)]/N\to K$; the absolute endpoints $0,KN$
are included when the corresponding initial weight is positive.
In the nonhomogeneous case, uniform convergence in
Theorem~\ref{thm:inhomogeneous_urn} suffices, with no additional rate or summability requirement.

\section{Optimal paths}
\label{sec:optimal_trajectories}

For every $\alpha\in[0,K]$, the minimum in
\eqref{eq:terminal_entropy_variational}
is attained. We call a path in $Q(\alpha)$ that attains
this minimum optimal. In this section we study
zero-cost optimal paths, following the method of Lemma~18 and
Corollaries~5--7 of
Ref.~\cite{Franchini_URNS}.
Their equation uses the mean function $\bar\pi$ defined in
\eqref{eq:mean_increment}; the frequencies
of the individual increments remain specified by the full vector $\pi$.

\subsection{Zero cost and fixed points}
\label{subsec:zero_cost_fixed_points}

The representation
\eqref{eq:path_rate_relative_entropy}
allows us to identify zero-cost paths through the equality
between local frequencies and transition probabilities.
This is the step corresponding to Lemma~18 of
Ref.~\cite{Franchini_URNS},
Subsection~3.2.2, Eqs.~(3.83)--(3.85).

\begin{lemma}[Zero-cost paths]
\label{lem:zero_cost_paths}
Let $\varphi\in Q$ and let
$\psi(\tau)=\varphi(\tau)/\tau$ for
$\tau>0$. We have $I(\varphi)=0$ if
and only if
\begin{equation}
\partial_\tau\varphi(\tau)
=\bar\pi\!\left(\frac{\varphi(\tau)}{\tau}\right)
\qquad\text{for almost every }\tau\in(0,1).
\label{eq:zero_cost_equation}
\end{equation}
In this case the minimizing count vector is given by
\begin{equation}
\gamma_k^*(\tau)
=\int_0^\tau\pi_k(\psi(s))\,ds,
\qquad k\in\Omega,\qquad 0\leq\tau\leq1.
\label{eq:zero_cost_counts}
\end{equation}
Moreover, $\varphi\in C^1([0,1])$,
the mean $\psi$ is continuously
differentiable on $(0,1]$ and satisfies
\begin{equation}
\tau\,\partial_\tau\psi(\tau)
=\bar\pi(\psi(\tau))-\psi(\tau),
\qquad 0<\tau\leq1.
\label{eq:zero_cost_average_equation}
\end{equation}
Derivatives at the endpoints are understood to be one-sided.
The mean $\psi$ is monotone, that is, either nondecreasing or
nonincreasing, and we have
\begin{equation}
\lim_{\tau\downarrow0}\psi(\tau)=\partial_\tau\varphi(0)\in C,
\label{eq:zero_cost_average_bounds}
\end{equation}
\begin{equation}
\inf C\leq\psi(\tau)\leq\sup C,\qquad0<\tau\leq1.
\label{eq:style_s6_mean_confinement}
\end{equation}
\end{lemma}

\begin{proof}
By Lemma~\ref{lem:local_frequency_minimum}, the integrand in
\eqref{eq:path_rate} is the minimum of the relative entropy with respect to
$\pi(\psi(\tau))$ among the frequencies with mean
$\partial_\tau\varphi(\tau)$. This minimum is nonnegative and vanishes
exactly when the vector $\pi(\psi(\tau))$ satisfies the mean constraint,
that is, when
$\partial_\tau\varphi(\tau)=\bar\pi(\psi(\tau))$.
An integral of a nonnegative function is zero if and only if its integrand
is zero almost everywhere: this proves
\eqref{eq:zero_cost_equation} in both directions.
In this case the minimizing frequencies are uniquely
$\pi_k(\psi(\tau))$ almost everywhere, and the construction
\eqref{eq:optimal_count_construction} becomes
\eqref{eq:zero_cost_counts}. The constraints
\eqref{eq:optimal_count_constraints}--\eqref{eq:style_s4_optimal_count_projection} and the uniqueness already proved in
Subsection~\ref{subsec:scalar_contraction} therefore give
admissibility, projection onto $\varphi$, and uniqueness of the counts.

On every interval $[\varepsilon,1]$,
with $\varepsilon>0$, the right-hand side
of \eqref{eq:zero_cost_equation} is
continuous. The absolute continuity of
$\varphi$ then implies that $\varphi$
is continuously differentiable and that
the equation holds at every point of
$(0,1]$. Differentiating the ratio
$\psi(\tau)=\varphi(\tau)/\tau$
yields
\eqref{eq:zero_cost_average_equation}.
This is the change of variables in Eqs.~(3.86)--(3.87) of
Ref.~\cite{Franchini_URNS}, expressed
in our notation for the mean.

Monotonicity follows from an identity that
does not require uniqueness for the differential equation. For
$0<\tau_1<\tau_2\leq1$, the chain rule
and \eqref{eq:zero_cost_average_equation}
give
\begin{equation}
\int_{\psi(\tau_1)}^{\psi(\tau_2)}\bigl(\bar\pi(u)-u\bigr)\,du
=\int_{\tau_1}^{\tau_2}\bigl(\bar\pi(\psi(s))-\psi(s)\bigr)\partial_s\psi(s)\,ds,
\label{eq:zero_cost_monotonicity_identity}
\end{equation}
\begin{equation}
\int_{\tau_1}^{\tau_2}\bigl(\bar\pi(\psi(s))-\psi(s)\bigr)\partial_s\psi(s)\,ds
=\int_{\tau_1}^{\tau_2}s\bigl(\partial_s\psi(s)\bigr)^2\,ds.
\label{eq:style_s6_monotonicity_energy}
\end{equation}
If $\psi(\tau_1)=\psi(\tau_2)$,
the last integral is zero. Its integrand
is continuous and nonnegative,
so $\psi$ is constant between those two
times. A continuous function that
changed direction would take the same
value at two times without being
constant between them. Hence $\psi$
is monotone on $(0,1]$.

Since $0\leq\psi\leq K$, monotonicity gives
a limit $\beta\in[0,K]$ as $\tau\downarrow0$.
Integrating \eqref{eq:zero_cost_equation} from zero to $\tau$,
we have
\begin{equation}
\psi(\tau)=\frac1\tau\int_0^\tau\bar\pi(\psi(s))\,ds.
\label{eq:display_014}
\end{equation}
The continuity of $\bar\pi$ makes this integral mean
converge to $\bar\pi(\beta)$; hence
$\beta=\bar\pi(\beta)$ and $\beta\in C$.

By the definition of $C$ and the signs
at the endpoints of the domain,
$\bar\pi(u)-u>0$ for $u<\inf C$
and $\bar\pi(u)-u<0$ for $u>\sup C$.
If $\psi(\tau)<\inf C$ at a positive
time, the derivative of the mean
would be positive at that point.
The monotonicity just proved would
therefore force $\psi$ to be nondecreasing,
and we would have
$\beta\leq\psi(\tau)<\inf C$,
contradicting $\beta\in C$.
Similarly, if $\psi(\tau)>\sup C$,
the derivative would be negative and the
mean nonincreasing, giving
$\beta\geq\psi(\tau)>\sup C$,
which is again impossible. This yields
the pointwise bounds in
\eqref{eq:style_s6_mean_confinement}.

Finally, $\varphi(\tau)/\tau\to\beta$
and $\varphi(0)=0$ give the right-hand
derivative $\partial_\tau\varphi(0)=\beta$.
By \eqref{eq:zero_cost_equation},
the derivative at positive times converges
to $\bar\pi(\beta)=\beta$.
Hence $\varphi\in C^1([0,1])$,
and the assertion in
\eqref{eq:zero_cost_average_bounds}
is also proved.
\end{proof}

In particular, the straight-line
path $\varphi(\tau)=\alpha\tau$
has zero cost if and only if
$\alpha\in C$. For this path,
the minimizing counts are
\begin{equation}
\gamma_k^*(\tau)=\tau\pi_k(\alpha),
\qquad k\in\Omega,\qquad\alpha\in C.
\label{eq:zero_cost_fixed_point_counts}
\end{equation}
The following result identifies all
terminal states with zero entropy.
It extends Corollary~5 of
Ref.~\cite{Franchini_URNS} to the strictly
positive class considered here.

\begin{theorem}[Zero-entropy interval]
\label{thm:zero_entropy_interval}
Under the assumptions of Theorem~\ref{thm:path_ldp},
we have
\begin{equation}
\phi(\alpha)=0
\qquad\Longleftrightarrow\qquad
\alpha\in[\inf C,\sup C].
\label{eq:zero_entropy_interval}
\end{equation}
For every $\alpha$ in this interval
there exists a zero-cost path in $Q(\alpha)$.
The optimal paths are precisely the solutions of
\eqref{eq:zero_cost_equation} with
$\varphi(1)=\alpha$.
For states outside the interval, we have
\begin{equation}
-\infty<\phi(\alpha)<0,
\qquad
\alpha\in[0,\inf C)\cup(\sup C,K].
\label{eq:strictly_negative_entropy}
\end{equation}
\end{theorem}

\begin{proof}
If $\phi(\alpha)=0$, attainment
of the minimum in
\eqref{eq:terminal_entropy_variational}
provides a path
$\varphi\in Q(\alpha)$ with zero cost.
Lemma~\ref{lem:zero_cost_paths}
then gives
$\alpha=\psi(1)\in[\inf C,\sup C]$.

Conversely, let
$\alpha\in[\inf C,\sup C]$.
If $\alpha\in C$, the path
$\varphi(\tau)=\alpha\tau$ has zero
cost, as just proved.
If $\alpha\notin C$, define
\begin{equation}
a:=\max\bigl(C\cap[0,\alpha]\bigr),
\label{eq:zero_cost_adjacent_fixed_points}
\end{equation}
\begin{equation}
b:=\min\bigl(C\cap[\alpha,K]\bigr).
\label{eq:style_s6_adjacent_fixed_point_b}
\end{equation}
Both sets are nonempty and compact.
Thus $a,b\in C$,
$a<\alpha<b$, and $(a,b)$ contains no
fixed points. The function
$\bar\pi(u)-u$ is continuous, nonzero,
and has constant sign on $(a,b)$.

We solve
\eqref{eq:zero_cost_average_equation}
by proceeding from the terminal condition
$\psi(1)=\alpha$ toward earlier times.
As long as $\psi$ remains in $(a,b)$,
separation of variables gives
\begin{equation}
\int_\alpha^{\psi(\tau)}
\frac{du}{\bar\pi(u)-u}
=\log\tau.
\label{eq:zero_cost_separated_equation}
\end{equation}
The antiderivative on the left has a
continuous, nonzero derivative on the interval,
so it is strictly monotone and
has an inverse. This constructs
the solution near $\tau=1$ and
allows it to be extended as long as its
value remains in $(a,b)$.
If $\bar\pi(u)-u>0$, as
$\tau$ decreases the mean decreases from
$\alpha$ toward $a$; if the sign is
negative, the mean increases from
$\alpha$ toward $b$.

The solution can be extended
to arbitrarily small times,
or it reaches the corresponding
endpoint of $(a,b)$ at a time
$\tau_*>0$. Indeed, stopping at a
positive time with a limiting value
inside $(a,b)$ would be ruled out
by local invertibility of the
antiderivative. In the second case,
we extend $\psi$ by keeping it
equal to the endpoint reached for
$0<\tau\leq\tau_*$.
Since that endpoint belongs to
$C$, the derivative in
\eqref{eq:zero_cost_average_equation}
tends to zero as $\tau$
approaches $\tau_*$ from above.
The extension is therefore
continuously differentiable and satisfies
the same equation at the joining point
as well. In both cases we obtain
a solution on $(0,1]$
with values in $[a,b]$.

Now set $\varphi(0)=0$ and
$\varphi(\tau)=\tau\psi(\tau)$
for $\tau>0$. The boundedness of
$\psi$ ensures continuity
at zero. At every positive time,
\begin{equation}
\partial_\tau\varphi(\tau)
=\psi(\tau)+\tau\partial_\tau\psi(\tau)
=\bar\pi(\psi(\tau))\in(0,K).
\label{eq:zero_cost_constructed_path}
\end{equation}
Integrating between two positive times
and passing to the limit if one of them
tends to zero yields the constraints
\eqref{eq:path_space}. Hence
$\varphi\in Q(\alpha)$ and, by
Lemma~\ref{lem:zero_cost_paths},
$I(\varphi)=0$. Nonnegativity
of the cost gives $\phi(\alpha)=0$.

For $\alpha$ outside the interval,
there are no zero-cost paths in $Q(\alpha)$.
The minimum is attained
and is therefore strictly positive;
it is also finite by
\eqref{eq:path_rate_bounds}.
This proves
\eqref{eq:strictly_negative_entropy}.
Finally, when $\phi(\alpha)=0$,
a path is optimal exactly
when it has zero cost, and its
characterization follows from
Lemma~\ref{lem:zero_cost_paths}.
\end{proof}

The proof uses only the compactness
of $C$ and, for each terminal state
outside $C$ but between its
endpoints, the single component $(a,b)$
containing it. It therefore also covers
infinite sets $C$ and sets with
intervals of fixed points.
If $\inf C=\sup C$, the zero-entropy
interval reduces to
that single fixed point.

The set of zero-cost cumulative
paths and the interval
\eqref{eq:zero_entropy_interval}
are determined by the mean function
$\bar\pi$ alone. The counts
\eqref{eq:zero_cost_counts} depend
on the individual $\pi_k$, as does
the value of positive costs in
\eqref{eq:path_rate}.
By Theorems~\ref{thm:finite_initial_history}
and~\ref{thm:inhomogeneous_urn}, the
initializations and the nonhomogeneous
laws already considered have the
same rate function. The same
characterizations therefore hold,
with the limiting urn functions.

The construction by separation of
variables produces at least one
optimal path for every
$\alpha\in[\inf C,\sup C]$.
Continuity alone does not imply
uniqueness when a solution
reaches a fixed point.
In Subsections~\ref{subsec:terminal_between_fixed_points}
and~\ref{subsec:terminal_at_fixed_points}
we will specify the formulas for the
paths and the conditions for
reaching or leaving a fixed
point at a positive time.

\subsection{Terminal states between fixed points}
\label{subsec:terminal_between_fixed_points}

Let $\alpha\in[\inf C,\sup C]\setminus C$
and let $a,b$ be the fixed points defined
in \eqref{eq:zero_cost_adjacent_fixed_points}--\eqref{eq:style_s6_adjacent_fixed_point_b}.
Then $a<\alpha<b$ and $(a,b)$ is the
component of the complement of $C$
containing $\alpha$.
We make the construction by separation of variables explicit,
following Corollary~6 of
Ref.~\cite{Franchini_URNS},
Eqs.~(2.14)--(2.18), and its
proof in Subsection~3.2.2,
Eqs.~(3.90)--(3.93).
We will specify both uniqueness when the
mean remains in $[a,b]$ and the
form of all optimal paths
with terminal state $\alpha$.

The function $\bar\pi(u)-u$ has a
constant, nonzero sign on $(a,b)$.
Set
\begin{equation}
\alpha^*:=
\begin{cases}
a,&\bar\pi(\alpha)>\alpha,\\
b,&\bar\pi(\alpha)<\alpha.
\end{cases}
\label{eq:component_initial_fixed_point}
\end{equation}
This is the fixed point toward which
the mean moves when solving the equation
\eqref{eq:zero_cost_average_equation}
backward from the terminal value $\alpha$,
as long as it remains in the component $(a,b)$.
For $\alpha,u\in(a,b)$ define,
using the integral notation of
Ref.~\cite{Franchini_URNS},
\begin{equation}
F_\pi(\alpha,u):=
\int_u^\alpha\frac{dz}{\bar\pi(z)-z}.
\label{eq:component_separation_integral}
\end{equation}
This function depends on the urn
functions through their mean.
We have $F_\pi(\alpha,\alpha)=0$ and
\begin{equation}
\partial_uF_\pi(\alpha,u)
=-\frac{1}{\bar\pi(u)-u},
\qquad u\in(a,b).
\label{eq:component_integral_derivative}
\end{equation}
On the interval $(a,\alpha]$ if
$\alpha^*=a$, or on $[\alpha,b)$
if $\alpha^*=b$, the values of
$F_\pi(\alpha,u)$ are nonnegative
and increase as $u$ approaches
$\alpha^*$. We define the improper
value and the corresponding time by
\begin{equation}
F_\pi(\alpha,\alpha^*)
:=\lim_{\substack{u\to\alpha^*\\u\in(a,b)}}F_\pi(\alpha,u)\in(0,+\infty],
\label{eq:component_departure_time}
\end{equation}
\begin{equation}
\tau_\alpha^*:=\exp\bigl[-F_\pi(\alpha,\alpha^*)\bigr]\in[0,1),
\label{eq:style_s6_departure_time}
\end{equation}
with the convention $e^{-\infty}=0$.

We denote by $F_{\pi,\alpha}^{-1}$
the inverse of $u\mapsto F_\pi(\alpha,u)$
restricted to $(a,\alpha]$ in the first
case and to $[\alpha,b)$ in the second.
Its domain is
$[0,F_\pi(\alpha,\alpha^*))$;
if the improper value is infinite,
the domain is $[0,+\infty)$.
The inverse takes the value $\alpha$
at $q=0$ and converges to $\alpha^*$
as $q$ tends to the upper endpoint
of its domain.

\begin{theorem}[Paths between fixed points]
\label{thm:component_optimal_paths}
With the preceding definitions, set
\begin{equation}
\psi^*(\tau):= \alpha^*,\qquad 0\leq\tau\leq\tau_\alpha^*,
\label{eq:component_optimal_average}
\end{equation}
\begin{equation}
\psi^*(\tau):= F_{\pi,\alpha}^{-1}\!\left(\log\frac1\tau\right),\qquad \tau_\alpha^*<\tau\leq1.
\label{eq:style_global_component_optimal_average_2}
\end{equation}
and
\begin{equation}
\varphi^*(\tau):=\tau\psi^*(\tau),
\qquad 0\leq\tau\leq1.
\label{eq:component_optimal_path}
\end{equation}
Then $\varphi^*\in Q(\alpha)$ and
$I(\varphi^*)=0$. It is the unique
optimal path in $Q(\alpha)$
whose mean belongs to $[a,b]$
at every positive time.

For every optimal path
$\varphi\in Q(\alpha)$, we have
\begin{equation}
\varphi(\tau)=\varphi^*(\tau),
\qquad \tau_\alpha^*\leq\tau\leq1.
\label{eq:component_common_terminal_path}
\end{equation}
If $\tau_\alpha^*=0$, the optimal
path is therefore unique throughout
$Q(\alpha)$.
If $\tau_\alpha^*>0$, all and only
the optimal paths are obtained
by choosing a path
$\varphi'\in Q(\alpha^*)$ with
$I(\varphi')=0$ and setting
\begin{equation}
\varphi(\tau)= \tau_\alpha^*\varphi'\!\left(\frac{\tau}{\tau_\alpha^*}\right),\qquad 0\leq\tau\leq\tau_\alpha^*,
\label{eq:component_all_optimal_paths}
\end{equation}
\begin{equation}
\varphi(\tau)= \varphi^*(\tau),\qquad \tau_\alpha^*<\tau\leq1.
\label{eq:style_global_component_all_optimal_paths_2}
\end{equation}
In this case uniqueness in $Q(\alpha)$
is equivalent to uniqueness of the
zero-cost path in $Q(\alpha^*)$.
\end{theorem}

\begin{proof}
Formula \eqref{eq:component_optimal_average}--\eqref{eq:style_global_component_optimal_average_2} makes explicit
the construction in Theorem~\ref{thm:zero_entropy_interval}, with
endpoint $\alpha^*$ and time fixed by
\eqref{eq:style_s6_departure_time}. That construction proves
$\varphi^*\in Q(\alpha)$ and $I(\varphi^*)=0$;
Lemma~\ref{lem:zero_cost_paths} gives its $C^1$ regularity,
including at the joining point.

Now consider any
optimal path
$\varphi\in Q(\alpha)$.
Theorem~\ref{thm:zero_entropy_interval}
ensures that its cost is zero.
The mean $\psi$ therefore satisfies
\eqref{eq:zero_cost_average_equation}
and is monotone.
Let $(t,1]$, with $0\leq t<1$,
be the maximal time interval
containing $1$ on which
$\psi(\tau)\in(a,b)$.
Separation of variables requires
\begin{equation}
F_\pi(\alpha,\psi(\tau))
=\log\frac1\tau,
\qquad t<\tau\leq1.
\label{eq:component_any_optimal_identity}
\end{equation}
The sign of $\bar\pi(u)-u$ determines
the direction of backward motion:
the mean remains between $\alpha$ and
$\alpha^*$ and is uniquely determined
by the inverse.

By \eqref{eq:component_any_optimal_identity}, the domain
$[0,F_\pi(\alpha,\alpha^*))$ of the inverse requires
$\tau>\tau_\alpha^*$ for every $\tau\in(t,1]$; hence
$t\geq\tau_\alpha^*$.
If $t>\tau_\alpha^*$, then $t>0$ and, by continuity,
\begin{equation}
\psi(t)=F_{\pi,\alpha}^{-1}\!\left(\log\frac1t\right)\in(a,b),
\label{eq:display_015}
\end{equation}
contradicting the maximality of $(t,1]$. Hence $t=\tau_\alpha^*$,
even when the improper integral is infinite and the joining time is zero.
The identity with the canonical path extends to the joining point by
continuity; at zero both cumulative paths are zero.
This gives \eqref{eq:component_common_terminal_path},
as well as global uniqueness when $\tau_\alpha^*=0$.

Now suppose that $\tau_\alpha^*>0$.
For the optimal path just
considered, define
\begin{equation}
\varphi'(s):=
\frac{\varphi(\tau_\alpha^*s)}{\tau_\alpha^*},
\qquad 0\leq s\leq1.
\label{eq:component_rescaled_initial_path}
\end{equation}
Rescaling preserves the increment
constraints in
\eqref{eq:path_space}, and the terminal
value is $\varphi'(1)=\alpha^*$.
Hence $\varphi'\in Q(\alpha^*)$.

Conversely, for any
$\varphi'\in Q(\alpha^*)$,
formula
\eqref{eq:component_all_optimal_paths}--\eqref{eq:style_global_component_all_optimal_paths_2}
defines a path in
$Q(\alpha)$: the two portions agree
at the joining point with value
$\tau_\alpha^*\alpha^*$, and each
satisfies the increment constraints.
On the initial portion, with
$s=\tau/\tau_\alpha^*$, the velocity
is $\partial_s\varphi'(s)$ and the mean
is $\varphi'(s)/s$.
The change of variables
$\tau=\tau_\alpha^*s$ in the integral
\eqref{eq:path_rate}, together with the
zero cost of the terminal portion, gives
the exact identity
\begin{equation}
I(\varphi)=\tau_\alpha^*I(\varphi'),
\qquad \varphi'\in Q(\alpha^*).
\label{eq:component_rescaling_cost}
\end{equation}
Since $\tau_\alpha^*>0$, the total
cost is zero if and only
if that of $\varphi'$ is zero.
This proves both directions
of the representation
\eqref{eq:component_all_optimal_paths}--\eqref{eq:style_global_component_all_optimal_paths_2}.
When $I(\varphi')=0$, admissibility and
\eqref{eq:component_rescaling_cost} allow us to apply
Lemma~\ref{lem:zero_cost_paths} to the complete path,
which is therefore $C^1$ at the joining point as well.
The correspondence between $\varphi'$ and
the complete path is one-to-one,
which gives the assertion
on global uniqueness.

It remains to prove uniqueness under the
constraint on the mean. If
$\bar\pi(\alpha)>\alpha$,
the mean of every optimal path
is nondecreasing and equals $a$ at
time $\tau_\alpha^*>0$.
At earlier times it is therefore at
most $a$; if it also lies in $[a,b]$,
it must be constantly equal
to $a$. In the case
$\bar\pi(\alpha)<\alpha$,
the mean is nonincreasing and the
same argument forces it to equal
$b$ before the joining time.
In both cases we obtain
\eqref{eq:component_optimal_average}--\eqref{eq:style_global_component_optimal_average_2}.
The case $\tau_\alpha^*=0$ is already
covered by global uniqueness.
\end{proof}

The criterion for a positive joining
time follows directly from the
improper integrals:
\begin{equation}
F_\pi(\alpha,\alpha^*)=
\begin{cases}
\displaystyle\int_a^\alpha
\frac{du}{\bar\pi(u)-u},
&\bar\pi(\alpha)>\alpha,\\[3mm]
\displaystyle\int_\alpha^b
\frac{du}{u-\bar\pi(u)},
&\bar\pi(\alpha)<\alpha.
\end{cases}
\label{eq:component_endpoint_integrals}
\end{equation}
If the relevant integral diverges,
$\tau_\alpha^*=0$ and the nonconstant
solution is defined at every
positive time.
If it is finite, $0<\tau_\alpha^*<1$;
the path
\eqref{eq:component_optimal_path}
coincides with $\alpha^*\tau$ up
to that time and then reaches the
terminal state $\alpha$.
For a general path in
\eqref{eq:component_all_optimal_paths}--\eqref{eq:style_global_component_all_optimal_paths_2},
the preceding portion may be
nonlinear: $\alpha^*$ fixes
its mean value at the joining point.

In particular, if $\bar\pi$ is
Lipschitz near $\alpha^*$ on the
side of the component, there is a
constant $c>0$ such that
$|\bar\pi(u)-u|\leq c|u-\alpha^*|$
for $u$ sufficiently close
to $\alpha^*$ in $(a,b)$.
The reciprocal has a divergent integral,
so $\tau_\alpha^*=0$.
This criterion is sufficient;
the exact condition is divergence
of the integral in
\eqref{eq:component_endpoint_integrals}.

The dependence of the time on the terminal
value is computed within the same
component. For $\alpha,\beta\in(a,b)$,
additivity of the integral gives
\begin{equation}
\tau_\beta^*
=\tau_\alpha^*\exp\bigl[-F_\pi(\beta,\alpha)\bigr].
\label{eq:component_time_comparison}
\end{equation}
The exponential factor is finite
and strictly positive.
Hence the time is zero for every
terminal value in the component,
or it is positive for all of them.
In the positive case, the integrand is continuous in the component:
the same relation shows that the time is continuously
differentiable, strictly decreasing if
$\bar\pi(u)-u>0$ and strictly increasing if
$\bar\pi(u)-u<0$.
Again in the positive case,
it tends to $1$ as the terminal
value approaches $\alpha^*$.
At the other endpoint of the component
its limit is
\begin{equation}
\exp\left[-\int_a^b
\frac{du}{|\bar\pi(u)-u|}\right],
\label{eq:component_full_interval_time}
\end{equation}
with the same convention for
the infinite integral. The limit
can therefore be positive when
the singularities at both
endpoints are integrable.
If instead $\tau_\alpha^*=0$ throughout
the component, both limits
of the time are zero.

To complete the description in the
case $K>1$, consider the counts.
If $\tau_\alpha^*>0$, let
$\gamma_k$ and $\gamma_k'$ be the
minimizing counts associated,
respectively, with a path
\eqref{eq:component_all_optimal_paths}--\eqref{eq:style_global_component_all_optimal_paths_2}
and with its rescaled initial
path $\varphi'$.
Formula \eqref{eq:zero_cost_counts}
and the change of variables give
\begin{equation}
\gamma_k(\tau)= \tau_\alpha^*\gamma_k'\!\left(\frac{\tau}{\tau_\alpha^*}\right),\qquad 0\leq\tau\leq\tau_\alpha^*,
\label{eq:component_optimal_counts}
\end{equation}
\begin{equation}
\gamma_k(\tau)= \tau_\alpha^*\gamma_k'(1)
+\int_{\tau_\alpha^*}^{\tau}\pi_k(\psi^*(s))\,ds,\qquad \tau_\alpha^*<\tau\leq1.
\label{eq:style_global_component_optimal_counts_2}
\end{equation}
for every $k\in\Omega$.
For the path
\eqref{eq:component_optimal_path},
we have $\varphi'(s)=\alpha^*s$ and
$\gamma_k'(s)=s\pi_k(\alpha^*)$.
If $\tau_\alpha^*=0$, the counts
are directly given by
\begin{equation}
\gamma_k(\tau)=\int_0^\tau\pi_k(\psi^*(s))\,ds.
\label{eq:style_s6_counts_without_prefix}
\end{equation}

On the common terminal portion, all
optimal paths have the same
local frequencies
$\pi_k(\psi^*(\tau))$.
When $\tau_\alpha^*>0$, the cumulative
counts may nevertheless retain
constant differences determined by the
values $\gamma_k'(1)$ from the preceding
portion. These satisfy
\begin{equation}
\sum_{k\in\Omega}\gamma_k'(1)=1,
\label{eq:component_initial_count_totals}
\end{equation}
\begin{equation}
\sum_{k\in\Omega}k\gamma_k'(1)=\alpha^*.
\label{eq:style_s6_prefix_count_mean}
\end{equation}
For $K=1$ the two relations determine
both counts, which therefore also
coincide on the terminal portion.
For $K>1$ the two relations do not
in general determine all the
components; formula
\eqref{eq:style_global_component_optimal_counts_2}
preserves this dependence on the initial
portion. The explicit example in
Subsection~\ref{subsec:nonlipschitz_example} and
Figure~\ref{fig:component_memory} exhibits these persistent differences.

The representation
\eqref{eq:component_all_optimal_paths}--\eqref{eq:style_global_component_all_optimal_paths_2}
reduces any alternatives among
optimal paths, when
$\tau_\alpha^*>0$, to the study of
zero-cost paths with terminal
state in $C$.
This will be the subject of
Subsection~\ref{subsec:terminal_at_fixed_points}.

\subsection{Terminal states coinciding with fixed points}
\label{subsec:terminal_at_fixed_points}

For every $\alpha\in C$, the path
$\varphi(\tau)=\alpha\tau$ has zero
cost. We study when other
optimal paths with the same
terminal state exist. We first consider
a single component of the complement
of $C$, following Corollary~7
of Ref.~\cite{Franchini_URNS},
Eqs.~(2.19)--(2.23), and its
proof in Subsection~3.2.2,
Eqs.~(3.94)--(3.95).
We will then give a general criterion that
also covers accumulation points
of $C$ and an integral characterization
of all zero-cost paths.

Let $a,b\in C$, with $a<b$ and
$C\cap(a,b)=\varnothing$.
We retain the point $\alpha^*$ from
\eqref{eq:component_initial_fixed_point}
and denote by $\alpha^\dagger$
the other endpoint of the component:
\begin{equation}
(\alpha^*,\alpha^\dagger)=
\begin{cases}
(a,b),&\bar\pi(u)>u\text{ on }(a,b),\\
(b,a),&\bar\pi(u)<u\text{ on }(a,b).
\end{cases}
\label{eq:fixed_point_component_endpoints}
\end{equation}
The direction of the equation for the mean
is therefore from $\alpha^*$ toward
$\alpha^\dagger$.
We extend $F_\pi$ in its first argument
by means of the improper integral
\begin{equation}
F_\pi(\alpha^\dagger,u)
:=\lim_{\substack{\alpha\to\alpha^\dagger\\\alpha\in(a,b)}}F_\pi(\alpha,u),
\qquad u\in(a,b),
\label{eq:fixed_point_arrival_integral}
\end{equation}
\begin{equation}
F_\pi(\alpha^\dagger,u)
=\int_u^{\alpha^\dagger}\frac{dz}{\bar\pi(z)-z}\in(0,+\infty],
\qquad u\in(a,b).
\label{eq:style_s6_arrival_integral_value}
\end{equation}
The orientation of the integral makes
it positive in both cases.
It is finite at one interior point
if and only if it is finite at all of them,
since between two interior points
the integrand is continuous and bounded.

When the value in
\eqref{eq:style_s6_arrival_integral_value}
is finite, set
\begin{equation}
F_\pi(\alpha^\dagger,\alpha^*):=\int_a^b\frac{du}{|\bar\pi(u)-u|},
\label{eq:fixed_point_full_transit_time}
\end{equation}
\begin{equation}
\theta^*:=\exp\bigl[-F_\pi(\alpha^\dagger,\alpha^*)\bigr]\in[0,1).
\label{eq:style_s6_full_transit_time}
\end{equation}
The total improper value may be
infinite because of the behavior
near $\alpha^*$; in that case
$\theta^*=0$.
The quantity $\theta^*$ coincides with
that in
\eqref{eq:component_full_interval_time}
and with the limit of $\tau_\alpha^*$
as $\alpha\to\alpha^\dagger$ within the
component.

Still in the case of finiteness at arrival,
$F_\pi(\alpha^\dagger,u)$ tends to
zero as $u\to\alpha^\dagger$.
Set $F_\pi(\alpha^\dagger,\alpha^\dagger)=0$ and denote the inverse with respect to $u$, on the side of the component, as follows:
\begin{equation}
F_{\pi,\alpha^\dagger}^{-1}\bigl(F_\pi(\alpha^\dagger,u)\bigr)=u,\qquad u\in(a,b).
\label{eq:style_s6_inverse_arrival_notation}
\end{equation} Its domain
is $[0,F_\pi(\alpha^\dagger,\alpha^*))$;
at zero it takes the value
$\alpha^\dagger$ and at the other endpoint
it tends to $\alpha^*$.
This inverse is used only
when the integral near the arrival
point is finite.

\begin{theorem}[Arrival at a fixed point from a component]
\label{thm:fixed_point_component_arrival}
Among optimal paths with
mean in $[a,b]$ at every positive
time, the only one with terminal state
$\alpha^*$ is
$\varphi(\tau)=\alpha^*\tau$.

If $F_\pi(\alpha^\dagger,u)=+\infty$
for $u\in(a,b)$, the only optimal
path with terminal state
$\alpha^\dagger$ and mean in $[a,b]$
is $\varphi(\tau)=\alpha^\dagger\tau$.

If instead the integral is finite,
the optimal paths with terminal
state $\alpha^\dagger$ and mean in
$[a,b]$ are the straight line just indicated
and, for every $0<t\leq1$, the
paths $\varphi^*(\tau)=\tau\psi^*(\tau)$
defined by
\begin{equation}
\psi^*(\tau)= \alpha^*,\qquad 0\leq\tau\leq\theta^*t,
\label{eq:fixed_point_component_family}
\end{equation}
\begin{equation}
\psi^*(\tau)= F_{\pi,\alpha^\dagger}^{-1}\!\left(\log\frac{t}{\tau}\right),\qquad \theta^*t<\tau<t,
\label{eq:style_global_fixed_point_component_family_2}
\end{equation}
\begin{equation}
\psi^*(\tau)= \alpha^\dagger,\qquad t\leq\tau\leq1.
\label{eq:style_global_fixed_point_component_family_3}
\end{equation}
The parameter $t$ is the first arrival
time of the mean at
$\alpha^\dagger$. The straight line
$\alpha^\dagger\tau$ corresponds
separately to the case $t=0$.

If $\theta^*>0$, all optimal
paths that cross $(a,b)$
before reaching
$\alpha^\dagger$, allowing values
of the mean outside $[a,b]$ at
earlier times, are obtained from
\eqref{eq:fixed_point_component_family}--\eqref{eq:style_global_fixed_point_component_family_3}
by replacing the initial portion with
a rescaled zero-cost
path. More precisely, for
$\varphi'\in Q(\alpha^*)$ with
$I(\varphi')=0$, set
\begin{equation}
\varphi(\tau)= \theta^*t\,\varphi'\!\left(\frac{\tau}{\theta^*t}\right),\qquad 0\leq\tau\leq\theta^*t,
\label{eq:fixed_point_general_component_prefix}
\end{equation}
\begin{equation}
\varphi(\tau)= \varphi^*(\tau),\qquad \theta^*t<\tau\leq1.
\label{eq:style_global_fixed_point_general_component_prefix_2}
\end{equation}
If $\theta^*=0$, for every $0<t\leq1$
formula
\eqref{eq:fixed_point_component_family}--\eqref{eq:style_global_fixed_point_component_family_3}
already gives the unique path that
reaches $\alpha^\dagger$ at
time $t$ by crossing the
component $(a,b)$.
\end{theorem}

\begin{proof}
By Theorem~\ref{thm:zero_entropy_interval}, the optimal
paths considered have zero cost.
Lemma~\ref{lem:zero_cost_paths} makes their mean monotone;
if it takes a value in $(a,b)$, its direction is
from $\alpha^*$ to $\alpha^\dagger$, as fixed in
\eqref{eq:fixed_point_component_endpoints}. Under the constraint
$\psi\in[a,b]$, the terminal state $\alpha^*$ therefore forces
the mean to be constantly $\alpha^*$.

For a nonconstant mean with terminal state
$\alpha^\dagger$, let $t\in(0,1]$ be the first arrival time.
Monotonicity forces $\psi=\alpha^\dagger$ on $[t,1]$.
The portion with mean in $(a,b)$ is an interval $(s,t)$;
separating variables and passing to the limit at arrival
yields, for both orientations,
\begin{equation}
F_\pi(\alpha^\dagger,\psi(\tau))=\log\frac{t}{\tau},
\qquad s<\tau<t.
\label{eq:fixed_point_arrival_identity}
\end{equation}
The right-hand side is finite, so an integral that diverges
at arrival rules out every nonconstant confined mean.

In the case of finiteness at arrival,
\begin{equation}
\psi(\tau)=F_{\pi,\alpha^\dagger}^{-1}\!\Bigl(\log\frac t\tau\Bigr),\qquad s<\tau<t.
\label{eq:style_s6_inverse_arrival_mean}
\end{equation}
The domain of the inverse requires $\tau>\theta^*t$ for every
$\tau\in(s,t)$, hence $s\geq\theta^*t$.
If $s>\theta^*t$, then $s>0$ and, by continuity,
\begin{equation}
\psi(s)=F_{\pi,\alpha^\dagger}^{-1}
\!\left(\log\frac ts\right)\in(a,b),
\label{eq:display_016}
\end{equation}
contradicting the maximality of $(s,t)$. Thus $s=\theta^*t$,
even when the total integral is infinite and $\theta^*=0$.
For $s>0$, the limit of the inverse gives $\psi(s)=\alpha^*$, so
\begin{equation}
s=\theta^*t,
\label{eq:fixed_point_two_matching_times}
\end{equation}
\begin{equation}
F_\pi(\alpha^\dagger,\alpha^*)=\log\frac ts,\qquad s>0.
\label{eq:style_s6_matching_integral_time}
\end{equation}
In this case monotonicity and the constraint $\psi\in[a,b]$ force
$\psi=\alpha^*$ before $s$. Together with the already proved permanence at
$\alpha^\dagger$ after $t$, this yields
\eqref{eq:fixed_point_component_family}--\eqref{eq:style_global_fixed_point_component_family_3}.

Conversely, the inverse formula satisfies
\eqref{eq:zero_cost_average_equation} on the intermediate portion;
the same holds on the constant portions. The mean is continuous at the
joining points and bounded, even when $s=0$. Hence
$\varphi^*(0)=0$ and $\varphi^*(\tau)=\tau\psi^*(\tau)$
define a continuous path with derivative
$\bar\pi(\psi^*)\in(0,K)$ away from the joining points.
Integrating on the individual portions, and down to zero by taking a limit,
yields $\varphi^*\in Q(\alpha^\dagger)$.
Lemma~\ref{lem:zero_cost_paths} now gives zero cost and
$C^1$ regularity of the entire path, including the joining points.

The same separation applies to a path that
crosses $(a,b)$ without being confined to it at earlier times:
by monotonicity, after arrival it remains at $\alpha^\dagger$.
The common portion begins at $s=\theta^*t$.
If $s>0$, the rescaled prefix
$\varphi'(r)=\varphi(sr)/s$, $0\leq r\leq1$,
belongs to $Q(\alpha^*)$, as in
\eqref{eq:component_rescaled_initial_path}.
Conversely, for every $\varphi'\in Q(\alpha^*)$,
\eqref{eq:fixed_point_general_component_prefix}--\eqref{eq:style_global_fixed_point_general_component_prefix_2} is admissible:
the portions agree at the value $s\alpha^*$ and preserve
the increment constraints. The change of variables
used in \eqref{eq:component_rescaling_cost} gives
\begin{equation}
I(\varphi)=\theta^*t\,I(\varphi'),
\qquad \varphi'\in Q(\alpha^*).
\label{eq:fixed_point_prefix_cost}
\end{equation}
The intermediate and terminal portions have zero cost.
Since $\theta^*t>0$, the complete path is optimal
if and only if $I(\varphi')=0$; only then does
Lemma~\ref{lem:zero_cost_paths} ensure a $C^1$ joining.
If $\theta^*=0$, the interior portion covers $(0,t)$ and
there is no prefix left to choose.
\end{proof}

To cover arrivals through infinitely many fixed points,
even without a final adjacent component, we use
extended spatial Lebesgue integrals below: the reciprocal
of zero is $+\infty$, and the integral over an interval
reduced to a point is zero.

\begin{theorem}[General uniqueness criterion at a fixed point]
\label{thm:fixed_point_global_uniqueness}
Let $\alpha\in C$.
There exists an optimal path
in $Q(\alpha)$ other than the straight line
$\varphi(\tau)=\alpha\tau$ if and
only if at least one of the
following two conditions holds: there exists
$u\in[\inf C,\alpha)$ such that
\begin{equation}
\bar\pi(v)-v\geq0,\qquad u\leq v\leq\alpha,
\label{eq:fixed_point_arrival_from_left}
\end{equation}
\begin{equation}
\int_u^\alpha\frac{dv}{\bar\pi(v)-v}<\infty;
\label{eq:style_s6_arrival_left_integrability}
\end{equation}
or there exists $u\in(\alpha,\sup C]$ such that
\begin{equation}
v-\bar\pi(v)\geq0,\qquad\alpha\leq v\leq u,
\label{eq:fixed_point_arrival_from_right}
\end{equation}
\begin{equation}
\int_\alpha^u\frac{dv}{v-\bar\pi(v)}<\infty.
\label{eq:style_s6_arrival_right_integrability}
\end{equation}
The straight line is therefore the unique
optimal path if neither
condition is satisfied.
The point $u$ need not
belong to $C$; the constraint
$u\in[\inf C,\sup C]$ is part
of the conditions.
\end{theorem}

\begin{proof}
Let $\varphi\in Q(\alpha)$ be optimal
and different from the straight line.
Its mean is monotone and
nonconstant by
Lemma~\ref{lem:zero_cost_paths}.
First suppose that it is
nondecreasing. Then there exists
$\tau_0\in(0,1)$ with
$u:=\psi(\tau_0)<\alpha$.
The confinement already proved gives
$u\in[\inf C,\alpha)$.
Every value in $[u,\alpha]$
is taken by the mean, and
$\bar\pi(\psi(\tau))-\psi(\tau)
=\tau\partial_\tau\psi(\tau)\geq0$;
thus the sign required in
\eqref{eq:fixed_point_arrival_from_left}
holds throughout the interval.
For $\varepsilon>0$, the chain
rule gives
\begin{equation}
\int_u^\alpha\frac{dv}{\bar\pi(v)-v+\varepsilon}
=\int_{\tau_0}^1\frac{\partial_\tau\psi(\tau)}
{\bar\pi(\psi(\tau))-\psi(\tau)+\varepsilon}\,d\tau,
\label{eq:fixed_point_regularized_arrival_bound}
\end{equation}
\begin{equation}
\int_u^\alpha\frac{dv}{\bar\pi(v)-v+\varepsilon}
=\int_{\tau_0}^1\frac{\bar\pi(\psi(\tau))-\psi(\tau)}
{\bar\pi(\psi(\tau))-\psi(\tau)+\varepsilon}\,\frac{d\tau}{\tau},
\label{eq:style_s6_regularized_arrival_equation}
\end{equation}
\begin{equation}
\int_u^\alpha\frac{dv}{\bar\pi(v)-v+\varepsilon}\leq\log\frac1{\tau_0}.
\label{eq:style_s6_regularized_arrival_bound}
\end{equation}
Monotone convergence as
$\varepsilon\downarrow0$ proves
finiteness of the extended integral.
If the mean is nonincreasing,
choose $u=\psi(\tau_0)>\alpha$,
and the same calculation, with
$v-\bar\pi(v)$, gives
\eqref{eq:fixed_point_arrival_from_right}--\eqref{eq:style_s6_arrival_right_integrability}.

For sufficiency, assume
\eqref{eq:fixed_point_arrival_from_left}--\eqref{eq:style_s6_arrival_left_integrability}
and set
\begin{equation}
\tau_0:=\exp\left[-\int_u^\alpha
\frac{dv}{\bar\pi(v)-v}\right]\in(0,1).
\label{eq:fixed_point_general_initial_time}
\end{equation}
Finiteness of the integral implies
that the zero set of the
denominator in $[u,\alpha]$ has
measure zero.
The function obtained by integrating
the reciprocal is absolutely
continuous and strictly monotone.
Indeed, for $u\leq v_1<v_2\leq\alpha$
we have
\begin{equation}
v_2-v_1\leq
\max_{v\in[v_1,v_2]}|\bar\pi(v)-v|
\int_{v_1}^{v_2}\frac{dv}{|\bar\pi(v)-v|}.
\label{eq:fixed_point_inverse_increment_bound}
\end{equation}
The maximum is positive on every
nondegenerate interval, since
the integral is finite.
We therefore define a mean on
$[\tau_0,1]$ by inverting the relation
\begin{equation}
\int_{\psi(\tau)}^\alpha
\frac{dv}{\bar\pi(v)-v}
=\log\frac1\tau,
\qquad \tau_0\leq\tau\leq1.
\label{eq:fixed_point_general_arrival_path}
\end{equation}
We have $\psi(\tau_0)=u$ and
$\psi(1)=\alpha$.

Equation \eqref{eq:fixed_point_inverse_increment_bound} makes
the inverse Lipschitz with respect to $\log(1/\tau)$.
Where $\bar\pi(\psi)-\psi>0$, differentiating the integral
identity gives \eqref{eq:zero_cost_average_equation}.
At a value in $C$, instead, the difference quotient of the inverse
is bounded by the local maximum in
\eqref{eq:fixed_point_inverse_increment_bound}, which tends
to zero by continuity. The derivative there is therefore zero:
the equation holds everywhere and the mean is $C^1$ on
$[\tau_0,1]$, with one-sided derivatives at the endpoints.
This argument also covers nonisolated zeros.

Since $u\in[\inf C,\sup C]$,
Theorem~\ref{thm:zero_entropy_interval}
provides a path
$\varphi'\in Q(u)$ with zero cost.
Define
\begin{equation}
\varphi(\tau)= \tau_0\varphi'\!\left(\frac{\tau}{\tau_0}\right),\qquad 0\leq\tau\leq\tau_0,
\label{eq:fixed_point_general_arrival_completion}
\end{equation}
\begin{equation}
\varphi(\tau)= \tau\psi(\tau),\qquad \tau_0<\tau\leq1.
\label{eq:style_global_fixed_point_general_arrival_completion_2}
\end{equation}
The portions agree at the value $\tau_0u$, and each
satisfies the increment constraints, so
$\varphi\in Q(\alpha)$. The terminal portion satisfies
\eqref{eq:zero_cost_equation} and therefore has a zero
cost integrand. The change of variables in
\eqref{eq:component_rescaling_cost}, valid for every
admissible prefix, gives $I(\varphi)=\tau_0I(\varphi')=0$.
Lemma~\ref{lem:zero_cost_paths} therefore gives $C^1$
regularity at the joining point as well. The mean moves from $u<\alpha$
to $\alpha$, so the path differs from the straight line.

In the case
\eqref{eq:fixed_point_arrival_from_right}--\eqref{eq:style_s6_arrival_right_integrability}
we instead use
\begin{equation}
\int_\alpha^{\psi(\tau)}
\frac{dv}{v-\bar\pi(v)}
=\log\frac1\tau,
\label{eq:fixed_point_general_right_arrival_path}
\end{equation}
with $\tau_0$ defined by the
corresponding integral from
$\alpha$ to $u$.
The mean is nonincreasing; the same local estimate
\eqref{eq:fixed_point_inverse_increment_bound}, with
$|\bar\pi(v)-v|$, and the rescaled prefix prove
admissibility and zero cost.
Lemma~\ref{lem:zero_cost_paths} then gives regularity.
\end{proof}

The constraint $u\in[\inf C,\sup C]$ guarantees the zero-cost
prefix down to zero; the arrival portion at positive
times alone is not enough for the variational problem.

In particular, if $\bar\pi$ is
locally Lipschitz near
$\alpha\in C$, the reciprocal of
$|\bar\pi(v)-v|$ has a divergent
integral on both sides:
the estimate
$|\bar\pi(v)-v|\leq c|v-\alpha|$
for a constant $c>0$ gives the
comparison with $1/|v-\alpha|$.
The straight line is then the unique
optimal path in $Q(\alpha)$.
The same conclusion holds when
$\alpha$ is in the interior of an interval
of fixed points, because every
nondegenerate interval arriving
there contains a portion on
which the denominator is zero.
Theorem~\ref{thm:fixed_point_global_uniqueness}
provides the exact criterion in
the other cases as well.

The following identity completes the
description of the paths without
requiring the fixed points crossed
to be isolated.
The indicator retains the notation
$\mathbb{I}(E)$ already used for
events.

\begin{lemma}[Integral identity for optimal paths]
\label{lem:fixed_point_integral_characterization}
Let $\alpha\in C$ and let
$\psi:(0,1]\to[0,K]$ be continuous,
with $\psi(1)=\alpha$.
Set $\varphi(0)=0$ and
$\varphi(\tau)=\tau\psi(\tau)$
for $\tau>0$.
Then $\varphi\in Q(\alpha)$
and $I(\varphi)=0$ if and only
if the following conditions hold.

The mean is nondecreasing and
$\bar\pi(\psi(\tau))-\psi(\tau)\geq0$
for every $\tau>0$, or it is
nonincreasing and
$\bar\pi(\psi(\tau))-\psi(\tau)\leq0$
for every $\tau>0$.
Moreover, for every
$0<\tau_1<\tau_2\leq1$,
\begin{equation}
\log\frac{\tau_2}{\tau_1}
=\int_{\min\{\psi(\tau_1),\psi(\tau_2)\}}
^{\max\{\psi(\tau_1),\psi(\tau_2)\}}\frac{du}{|\bar\pi(u)-u|}
{}+\int_{\tau_1}^{\tau_2}\mathbb I(\psi(s)\in C)\,\frac{ds}{s}.
\label{eq:fixed_point_complete_time_identity}
\end{equation}
\end{lemma}

\begin{proof}
If $I(\varphi)=0$,
Lemma~\ref{lem:zero_cost_paths}
gives monotonicity, the required sign,
and the equation for the mean.
Regularizing the spatial
denominator with
$|\bar\pi(u)-u|+\varepsilon$,
the chain rule gives the same
calculation as in
\eqref{eq:fixed_point_regularized_arrival_bound}--\eqref{eq:style_s6_regularized_arrival_bound},
for both directions of monotonicity.
Monotone convergence for
the spatial integral and
dominated convergence for the time integral
yield
\begin{equation}
\int_{\min\{\psi(\tau_1),\psi(\tau_2)\}}
^{\max\{\psi(\tau_1),\psi(\tau_2)\}}
\frac{du}{|\bar\pi(u)-u|}
=\int_{\tau_1}^{\tau_2}
\mathbb{I}(\psi(s)\notin C)\,\frac{ds}{s}.
\label{eq:fixed_point_moving_time_identity}
\end{equation}
Adding the integral over the times
at which $\psi(s)\in C$ gives
\eqref{eq:fixed_point_complete_time_identity}.
If the mean is constant between the two
times, its value belongs
to $C$: the spatial integral is
zero and the identity remains valid.

Conversely, assume the conditions
of the statement. For two fixed
times, let $u$ and $v$ be the minimum
and maximum of the corresponding
values of the mean.
The identity and the estimate
\eqref{eq:fixed_point_inverse_increment_bound},
applied with the absolute value
of the denominator, give
\begin{equation}
v-u\leq
\max_{z\in[u,v]}|\bar\pi(z)-z|
\log\frac{\tau_2}{\tau_1}.
\label{eq:fixed_point_mean_increment_bound}
\end{equation}
For $u=v$ the inequality is
immediate. In the other cases
the spatial integral is finite by
\eqref{eq:fixed_point_complete_time_identity},
so the estimate applies.
In particular, $\psi$ is Lipschitz
on every compact time interval
contained in $(0,1]$.

Fix $\tau_0>0$.
If $\psi(\tau_0)\in C$, divide
\eqref{eq:fixed_point_mean_increment_bound}
by the distance between $\tau_0$ and
a nearby time.
The continuity of $\psi$ and of
$\bar\pi$ makes the maximum
on the right-hand side tend to zero, while
the difference quotient of the
logarithm tends to $1/\tau_0$.
Hence $\partial_\tau\psi(\tau_0)=0$.

If $\psi(\tau_0)\notin C$,
there exists a time neighborhood in which
the mean does not meet $C$.
The last integral in
\eqref{eq:fixed_point_complete_time_identity}
is zero on this neighborhood.
Compatibility between the sign and
the direction of monotonicity transforms
the identity into the oriented relation
\begin{equation}
\int_{\psi(\tau_1)}^{\psi(\tau_2)}
\frac{du}{\bar\pi(u)-u}
=\log\frac{\tau_2}{\tau_1}.
\label{eq:fixed_point_local_oriented_identity}
\end{equation}
The antiderivative has a continuous,
nonzero derivative. Inverting it
yields
$\tau_0\partial_\tau\psi(\tau_0)
=\bar\pi(\psi(\tau_0))-\psi(\tau_0)$.
The equation therefore holds at every
positive time, and its right-hand
side shows that $\psi$ is $C^1$
on $(0,1]$.

It follows that
$\partial_\tau\varphi(\tau)
=\bar\pi(\psi(\tau))\in(0,K)$.
Boundedness of the mean gives
continuity of $\varphi$ at zero;
integrating the derivative yields
the constraints of $Q$.
Since $\varphi(1)=\alpha$,
we have $\varphi\in Q(\alpha)$,
and Lemma~\ref{lem:zero_cost_paths}
concludes that $I(\varphi)=0$.
\end{proof}

In \eqref{eq:fixed_point_complete_time_identity},
the spatial integral measures the logarithmic transit
time and the time integral measures the time spent at
fixed points, even if the times spent there do not form
a union of intervals or infinitely many fixed points are
crossed. On a single component, one recovers the
separation of variables and the earlier constant portions.

Finally, for an optimal path
other than the straight line with terminal state
$\alpha\in C$, monotonicity
ensures a first arrival time
$t\in(0,1]$ at $\alpha$.
For $t\leq\tau\leq1$ the mean is
constant and $\varphi(\tau)=\alpha\tau$.
The counts, however, satisfy
\begin{equation}
\gamma_k(\tau)
=\gamma_k(t)+(\tau-t)\pi_k(\alpha),
\qquad k\in\Omega,\qquad t\leq\tau\leq1.
\label{eq:fixed_point_counts_after_arrival}
\end{equation}
The formula follows from
\eqref{eq:zero_cost_counts}.
For $K>1$ the value $\gamma_k(t)$
retains information from the preceding
portion, as in
\eqref{eq:style_global_component_optimal_counts_2};
it is not in general individually
determined by the constraints on
the sum and the mean.
The three subsections thus complete
the characterization of
zero-cost paths in the
class of strictly positive urn
functions considered here.

\section{Cumulant generating function}
\label{sec:cumulants}

The terminal probabilities studied in
Section~\ref{sec:terminal_entropy}
determine the exponential behavior
of the moments of the cumulative sum $M_N$.
The cumulant generating function
allows this behavior to be expressed
through a variational problem and
connected to the recursion of the process.

\subsection{Existence and variational representation}
\label{subsec:cumulant_variational}

For $N\geq1$ and $\lambda\in\mathbb{R}$,
we define the cumulant generating
function, normalized by the number
of steps, by
\begin{equation}
\zeta_N(\lambda):=\frac1N\log\mathbb E\bigl[e^{\lambda M_N}\bigr]
=\frac1N\log\sum_{r=0}^{KN}e^{\lambda r}\mathbb P(M_N=r).
\label{eq:cumulant_finite}
\end{equation}
The parameter $\lambda$ multiplies
$M_N=N\psi_N$, so the exponent
has the same scale $N$ as the LDP.
The definition corresponds to
Ref.~\cite{Franchini_URNS},
Eqs.~(2.29) and~(3.96), with the domain
of the mean extended from $[0,1]$ to $[0,K]$.

Since $0\leq M_N\leq KN$,
the function in
\eqref{eq:cumulant_finite} is finite
for every $\lambda$ and satisfies
\begin{equation}
\min\{0,K\lambda\}\leq\zeta_N(\lambda)\leq\max\{0,K\lambda\},
\label{eq:cumulant_finite_bounds}
\end{equation}
\begin{equation}
\zeta_N(0)=0.
\label{eq:style_s7_001}
\end{equation}
Moreover, for $\lambda\leq\mu$,
the inequality
$1\leq e^{(\mu-\lambda)M_N}
\leq e^{NK(\mu-\lambda)}$
gives
\begin{equation}
0\leq\zeta_N(\mu)-\zeta_N(\lambda)
\leq K(\mu-\lambda).
\label{eq:cumulant_finite_lipschitz}
\end{equation}
These estimates hold for all the
laws considered in
Subsection~\ref{subsec:initial_and_inhomogeneous},
since they depend only on the
possible values of the increments.

\begin{theorem}[Limit of the cumulant generating function]
\label{thm:cumulant_variational}
Under the assumptions of
Theorem~\ref{thm:path_ldp},
for every $\lambda\in\mathbb{R}$
the following finite limit exists:
\begin{equation}
\zeta(\lambda)
:=\lim_{N\to\infty}\zeta_N(\lambda)
=\max_{\alpha\in[0,K]}
\{\lambda\alpha+\phi(\alpha)\}.
\label{eq:cumulant_variational_terminal}
\end{equation}
The equivalent representations
on scalar paths and on counts
are, respectively,
\begin{equation}
\zeta(\lambda)
=\max_{\varphi\in Q}
\{\lambda\varphi(1)-I(\varphi)\}
\label{eq:cumulant_variational_paths}
\end{equation}
and
\begin{equation}
\zeta(\lambda)=\max_{\gamma\ \mathrm{admissible}}\biggl\{
\lambda\sum_{k\in\Omega}k\gamma_k(1)-I_{\mathrm c}(\gamma)\biggr\}.
\label{eq:cumulant_variational_counts}
\end{equation}
All maxima are attained.
In the count formula, $I_{\mathrm c}$ is defined by
\eqref{eq:count_rate}, with the convention $0\log0=0$.

The function $\zeta$ is convex,
nondecreasing, and $K$-Lipschitz,
and satisfies $\zeta(0)=0$.
The convergence in
\eqref{eq:cumulant_variational_terminal}
is uniform on every compact
interval:
\begin{equation}
\lim_{N\to\infty}
\sup_{|\lambda|\leq R}
|\zeta_N(\lambda)-\zeta(\lambda)|=0,
\qquad R>0.
\label{eq:cumulant_compact_convergence}
\end{equation}
The same conclusions, with the
same function $\zeta$, hold
for the initial conditions of
Theorem~\ref{thm:finite_initial_history}
and for the nonhomogeneous urn
functions of
Theorem~\ref{thm:inhomogeneous_urn},
using in each case the expectation
under the law
considered.
\end{theorem}

\begin{proof}
Theorem~\ref{thm:terminal_entropy} provides the LDP on
$[0,K]$ with good rate function $-\phi$.
For every fixed $\lambda$, $\alpha\mapsto\lambda\alpha$
is continuous and bounded on this compact set; the tail
condition of the Laplace principle is therefore automatic.
Theorem~4.3.1 of Ref.~\cite{Dembo_Zeitouni} gives
\eqref{eq:cumulant_variational_terminal}, with the maximum
attained by continuity of $\phi$.
This is the step recalled in Ref.~\cite{Franchini_URNS}
after Eq.~(2.29).

For every $\alpha$, definition
\eqref{eq:terminal_entropy_variational}
gives
\begin{equation}
\lambda\alpha+\phi(\alpha)
=\max_{\varphi\in Q(\alpha)}
\{\lambda\varphi(1)-I(\varphi)\}.
\label{eq:cumulant_terminal_to_path}
\end{equation}
Every $\varphi\in Q$ belongs
to $Q(\varphi(1))$; taking
the maximum over $\alpha$ gives
\eqref{eq:cumulant_variational_paths}.
A terminal maximizer
and a path attaining
the minimum of $I$ with that
terminal value provide
a maximizing path.
Similarly,
\eqref{eq:terminal_entropy_counts}
gives
\eqref{eq:cumulant_variational_counts}:
as $\alpha$ varies, the constraint
on the weighted sum of the terminal
counts ranges over all admissible
vectors. This maximum
is also attained, because
both the terminal maximum
and the maximum over counts with
fixed terminal value are attained.

The first representation
expresses $\zeta$ as a maximum
of affine functions of
$\lambda$, so $\zeta$
is convex.
Passing to the limit in
\eqref{eq:cumulant_finite_bounds}
and
\eqref{eq:cumulant_finite_lipschitz}
yields $\zeta(0)=0$,
monotonicity, and the Lipschitz
constant $K$.
To verify
\eqref{eq:cumulant_compact_convergence},
fix $R>0$ and choose
a finite grid
$\lambda_1,\ldots,\lambda_m$
in $[-R,R]$ such that every
point of the interval is at
distance at most $\delta$ from a point
of the grid. The common
Lipschitz constant gives
\begin{equation}
\sup_{|\lambda|\leq R}
|\zeta_N(\lambda)-\zeta(\lambda)|
\leq
\max_{1\leq j\leq m}
|\zeta_N(\lambda_j)-\zeta(\lambda_j)|
+2K\delta.
\label{eq:cumulant_finite_grid}
\end{equation}
First let $N$
tend to infinity, using
convergence at the grid
points, and then let $\delta$ tend to zero.

The laws in Theorems~\ref{thm:finite_initial_history}
and~\ref{thm:inhomogeneous_urn} have the same terminal LDP
with rate $-\phi$, by
Subsection~\ref{subsec:terminal_variational}; the proof therefore gives
the same limiting function.
\end{proof}

The finite initializations and nonhomogeneous urns of
Subsection~\ref{subsec:initial_and_inhomogeneous} thus have
the same functions $I$, $\phi$, and $\zeta$. The properties of these
limiting functions and the corresponding variational problems established
below hold for all these laws. Injectivity and affinity
refer to the limiting vector $\pi$; in the inverse constructions,
this is the vector to be prescribed. The finite-time identities
and restrictions to the respective reachable supports remain distinct.

The terminal representation
also allows us to identify
the one-sided derivatives of
$\zeta$ without introducing
further assumptions on the urn
functions. We denote by
$\zeta'_-(\lambda)$ and
$\zeta'_+(\lambda)$ the
left and right derivatives;
they exist and are finite
by convexity and the
Lipschitz property.

\begin{lemma}[Maximizers and one-sided derivatives]
\label{lem:cumulant_one_sided_derivatives}
For every $\lambda\in\mathbb{R}$,
we have
\begin{equation}
\zeta'_-(\lambda)=\min\bigl\{\alpha\in[0,K]:
\lambda\alpha+\phi(\alpha)=\zeta(\lambda)\bigr\},
\label{eq:cumulant_one_sided_maximizers}
\end{equation}
\begin{equation}
\zeta'_+(\lambda)=\max\bigl\{\alpha\in[0,K]:
\lambda\alpha+\phi(\alpha)=\zeta(\lambda)\bigr\}.
\label{eq:style_s7_002}
\end{equation}
In particular, $\zeta$ is
differentiable at $\lambda$
if and only if the maximum
in
\eqref{eq:cumulant_variational_terminal}
is attained at a unique
value of $\alpha$; that
value is $\zeta'(\lambda)$.
\end{lemma}

\begin{proof}
Let $A_\lambda$ be the nonempty compact set of maximizers
in \eqref{eq:cumulant_variational_terminal}, and set
$a_-=\min A_\lambda$, $a_+=\max A_\lambda$.
For $h>0$ choose $a_h^+\in A_{\lambda+h}$ and
$a_h^-\in A_{\lambda-h}$. Evaluating each maximum at a
maximizer for the other parameter yields
\begin{equation}
a_+\leq\frac{\zeta(\lambda+h)-\zeta(\lambda)}h\leq a_h^+,
\label{eq:display_017}
\end{equation}
\begin{equation}
a_h^-\leq\frac{\zeta(\lambda)-\zeta(\lambda-h)}h\leq a_-.
\label{eq:style_s7_003}
\end{equation}
Every subsequential limit of $a_h^\pm$ as
$h\downarrow0$ belongs to $A_\lambda$: indeed,
\begin{equation}
\lambda a_h^\pm+\phi(a_h^\pm)
=\zeta(\lambda\pm h)\mp h a_h^\pm\longrightarrow\zeta(\lambda),
\label{eq:display_018}
\end{equation}
and $\phi$ is continuous. By compactness,
$\limsup_{h\downarrow0}a_h^+\leq a_+$ and
$\liminf_{h\downarrow0}a_h^-\geq a_-$.
The two secant estimates therefore give
$\zeta'_+(\lambda)=a_+$ and $\zeta'_-(\lambda)=a_-$.
The derivatives coincide exactly when $A_\lambda$ contains
a single point.
\end{proof}

At $\lambda=0$ the maximizers
are the points where
$\phi$ vanishes.
Theorem~\ref{thm:zero_entropy_interval}
identifies them with
$[\inf C,\sup C]$.
The lemma therefore gives
\begin{equation}
\zeta'_-(0)=\inf C,
\label{eq:cumulant_derivatives_at_zero}
\end{equation}
\begin{equation}
\zeta'_+(0)=\sup C.
\label{eq:style_s7_004}
\end{equation}
Consequently, $\zeta$ is
differentiable at zero
if and only if $C$
contains a single point.

The location of the maximizers
is also constrained by the sign
of $\lambda$.
If $\lambda<0$ and
$\alpha>\inf C$, then
$\lambda\alpha+\phi(\alpha)
\leq\lambda\alpha
<\lambda\inf C$,
while the value
$\lambda\inf C$ is attained
at $\alpha=\inf C$.
Hence every maximizer
belongs to $[0,\inf C]$.
For $\lambda>0$, the same
comparison with $\alpha=\sup C$
shows that all
maximizers belong
to $[\sup C,K]$.
In particular,
\begin{equation}
0\leq\zeta'_-(\lambda)\leq\zeta'_+(\lambda)\leq\inf C,
\qquad\lambda<0,
\label{eq:cumulant_derivative_branches}
\end{equation}
\begin{equation}
\sup C\leq\zeta'_-(\lambda)\leq\zeta'_+(\lambda)\leq K,
\qquad\lambda>0.
\label{eq:style_s7_005}
\end{equation}

Finally,
\eqref{eq:cumulant_variational_terminal}
gives, for every
$\alpha\in[0,K]$, the
inequality
\begin{equation}
\phi(\alpha)
\leq\inf_{\lambda\in\mathbb{R}}
\{\zeta(\lambda)-\lambda\alpha\}.
\label{eq:cumulant_entropy_dual_bound}
\end{equation}
Equality requires identifying the states that attain the
maximum in the variational formula. In
Subsection~\ref{subsec:entropy_shape}, dual contact and the range
of slopes will provide this identification; the concavity of $\phi$
will follow as a consequence.

\subsection{Recursion and limiting equation}
\label{subsec:cumulant_recursion}

The relation between the
cumulant generating function and
the urn functions is obtained
by combining the recursion of the
process with an exponential
change of measure.
We follow the procedure
of Ref.~\cite{Franchini_URNS},
Subsections~3.3 and~3.3.1, Eq.~(3.97),
making explicit the concentration
needed to pass to the limit.

For $N\geq1$, $\lambda\in\mathbb{R}$,
and $E\subseteq\Omega^S$,
define
\begin{equation}
\mathbb{P}_{N,\lambda}(E)
:=\mathbb{E}\left[
\mathbb{I}(\sigma\in E)
e^{\lambda M_N-N\zeta_N(\lambda)}
\right].
\label{eq:cumulant_changed_measure}
\end{equation}
Normalization follows
from the definition of
$\zeta_N$.
We denote by
$\mathbb{E}_{N,\lambda}$
the corresponding expectation.
The first index specifies
the number of steps;
in particular,
$\mathbb{P}_{N,0}=\mathbb{P}$,
whereas $\mathbb{P}_0$ retains
its meaning as the uniform
reference law.
Differentiating the finite sum
in \eqref{eq:cumulant_finite}
also yields
\begin{equation}
\zeta_N'(\lambda)
=\mathbb{E}_{N,\lambda}[\psi_N]
\in[0,K].
\label{eq:cumulant_finite_derivative}
\end{equation}

The following lemma identifies
the terminal law after
the change of measure
and justifies the limit
of the expectations of the
urn functions.

\begin{lemma}[Terminal law after the change of measure]
\label{lem:cumulant_changed_ldp}
Under the assumptions of
Theorem~\ref{thm:path_ldp},
for every fixed $\lambda\in\mathbb{R}$,
the laws of $\psi_N$
under $\mathbb{P}_{N,\lambda}$
satisfy an LDP on $[0,K]$
with speed $N$ and good
rate function
\begin{equation}
\alpha\longmapsto
\zeta(\lambda)-\lambda\alpha-\phi(\alpha).
\label{eq:cumulant_changed_rate}
\end{equation}
In particular, for every
Borel set $B\subseteq[0,K]$
fixed as $N$ varies,
\begin{equation}
\liminf_{N\to\infty}\frac{1}{N}
\log\mathbb{P}_{N,\lambda}(\psi_N\in B)
\geq\sup_{\alpha\in B^\circ}
\{\lambda\alpha+\phi(\alpha)-\zeta(\lambda)\}
\label{eq:cumulant_changed_ldp_lower}
\end{equation}
and
\begin{equation}
\limsup_{N\to\infty}\frac{1}{N}
\log\mathbb{P}_{N,\lambda}(\psi_N\in B)
\leq\sup_{\alpha\in\overline B}
\{\lambda\alpha+\phi(\alpha)-\zeta(\lambda)\}.
\label{eq:cumulant_changed_ldp_upper}
\end{equation}
The interior and closure
are relative to $[0,K]$;
the convention
$\sup\varnothing=-\infty$ is retained.

If $\zeta$ is differentiable
at $\lambda$, for every
$\varepsilon>0$ we have
\begin{equation}
\limsup_{N\to\infty}\frac{1}{N}
\log\mathbb{P}_{N,\lambda}
\left(|\psi_N-\zeta'(\lambda)|\geq\varepsilon\right)
<0.
\label{eq:cumulant_changed_concentration}
\end{equation}
Consequently, for every
$f\in C([0,K])$,
\begin{equation}
\lim_{N\to\infty}
\mathbb{E}_{N,\lambda}[f(\psi_N)]
=f\bigl(\zeta'(\lambda)\bigr).
\label{eq:cumulant_changed_expectation}
\end{equation}
The conclusions also hold
for the laws in
Theorems~\ref{thm:finite_initial_history}
and \ref{thm:inhomogeneous_urn},
with the change of measure
defined relative to the
law considered.
\end{lemma}

\begin{proof}
The terminal LDP in Theorem~\ref{thm:terminal_entropy}
has good rate function $-\phi$ on the compact set $[0,K]$.
The change \eqref{eq:cumulant_changed_measure} uses
the continuous bounded potential $\alpha\mapsto\lambda\alpha$;
the logarithm of its normalization, divided by $N$,
converges to $\zeta(\lambda)$ by
Theorem~\ref{thm:cumulant_variational}.
The LDP under an exponential change of measure,
Ref.~\cite{Dembo_Zeitouni}, Exercise~4.3.11,
therefore gives \eqref{eq:cumulant_changed_rate} and the
stated estimates. The cost is continuous and nonnegative,
with minimum zero; it is good by compactness.

Let $U$ be a relatively open set containing all the
maximizers of \eqref{eq:cumulant_variational_terminal}.
On the compact complement, if nonempty, the cost has
a strictly positive minimum. The upper bound gives
\begin{equation}
\limsup_{N\to\infty}\frac1N
\log\mathbb{P}_{N,\lambda}(\psi_N\notin U)
\leq\sup_{\alpha\in[0,K]\setminus U}
\{\lambda\alpha+\phi(\alpha)-\zeta(\lambda)\}<0.
\label{eq:cumulant_changed_maximizer_concentration}
\end{equation}
The convention $\sup\varnothing=-\infty$ covers
the empty complement. At points of differentiability,
Lemma~\ref{lem:cumulant_one_sided_derivatives}
identifies the unique maximizer with $\zeta'(\lambda)$;
choosing the window of radius $\varepsilon$ yields
\eqref{eq:cumulant_changed_concentration}.

For $f\in C([0,K])$, split the expectation over the
window and its complement: the error on the former is
at most the oscillation of $f$ within distance $\varepsilon$
of $\zeta'(\lambda)$, while on the latter it is at most
$2\|f\|_\infty$ times the probability of the complement.
Letting first $N\to\infty$ and then $\varepsilon\downarrow0$
yields \eqref{eq:cumulant_changed_expectation}.
In particular, \eqref{eq:cumulant_finite_derivative}
gives $\zeta_N'(\lambda)\to\zeta'(\lambda)$ at points
of differentiability of $\zeta$.
The other laws considered have the same terminal LDP
and the same $\zeta$, so they satisfy all the assumptions
just verified.
\end{proof}

We turn to the recursion.
The transition probability
\eqref{eq:transition_probability}
gives exactly
\begin{equation}
\mathbb{E}\left[
e^{\lambda M_{N+1}}
\mid\sigma_1,\ldots,\sigma_N
\right]
=e^{\lambda M_N}
\sum_{k\in\Omega}\pi_k(\psi_N)e^{k\lambda}.
\label{eq:cumulant_conditional_recursion}
\end{equation}
Taking expectations
and dividing by
$\mathbb{E}[e^{\lambda M_N}]$,
we obtain
\begin{equation}
e^{(N+1)\zeta_{N+1}(\lambda)-N\zeta_N(\lambda)}
=\mathbb{E}_{N,\lambda}
\left[\sum_{k\in\Omega}
\pi_k(\psi_N)e^{k\lambda}\right].
\label{eq:cumulant_exact_recursion}
\end{equation}
Equivalently,
\begin{equation}
(N+1)\bigl\{\zeta_{N+1}(\lambda)-\zeta_N(\lambda)\bigr\}
=-\zeta_N(\lambda) 
+\log\mathbb E_{N,\lambda}\Bigl[\sum_{k\in\Omega}\pi_k(\psi_N)e^{k\lambda}\Bigr].
\label{eq:cumulant_exact_difference}
\end{equation}
These are the forms of the
recursion corresponding
to Ref.~\cite{Franchini_URNS},
Eqs.~(3.99)--(3.100).
The sum inside the
expectation is a
convex combination of
$e^{k\lambda}$.
The logarithm of its
expectation and
$\zeta_N(\lambda)$
therefore both belong
to the interval
$[\min\{0,K\lambda\},
\max\{0,K\lambda\}]$.
This gives the estimate
\begin{equation}
\left|(N+1)
\{\zeta_{N+1}(\lambda)-\zeta_N(\lambda)\}
\right|\leq K|\lambda|.
\label{eq:cumulant_difference_bound}
\end{equation}

\begin{theorem}[Limiting equation and strict convexity]
\label{thm:cumulant_limit_equation}
Under the assumptions of
Theorem~\ref{thm:path_ldp},
at every point $\lambda$
where $\zeta$ is
differentiable, we have
\begin{equation}
\zeta(\lambda)
=\log\left[
\sum_{k\in\Omega}
\pi_k\bigl(\zeta'(\lambda)\bigr)e^{k\lambda}
\right],
\label{eq:cumulant_implicit_equation}
\end{equation}
and moreover
\begin{equation}
\lim_{N\to\infty}
(N+1)\{\zeta_{N+1}(\lambda)-\zeta_N(\lambda)\}
=0.
\label{eq:cumulant_difference_limit}
\end{equation}
For every $\lambda\in\mathbb{R}$,
the two one-sided derivatives
separately satisfy
the same equation:
\begin{equation}
e^{\zeta(\lambda)}=\sum_{k\in\Omega}\pi_k\bigl(\zeta'_-(\lambda)\bigr)e^{k\lambda},
\label{eq:cumulant_one_sided_equation}
\end{equation}
\begin{equation}
e^{\zeta(\lambda)}=\sum_{k\in\Omega}\pi_k\bigl(\zeta'_+(\lambda)\bigr)e^{k\lambda}.
\label{eq:style_s7_006}
\end{equation}
The function $\zeta$ is
strictly convex.
In particular,
\begin{equation}
0<\zeta'_-(\lambda)
\leq\zeta'_+(\lambda)<K,
\qquad \lambda\in\mathbb{R},
\label{eq:cumulant_interior_derivatives}
\end{equation}
and the localizations
\eqref{eq:cumulant_derivative_branches}--\eqref{eq:style_s7_005}
strengthen to
\begin{equation}
\zeta'_+(\lambda)<\inf C,\qquad\lambda<0,
\label{eq:cumulant_strict_derivative_branches}
\end{equation}
\begin{equation}
\sup C<\zeta'_-(\lambda),\qquad\lambda>0.
\label{eq:style_s7_007}
\end{equation}
The conclusions also hold
for the initial conditions
and nonhomogeneous urn
functions of
Subsection~\ref{subsec:initial_and_inhomogeneous},
with the limiting urn
functions in equations
\eqref{eq:cumulant_implicit_equation}
and
\eqref{eq:cumulant_one_sided_equation}--\eqref{eq:style_s7_006}.
\end{theorem}

\begin{proof}
Fix a point of
differentiability $\lambda$.
The function
\begin{equation}
\alpha\longmapsto\sum_{k\in\Omega}\pi_k(\alpha)e^{k\lambda}
\label{eq:style_s7_049}
\end{equation}
is continuous on $[0,K]$.
Lemma~\ref{lem:cumulant_changed_ldp}
therefore gives
\begin{equation}
\lim_{N\to\infty}
\mathbb{E}_{N,\lambda}
\left[\sum_{k\in\Omega}
\pi_k(\psi_N)e^{k\lambda}\right]
=\sum_{k\in\Omega}
\pi_k\bigl(\zeta'(\lambda)\bigr)e^{k\lambda}>0.
\label{eq:cumulant_recursion_factor_limit}
\end{equation}
To identify the
limit, sum the
logarithmic form of
\eqref{eq:cumulant_exact_recursion}
from time $1$
to time $N-1$:
\begin{equation}
\zeta_N(\lambda)
=\frac{\zeta_1(\lambda)}{N}
+\frac{1}{N}\sum_{n=1}^{N-1}
\log\mathbb{E}_{n,\lambda}
\left[\sum_{k\in\Omega}
\pi_k(\psi_n)e^{k\lambda}\right].
\label{eq:cumulant_recursion_telescope}
\end{equation}
Each term of the sum,
as $n\to\infty$,
converges to the logarithm
of the right-hand side of
\eqref{eq:cumulant_recursion_factor_limit}.
The Cesàro mean has
the same limit;
Theorem~\ref{thm:cumulant_variational}
identifies the limit
of the left-hand side
with $\zeta(\lambda)$.
This yields
\eqref{eq:cumulant_implicit_equation}.
Substituting this
identity into
\eqref{eq:cumulant_exact_difference}
and using
\eqref{eq:cumulant_recursion_factor_limit}
also yields
\eqref{eq:cumulant_difference_limit}.
This step corresponds to Lemma~20 of
Ref.~\cite{Franchini_URNS}, Eqs.~(3.102)--(3.103):
the concentration in Lemma~\ref{lem:cumulant_changed_ldp}
justifies the limit of expectations, and telescoping
identifies the equation before deducing
\eqref{eq:cumulant_difference_limit}.

To extend the identity
to one-sided derivatives,
now fix any
$\lambda\in\mathbb{R}$.
A convex function
finite on the real line is
differentiable on a
dense set.
We can therefore choose
points of differentiability
$\lambda_n\uparrow\lambda$
and $\mu_n\downarrow\lambda$.
By monotonicity of the
derivatives of a convex
function,
\begin{equation}
\zeta'(\lambda_n)\longrightarrow\zeta'_-(\lambda),
\label{eq:cumulant_one_sided_approximation}
\end{equation}
\begin{equation}
\zeta'(\mu_n)\longrightarrow\zeta'_+(\lambda).
\label{eq:style_s7_008}
\end{equation}
Applying
\eqref{eq:cumulant_implicit_equation}
at these points and using
the continuity of
$\zeta$ and the $\pi_k$
yields
\eqref{eq:cumulant_one_sided_equation}--\eqref{eq:style_s7_006}.
At $\lambda=0$ the two
identities reduce
to normalization of
the probabilities;
the slopes at that
point are determined
by
\eqref{eq:cumulant_derivatives_at_zero}--\eqref{eq:style_s7_004}.

We prove strict
convexity.
If $\zeta$ were affine
on a nonempty open
interval, its
derivative would be
constant on that
interval, with a
value $\alpha\in[0,K]$.
Equation
\eqref{eq:cumulant_implicit_equation}
would then give, throughout
the interval,
\begin{equation}
\zeta(\lambda)=\log\Bigl(\sum_{k\in\Omega}\pi_k(\alpha)e^{k\lambda}\Bigr),
\label{eq:style_s7_050}
\end{equation}
with $\alpha$ fixed.
However,
\begin{equation}
\frac{d^2}{d\lambda^2}\log\Bigl(\sum_{k\in\Omega}\pi_k(\alpha)e^{k\lambda}\Bigr) 
=\frac{\displaystyle\sum_{0\leq j<k\leq K}(k-j)^2\pi_j(\alpha)\pi_k(\alpha)e^{(j+k)\lambda}}
{\displaystyle\Bigl(\sum_{k\in\Omega}\pi_k(\alpha)e^{k\lambda}\Bigr)^2}>0.
\label{eq:cumulant_fixed_state_strict_convexity}
\end{equation}
Positivity already follows
from the term $j=0$, $k=K$,
since $\pi_0(\alpha)$
and $\pi_K(\alpha)$
are strictly positive.
This is a contradiction;
a convex function
without affine intervals
is strictly
convex.

For $\lambda<\mu$,
convexity gives
$\zeta'_+(\lambda)
\leq\zeta'_-(\mu)$.
Equality would force
$\zeta$ to be affine
on $[\lambda,\mu]$;
hence
$\zeta'_+(\lambda)
<\zeta'_-(\mu)$.
Comparing each
point with a point
to its left
and one to its right,
and using
$0\leq\zeta'_-\leq\zeta'_+\leq K$,
yields
\eqref{eq:cumulant_interior_derivatives}.
Comparing instead
with zero
and using
\eqref{eq:cumulant_derivatives_at_zero}--\eqref{eq:style_s7_004}
yields
\eqref{eq:cumulant_strict_derivative_branches}--\eqref{eq:style_s7_007}.

For a fixed law on the first $m$ increments, including a mixture,
with $m$ and the law independent of $N$, the recursion
\eqref{eq:cumulant_exact_recursion} holds for $N\geq m$;
in the nonhomogeneous case it involves the $\pi_{N,k}$.
For every fixed $\lambda$, the difference between the recursive factor
and the one constructed with the limiting functions is bounded,
uniformly in the state and the tilted law, by
\begin{equation}
\sum_{k\in\Omega}e^{k\lambda}
\sup_{\alpha\in[0,K]}
|\pi_{N,k}(\alpha)-\pi_k(\alpha)|
\longrightarrow0.
\label{eq:cumulant_inhomogeneous_factor_error}
\end{equation}
In the homogeneous case the error is zero.
Lemma~\ref{lem:cumulant_changed_ldp} gives concentration
for the same limiting function $\zeta$.
At its points of differentiability, the factor constructed with
the limiting probabilities therefore converges to
\begin{equation}
\sum_{k\in\Omega}\pi_k\bigl(\zeta'(\lambda)\bigr)e^{k\lambda}
=e^{\zeta(\lambda)}>0,
\label{eq:display_019}
\end{equation}
by \eqref{eq:cumulant_implicit_equation}, already proved for
the base law. The logarithm of the recursive factor therefore converges
to $\zeta(\lambda)$; the difference form
\eqref{eq:cumulant_exact_difference} and $\zeta_N\to\zeta$
give \eqref{eq:cumulant_difference_limit} directly.
The other conclusions concern the same $\zeta$ and the limiting
urn functions.
\end{proof}

For $\lambda\neq0$,
normalization of the
probabilities and the identity
\begin{equation}
e^{k\lambda}-1=(e^\lambda-1)\sum_{j=0}^{k-1}e^{j\lambda}
\label{eq:style_s7_051}
\end{equation}
give, at points of
differentiability,
\begin{equation}
\frac{e^{\zeta(\lambda)}-1}{e^\lambda-1}
=\sum_{k=1}^{K}
\pi_k\bigl(\zeta'(\lambda)\bigr)
\sum_{j=0}^{k-1}e^{j\lambda}.
\label{eq:cumulant_generalized_binary_equation}
\end{equation}
The same formula holds
separately with
$\zeta'_-$ and $\zeta'_+$.
For $K=1$, the right-hand
side reduces to
$\pi_1(\zeta'(\lambda))$,
recovering
Ref.~\cite{Franchini_URNS},
Eq.~(2.31).
For $K>1$ all
components enter
with weights depending
on $\lambda$.
For example, for $K=2$
and $\alpha\in[0,2]$,
\begin{equation}
\pi_1(\alpha)
+(1+e^\lambda)\pi_2(\alpha)
=\bar\pi(\alpha)
+(e^\lambda-1)\pi_2(\alpha).
\label{eq:cumulant_two_increment_dependence}
\end{equation}
Invertibility of the
mean function
$\bar\pi$ alone therefore
does not determine the inversion
needed to obtain
$\zeta'$.
For fixed
$\lambda$, one must study
the dependence on
$\alpha$ of the exponential sum in
\eqref{eq:cumulant_map_in_mean}.
The conditions that
allow the equation
to be made explicit and
its solution to be
characterized will
be examined in
Subsection~\ref{subsec:cumulant_boundary_uniqueness}.

\subsection{Boundary conditions and uniqueness}
\label{subsec:cumulant_boundary_uniqueness}

The boundary conditions
for the limiting equation are
obtained from the variational
representation and the values
of the entropy at $0$ and $K$.
To characterize the solution
uniquely, we will then use
an injectivity assumption
with respect to the mean,
formulated separately
on the two half-lines of the
parameter $\lambda$.

\begin{lemma}[Boundary conditions]
\label{lem:cumulant_boundary_values}
Under the assumptions of
Theorem~\ref{thm:path_ldp},
the function $\zeta$ satisfies
\begin{equation}
\lim_{\lambda\to-\infty}\zeta(\lambda)=\phi(0)=\log\pi_0(0),
\label{eq:cumulant_infinite_boundary_values}
\end{equation}
\begin{equation}
\lim_{\lambda\to+\infty}\bigl\{\zeta(\lambda)-K\lambda\bigr\}=\phi(K)=\log\pi_K(K).
\label{eq:style_s7_009}
\end{equation}
The one-sided derivatives
satisfy
\begin{equation}
\lim_{\lambda\to-\infty}\zeta'_-(\lambda)=\lim_{\lambda\to-\infty}\zeta'_+(\lambda)=0,
\label{eq:cumulant_infinite_boundary_slopes}
\end{equation}
\begin{equation}
\lim_{\lambda\to+\infty}\zeta'_-(\lambda)=\lim_{\lambda\to+\infty}\zeta'_+(\lambda)=K.
\label{eq:style_s7_010}
\end{equation}
and
\begin{equation}
\lim_{\lambda\uparrow0}\zeta'_-(\lambda)=\lim_{\lambda\uparrow0}\zeta'_+(\lambda)=\inf C,
\label{eq:cumulant_zero_boundary_slopes}
\end{equation}
\begin{equation}
\lim_{\lambda\downarrow0}\zeta'_-(\lambda)=\lim_{\lambda\downarrow0}\zeta'_+(\lambda)=\sup C.
\label{eq:style_s7_011}
\end{equation}
Moreover,
$\lim_{\lambda\to0^\pm}\zeta(\lambda)=0$.
\end{lemma}

\begin{proof}
For every $\lambda$,
choose any
maximizer $\alpha_\lambda$
in
\eqref{eq:cumulant_variational_terminal}.
If $\lambda<0$, comparison
with $\alpha=0$ and the
inequality $\phi\leq0$
give
$\phi(0)\leq
\lambda\alpha_\lambda+
\phi(\alpha_\lambda)
\leq-|\lambda|\alpha_\lambda$.
If $\lambda>0$, comparison
with $\alpha=K$ gives
$K\lambda+\phi(K)\leq
\lambda\alpha_\lambda+
\phi(\alpha_\lambda)
\leq\lambda\alpha_\lambda$.
Hence
\begin{equation}
0\leq|\lambda|\alpha_\lambda\leq-\phi(0),\qquad\lambda<0,
\label{eq:cumulant_endpoint_maximizer_bounds}
\end{equation}
\begin{equation}
0\leq\lambda(K-\alpha_\lambda)\leq-\phi(K),\qquad\lambda>0.
\label{eq:style_s7_012}
\end{equation}
All maximizers
therefore tend to zero
as $\lambda\to-\infty$
and to $K$ as
$\lambda\to+\infty$.
On the other hand,
\begin{equation}
\phi(0)\leq\zeta(\lambda)\leq\phi(\alpha_\lambda),\qquad\lambda<0,
\label{eq:cumulant_endpoint_value_bounds}
\end{equation}
\begin{equation}
\phi(K)\leq\zeta(\lambda)-K\lambda\leq\phi(\alpha_\lambda),\qquad\lambda>0.
\label{eq:style_s7_013}
\end{equation}
The continuity of $\phi$
and
\eqref{eq:terminal_endpoint_values}--\eqref{eq:style_s5_endpoint_value_K}
prove
\eqref{eq:cumulant_infinite_boundary_values}--\eqref{eq:style_s7_009}.
By
Lemma~\ref{lem:cumulant_one_sided_derivatives},
both one-sided
derivatives are maximizers.
The estimates
\eqref{eq:cumulant_endpoint_maximizer_bounds}--\eqref{eq:style_s7_012}
therefore also give
\eqref{eq:cumulant_infinite_boundary_slopes}--\eqref{eq:style_s7_010}.

By convexity, as $\lambda$ tends to the origin from either
half-line, the one-sided derivatives have the limit given by the corresponding
one-sided derivative at zero. The values already identified in
\eqref{eq:cumulant_derivatives_at_zero}--\eqref{eq:style_s7_004} are $\inf C$ and $\sup C$, giving
\eqref{eq:cumulant_zero_boundary_slopes}--\eqref{eq:style_s7_011}.
Continuity and $\zeta(0)=0$ give the limit of the values.
\end{proof}

The assumption that allows
the derivative to be determined
in the limiting equation
concerns, for every
$\lambda\neq0$, the map
\begin{equation}
\alpha\longmapsto
\sum_{k\in\Omega}\pi_k(\alpha)e^{k\lambda}.
\label{eq:cumulant_map_in_mean}
\end{equation}
The relevant intervals
are $[0,\inf C]$ for
$\lambda<0$ and
$[\sup C,K]$ for
$\lambda>0$.
Injectivity on these
intervals is an
additional condition,
to be checked for the
urn functions considered.

\begin{theorem}[Global characterization on the two half-lines]
\label{thm:cumulant_boundary_uniqueness}
Under the assumptions of Theorem~\ref{thm:path_ldp}, the following
two independent conclusions hold.
If, for every $\lambda<0$, the map
\eqref{eq:cumulant_map_in_mean} is injective on $[0,\inf C]$,
then $\zeta\in C^1((-\infty,0))$ and its restriction is the unique
function $f\in C^1((-\infty,0))$ satisfying
\begin{equation}
e^{f(\lambda)}=\sum_{k\in\Omega}\pi_k\bigl(f'(\lambda)\bigr)e^{k\lambda},\qquad\lambda<0,
\label{eq:cumulant_negative_boundary_problem}
\end{equation}
\begin{equation}
0\leq f'(\lambda)\leq\inf C,\qquad\lambda<0.
\label{eq:style_s7_014}
\end{equation}
If, for every $\lambda>0$, the same map is injective on
$[\sup C,K]$, then $\zeta\in C^1((0,\infty))$ and its
restriction is the unique function $f\in C^1((0,\infty))$ satisfying
\begin{equation}
e^{f(\lambda)}=\sum_{k\in\Omega}\pi_k\bigl(f'(\lambda)\bigr)e^{k\lambda},\qquad\lambda>0,
\label{eq:cumulant_positive_boundary_problem}
\end{equation}
\begin{equation}
\sup C\leq f'(\lambda)\leq K,\qquad\lambda>0.
\label{eq:style_s7_015}
\end{equation}
In each class, the equation and the global slope constraints
already imply the zero limit at the origin and the corresponding asymptote
in \eqref{eq:cumulant_infinite_boundary_values}--\eqref{eq:style_s7_009}; convexity of the
candidate is not assumed. After identification with $\zeta$, the
slope limits in
\eqref{eq:cumulant_infinite_boundary_slopes}--\eqref{eq:style_s7_010} and
\eqref{eq:cumulant_zero_boundary_slopes}--\eqref{eq:style_s7_011} also hold.
On each half-line for which the injectivity assumption holds, the derivative
$\zeta'$ is continuous and strictly increasing, with range
\begin{equation}
\{\zeta'(\lambda):\lambda<0\}=(0,\inf C),
\label{eq:cumulant_derivative_images}
\end{equation}
\begin{equation}
\{\zeta'(\lambda):\lambda>0\}=(\sup C,K).
\label{eq:style_s7_016}
\end{equation}
\end{theorem}

\begin{proof}
The one-sided derivatives of $\zeta$ belong to the respective interval
of means and satisfy \eqref{eq:cumulant_one_sided_equation}--\eqref{eq:style_s7_006}.
Injectivity identifies them; convexity gives continuity of the
derivative, and strict convexity gives its strict increase.
Theorem~\ref{thm:cumulant_limit_equation} therefore ensures
existence in each of the two classes;
Lemma~\ref{lem:cumulant_boundary_values} gives
\eqref{eq:cumulant_derivative_images}--\eqref{eq:style_s7_016}.

Before making the comparison, we derive the boundary data for every candidate.
Normalization of the probabilities and $0\leq f'\leq K$ give
\begin{equation}
\min\{0,K\lambda\}\leq f(\lambda)\leq\max\{0,K\lambda\},
\label{eq:cumulant_candidate_bounds}
\end{equation}
so $f(\lambda)\to0$ at the origin from the half-line considered.
On the negative half-line we also have
\begin{equation}
\log\min_{u\in[0,K]}\pi_0(u)\leq f(\lambda)\leq0.
\label{eq:style_s7_052}
\end{equation}
Since $f'\geq0$, the limit at $-\infty$ exists and is finite.
Integrability of $f'$ on every left tail provides a sequence
$\lambda_n\to-\infty$ with $f'(\lambda_n)\to0$.
The equation along this sequence gives
\begin{equation}
e^{f(\lambda_n)}=\pi_0\bigl(f'(\lambda_n)\bigr)
+\sum_{k=1}^K\pi_k\bigl(f'(\lambda_n)\bigr)e^{k\lambda_n}.
\label{eq:style_ref_s7_negative_sequence}
\end{equation}
Passing to the limit yields
\begin{equation}
\lim_{\lambda\to-\infty}f(\lambda)=\log\pi_0(0).
\label{eq:cumulant_candidate_negative_asymptote}
\end{equation}
On the positive half-line,
\begin{equation}
\log\min_{u\in[0,K]}\pi_K(u)\leq f(\lambda)-K\lambda\leq0.
\label{eq:style_s7_053}
\end{equation}
The function $f-K\lambda$ is nonincreasing and bounded;
integrability of $K-f'$ on every right tail gives a sequence
$\lambda_n\to+\infty$ with $f'(\lambda_n)\to K$.
Dividing the equation by $e^{K\lambda_n}$ yields
\begin{equation}
e^{f(\lambda_n)-K\lambda_n}=\pi_K\bigl(f'(\lambda_n)\bigr)
+\sum_{k=0}^{K-1}\pi_k\bigl(f'(\lambda_n)\bigr)e^{-(K-k)\lambda_n}.
\label{eq:style_global_candidate_positive_sequence}
\end{equation}
Passing to the limit gives
\begin{equation}
\lim_{\lambda\to+\infty}\bigl\{f(\lambda)-K\lambda\bigr\}=\log\pi_K(K).
\label{eq:cumulant_candidate_positive_asymptote}
\end{equation}
Full limits of the derivatives of the candidates have not been assumed.

Let $f_1,f_2$ be two candidates on the same half-line and
$w=f_1-f_2$. The boundary data just derived give $w\to0$
at both endpoints; on the positive half-line the terms
$K\lambda$ cancel. If it were not zero, $w$ would attain
a nonzero interior extremum at some point $\lambda_*$.
At that point $w'(\lambda_*)=0$, so the equations for the two
candidates give
\begin{equation}
e^{f_1(\lambda_*)}=\sum_{k\in\Omega}\pi_k\bigl(f_1'(\lambda_*)\bigr)e^{k\lambda_*} 
=\sum_{k\in\Omega}\pi_k\bigl(f_2'(\lambda_*)\bigr)e^{k\lambda_*}=e^{f_2(\lambda_*)}.
\label{eq:cumulant_solution_comparison}
\end{equation}
This would imply $w(\lambda_*)=0$, a contradiction.
Hence $f_1=f_2$.
Identification with $\zeta$ finally transfers to the
candidates the slope limits already proved.
\end{proof}

The comparison uses the boundary data and the equation at a
stationary point, without Lipschitz conditions on the inverse. For $K=1$ the exponential sum is
$1+(e^\lambda-1)\pi_1(\alpha)$: injectivity is equivalent to that of
$\pi_1$ on the corresponding interval. This gives the characterization
for general increments corresponding to Ref.~\cite{Franchini_URNS},
Theorem~9, Eqs.~(2.32)--(2.33).

A sufficient condition
for verifying injectivity
can be expressed directly
in terms of the increment
probabilities.

\begin{lemma}[Condition on probability sums]
\label{lem:cumulant_probability_sums}
On one of the intervals
$[0,\inf C]$ or
$[\sup C,K]$,
suppose that for
every pair
$\alpha<\beta$ we have
\begin{equation}
\sum_{k=j}^{K}\pi_k(\alpha)
\leq
\sum_{k=j}^{K}\pi_k(\beta),
\qquad j=1,\ldots,K,
\label{eq:cumulant_probability_sum_order}
\end{equation}
with at least one
strict inequality
for every pair.
Then, on that
interval, the map
\eqref{eq:cumulant_map_in_mean}
is strictly increasing
for every $\lambda>0$
and strictly
decreasing for every
$\lambda<0$.
If all the
inequalities in
\eqref{eq:cumulant_probability_sum_order}
are reversed,
with the same
strict inequality
condition for
every pair,
the directions of monotonicity
are interchanged.
\end{lemma}

\begin{proof}
Normalization of
the probabilities
and interchange of
two finite sums
give
\begin{equation}
\sum_{k\in\Omega}\pi_k(\alpha)e^{k\lambda}
=
1+(e^\lambda-1)
\sum_{j=1}^{K}e^{(j-1)\lambda}
\sum_{k=j}^{K}\pi_k(\alpha).
\label{eq:cumulant_probability_sum_identity}
\end{equation}
All coefficients
$e^{(j-1)\lambda}$
are strictly
positive.
The assumptions therefore
imply that the sum
multiplying
$e^\lambda-1$
is strictly
increasing in
$\alpha$.
The sign of
$e^\lambda-1$
gives the two stated
directions.
Reversing the
inequalities on the
probability sums
also reverses
these directions.
\end{proof}

For $K=1$ the criterion
reduces to
strict monotonicity
of $\pi_1$.
For $K=2$, the two
sums to compare
are
\begin{equation}
\pi_1(\alpha)+\pi_2(\alpha)=1-\pi_0(\alpha),
\qquad
\pi_2(\alpha).
\label{eq:cumulant_two_increment_order}
\end{equation}
It is therefore sufficient
that $\pi_0$ be
nonincreasing and
$\pi_2$ nondecreasing,
and that, for every
$\alpha<\beta$
in the interval,
at least one of
$\pi_0(\alpha)>\pi_0(\beta)$
and
$\pi_2(\alpha)<\pi_2(\beta)$
be satisfied.
The component
$\pi_1$ is determined
by normalization
and does not require a
separate monotonicity
condition.

\subsection{Shape and regularity of the entropy}
\label{subsec:entropy_shape}

The inequality
\eqref{eq:cumulant_entropy_dual_bound}
provides an upper
bound for $\phi$
through the cumulant
generating function.
To obtain equality,
one must identify
the terminal states
that attain the maximum
in
\eqref{eq:cumulant_variational_terminal}.
The ranges of the derivative
determined in
Theorem~\ref{thm:cumulant_boundary_uniqueness}
will allow us to cover
in full each
of the two intervals
outside $[\inf C,\sup C]$.

\begin{lemma}[Equality in the inverse formula]
\label{lem:entropy_dual_equality}
Under the assumptions of
Theorem~\ref{thm:path_ldp},
the identity
\begin{equation}
\phi(\alpha)
=\inf_{\lambda\in\mathbb{R}}
\{\zeta(\lambda)-\lambda\alpha\}
\label{eq:entropy_dual_known_states}
\end{equation}
holds for
$\alpha\in\{0,K\}\cup[\inf C,\sup C]$
and for every value
$\alpha=\zeta'(\lambda_0)$
corresponding to a
point of differentiability
$\lambda_0$ of $\zeta$.
For $\alpha\in[\inf C,\sup C]$,
the minimum is attained
uniquely at $\lambda=0$.
For $\alpha=\zeta'(\lambda_0)$,
it is attained uniquely
at $\lambda=\lambda_0$.
At the endpoints $\alpha=0$
and $\alpha=K$, the infimum
is not attained at
a finite parameter;
it is obtained as a limit
as $\lambda\to-\infty$
and $\lambda\to+\infty$,
respectively.
\end{lemma}

\begin{proof}
At a point of
differentiability,
Lemma~\ref{lem:cumulant_one_sided_derivatives}
identifies the unique
maximizer
in
\eqref{eq:cumulant_variational_terminal}
with the value
$\zeta'(\lambda)$.
Hence
\begin{equation}
\phi\bigl(\zeta'(\lambda)\bigr)
=\zeta(\lambda)-\lambda\zeta'(\lambda).
\label{eq:entropy_parametric_value}
\end{equation}
The inequality
\eqref{eq:cumulant_entropy_dual_bound}
and evaluation of
the infimum at that
parameter give
\eqref{eq:entropy_dual_known_states}
at the point considered.

For $\alpha\in[\inf C,\sup C]$,
we have $\phi(\alpha)=0$,
and the variational
representation implies
$\zeta(\lambda)-\lambda\alpha\geq0$.
At $\lambda=0$,
equality holds.
Uniqueness of each
attained minimum follows
from strict convexity
of
$\lambda\mapsto\zeta(\lambda)-\lambda\alpha$.

Finally,
Lemma~\ref{lem:cumulant_boundary_values}
gives
\begin{equation}
\inf_{\lambda\in\mathbb R}\zeta(\lambda)=\phi(0),
\label{eq:entropy_dual_endpoints}
\end{equation}
\begin{equation}
\inf_{\lambda\in\mathbb R}\bigl\{\zeta(\lambda)-K\lambda\bigr\}=\phi(K).
\label{eq:style_s7_017}
\end{equation}
The inequalities
\eqref{eq:cumulant_interior_derivatives}
imply that $\zeta$
is strictly increasing
and that
$\lambda\mapsto\zeta(\lambda)-K\lambda$
is strictly
decreasing.
The two infima
therefore cannot
be attained
at a finite parameter.
\end{proof}

The injectivity
assumptions on the two
half-lines remain
independent for the
conclusions about
entropy as well.

\begin{theorem}[Shape and regularity of the entropy]
\label{thm:entropy_shape}
Under the assumptions of
Theorem~\ref{thm:path_ldp},
suppose that the
injectivity assumption
of
Theorem~\ref{thm:cumulant_boundary_uniqueness}
for the negative
half-line holds.
Then
\begin{equation}
\phi(\alpha)
=\inf_{\lambda\leq0}
\{\zeta(\lambda)-\lambda\alpha\},
\qquad 0\leq\alpha\leq\inf C.
\label{eq:entropy_negative_duality}
\end{equation}
The restriction of
$\phi$ to $[0,\inf C]$
is absolutely
continuous, strictly
concave, and strictly
increasing from
$\log\pi_0(0)$ to zero.
On the interval
$(0,\inf C)$ it is
of class $C^1$,
with
\begin{equation}
\phi'(\alpha)>0,
\label{eq:entropy_negative_derivative_limits}
\end{equation}
\begin{equation}
\lim_{\alpha\downarrow0}\phi'(\alpha)=+\infty,
\label{eq:style_s7_018}
\end{equation}
\begin{equation}
\lim_{\alpha\uparrow\inf C}\phi'(\alpha)=0.
\label{eq:style_s7_019}
\end{equation}

Under only the
injectivity assumption for
the positive half-line,
we instead have
\begin{equation}
\phi(\alpha)
=\inf_{\lambda\geq0}
\{\zeta(\lambda)-\lambda\alpha\},
\qquad \sup C\leq\alpha\leq K.
\label{eq:entropy_positive_duality}
\end{equation}
The restriction of
$\phi$ to $[\sup C,K]$
is absolutely
continuous, strictly
concave, and strictly
decreasing from zero
to $\log\pi_K(K)$.
On the interval
$(\sup C,K)$ it is
of class $C^1$,
with
\begin{equation}
\phi'(\alpha)<0,
\label{eq:entropy_positive_derivative_limits}
\end{equation}
\begin{equation}
\lim_{\alpha\downarrow\sup C}\phi'(\alpha)=0,
\label{eq:style_s7_020}
\end{equation}
\begin{equation}
\lim_{\alpha\uparrow K}\phi'(\alpha)=-\infty.
\label{eq:style_s7_021}
\end{equation}

On each half-line
for which the
injectivity assumption is
satisfied, the relation
\begin{equation}
\phi'(\zeta'(\lambda))=-\lambda.
\label{eq:entropy_parametric_derivative}
\end{equation}
holds.
For every $\alpha$
in the corresponding
open interval,
the dual minimum
is attained at the unique
parameter $\lambda$
such that
$\zeta'(\lambda)=\alpha$.

If both injectivity
assumptions hold,
then
\begin{equation}
\phi(\alpha)
=\inf_{\lambda\in\mathbb{R}}
\{\zeta(\lambda)-\lambda\alpha\},
\qquad 0\leq\alpha\leq K.
\label{eq:entropy_full_duality}
\end{equation}
The function $\phi$
is concave and
absolutely continuous
on $[0,K]$, and is
of class $C^1$
on $(0,K)$.
In particular,
\begin{equation}
\phi'(\alpha)=0,
\qquad \inf C\leq\alpha\leq\sup C.
\label{eq:entropy_zero_interval_derivative}
\end{equation}
For $0<\alpha<K$,
the infimum in
\eqref{eq:entropy_full_duality}
is a unique minimum;
at the endpoints the
qualification in
Lemma~\ref{lem:entropy_dual_equality} applies.
Concavity is
strict on all of
$[0,K]$ if and
only if $C$
contains a single
point.
\end{theorem}

\begin{proof}
On the negative half-line,
Theorem~\ref{thm:cumulant_boundary_uniqueness} makes $\zeta'$ continuous
and strictly increasing, with range $(0,\inf C)$.
Every point of this interval is therefore $\zeta'(\lambda)$ for a unique
$\lambda<0$; Lemma~\ref{lem:entropy_dual_equality} gives dual
equality, with the minimum at that parameter. The same lemma covers
$\alpha=\inf C$ and $\alpha=0$, with parameter zero and limit at
$-\infty$, respectively. This proves
\eqref{eq:entropy_negative_duality} on the closed interval.
On the positive half-line the range $(\sup C,K)$ similarly gives
\eqref{eq:entropy_positive_duality}, including the boundaries.
Only now, as an infimum of affine functions of $\alpha$, does each
representation imply concavity on the corresponding interval.

\phantomsection\label{par:entropy_secants_and_integrals}
Secants and integrals suffice for regularity. If
$\alpha=\zeta'(\lambda)<\beta=\zeta'(\mu)$ belong to the same
open branch, then $\lambda<\mu$ and the two evaluations of the dual
formula at $\lambda$ and $\mu$ give
\begin{equation}
-\mu\leq\frac{\phi(\beta)-\phi(\alpha)}{\beta-\alpha}\leq-\lambda.
\label{eq:entropy_secant_bounds}
\end{equation}
The inverse of $\zeta'$ is continuous. Letting $\beta$ tend to $\alpha$
and using the same comparison from the left yields
$\phi'(\zeta'(\lambda))=-\lambda$, without assuming existence of
$\zeta''$. The derivative of $\phi$ is continuous and strictly decreasing.
Its sign and the limits of the parameter give the monotonicity properties and
\eqref{eq:entropy_negative_derivative_limits}--\eqref{eq:style_s7_019} and \eqref{eq:entropy_positive_derivative_limits}--\eqref{eq:style_s7_021}.
Continuity of $\phi$ extends these monotonicity properties to the closed branches;
strict increase of the parameter rules out every nontrivial affine
portion and gives strict concavity on each of them.

The fundamental theorem of calculus holds on interior compact intervals.
Letting the endpoints tend to the boundaries and using the sign of $\phi'$
and the continuity of $\phi$ yields, under the respective assumption,
\begin{equation}
\int_0^{\inf C}\phi'(\alpha)\,d\alpha=-\log\pi_0(0)<\infty,
\label{eq:entropy_integral_representation}
\end{equation}
\begin{equation}
\phi(x)=\phi(0)+\int_0^x\phi'(\alpha)\,d\alpha,\qquad0\leq x\leq\inf C,
\label{eq:style_s7_022}
\end{equation}
\begin{equation}
\int_{\sup C}^K|\phi'(\alpha)|\,d\alpha=-\log\pi_K(K)<\infty,
\label{eq:style_s7_023}
\end{equation}
\begin{equation}
\phi(x)=\int_{\sup C}^x\phi'(\alpha)\,d\alpha,\qquad\sup C\leq x\leq K.
\label{eq:style_s7_024}
\end{equation}
Integrability and integral reconstruction prove absolute
continuity up to the physical boundaries.

If both assumptions hold, the two branches and the interval
$[\inf C,\sup C]$ cover $[0,K]$.
Lemma~\ref{lem:entropy_dual_equality} completes
\eqref{eq:entropy_full_duality}; global concavity follows from this.
On the central interval $\phi=0$, and the limits of the derivative from the outer branches
are zero. Continuity and the mean value theorem give a zero derivative
at the joining points as well, including the case $\inf C=\sup C$.
Thus $\phi\in C^1((0,K))$ and
\eqref{eq:entropy_zero_interval_derivative} holds.
The integral reconstructions and the constant portion give global
absolute continuity and
\begin{equation}
\int_0^K|\phi'(\alpha)|\,d\alpha
=-\log\pi_0(0)-\log\pi_K(K).
\label{eq:entropy_total_variation}
\end{equation}
If $\inf C<\sup C$, the constant portion prevents global strict
concavity. If they coincide, the derivative, strictly decreasing on the branches,
positive on the left and negative on the right, rules out every nontrivial affine
interval: concavity is strict on $[0,K]$.
\end{proof}

Theorem~\ref{thm:entropy_shape} extends Corollary~11 of
Ref.~\cite{Franchini_URNS}. Its proof separates dual contact
from the secant bound
\eqref{eq:entropy_secant_bounds} and the integral representations
\eqref{eq:entropy_integral_representation}--\eqref{eq:style_s7_024}:
the latter will be reused in the affine class after verifying there
the regularity of the cumulant generating function.

For $\alpha=\zeta'(\lambda)$, the parametric identities give
\begin{equation}
\lambda=-\phi'(\alpha),
\label{eq:entropy_cumulant_substitution}
\end{equation}
\begin{equation}
\zeta(\lambda)=\phi(\alpha)-\alpha\phi'(\alpha).
\label{eq:style_s7_025}
\end{equation}
Substituting into \eqref{eq:cumulant_implicit_equation} yields
the entropy equation in the two equivalent forms
\begin{equation}
e^{\phi(\alpha)-\alpha\phi'(\alpha)}=\sum_{k\in\Omega}\pi_k(\alpha)e^{-k\phi'(\alpha)},
\label{eq:entropy_equation}
\end{equation}
\begin{equation}
e^{\phi(\alpha)}=\sum_{k\in\Omega}\pi_k(\alpha)e^{(\alpha-k)\phi'(\alpha)}.
\label{eq:style_s7_026}
\end{equation}
The identity holds on $(0,\inf C)$ under the negative assumption and on
$(\sup C,K)$ under the positive one. If both hold, it holds for every
$0<\alpha<K$: on the central interval $\phi=\phi'=0$, leaving
normalization. The first form will be used in the inverse problem.

Still with $\alpha=\zeta'(\lambda)$,
\begin{equation}
\alpha\phi'(\alpha)=\phi(\alpha)-\zeta(\lambda),
\label{eq:entropy_endpoint_product_identities}
\end{equation}
\begin{equation}
(K-\alpha)\phi'(\alpha)=\zeta(\lambda)-K\lambda-\phi(\alpha).
\label{eq:style_s7_027}
\end{equation}
Lemma~\ref{lem:cumulant_boundary_values} and continuity of $\phi$
give, under the respective branch assumption,
\begin{equation}
\lim_{\alpha\downarrow0}\alpha\phi'(\alpha)=0,
\label{eq:entropy_endpoint_product_limits}
\end{equation}
\begin{equation}
\lim_{\alpha\uparrow K}(K-\alpha)\phi'(\alpha)=0.
\label{eq:style_s7_028}
\end{equation}
In the equivalent form \eqref{eq:style_s7_026}, only the term $k=0$
survives at the left boundary and only the term $k=K$ at the right boundary;
we recover $e^{\phi(0)}=\pi_0(0)$ and $e^{\phi(K)}=\pi_K(K)$.

\subsection{Inverse problem}
\label{subsec:inverse_problem}

Equation \eqref{eq:entropy_equation} suggests starting
from a function $f$ and seeking urn probabilities that
realize it as the entropy of the terminal state. For
binary increments, the corresponding formula is given in
Ref.~\cite{Franchini_URNS}, Corollary~10, Eq.~(2.34).
Here we consider functions that are concave and continuous on the
closed interval. A verification criterion based on the functional
\eqref{eq:path_rate} will allow us to characterize those
realizable as entropies, explicitly construct
at least one urn, and determine the freedom in the probabilities
when $K>1$. Interior $C^1$ regularity will turn out to be
a necessary condition for realizability in
this class.

The following verification result uses the regularity
of the prescribed function and strict positivity of the
probabilities. It does not require the injectivity assumptions
of Theorem~\ref{thm:cumulant_boundary_uniqueness}.

\begin{theorem}[Verification of the prescribed entropy]
\label{thm:inverse_verification}
Let the functions $\pi_k$ satisfy the assumptions of
Theorem~\ref{thm:path_ldp}. Let
$f:[0,K]\to\mathbb{R}$ be continuous and concave, of class
$C^1$ on $(0,K)$, with
\begin{equation}
\lim_{\alpha\downarrow0}f'(\alpha)=+\infty,
\label{eq:inverse_slope_conditions}
\end{equation}
\begin{equation}
\lim_{\alpha\uparrow K}f'(\alpha)=-\infty.
\label{eq:style_s7_029}
\end{equation}
Suppose that
\begin{equation}
e^{f(\alpha)-\alpha f'(\alpha)}
=
\sum_{k\in\Omega}\pi_k(\alpha)e^{-k f'(\alpha)},
\qquad 0<\alpha<K.
\label{eq:inverse_equation}
\end{equation}
Then
\begin{equation}
\phi(\alpha)=f(\alpha)
=-\min_{\varphi\in Q(\alpha)}I(\varphi),
\qquad 0\leq\alpha\leq K.
\label{eq:inverse_entropy_identification}
\end{equation}
For every $\alpha$, the minimum can be attained by a
path constructed through the boundary problem
\eqref{eq:inverse_optimal_mean_equation}--\eqref{eq:style_s7_034} below.
\end{theorem}

\begin{proof}
The derivative of $f$ is continuous, nonincreasing, and, by
\eqref{eq:inverse_slope_conditions}--\eqref{eq:style_s7_029}, vanishes at at least one
interior point. At each of its zeros, a maximum of $f$, the inverse equation gives
\begin{equation}
e^f=\sum_k\pi_k=1.
\label{eq:style_s7_054}
\end{equation}
 Hence
\begin{equation}
\max_{\alpha\in[0,K]}f(\alpha)=0,
\label{eq:inverse_normalization}
\end{equation}
\begin{equation}
f\leq0.
\label{eq:style_s7_030}
\end{equation}
Before identifying $f$ with the entropy, concavity gives, near the
respective boundaries,
\begin{equation}
0\leq\alpha f'(\alpha)\leq f(\alpha)-f(0),
\label{eq:inverse_boundary_product_bounds}
\end{equation}
\begin{equation}
0\leq(K-\alpha)|f'(\alpha)|\leq f(\alpha)-f(K).
\label{eq:style_s7_031}
\end{equation}
By continuity, these imply
\begin{equation}
\lim_{\alpha\downarrow0}\alpha f'(\alpha)=0,
\label{eq:inverse_boundary_products}
\end{equation}
\begin{equation}
\lim_{\alpha\uparrow K}(K-\alpha)f'(\alpha)=0.
\label{eq:style_s7_032}
\end{equation}
Multiplying \eqref{eq:inverse_equation} by $e^{\alpha f'(\alpha)}$
and passing to the boundaries, the term $k=0$
and the term $k=K$ survive, respectively, so
\begin{equation}
f(0)=\log\pi_0(0),
\label{eq:inverse_boundary_values}
\end{equation}
\begin{equation}
f(K)=\log\pi_K(K).
\label{eq:style_s7_033}
\end{equation}

For every $\varphi\in Q(\alpha)$ set
$\psi(\tau)=\varphi(\tau)/\tau$. Where $0<\psi(\tau)<K$, evaluating
the local supremum at the parameter $-f'(\psi(\tau))$ gives almost everywhere
\begin{equation}
\sup_{\lambda\in\mathbb R}\Bigl\{\lambda\partial_\tau\varphi
-\log\sum_{k\in\Omega}\pi_k(\psi)e^{k\lambda}\Bigr\}
\geq-f'(\psi)\partial_\tau\varphi-f(\psi)+\psi f'(\psi),
\label{eq:inverse_path_inequality}
\end{equation}
\begin{equation}
-f'(\psi)\partial_\tau\varphi-f(\psi)+\psi f'(\psi)
=-\frac{d}{d\tau}\bigl\{\tau f(\psi(\tau))\bigr\}.
\label{eq:style_s7_055}
\end{equation}
The chain rule applies on every compact time interval on which
$\psi$ stays away from the boundaries. Since $\varphi$ and $K\tau-\varphi$ are
nondecreasing, each physical boundary can be followed only
on an initial portion; after entering the interior, the path remains there.
If the path is always interior for $\tau>0$, we integrate from
$\varepsilon$ to one: the initial term tends to zero because
$|\varepsilon f(\psi(\varepsilon))|\leq\varepsilon\|f\|_\infty$.
If it follows the zero boundary up to $\tau_0$, the initial cost is
$-\tau_0\log\pi_0(0)=-\tau_0f(0)$, by
Lemma~\ref{lem:local_frequency_minimum}; the limit of the integral
on the interior portion produces $\tau_0f(0)$ and cancels it.
On the boundary $K\tau$, the same cancellation holds with $f(K)$.
For paths lying entirely on a boundary, equality follows directly
from \eqref{eq:inverse_boundary_values}--\eqref{eq:style_s7_033}. In all cases,
\begin{equation}
I(\varphi)\geq-f(\alpha),\qquad\varphi\in Q(\alpha).
\label{eq:inverse_action_lower_bound}
\end{equation}

To obtain equality, we construct directly in $\tau$ a
solution of
\begin{equation}
\tau\partial_\tau\psi^*(\tau)
=\frac{\displaystyle\sum_{k\in\Omega}k\pi_k\bigl(\psi^*(\tau)\bigr)e^{-kf'\bigl(\psi^*(\tau)\bigr)}}
{\displaystyle\sum_{k\in\Omega}\pi_k\bigl(\psi^*(\tau)\bigr)e^{-kf'\bigl(\psi^*(\tau)\bigr)}}-\psi^*(\tau),
\label{eq:inverse_optimal_mean_equation}
\end{equation}
\begin{equation}
\psi^*(1)=\alpha.
\label{eq:style_s7_034}
\end{equation}
The ratio is continuous and belongs to $(0,K)$ for interior means.
The extreme slopes of $f$ extend it with value zero at zero and
$K$ at $K$: at the right boundary it suffices to divide the sums by $e^{-Kf'}$.
Denote by $h(u)$ the right-hand side extended in this way. It is continuous
on $[0,K]$ and vanishes at the boundaries. Its zeros are those of the tilted
field and are not identified with $C$.

If $h(\alpha)=0$, choose $\psi^*\equiv\alpha$, including at the physical
endpoints. Otherwise, let $(a,b)$ be the component of $\{h\neq0\}$ that
contains $\alpha$, and let $c=a$ if $h>0$, $c=b$ if $h<0$ on it.
Define the arrival time by the extended integral
\begin{equation}
\tau_*:=\exp\Bigl(\int_\alpha^c\frac{du}{h(u)}\Bigr)\in[0,1),
\label{eq:display_020}
\end{equation}
\begin{equation}
e^{-\infty}:=0.
\label{eq:style_s7_035}
\end{equation}
For $\tau_*<\tau\leq1$, the identity
\begin{equation}
\int_\alpha^{\psi^*(\tau)}\frac{du}{h(u)}=\log\tau
\label{eq:display_021}
\end{equation}
defines a unique mean between $\alpha$ and $c$, with
$\tau\partial_\tau\psi^*=h(\psi^*)$ and limit $c$ at time $\tau_*$.
If $\tau_*>0$, extend it by $\psi^*(\tau)=c$ for
$0<\tau\leq\tau_*$: since $h(c)=0$, the joining is $C^1$.
If $\tau_*=0$, it is already defined at all positive times.
This construction uses the sign of the tilted field even when
its zero set is infinite.

We reconstruct the path and its constraints together:
\begin{equation}
\varphi^*(\tau)=\tau\psi^*(\tau),\qquad\tau>0,
\label{eq:inverse_optimal_path}
\end{equation}
\begin{equation}
\varphi^*(0)=0,
\label{eq:style_s7_036}
\end{equation}
\begin{equation}
0\leq\partial_\tau\varphi^*=\psi^*+\tau\partial_\tau\psi^*\leq K.
\label{eq:style_s7_037}
\end{equation}
The last derivative is the ratio in
\eqref{eq:inverse_optimal_mean_equation}, with the boundary extensions.
The constraint $0\leq\varphi^*(\tau)\leq K\tau$ and the derivative constraint
extend the path absolutely continuously to the origin;
hence $\varphi^*\in Q(\alpha)$.

Where $0<\psi^*<K$, the velocity is strictly between zero and $K$, and
Lemma~\ref{lem:local_frequency_minimum}, through
\eqref{eq:local_frequency_parameter}, identifies the parameter
\begin{equation}
\lambda\!\left(\partial_\tau\varphi^*,\psi^*\right)=-f'(\psi^*).
\label{eq:inverse_local_parameter}
\end{equation}
It attains the local supremum, so
\eqref{eq:inverse_path_inequality} is an equality along the path.
The initial terms and any boundary portions are treated as
above, yielding, for every $\alpha\in[0,K]$,
\begin{equation}
I(\varphi^*)=-f(\alpha).
\label{eq:inverse_action_equality}
\end{equation}
Together with \eqref{eq:inverse_action_lower_bound} and definition
\eqref{eq:terminal_entropy_variational}, this proves
\eqref{eq:inverse_entropy_identification}.
\end{proof}

The identification just proved and
Theorem~\ref{thm:zero_entropy_interval} automatically
give
\begin{equation}
\{\alpha\in[0,K]:f(\alpha)=0\}
=[\inf C,\sup C].
\label{eq:inverse_zero_interval}
\end{equation}
The endpoints of the zero portion are therefore
compatible with the fixed points of the constructed urn.

We can now characterize the functions that
admit at least one realization in the class
considered. The following construction uses a
scalar parameter $\pi_f$, in agreement with the
notation of the binary inverse formula in
Ref.~\cite{Franchini_URNS}, Eq.~(2.34).

\begin{theorem}[Characterization and construction of concave entropies]
\label{thm:inverse_binomial_construction}
Let $f:[0,K]\to\mathbb{R}$ be continuous and concave.
There exist continuous, strictly
positive urn functions, normalized on $\Omega$, whose
terminal-state entropy is $f$, if and
only if
\begin{equation}
f\in C^1((0,K)),
\label{eq:inverse_realizability_conditions}
\end{equation}
\begin{equation}
\max_{\alpha\in[0,K]}f(\alpha)=0,
\label{eq:style_s7_038}
\end{equation}
\begin{equation}
\lim_{\alpha\downarrow0}f'(\alpha)=+\infty,
\label{eq:style_s7_039}
\end{equation}
\begin{equation}
\lim_{\alpha\uparrow K}f'(\alpha)=-\infty.
\label{eq:style_s7_040}
\end{equation}
When these conditions are satisfied,
an explicit choice is
\begin{equation}
\pi_k(\alpha)
=\binom{K}{k}\pi_f(\alpha)^k
\bigl(1-\pi_f(\alpha)\bigr)^{K-k},
\qquad k\in\Omega,
\label{eq:inverse_binomial_probabilities}
\end{equation}
where
\begin{equation}
\pi_f(\alpha):=\frac{e^{[f(\alpha)-\alpha f'(\alpha)]/K}-1}{e^{-f'(\alpha)}-1},
\qquad0<\alpha<K,\qquad f'(\alpha)\neq0,
\label{eq:inverse_binomial_parameter}
\end{equation}
\begin{equation}
\pi_f(\alpha):=\frac\alpha K,\qquad0<\alpha<K,\qquad f'(\alpha)=0,
\label{eq:style_s7_056}
\end{equation}
\begin{equation}
\pi_f(0):=1-e^{f(0)/K},
\label{eq:style_s7_057}
\end{equation}
\begin{equation}
\pi_f(K):=e^{f(K)/K}.
\label{eq:style_s7_058}
\end{equation}
The function $\pi_f$ is continuous and takes values
strictly between zero and one.
For this choice we also have
\begin{equation}
C=\{\alpha\in[0,K]:f(\alpha)=0\}.
\label{eq:inverse_binomial_fixed_points}
\end{equation}
\end{theorem}

\begin{proof}
If $f$ is the entropy of an urn in the class, its maximum is zero by
\eqref{eq:style_s5_entropy_maximum}. The other necessary conditions
all follow from strict convexity of $\zeta$, proved in
Theorem~\ref{thm:cumulant_limit_equation}.
At an interior point, concavity makes $\alpha$ a maximizer of
$\lambda\alpha+f(\alpha)$ for every
$\lambda\in[-f'_-(\alpha),-f'_+(\alpha)]$, so
\begin{equation}
\zeta(\lambda)=\lambda\alpha+f(\alpha),
\qquad -f'_-(\alpha)\leq\lambda\leq-f'_+(\alpha).
\label{eq:inverse_corner_affine_cumulant}
\end{equation}
A nondegenerate interval would contradict strict convexity:
$f$ is therefore differentiable in the interior and, being concave,
has a continuous derivative. A finite limit of $f'$ at zero would give,
by the tangent inequality, $\zeta(\lambda)=f(0)$ on a
left half-line; a finite limit at $K$ would give
$\zeta(\lambda)=K\lambda+f(K)$ on a right half-line.
Both would be affine portions, again ruled out by strict
convexity. Monotonicity of $f'$ therefore forces the two infinite limits
in \eqref{eq:style_s7_039}--\eqref{eq:style_s7_040}.

Now assume these conditions. The set of maxima is a nonempty
interval contained in $(0,K)$, on which $f=f'=0$; outside it,
$f<0$, with positive derivative on the left and negative derivative on the right.
In particular, $f(0),f(K)<0$.
For a maximum point $c$, the tangent inequality gives
\begin{equation}
0=f(c)\leq f(\alpha)+(c-\alpha)f'(\alpha).
\label{eq:inverse_tangent_at_maximum}
\end{equation}
Since $0<c<K$, these imply
\begin{equation}
-Kf'(\alpha)<f(\alpha)-\alpha f'(\alpha)<0,\qquad f'(\alpha)>0,
\label{eq:inverse_binomial_parameter_bounds}
\end{equation}
\begin{equation}
0<f(\alpha)-\alpha f'(\alpha)<-Kf'(\alpha),\qquad f'(\alpha)<0.
\label{eq:style_s7_041}
\end{equation}
After division by $K$, these inequalities give
$0<\pi_f(\alpha)<1$. Where $f'=0$, instead,
$\pi_f(\alpha)=\alpha/K\in(0,1)$.

The same tangent, for $f(c)=0$ and $f'(\alpha)\neq0$, implies
\begin{equation}
\left|\frac{f(\alpha)}{f'(\alpha)}\right|\leq|\alpha-c|.
\label{eq:inverse_zero_slope_ratio}
\end{equation}
If $c$ is approached by points with nonzero derivative, using
$(e^x-1)/x\to1$ yields the joining relation
\begin{equation}
\lim_{\substack{\alpha\to c\\f'(\alpha)\neq0}}\pi_f(\alpha)
=\lim_{\substack{\alpha\to c\\f'(\alpha)\neq0}}
\frac{\alpha-f(\alpha)/f'(\alpha)}{K}=\frac cK.
\label{eq:inverse_binomial_join}
\end{equation}
In the interior of the interval of maxima, the joining relation is already given by the
definition. At the physical boundaries,
\eqref{eq:inverse_boundary_product_bounds}--\eqref{eq:style_s7_032}
apply, since they depend only on the properties of $f$. The limit at zero is
direct; for the limit at $K$, we rewrite
\begin{equation}
\pi_f(\alpha)=
\frac{e^{[f(\alpha)+(K-\alpha)f'(\alpha)]/K}-e^{f'(\alpha)}}
{1-e^{f'(\alpha)}}.
\label{eq:inverse_binomial_right_rewriting}
\end{equation}
Hence
\begin{equation}
\lim_{\alpha\downarrow0}\pi_f(\alpha)=1-e^{f(0)/K},
\label{eq:inverse_binomial_boundary_limits}
\end{equation}
\begin{equation}
\lim_{\alpha\uparrow K}\pi_f(\alpha)=e^{f(K)/K}.
\label{eq:style_s7_042}
\end{equation}
The limits belong to $(0,1)$, so $\pi_f$ is continuous and stays
uniformly away from zero and one on the compact set $[0,K]$.
The probabilities \eqref{eq:inverse_binomial_probabilities} are therefore
continuous, uniformly positive, and normalized.

For $f'(\alpha)\neq0$, the binomial formula gives
\begin{equation}
\sum_{k\in\Omega}\pi_k(\alpha)e^{-kf'(\alpha)}
=\bigl(1-\pi_f(\alpha)+\pi_f(\alpha)e^{-f'(\alpha)}\bigr)^K
=e^{f(\alpha)-\alpha f'(\alpha)}.
\label{eq:inverse_binomial_identity}
\end{equation}
Where $f'=0$, we also have $f=0$, and the identity follows from normalization.
Theorem~\ref{thm:inverse_verification} identifies the entropy with $f$.
Finally, \eqref{eq:inverse_zero_interval} gives
$C\subseteq\{f=0\}=[\inf C,\sup C]$; conversely, on $\{f=0\}$
the construction gives $\bar\pi(\alpha)=K\pi_f(\alpha)=\alpha$.
The two inclusions prove \eqref{eq:inverse_binomial_fixed_points}.
\end{proof}

For $K=1$, at points where $f'\neq0$,
formula
\eqref{eq:inverse_binomial_parameter}
becomes
\begin{equation}
\pi_1(\alpha)=\pi_f(\alpha)=\frac{e^{f(\alpha)-\alpha f'(\alpha)}-1}{e^{-f'(\alpha)}-1},
\label{eq:inverse_binary_formula}
\end{equation}
\begin{equation}
\pi_0(\alpha)=1-\pi_f(\alpha).
\label{eq:style_s7_043}
\end{equation}
This recovers the formula in
Ref.~\cite{Franchini_URNS}, Corollary~10,
Eq.~(2.34), with a direct verification
in the class $C^1((0,1))$ considered here.
The definition on the zero portion and the boundary
limits complete the quotient even
where it cannot be used directly.

For $K>1$, the binomial choice provides
an explicit realization, while
the inverse equation admits other
families of probabilities.
The following lemma describes all its
strictly positive solutions
at a point where $f'\neq0$.

\begin{lemma}[Probabilities compatible with the inverse equation]
\label{lem:inverse_probability_freedom}
Let $f$ be continuous and concave on $[0,K]$,
with conditions
\eqref{eq:inverse_realizability_conditions}--\eqref{eq:style_s7_040},
and fix $0<\alpha<K$ with
$f'(\alpha)\neq0$.
Once the intermediate probabilities
$\pi_j(\alpha)>0$, $1\leq j\leq K-1$,
are chosen, normalization and
\eqref{eq:inverse_equation} require
\begin{equation}
\pi_K(\alpha)=\frac{e^{f(\alpha)-\alpha f'(\alpha)}-1
-\displaystyle\sum_{j=1}^{K-1}\pi_j(\alpha)\bigl(e^{-jf'(\alpha)}-1\bigr)}{e^{-Kf'(\alpha)}-1},
\label{eq:inverse_general_probabilities}
\end{equation}
\begin{equation}
\pi_0(\alpha)=1-\sum_{j=1}^{K-1}\pi_j(\alpha)-\pi_K(\alpha).
\label{eq:style_s7_044}
\end{equation}
These two probabilities are strictly
positive if and only if
\begin{multline}
\min\bigl\{1,e^{-Kf'(\alpha)}\bigr\}\Bigl(1-\sum_{j=1}^{K-1}\pi_j(\alpha)\Bigr)\\
<e^{f(\alpha)-\alpha f'(\alpha)}-\sum_{j=1}^{K-1}\pi_j(\alpha)e^{-jf'(\alpha)}\\
<\max\bigl\{1,e^{-Kf'(\alpha)}\bigr\}\Bigl(1-\sum_{j=1}^{K-1}\pi_j(\alpha)\Bigr).
\label{eq:inverse_probability_positivity}
\end{multline}
The set of positive choices is
nonempty and has dimension $K-1$
in the affine set determined
by the two constraints.
For $K=1$ the empty sums are zero
and the solution is unique.
\end{lemma}

\begin{proof}
Normalization expresses
$\pi_0$ in terms of the other probabilities.
Substituting into
\eqref{eq:inverse_equation} yields
\eqref{eq:inverse_general_probabilities};
the denominator is nonzero because
$f'(\alpha)\neq0$.
The contribution of the two endpoint
probabilities to the exponential sum is
\begin{equation}
\pi_0(\alpha)+\pi_K(\alpha)e^{-Kf'(\alpha)}
=e^{f(\alpha)-\alpha f'(\alpha)}-\sum_{j=1}^{K-1}\pi_j(\alpha)e^{-jf'(\alpha)},
\label{eq:inverse_residual_constraints}
\end{equation}
\begin{equation}
\pi_0(\alpha)+\pi_K(\alpha)=1-\sum_{j=1}^{K-1}\pi_j(\alpha).
\label{eq:style_s7_045}
\end{equation}
The double inequality \eqref{eq:inverse_probability_positivity}
already implies that the residual mass is positive: otherwise its
lower bound would be greater than or equal to its upper bound.
For a positive residual mass, the exponential moment is a combination
with strictly positive coefficients
of $1$ and $e^{-K f'(\alpha)}$ if and
only if it lies strictly between
the two values indicated in
\eqref{eq:inverse_probability_positivity}.
This proves the equivalence.

The constraints on the sum and on the
exponential moment are linearly
independent: the coefficients
$e^{-k f'(\alpha)}$ are not all
equal. They therefore define
an affine set of dimension
$(K+1)-2=K-1$.
The construction
\eqref{eq:inverse_binomial_probabilities}
provides a point with all
components positive, and positivity
remains valid in a neighborhood of it
relative to that affine set.
The stated nonemptiness and
dimension follow.
For $K=1$, the formula coincides with
\eqref{eq:inverse_binary_formula}.
\end{proof}

The pointwise choices in Lemma~\ref{lem:inverse_probability_freedom}
must extend to continuous, strictly positive functions
on the entire compact interval, including at the boundaries and where $f'=0$.
On the portion where $f=f'=0$, the inverse equation reduces to normalization.
If a point $c$ of this portion is approached by states with $f'\neq0$, dividing
the inverse equation minus one by $-f'(\alpha)$ and using
\eqref{eq:inverse_zero_slope_ratio}, continuity forces
$\bar\pi(c)=c$. This condition is not imposed in the interior
of the interval of maxima.

For $K>1$, a global perturbation of the binomial probabilities
explicitly provides different urns with the same entropy. Set
\begin{equation}
r(\alpha):=\frac1{1+e^{-f'(\alpha)}},\qquad0<\alpha<K,
\label{eq:display_022}
\end{equation}
\begin{equation}
r(0)=1,
\label{eq:style_s7_046}
\end{equation}
\begin{equation}
r(K)=0.
\label{eq:style_s7_047}
\end{equation}
The slope limits make $r$ continuous on $[0,K]$;
where $f'=0$, we have $r=1/2$. To components $0,1,2$, distinct
for $K>1$, we add the variations
\begin{equation}
\bigl(\delta\pi_0,\delta\pi_1,\delta\pi_2\bigr)=\varepsilon\bigl(-(1-r),\,1,\,-r\bigr),
\label{eq:inverse_global_perturbation}
\end{equation}
\begin{equation}
\delta\pi_k=0,\qquad k\geq3.
\label{eq:style_s7_048}
\end{equation}
Their sum is zero. To verify the exponential moment,
write locally $q=e^{-f'(\alpha)}$; then
\begin{equation}
\sum_{k\in\Omega}q^k\delta\pi_k
=\varepsilon[-(1-r)+q-rq^2]=0,
\label{eq:display_023}
\end{equation}
even where $f'=0$. Each variation has absolute value at most
$|\varepsilon|$. Choosing
\begin{equation}
|\varepsilon|<\min_{k\in\Omega}\min_{\alpha\in[0,K]}\pi_k(\alpha),
\label{eq:display_024}
\end{equation}
with the positive minimum computed on the binomial vector before the
perturbation preserves continuity and positivity on the entire
compact interval. Moreover, $\delta\pi_0(0)=0$; for $K=2$ we have
$\delta\pi_K(K)=-\varepsilon r(K)=0$, while for $K>2$ the
component $K$ does not change. Theorem~\ref{thm:inverse_verification}
again identifies the entropy with $f$. Distinct parameters give
distinct urns because component $1$ changes by $\varepsilon$.

The dimension $K-1$ in the lemma concerns the choices at each point
with $f'\neq0$, not the dimension of the space of global functions.
For $K=1$ the solution is unique at such points, whereas in the interior of
a nondegenerate interval of maxima some freedom may remain:
the two probabilities can be continuously perturbed with
opposite variations, compactly supported in that interior and small
enough to preserve positivity. The inverse equation remains valid,
and the verification theorem again gives the same entropy.

The next section will study
the generating function directly
for affine urn functions.

\section{Affine urn functions}
\label{sec:affine_urns}

For affine urn functions, the equation for the cumulant generating function becomes
a linear equation, following the method of
Ref.~\cite{Franchini_URNS}, Subsection~3.3.2,
Eq.~(3.109)--(3.111). For $K>1$, the coefficient of the derivative may
also vanish away from the origin: we first determine the matching conditions and
uniqueness, then the quadratures and the regularity of the entropy.

\subsection{Definition and equation for the generating function}
\label{subsec:affine_definition_equation}

We consider
\begin{equation}
\pi_k(\alpha)=a_k+b_k\alpha,
\qquad k\in\Omega,\qquad 0\leq\alpha\leq K,
\label{eq:affine_urn_functions}
\end{equation}
with real coefficients such that
\begin{equation}
\sum_{k\in\Omega}a_k=1,
\label{eq:affine_admissibility}
\end{equation}
\begin{equation}
\sum_{k\in\Omega}b_k=0,
\label{eq:style_s8_slope_normalization}
\end{equation}
\begin{equation}
a_k>0,\qquad k\in\Omega,
\label{eq:style_s8_left_probabilities_positive}
\end{equation}
\begin{equation}
a_k+Kb_k>0,\qquad k\in\Omega.
\label{eq:style_s8_right_probabilities_positive}
\end{equation}
These are precisely normalization and strict positivity on $[0,K]$:
each probability linearly interpolates its values at the two boundaries.
Set
\begin{equation}
a:=\sum_{k\in\Omega}k a_k,
\label{eq:affine_mean_coefficients}
\end{equation}
\begin{equation}
b:=\sum_{k\in\Omega}k b_k.
\label{eq:style_s8_mean_slope}
\end{equation}
The mean $\bar\pi(\alpha)=a+b\alpha$ satisfies
\begin{equation}
0<a<K,
\label{eq:affine_mean_constraints}
\end{equation}
\begin{equation}
0<a+Kb<K,
\label{eq:style_s8_right_mean_bounds}
\end{equation}
\begin{equation}
-\frac aK<b<1-\frac aK.
\label{eq:style_s8_mean_slope_bounds}
\end{equation}
Thus $-1<b<1$ and
\begin{equation}
C=\left\{\frac{a}{1-b}\right\}\subset(0,K).
\label{eq:affine_fixed_point}
\end{equation}
Define
\begin{equation}
A(\lambda):=\sum_{k\in\Omega}a_ke^{k\lambda},
\label{eq:affine_exponential_sums}
\end{equation}
\begin{equation}
B(\lambda):=\sum_{k\in\Omega}b_ke^{k\lambda}.
\label{eq:style_s8_exponential_slope_sum}
\end{equation}
We have $A>0$ on $\mathbb R$ and
\begin{equation}
A(0)=1,
\label{eq:affine_sums_at_zero}
\end{equation}
\begin{equation}
B(0)=0,
\label{eq:style_s8_slope_sum_at_zero}
\end{equation}
\begin{equation}
A'(0)=a,
\label{eq:style_s8_sum_derivative_at_zero}
\end{equation}
\begin{equation}
B'(0)=b.
\label{eq:style_s8_slope_sum_derivative_at_zero}
\end{equation}
The identity
\begin{equation}
\sum_{k\in\Omega}\pi_k(\alpha)e^{k\lambda}
=A(\lambda)+\alpha B(\lambda)
\label{eq:affine_exponential_decomposition}
\end{equation}
allows us to apply the limiting equation of
Theorem~\ref{thm:cumulant_limit_equation} directly.

\begin{lemma}[Zeros of the coefficient of the derivative]
\label{lem:affine_coefficient_zeros}
Under \eqref{eq:affine_admissibility}--\eqref{eq:style_s8_right_probabilities_positive},
we have $B\equiv0$ if and only if
$b_k=0$ for every $k\in\Omega$.
Otherwise, $B$ has at most $K$ real
zeros, counting multiplicities,
and one of them is zero.
For every $\lambda_0$ such that
$B(\lambda_0)=0$, we have
\begin{equation}
A(\lambda_0)-B'(\lambda_0)
>
\frac{A'(\lambda_0)}{K}>0.
\label{eq:affine_root_denominator}
\end{equation}
In particular,
\begin{equation}
0<
\frac{A'(\lambda_0)}
{A(\lambda_0)-B'(\lambda_0)}
<K.
\label{eq:affine_root_slope_bounds}
\end{equation}
\end{lemma}

\begin{proof}
The substitution $x=e^\lambda$ gives the polynomial representation
\begin{equation}
B(\log x)=\sum_{k\in\Omega}b_kx^k,\qquad x>0,
\label{eq:style_s8_polynomial_representation}
\end{equation}
which yields the alternative
$B\equiv0$ and the bound on the number of zeros, counting multiplicities.
At a zero $\lambda_0$, positivity of $a_1$ gives $A'(\lambda_0)>0$
and $B(\lambda_0)=0$ gives the identity
\begin{equation}
K\bigl\{A(\lambda_0)-B'(\lambda_0)\bigr\}-A'(\lambda_0)
=\sum_{k\in\Omega}(K-k)(a_k+Kb_k)e^{k\lambda_0}>0.
\label{eq:display_025}
\end{equation}
The sum is strictly positive already because of the term $k=0$.
Equations \eqref{eq:affine_root_denominator} and
\eqref{eq:affine_root_slope_bounds} follow.
\end{proof}

In the nonconstant binary case,
$B(\lambda)=b(e^\lambda-1)$,
so the only zero is the origin.
For $K>1$, other zeros may
appear. For example, for $K=2$,
\begin{equation}
\Bigl(a_0,a_1,a_2\Bigr)=\Bigl(\frac13,\frac13,\frac13\Bigr),
\label{eq:affine_extra_root_example}
\end{equation}
\begin{equation}
\Bigl(b_0,b_1,b_2\Bigr)=\Bigl(\frac1{12},-\frac18,\frac1{24}\Bigr)
\label{eq:style_s8_extra_root_slopes}
\end{equation}
satisfy
\eqref{eq:affine_admissibility}--\eqref{eq:style_s8_right_probabilities_positive}
and give
\begin{equation}
B(\lambda)
=\frac1{24}(e^\lambda-1)(e^\lambda-2).
\label{eq:affine_extra_root_factorization}
\end{equation}
Thus the zero
$\lambda=\log2$ also appears.
Moreover, the scalar condition
$b=0$ only makes $B'(0)$ vanish:
for $K>1$, it does not imply that all
the coefficients $b_k$ vanish.
If $B\not\equiv0$ and $b=0$,
the origin is a multiple zero.

\begin{theorem}[Equation for the cumulant generating function of affine urns]
\label{thm:affine_cumulant_equation}
Under \eqref{eq:affine_admissibility}--\eqref{eq:style_s8_right_probabilities_positive},
the cumulant generating function
is of class $C^1$ on $\mathbb{R}$
and satisfies
\begin{equation}
e^{\zeta(\lambda)}
=A(\lambda)+B(\lambda)\zeta'(\lambda),
\qquad \lambda\in\mathbb{R}.
\label{eq:affine_cumulant_equation}
\end{equation}
It is strictly convex, its
derivative is continuous and strictly
increasing, and
\begin{equation}
\{\zeta'(\lambda):\lambda\in\mathbb{R}\}
=(0,K).
\label{eq:affine_derivative_image}
\end{equation}
At every zero $\lambda_0$ of $B$,
we have
\begin{equation}
\zeta(\lambda_0)=\log A(\lambda_0),
\label{eq:affine_root_matching}
\end{equation}
\begin{equation}
\zeta'(\lambda_0)=\frac{A'(\lambda_0)}{A(\lambda_0)-B'(\lambda_0)}.
\label{eq:style_s8_root_matching_derivative}
\end{equation}
In particular,
\begin{equation}
\zeta(0)=0,
\label{eq:affine_cumulant_at_zero}
\end{equation}
\begin{equation}
\zeta'(0)=\frac{a}{1-b}.
\label{eq:style_s8_cumulant_derivative_at_zero}
\end{equation}
The boundary conditions are
\begin{equation}
\lim_{\lambda\to-\infty}\zeta(\lambda)=\log a_0,
\label{eq:affine_cumulant_boundaries}
\end{equation}
\begin{equation}
\lim_{\lambda\to+\infty}\bigl\{\zeta(\lambda)-K\lambda\bigr\}=\log(a_K+Kb_K).
\label{eq:style_s8_cumulant_right_boundary}
\end{equation}
The function $\zeta$ is the unique
solution of class $C^1(\mathbb{R})$
of \eqref{eq:affine_cumulant_equation}
with the conditions
\eqref{eq:affine_cumulant_boundaries}--\eqref{eq:style_s8_cumulant_right_boundary}.
If $b_k=0$ for every $k$,
the solution is
\begin{equation}
\zeta(\lambda)=\log A(\lambda).
\label{eq:affine_constant_probabilities}
\end{equation}
\end{theorem}

\begin{proof}
The general equation for the one-sided derivatives gives
\begin{equation}
e^{\zeta(\lambda)}=A(\lambda)+B(\lambda)\zeta'_-(\lambda),
\label{eq:affine_one_sided_equations}
\end{equation}
\begin{equation}
e^{\zeta(\lambda)}=A(\lambda)+B(\lambda)\zeta'_+(\lambda).
\label{eq:style_s8_cumulant_right_derivative_equation}
\end{equation}
If $B\equiv0$, we have $\zeta=\log A$, and
\eqref{eq:affine_root_matching}--\eqref{eq:style_s8_root_matching_derivative} holds everywhere with $B'=0$.
Otherwise, the one-sided derivatives coincide where $B\neq0$.
Fix a zero $\lambda_0$ and a neighborhood containing no other zeros.
Its value is $e^{\zeta(\lambda_0)}=A(\lambda_0)$; for sufficiently small $h\neq0$,
\begin{equation}
\frac{e^{\zeta(\lambda_0+h)}-e^{\zeta(\lambda_0)}}{h}
=\frac{A(\lambda_0+h)-A(\lambda_0)}{h}
 +\frac{B(\lambda_0+h)}{h}\zeta'(\lambda_0+h).
\label{eq:affine_root_difference_quotient}
\end{equation}
By convexity, $\zeta'(\lambda_0+h)\to\zeta'_\pm(\lambda_0)$
as $h\to0^\pm$; the limit of the left-hand side is
$A(\lambda_0)\zeta'_\pm(\lambda_0)$. Hence
\begin{equation}
\bigl\{A(\lambda_0)-B'(\lambda_0)\bigr\}\zeta'_\pm(\lambda_0)
=A'(\lambda_0).
\label{eq:affine_root_one_sided_slopes}
\end{equation}
The denominator is positive by
\eqref{eq:affine_root_denominator}; the two slopes coincide and
satisfy \eqref{eq:style_s8_root_matching_derivative}, including at multiple zeros.
No second derivative has been used, nor has there been division by $B'(\lambda_0)$.

Convexity and differentiability everywhere imply
$\zeta\in C^1(\mathbb R)$. The strict convexity in
Theorem~\ref{thm:cumulant_limit_equation} makes $\zeta'$ strictly
increasing; the limits of the slopes in
Lemma~\ref{lem:cumulant_boundary_values} give
\eqref{eq:affine_derivative_image}. The same lemma gives the asymptotes
\eqref{eq:affine_cumulant_boundaries}--\eqref{eq:style_s8_cumulant_right_boundary}, since
$\pi_0(0)=a_0$ and $\pi_K(K)=a_K+Kb_K$.
At the origin, \eqref{eq:affine_sums_at_zero}--\eqref{eq:style_s8_slope_sum_derivative_at_zero} gives
\eqref{eq:affine_cumulant_at_zero}--\eqref{eq:style_s8_cumulant_derivative_at_zero}.

For the quadratures, we introduce the transformation:
\begin{equation}
y(\lambda):=e^{-\zeta(\lambda)}-1.
\label{eq:affine_linearizing_function}
\end{equation}
We have $y\in C^1(\mathbb R)$, $1+y>0$, and the equation for the cumulant generating function
is equivalent to
\begin{equation}
B(\lambda)y'(\lambda)=A(\lambda)y(\lambda)+A(\lambda)-1.
\label{eq:affine_linear_equation}
\end{equation}
Where $B\neq0$, this becomes
\begin{equation}
y'(\lambda)-\frac{A(\lambda)}{B(\lambda)}y(\lambda)
=\frac{A(\lambda)-1}{B(\lambda)}.
\label{eq:affine_explicit_linear_equation}
\end{equation}
For uniqueness, let $f_1,f_2\in C^1(\mathbb R)$ be two real-valued candidates
satisfying the equation and the two stated asymptotes. The difference
$w=f_1-f_2$ tends to zero at both infinities, since on the right the
terms $K\lambda$ cancel. If it were nonzero, it would attain
a nonzero interior extremum $\lambda_*$, with $w'(\lambda_*)=0$.
Subtracting the equations at that point would give
\begin{equation}
e^{f_1(\lambda_*)}-e^{f_2(\lambda_*)}
=B(\lambda_*)w'(\lambda_*)=0,
\label{eq:display_026}
\end{equation}
contradicting $w(\lambda_*)\neq0$. Thus $f_1=f_2$. The comparison does not
require convexity or constraints on the slopes of the candidates and does not
divide by $B$: it includes its multiple zeros and the case $B\equiv0$.
\end{proof}

For $K=1$, the coefficients are $a_1=a$, $b_1=b$, $a_0=1-a$,
$b_0=-b$, with $0<a<1$ and $0<a+b<1$. Thus
\begin{equation}
A(\lambda)=1+a(e^\lambda-1),
\label{eq:affine_binary_sums}
\end{equation}
\begin{equation}
B(\lambda)=b(e^\lambda-1),
\label{eq:style_s8_binary_slope_sum}
\end{equation}
and, for $\lambda\neq0$,
\begin{equation}
a+b\zeta'(\lambda)
=\frac{e^{\zeta(\lambda)}-1}{e^\lambda-1}.
\label{eq:affine_binary_cumulant_equation}
\end{equation}
This binary form is also the equation in
Ref.~\cite{FB}, Section~VIII, Eq.~(90)--(92).

For the variable $y$ in \eqref{eq:affine_linearizing_function},
the matching data of the theorem become
\begin{equation}
y(\lambda_0)=\frac{1}{A(\lambda_0)}-1,\qquad B(\lambda_0)=0,
\label{eq:affine_linear_root_data}
\end{equation}
\begin{equation}
y'(\lambda_0)=-\frac{A'(\lambda_0)}{A(\lambda_0)\bigl\{A(\lambda_0)-B'(\lambda_0)\bigr\}},\qquad B(\lambda_0)=0.
\label{eq:style_s8_linear_root_derivative}
\end{equation}
At the origin, they are $y(0)=0$ and $y'(0)=-a/(1-b)$.
At the endpoints of the real line, we have
\begin{equation}
\lim_{\lambda\to-\infty}y(\lambda)=\frac{1}{a_0}-1,
\label{eq:affine_linear_boundary_data}
\end{equation}
\begin{equation}
\lim_{\lambda\to+\infty}y(\lambda)=-1,
\label{eq:style_s8_linear_right_limit}
\end{equation}
\begin{equation}
\lim_{\lambda\to+\infty}e^{K\lambda}\bigl\{1+y(\lambda)\bigr\}=\frac{1}{a_K+Kb_K}.
\label{eq:style_s8_linear_right_asymptote}
\end{equation}
The quadratures in the next subsection will use
\eqref{eq:affine_explicit_linear_equation} and these data.

\subsection{Integral representation}
\label{subsec:affine_integral_representation}

The integrating factor of \eqref{eq:affine_explicit_linear_equation}
selects the right endpoint when $B>0$ and the left endpoint when
$B<0$. Boundedness near that endpoint determines the solution
on each component of $\{B\neq0\}$.

\begin{theorem}[Integral representation for affine urns]
\label{thm:affine_integral_representation}
Assume
\eqref{eq:affine_admissibility}--\eqref{eq:style_s8_right_probabilities_positive}.
If $B\not\equiv0$, let
$J=(\lambda_-,\lambda_+)$
be a maximal interval on which
$B$ does not vanish, with
$-\infty\leq\lambda_-<\lambda_+\leq+\infty$.
If $B>0$ on $J$, for every
$\lambda\in J$ we have
\begin{equation}
1+y(\lambda)=e^{-\zeta(\lambda)}
=
\int_\lambda^{\lambda_+}
\frac{1}{B(t)}
\exp\left\{-\int_\lambda^t
\frac{A(u)}{B(u)}\,du\right\}\,dt.
\label{eq:affine_integral_positive_coefficient}
\end{equation}
If $B<0$ on $J$, we instead have
\begin{equation}
1+y(\lambda)=e^{-\zeta(\lambda)}
=
\int_{\lambda_-}^\lambda
\frac{1}{-B(t)}
\exp\left\{-\int_\lambda^t
\frac{A(u)}{B(u)}\,du\right\}\,dt.
\label{eq:affine_integral_negative_coefficient}
\end{equation}
The improper integrals are
absolutely convergent and
\begin{equation}
0<1+y(\lambda)\leq\frac{1}{a_0},
\qquad \lambda\in J.
\label{eq:affine_integral_positive_bound}
\end{equation}

The function $y$ thus
represented is the unique
solution of class $C^1(J)$
of
\eqref{eq:affine_explicit_linear_equation}
that remains bounded in a
one-sided neighborhood of the boundary
$\lambda_+$ if $B>0$, or
of the boundary $\lambda_-$ if $B<0$.
The formulas also hold
when that boundary is a multiple
zero of $B$ or is infinite.

Together with the values
\eqref{eq:affine_linear_root_data}--\eqref{eq:style_s8_linear_root_derivative}
at the zeros of $B$, these
representations determine
the cumulant generating function on the whole real line.
The matching is of class
$C^1$ and the conditions at the
endpoints are
\eqref{eq:affine_linear_boundary_data}--\eqref{eq:style_s8_linear_right_asymptote}.
If $B\equiv0$, the
representation is instead
the algebraic one
\begin{equation}
1+y(\lambda)=\frac{1}{A(\lambda)},
\label{eq:affine_integral_constant_case}
\end{equation}
\begin{equation}
\zeta(\lambda)=\log A(\lambda).
\label{eq:style_s8_constant_cumulant}
\end{equation}
\end{theorem}

\begin{proof}
If $B\equiv0$, \eqref{eq:affine_linear_equation} suffices.
Otherwise, the zeros of $B$ are finite in number. Fix
$\lambda\in J$ and set
\begin{equation}
F(t):=\exp\left\{-\int_\lambda^t\frac{A(u)}{B(u)}\,du\right\},
\qquad t\in J.
\label{eq:affine_local_integrating_factor}
\end{equation}
Then
\begin{equation}
F(\lambda)=1,
\label{eq:affine_integrating_factor_derivative}
\end{equation}
\begin{equation}
F(t)>0,
\label{eq:style_s8_integrating_factor_positive}
\end{equation}
\begin{equation}
F'(t)=-\frac{A(t)}{B(t)}F(t).
\label{eq:style_s8_integrating_factor_derivative}
\end{equation}
The selected endpoint is the right one if $B>0$, the left one
if $B<0$. The inequality
\begin{equation}
A(t)\geq a_0>0,
\label{eq:affine_integrating_factor_global_bound}
\end{equation}
\begin{equation}
|B(t)|\leq\left(\max_{k\in\Omega}\frac{|b_k|}{a_k}\right)A(t)
\label{eq:style_s8_integrating_factor_global_ratio_bound}
\end{equation}
shows that $A/|B|$ is bounded away from zero and has a divergent integral
toward every infinite endpoint. At a finite endpoint $\lambda_0$,
the regularity of $B$ gives, for some $c>0$,
\begin{equation}
|B(t)|\leq c|t-\lambda_0|,
\label{eq:affine_integrating_factor_finite_bound}
\end{equation}
\begin{equation}
\frac{A(t)}{|B(t)|}\geq\frac{a_0}{c|t-\lambda_0|}
\label{eq:style_s8_integrating_factor_finite_ratio_bound}
\end{equation}
in the relevant neighborhood, including for a multiple zero.
In the selected direction, the oriented integral of $A/B$ is positive
and diverges; hence
\begin{equation}
\lim_{t\uparrow\lambda_+}F(t)=0,\qquad B>0,
\label{eq:affine_integrating_factor_selected_limit}
\end{equation}
\begin{equation}
\lim_{t\downarrow\lambda_-}F(t)=0,\qquad B<0.
\label{eq:style_s8_integrating_factor_left_limit}
\end{equation}
Moreover,
\begin{equation}
\frac{F(t)}{|B(t)|}\leq\frac{|F'(t)|}{a_0},
\label{eq:affine_integral_absolute_bounds}
\end{equation}
\begin{equation}
\frac{F(t)|A(t)-1|}{|B(t)|}\leq\left(1+\frac1{a_0}\right)|F'(t)|.
\label{eq:style_s8_shifted_integral_absolute_bound}
\end{equation}
Since $F$ is monotone, its variation between $\lambda$ and
the selected endpoint is one. Equivalently,
\begin{equation}
\int_\lambda^{\lambda_+}\frac{A(t)F(t)}{B(t)}\,dt=1,\qquad B>0,
\label{eq:affine_integrating_factor_total_variation}
\end{equation}
\begin{equation}
\int_{\lambda_-}^{\lambda}\frac{A(t)F(t)}{-B(t)}\,dt=1,\qquad B<0.
\label{eq:style_s8_integrating_factor_left_variation}
\end{equation}
The preceding estimates prove absolute convergence of the integrals
for both variables $1+y$ and $y$.

Now let $v\in C^1(J)$ be any solution
of \eqref{eq:affine_explicit_linear_equation}, bounded near
the selected endpoint. The only product identity needed is
\begin{equation}
\Bigl\{F(t)\bigl[1+v(t)\bigr]\Bigr\}'=-\frac{F(t)}{B(t)}.
\label{eq:display_027}
\end{equation}
Integrating it between $\lambda$ and $t$ gives
\begin{equation}
1+v(\lambda)=F(t)\bigl[1+v(t)\bigr]+
\int_\lambda^t\frac{F(u)}{B(u)}\,du.
\label{eq:display_028}
\end{equation}
The product at the boundary tends to zero by boundedness of $v$.
Passing to the right endpoint for $B>0$, and to the left endpoint
for $B<0$, gives
\eqref{eq:affine_integral_positive_coefficient} and
\eqref{eq:affine_integral_negative_coefficient}, respectively, with $v$ in place of $y$.
The formula determines every candidate and also proves local uniqueness.
The integrand is positive; \eqref{eq:affine_integral_absolute_bounds}
and \eqref{eq:affine_integrating_factor_total_variation}--\eqref{eq:style_s8_integrating_factor_left_variation} give
\eqref{eq:affine_integral_positive_bound}.

The function $e^{-\zeta}-1$ from
Theorem~\ref{thm:affine_cumulant_equation} is such a candidate:
it is bounded at finite endpoints by continuity and at infinite ones
by \eqref{eq:affine_linear_boundary_data}--\eqref{eq:style_s8_linear_right_asymptote}. Thus the represented
solution exists, and coincides with this function throughout $J$.
The $C^1$ matching, including at the boundary opposite the selected one,
and the asymptotes follow from
\eqref{eq:affine_linear_root_data}--\eqref{eq:style_s8_linear_root_derivative} and
\eqref{eq:affine_linear_boundary_data}--\eqref{eq:style_s8_linear_right_asymptote}.
\end{proof}

Subtracting \eqref{eq:affine_integrating_factor_total_variation}--\eqref{eq:style_s8_integrating_factor_left_variation}
from the representations of $1+y$ directly gives the formulas
for $y$ in \eqref{eq:affine_linearizing_function}.
If $B>0$ on $J$,
\begin{equation}
y(\lambda)
=-\int_\lambda^{\lambda_+}
\frac{A(t)-1}{B(t)}
\exp\left\{-\int_\lambda^t
\frac{A(u)}{B(u)}\,du\right\}\,dt;
\label{eq:affine_shifted_integral_positive_coefficient}
\end{equation}
if $B<0$ on $J$,
\begin{equation}
y(\lambda)
=\int_{\lambda_-}^\lambda
\frac{A(t)-1}{B(t)}
\exp\left\{-\int_\lambda^t
\frac{A(u)}{B(u)}\,du\right\}\,dt.
\label{eq:affine_shifted_integral_negative_coefficient}
\end{equation}
The estimate
\eqref{eq:style_s8_shifted_integral_absolute_bound}
justifies their absolute
convergence.
The forms
\eqref{eq:affine_integral_positive_coefficient}--\eqref{eq:affine_integral_negative_coefficient}
make the positivity of $1+y$
immediately visible and
allow us to recover
$\zeta=-\log(1+y)$.
The procedure does not require
$b_0$ or $b_K$ to be
nonzero: the estimate
\eqref{eq:style_s8_integrating_factor_global_ratio_bound}
also covers these cases.
Strict positivity of the
probabilities in
\eqref{eq:style_s8_left_probabilities_positive}--\eqref{eq:style_s8_right_probabilities_positive}
always remains in force.

For $K=1$ and $b\neq0$,
\eqref{eq:affine_binary_sums}--\eqref{eq:style_s8_binary_slope_sum}
gives
\begin{equation}
\frac{A(\lambda)}{B(\lambda)}=\frac{a}{b}+\frac{1}{b(e^\lambda-1)},\qquad \lambda\neq0,
\label{eq:affine_binary_integrating_factor_coefficients}
\end{equation}
\begin{equation}
\frac{A(\lambda)-1}{B(\lambda)}=\frac{a}{b},\qquad \lambda\neq0.
\label{eq:style_s8_binary_inhomogeneous_coefficient}
\end{equation}
Up to a positive
multiplicative constant,
the integrating factor is
therefore
\begin{equation}
e^{-a\lambda/b}
|1-e^{-\lambda}|^{-1/b},
\qquad \lambda\neq0.
\label{eq:affine_binary_integrating_factor}
\end{equation}
For $\lambda>0$, the change
of variable
$z=1-e^{-\lambda}$
expresses it as
$(1-z)^{a/b}z^{-1/b}$,
which is the factor used in
Ref.~\cite{Franchini_URNS},
Eq.~(3.111).
If $b>0$, the selected boundary
is $+\infty$ on the positive
half-line and $-\infty$ on
the negative one; if $b<0$,
the origin is used from
both sides.
For $\lambda>0$, these
choices coincide with
those that determine
the constant in
Ref.~\cite{Franchini_URNS},
Eq.~(3.113) and~(3.116).
The evaluation of the integrals
in the binary case will
be taken up again in the section
devoted to examples.

\subsection{Entropy and regularity}
\label{subsec:affine_entropy_regularity}

The global regularity of the
cumulant generating function obtained in
Theorem~\ref{thm:affine_cumulant_equation}
allows us to apply
Lemma~\ref{lem:entropy_dual_equality}
to every interior terminal
state. The integral
representation in the preceding
subsection thus also determines
the entropy. We will then
study the second derivative
of the cumulant generating function at the points
where $B$ vanishes and
the corresponding regularity
of the entropy.

\begin{theorem}[Entropy for affine urns]
\label{thm:affine_entropy}
Under
\eqref{eq:affine_admissibility}--\eqref{eq:style_s8_right_probabilities_positive},
we have
\begin{equation}
\phi(\alpha)
=\inf_{\lambda\in\mathbb{R}}
\bigl\{\zeta(\lambda)-\lambda\alpha\bigr\},
\qquad 0\leq\alpha\leq K.
\label{eq:affine_entropy_duality}
\end{equation}
For $0<\alpha<K$, the infimum
is a unique minimum,
attained at the parameter
$\lambda$ characterized by
$\zeta'(\lambda)=\alpha$.
For $\alpha=0$ and $\alpha=K$,
the infimum is obtained in the limit
as $\lambda\to-\infty$ and
$\lambda\to+\infty$,
respectively, and is not
attained at a finite
parameter.

The function $\phi$ is
strictly concave and
absolutely continuous
on $[0,K]$, is of class
$C^1$ on $(0,K)$, and satisfies
\begin{equation}
\phi\bigl(\zeta'(\lambda)\bigr)=\zeta(\lambda)-\lambda\zeta'(\lambda),\qquad \lambda\in\mathbb{R},
\label{eq:affine_entropy_parametric_relations}
\end{equation}
\begin{equation}
\phi'\bigl(\zeta'(\lambda)\bigr)=-\lambda,\qquad \lambda\in\mathbb{R}.
\label{eq:style_s8_entropy_parametric_derivative}
\end{equation}
Its unique maximum is
\begin{equation}
\phi\left(\frac{a}{1-b}\right)=0.
\label{eq:affine_entropy_unique_maximum}
\end{equation}
It is strictly increasing
on $[0,a/(1-b)]$ and
strictly decreasing
on $[a/(1-b),K]$.
The boundary values and
slopes are
\begin{equation}
\phi(0)=\log a_0,
\label{eq:affine_entropy_boundary_values}
\end{equation}
\begin{equation}
\phi(K)=\log(a_K+Kb_K),
\label{eq:style_s8_entropy_right_boundary}
\end{equation}
\begin{equation}
\lim_{\alpha\downarrow0}\phi'(\alpha)=+\infty,
\label{eq:style_s8_entropy_left_slope}
\end{equation}
\begin{equation}
\lim_{\alpha\uparrow K}\phi'(\alpha)=-\infty.
\label{eq:style_s8_entropy_right_slope}
\end{equation}
Moreover,
\begin{equation}
\phi'\left(\frac{a}{1-b}\right)=0,
\label{eq:affine_entropy_derivative_integrability}
\end{equation}
\begin{equation}
\int_0^K|\phi'(\alpha)|\,d\alpha=-\log a_0-\log(a_K+Kb_K).
\label{eq:style_s8_entropy_total_variation}
\end{equation}
\end{theorem}

\begin{proof}
Theorem~\ref{thm:affine_cumulant_equation} gives
$\zeta\in C^1(\mathbb R)$ and a strictly increasing homeomorphism
$\zeta':\mathbb R\to(0,K)$. Thus
Lemma~\ref{lem:entropy_dual_equality} provides duality
at every interior state and at the physical boundaries, with the unique minimum and
the stated limits. Only after this identification does the infimum of
affine functions prove concavity of $\phi$.

We now reuse the secant argument in the proof of
Theorem~\ref{thm:entropy_shape},
\eqref{eq:entropy_secant_bounds}: continuity of the inverse of
$\zeta'$ identifies the derivative as $-\lambda$ and gives
\eqref{eq:style_s8_entropy_parametric_derivative}.
The argument requires only the duality and the image of $\zeta'$
just proved, and therefore also holds at the states corresponding
to the nonzero zeros of $B$, or when $B\equiv0$.
The derivative is continuous and strictly decreasing; together with continuity
at the boundaries, this gives strict concavity on $[0,K]$.

The data $\zeta(0)=0$, $\zeta'(0)=a/(1-b)$ identify the maximum and
the two monotonicity properties. The asymptotes
\eqref{eq:affine_cumulant_boundaries}--\eqref{eq:style_s8_cumulant_right_boundary} and the image $(0,K)$ give
\eqref{eq:affine_entropy_boundary_values}--\eqref{eq:style_s8_entropy_right_slope}.
For absolute continuity, the integral argument
in the same proof, \eqref{eq:entropy_integral_representation}--\eqref{eq:style_s7_024},
is applied to the two branches separated by $a/(1-b)$: we first integrate over
interior compact sets and pass to the boundaries using the sign of the derivative and
continuity of $\phi$. This yields
$\phi'\in L^1(0,K)$, the total variation
\eqref{eq:style_s8_entropy_total_variation}, and the reconstruction
\begin{equation}
\phi(\alpha)=\phi(0)+\int_0^\alpha\phi'(u)\,du,
\qquad 0\leq\alpha\leq K.
\label{eq:display_029}
\end{equation}
The two arguments from the general proof have been reused, without
applying its statement, which is subject to injectivity.
\end{proof}

For the entire positive
affine class, the theorem
gives the entropy properties
of Corollary~11 in
Ref.~\cite{Franchini_URNS}.
Here dual equality
follows from the global
image of $\zeta'$ in
Theorem~\ref{thm:affine_cumulant_equation};
the injectivity hypotheses
of
Theorem~\ref{thm:entropy_shape}
are unnecessary.

The substitution
$\lambda=-\phi'(\alpha)$
into the affine equation for the
cumulant generating function gives
\begin{equation}
e^{\phi(\alpha)-\alpha\phi'(\alpha)}
=A\bigl(-\phi'(\alpha)\bigr)
+\alpha B\bigl(-\phi'(\alpha)\bigr),
\qquad 0<\alpha<K.
\label{eq:affine_entropy_equation}
\end{equation}
This is the specialization of \eqref{eq:entropy_equation}, which also holds
at the states corresponding to the zeros of $B$.
The parametric identities just proved and continuity at the boundaries
also allow the same passage to the limit as in
\eqref{eq:entropy_endpoint_product_identities}--\eqref{eq:style_s7_027}:
$\alpha\phi'(\alpha)\to0$ as $\alpha\downarrow0$ and
$(K-\alpha)\phi'(\alpha)\to0$ as $\alpha\uparrow K$.

The formulas in Subsection~\ref{subsec:affine_integral_representation}
give the entropy through the parametrization
\begin{equation}
\alpha=\frac{1-A(\lambda)\bigl[1+y(\lambda)\bigr]}{B(\lambda)\bigl[1+y(\lambda)\bigr]},\qquad B(\lambda)\neq0,
\label{eq:affine_entropy_integral_parametrization}
\end{equation}
\begin{equation}
\phi(\alpha)=-\log\bigl[1+y(\lambda)\bigr]-\lambda\alpha,\qquad B(\lambda)\neq0.
\label{eq:style_s8_entropy_integral_value}
\end{equation}
The first identity is
\eqref{eq:affine_cumulant_equation}
solved for
$\zeta'(\lambda)$.
At the zeros of $B$,
we instead use
the matching data:
\begin{equation}
\alpha=\frac{A'(\lambda_0)}{A(\lambda_0)-B'(\lambda_0)},\qquad B(\lambda_0)=0,
\label{eq:affine_entropy_root_parametrization}
\end{equation}
\begin{equation}
\phi(\alpha)=\log A(\lambda_0)-\lambda_0\alpha,\qquad B(\lambda_0)=0.
\label{eq:style_s8_entropy_root_value}
\end{equation}
If $B\equiv0$, the formulas
become
\begin{equation}
\alpha=\frac{A'(\lambda)}{A(\lambda)},
\label{eq:affine_entropy_constant_parametrization}
\end{equation}
\begin{equation}
\phi(\alpha)=\log A(\lambda)-\lambda\alpha.
\label{eq:style_s8_entropy_constant_value}
\end{equation}
In every case, as
$\lambda$ varies over the real line,
the state $\alpha$ traverses
$(0,K)$ exactly once.

Regularity of the second
derivative requires examining
the equation near each
zero of $B$.
The following result
also provides the criterion
for the curvature of the
entropy to vanish.

\begin{theorem}[Second derivatives and regularity of the entropy]
\label{thm:affine_second_derivatives}
Under
\eqref{eq:affine_admissibility}--\eqref{eq:style_s8_right_probabilities_positive},
if $B\not\equiv0$, the
cumulant generating function is real analytic
on every interval on which
$B$ does not vanish, and
\begin{equation}
\zeta''(\lambda)>0,
\qquad B(\lambda)\neq0.
\label{eq:affine_strict_second_derivative}
\end{equation}
Let $\lambda_0$ be a zero
of $B$ and set
\begin{equation}
\alpha_0:=\zeta'(\lambda_0)
=\frac{A'(\lambda_0)}
{A(\lambda_0)-B'(\lambda_0)}.
\label{eq:affine_root_terminal_state}
\end{equation}
If $2B'(\lambda_0)<A(\lambda_0)$,
the function $\zeta$ is
of class $C^2$ in a
neighborhood of $\lambda_0$ and
\begin{equation}
\zeta''(\lambda_0)
=
\frac{
\displaystyle\sum_{k\in\Omega}
(k-\alpha_0)^2
(a_k+b_k\alpha_0)e^{k\lambda_0}
}{
A(\lambda_0)-2B'(\lambda_0)
}>0.
\label{eq:affine_root_finite_curvature}
\end{equation}
If $2B'(\lambda_0)\geq A(\lambda_0)$,
we instead have
\begin{equation}
\lim_{\substack{\lambda\to\lambda_0\\\lambda\neq\lambda_0}}\zeta''(\lambda)=+\infty,
\label{eq:affine_root_infinite_curvature}
\end{equation}
\begin{equation}
\lim_{h\to0}\frac{\zeta'(\lambda_0+h)-\zeta'(\lambda_0)}{h}=+\infty.
\label{eq:style_s8_root_infinite_difference_quotient}
\end{equation}
More precisely, if
$2B'(\lambda_0)=A(\lambda_0)$,
\begin{equation}
\zeta''(\lambda_0+h)
\sim
\frac{
\displaystyle\sum_{k\in\Omega}
(k-\alpha_0)^2
(a_k+b_k\alpha_0)e^{k\lambda_0}
}{B'(\lambda_0)}
\log\frac{1}{|h|},
\qquad h\to0^\pm.
\label{eq:affine_root_logarithmic_curvature}
\end{equation}
If $2B'(\lambda_0)>A(\lambda_0)$,
there exist constants
$c_+,c_->0$ such that
\begin{equation}
\zeta''(\lambda_0+h)
\sim c_\pm
|h|^{A(\lambda_0)/B'(\lambda_0)-2},
\qquad h\to0^\pm.
\label{eq:affine_root_power_curvature}
\end{equation}

For every affine urn in
the class considered,
the entropy belongs to
$C^2((0,K))$.
If $B\not\equiv0$, it is
real analytic on the
images, under $\zeta'$,
of the intervals containing
no zeros of $B$, and
\begin{equation}
\phi''\bigl(\zeta'(\lambda)\bigr)
=-\frac{1}{\zeta''(\lambda)}<0,
\qquad B(\lambda)\neq0.
\label{eq:affine_entropy_curvature_reciprocal}
\end{equation}
At the states in
\eqref{eq:affine_root_terminal_state},
we have
\begin{equation}
\phi''(\alpha_0)=-\frac{1}{\zeta''(\lambda_0)},\qquad 2B'(\lambda_0)<A(\lambda_0),
\label{eq:affine_entropy_root_curvature}
\end{equation}
\begin{equation}
\phi''(\alpha_0)=0,\qquad 2B'(\lambda_0)\geq A(\lambda_0).
\label{eq:style_s8_entropy_root_zero_curvature}
\end{equation}
If $B\equiv0$, the cumulant generating function
is real analytic on
$\mathbb{R}$, with
$\zeta''>0$, and the entropy
is real analytic on
$(0,K)$, with $\phi''<0$.
\end{theorem}

\begin{proof}
If $B\equiv0$, we have
$\zeta=\log A$ and
\begin{equation}
\zeta''(\lambda)
=
\frac{1}{A(\lambda)}
\sum_{k\in\Omega}
\left(k-\frac{A'(\lambda)}{A(\lambda)}\right)^2
a_ke^{k\lambda}>0.
\label{eq:affine_constant_curvature}
\end{equation}
Analyticity of $\zeta$
follows from that of $A>0$.
The analytic inverse of
$\zeta'$ and
\eqref{eq:style_s8_entropy_parametric_derivative}
then give analyticity
of $\phi$ and the relation
between the second derivatives.

Assume $B\not\equiv0$.
Where $B$ does not vanish,
the linear equation
\eqref{eq:affine_explicit_linear_equation}
has real analytic
coefficients. The local
integrating factor formula
shows that
$y$ is real analytic,
and the same holds for
$\zeta=-\log(1+y)$,
since $1+y>0$.
Set $\alpha(t)=\zeta'(t)$ and, where $B(t)\neq0$,
$w(t)=\zeta''(t)$. Define on the whole real line the continuous function
\begin{equation}
V(t):=\sum_{k\in\Omega}\bigl(k-\alpha(t)\bigr)^2
\bigl\{a_k+b_k\alpha(t)\bigr\}e^{kt}>0.
\label{eq:affine_local_curvature_forcing}
\end{equation}
Positivity follows from $\alpha(t)\in(0,K)$, positivity of the
probabilities, and the presence of at least two distinct increments.
Differentiating \eqref{eq:affine_cumulant_equation} only where
$B\neq0$ gives
\begin{equation}
Bw=(e^\zeta-B')\alpha-A'.
\label{eq:affine_differentiated_cumulant_equations}
\end{equation}
This identity and $e^\zeta=A+B\alpha$ give
\begin{equation}
A'+B'\alpha=e^\zeta\alpha-Bw,
\label{eq:display_030}
\end{equation}
\begin{equation}
V=A''+B''\alpha-e^\zeta\alpha^2+2B\alpha w.
\label{eq:style_s8_local_variance_identity}
\end{equation}
Differentiating the equation for the cumulant generating function once more gives
\begin{equation}
B(t)w'(t)=c(t)w(t)-V(t),\qquad B(t)\neq0,
\label{eq:affine_local_curvature_equation}
\end{equation}
\begin{equation}
c(t):=3e^{\zeta(t)}-2A(t)-2B'(t).
\label{eq:style_s8_local_curvature_coefficient}
\end{equation}
Indeed, $c=e^\zeta-2B'+2B\alpha$. Convexity gives $w\geq0$;
an interior zero of this function would be a minimum with $w'=0$, contradicting
\eqref{eq:affine_local_curvature_equation} and $V>0$. This proves
\eqref{eq:affine_strict_second_derivative}. The analytic inverse of
$\zeta'$ gives \eqref{eq:affine_entropy_curvature_reciprocal} and
analyticity of $\phi$ on the images of the intervals considered.

Now fix a zero $\lambda_0$ of $B$ and a neighborhood containing no other
zeros. The function $c$ is $C^1$ also at the zero, and
$c(\lambda_0)=A(\lambda_0)-2B'(\lambda_0)$;
$V(\lambda_0)>0$ is precisely the numerator of
\eqref{eq:affine_root_finite_curvature}. Moreover, $w>0$ is integrable
on each side: the fundamental theorem of calculus on interior
intervals and continuity of $\zeta'$ give, for sufficiently small $\varepsilon>0$,
\begin{equation}
\int_{\lambda_0}^{\lambda_0+\varepsilon}w(t)\,dt=\zeta'(\lambda_0+\varepsilon)-\alpha_0<\infty,
\label{eq:affine_local_curvature_integrability}
\end{equation}
\begin{equation}
\int_{\lambda_0-\varepsilon}^{\lambda_0}w(t)\,dt=\alpha_0-\zeta'(\lambda_0-\varepsilon)<\infty.
\label{eq:style_s8_local_curvature_left_integrability}
\end{equation}

On each side, choose a point $t_*$ and set
\begin{equation}
H(t):=\exp\left\{\int_{t_*}^t\frac{c(u)}{B(u)}\,du\right\}>0.
\label{eq:affine_curvature_homogeneous_factor}
\end{equation}
The equation for $w$ gives
\begin{equation}
\left(\frac wH\right)'=-\frac{V}{BH},
\label{eq:affine_curvature_ratio_derivatives}
\end{equation}
\begin{equation}
\left(\frac1H\right)'=-\frac{c}{BH}.
\label{eq:style_s8_curvature_inverse_factor_derivative}
\end{equation}
At every simple zero, the $C^1$ regularity of $c$ and analyticity
of $B$ give
\begin{equation}
\frac{c(t)}{B(t)}
=\frac{A(\lambda_0)/B'(\lambda_0)-2}{t-\lambda_0}+O(1).
\label{eq:affine_curvature_simple_root_coefficient}
\end{equation}
Thus the ratio
\begin{equation}
\frac{H(t)}{|t-\lambda_0|^{A(\lambda_0)/B'(\lambda_0)-2}}
\label{eq:style_s8_homogeneous_factor_power_ratio}
\end{equation}
tends, from each side, to a finite positive constant.
When $c(\lambda_0)\neq0$, we also have
\begin{equation}
\frac{(w/H)'(t)}{(1/H)'(t)}
=\frac{V(t)}{c(t)}
\longrightarrow
\frac{V(\lambda_0)}{A(\lambda_0)-2B'(\lambda_0)}.
\label{eq:affine_curvature_ratio_limit}
\end{equation}
In this case, $1/H$ is monotone in a sufficiently small one-sided neighborhood.
If \mbox{$H\to\infty$}, the identity
\begin{equation}
\left(\frac wH\right)'=\frac Vc\left(\frac1H\right)',
\label{eq:style_s8_curvature_ratio_identity}
\end{equation}
with $V/c$ bounded and $1/H$ of finite variation, shows that $w/H$ has a
finite nonnegative limit.

Assume $A(\lambda_0)-2B'(\lambda_0)>0$.
At simple zeros, \eqref{eq:affine_curvature_simple_root_coefficient}
gives $H\to0$ if $B'(\lambda_0)>0$, and
\begin{equation}
H(t)\asymp |t-\lambda_0|^{A(\lambda_0)/B'(\lambda_0)-2},\qquad B'(\lambda_0)<0,
\label{eq:style_s8_homogeneous_factor_negative_slope_power}
\end{equation}
with exponent less than $-2$ in the latter case.
At multiple zeros $B'(\lambda_0)=0$, instead,
\begin{equation}
\left|\frac{(t-\lambda_0)c(t)}{B(t)}\right|\longrightarrow\infty.
\label{eq:affine_curvature_multiple_root_coefficient}
\end{equation}
The quotient has constant sign on each side sufficiently close
to the zero, since $B$ is analytic and $c(\lambda_0)=A(\lambda_0)>0$.
Thus $H\to0$ or $H\to\infty$; in the latter case,
\begin{equation}
\frac{d\log H}{d\log|t-\lambda_0|}\leq-2
\label{eq:style_s8_homogeneous_factor_log_slope}
\end{equation}
near the zero and $H(t)\geq c_*|t-\lambda_0|^{-2}$ for some $c_*>0$.

If $H\to0$, the positive limit in
\eqref{eq:affine_curvature_ratio_limit} and integration of the
ratio identities give $w/H\to\infty$; l'Hôpital's rule
applied to
\begin{equation}
w=\frac{w/H}{1/H}
\label{eq:style_s8_curvature_ratio_form}
\end{equation}
determines its limit.
If \mbox{$H\to\infty$}, the nonnegative limit of $w/H$ must be
zero: a positive limit would make $w$ nonintegrable, by the
estimates just obtained, contradicting
\eqref{eq:affine_local_curvature_integrability}--\eqref{eq:style_s8_local_curvature_left_integrability}.
The same ratio is then of the form $0/0$, and l'Hôpital's rule gives
\begin{equation}
\lim_{\substack{t\to\lambda_0\\t\neq\lambda_0}}w(t)
=\frac{V(\lambda_0)}{A(\lambda_0)-2B'(\lambda_0)}>0.
\label{eq:affine_curvature_finite_limit}
\end{equation}
The mean value theorem applied to $\zeta'$, which is continuous also at the zero,
identifies this limit with $\zeta''(\lambda_0)$ and proves its
continuity. This includes all multiple zeros.

If $2B'(\lambda_0)=A(\lambda_0)$, the zero is simple and
\eqref{eq:affine_curvature_simple_root_coefficient} gives a
finite positive limit of $H$ on each side. Integrating
\eqref{eq:affine_curvature_ratio_derivatives}, with
$V(t)\to V(\lambda_0)>0$ and
$B(t)/(t-\lambda_0)\to B'(\lambda_0)>0$, gives
\begin{equation}
\frac{w(t)}{\log(1/|t-\lambda_0|)}
\longrightarrow\frac{V(\lambda_0)}{B'(\lambda_0)}.
\label{eq:affine_curvature_logarithmic_limit}
\end{equation}
The one-sided limit of $H$ cancels: this is the coefficient in
\eqref{eq:affine_root_logarithmic_curvature}.

If $2B'(\lambda_0)>A(\lambda_0)$,
Lemma~\ref{lem:affine_coefficient_zeros} gives
$B'(\lambda_0)<A(\lambda_0)$ and hence
\begin{equation}
-1<\frac{A(\lambda_0)}{B'(\lambda_0)}-2<0.
\label{eq:affine_curvature_power_range}
\end{equation}
The factor $H$ diverges with this power, and $w/H$ has a
finite nonnegative limit. If it were zero, l'Hôpital's rule and
\eqref{eq:affine_curvature_ratio_limit} would give a
negative limit of $w$, which is impossible. Thus the limit is positive from
both sides: this gives the constants $c_\pm>0$ in
\eqref{eq:affine_root_power_curvature}, without imposing their equality.
In the two divergent regimes, the mean value theorem applied to $\zeta'$ also gives
the limit of the difference quotient in
\eqref{eq:style_s8_root_infinite_difference_quotient}.

It remains to establish regularity
of the entropy at the
states $\alpha_0$.
For sufficiently small
$h\neq0$,
\eqref{eq:style_s8_entropy_parametric_derivative}
gives exactly
\begin{equation}
\frac{
\phi'\bigl(\zeta'(\lambda_0+h)\bigr)
-\phi'(\alpha_0)
}{
\zeta'(\lambda_0+h)-\alpha_0
}
=
-\frac{h}{
\zeta'(\lambda_0+h)-\zeta'(\lambda_0)
}.
\label{eq:affine_entropy_root_difference_quotient}
\end{equation}
Since $\zeta'$ is
continuous and strictly
increasing, these
increments traverse
a neighborhood of $\alpha_0$.
If the curvature of the
cumulant generating function has a finite
limit, the ratio
tends to
$-1/\zeta''(\lambda_0)$;
if it diverges, the ratio tends to
zero.
This proves existence
of the second derivative
of $\phi$ and
\eqref{eq:affine_entropy_root_curvature}--\eqref{eq:style_s8_entropy_root_zero_curvature}.
The formula
\eqref{eq:affine_entropy_curvature_reciprocal}
also shows that
$\phi''$ tends to the
same value as
one approaches
$\alpha_0$.
The zeros of $B$ are
finite in number:
the continuity thus
obtained, together with
that on the complementary
intervals, proves
$\phi\in C^2((0,K))$.
\end{proof}

Multiple zeros satisfy $B'=0$ and have finite curvature.
The curvature of the entropy vanishes only at the states of
simple zeros with $2B'\geq A$; these points are finite in number,
in agreement with strict concavity.
At the origin, $A(0)=1$ and $B'(0)=b$. For $b<1/2$,
\begin{equation}
\zeta''(0)=
\frac{\displaystyle\sum_{k\in\Omega}
\left(k-\frac{a}{1-b}\right)^2
\pi_k\left(\frac{a}{1-b}\right)}{1-2b}.
\label{eq:affine_origin_finite_curvature}
\end{equation}
For $b=1/2$, the logarithmic divergence of
\eqref{eq:affine_root_logarithmic_curvature} holds; for $b>1/2$, the
power $1/b-2$ in \eqref{eq:affine_root_power_curvature} holds.
This recovers the threshold discussed in Ref.~\cite{Franchini_URNS},
after Corollary~12, Eq.~(2.43)--(2.48).
Since $\phi\in C^2((0,K))$, at the maximum $\alpha_0=a/(1-b)$
we have
\begin{equation}
\phi(\alpha_0+h)=\frac12\phi''(\alpha_0)h^2+o(h^2).
\label{eq:style_s8_entropy_maximum_expansion}
\end{equation}
The quadratic coefficient is negative for $b<1/2$ and zero for
$b\geq1/2$. The theorem's criterion also holds for the nonzero
zeros of $B$ when $K>1$.

\section{Examples}
\label{sec:examples}

We first specialize Theorems~\ref{thm:intro_path_ldp}--\ref{thm:intro_inverse_entropy}
to arbitrary continuous, strictly positive urn functions for $K=1$ and
$K=2$. These give explicit local costs, variational formulas for the
terminal entropy, equations for the cumulant generating function,
and inverse constructions. The subsequent affine examples apply
Theorem~\ref{thm:intro_affine}: the binary case recovers the quadratures
of Corollary~12 in Ref.~\cite{Franchini_URNS} and specifies
regularity at the origin; the cases with three increments distinguish
the role of the mean from that of the individual probabilities.
Throughout, the process has the initialization in
Subsection~\ref{subsec:main_model}. For a different fixed initial
history, the uniform point estimates below are restricted to its
attainable terminal states, as in Theorem~\ref{thm:point_probabilities}.

The examples and figures in Section~\ref{sec:main_results} accompany
the main results they illustrate. We give here the formulas and proofs,
as well as the general specializations to $K=1$ and $K=2$.
The following construction also covers optimal paths with positive cost,
which are needed for the affine examples.

\begin{lemma}[Optimal paths for the affine examples]
\label{lem:affine_plot_paths}
Assume \eqref{eq:affine_admissibility}--\eqref{eq:style_s8_right_probabilities_positive}
and $\zeta\in C^2(\mathbb R)$.
Fix $\alpha\in(0,K)$, and let $\lambda_1$ be the unique real number
such that $\zeta'(\lambda_1)=\alpha$.
There is a unique solution on $(0,1]$ of
\begin{equation}
\tau\lambda_\alpha'(\tau)
=-B\bigl(\lambda_\alpha(\tau)\bigr)e^{-\zeta(\lambda_\alpha(\tau))},
\qquad\lambda_\alpha(1)=\lambda_1.
\label{eq:affine_optimal_tilt}
\end{equation}
It gives an optimal path in $Q(\alpha)$ by
\begin{equation}
\psi_\alpha^*(\tau)=\zeta'\bigl(\lambda_\alpha(\tau)\bigr),
\qquad\varphi_\alpha^*(\tau)=\tau\psi_\alpha^*(\tau),
\qquad0<\tau\leq1,
\label{eq:affine_optimal_mean}
\end{equation}
with $\varphi_\alpha^*(0)=0$ and $I(\varphi_\alpha^*)=-\phi(\alpha)$.
If $B(\lambda_1)=0$, the mean is constant; its cost is zero
if and only if $\lambda_1=0$.
\end{lemma}

\begin{proof}
Set $s=-\log\tau$. In this variable the equation has the locally
Lipschitz field $B e^{-\zeta}$, which is bounded: by
\eqref{eq:affine_cumulant_equation} and strict positivity,
its absolute value is at most $\max_k|b_k|/\varepsilon$, where
$\varepsilon=\min_{k,\beta}\pi_k(\beta)>0$.
Thus the solution exists for every $s\geq0$.
Define the local increment probabilities along the proposed path by
\begin{equation}
r_k(\tau)=\pi_k\bigl(\psi_\alpha^*(\tau)\bigr)
e^{k\lambda_\alpha(\tau)-\zeta(\lambda_\alpha(\tau))},
\qquad k\in\Omega.
\label{eq:affine_optimal_local_counts}
\end{equation}
The cumulant equation normalizes this vector. Differentiating that
equation and using \eqref{eq:affine_optimal_tilt} gives
\begin{equation}
\partial_\tau\varphi_\alpha^*(\tau)=\sum_{k\in\Omega}kr_k(\tau)\in[0,K].
\label{eq:affine_optimal_velocity}
\end{equation}
Since the mean lies in $(0,K)$, extension by zero at the initial time
therefore gives an admissible path. For $\beta,v\in[0,K]$, let $L(\beta,v)$ denote the supremum
on the right-hand side of \eqref{eq:local_frequency_duality}
with $\alpha=v$. The exponential tilt attains this supremum, so
\begin{equation}
L\bigl(\psi_\alpha^*,\partial_\tau\varphi_\alpha^*\bigr)
=\lambda_\alpha\partial_\tau\varphi_\alpha^*-
\zeta(\lambda_\alpha).
\label{eq:affine_optimal_local_cost}
\end{equation}
Affine duality identifies the right-hand side as
\begin{equation}
\lambda_\alpha\partial_\tau\varphi_\alpha^*-
\zeta(\lambda_\alpha)
=-\frac{d}{d\tau}\bigl\{\tau\phi(\psi_\alpha^*(\tau))\bigr\}.
\label{eq:affine_optimal_calibration}
\end{equation}
Integrating from a positive time to $1$ and letting that time decrease
to zero proves the cost identity, because $\phi$ is bounded on $[0,K]$.
The variational representation of the entropy proves optimality.
At a zero of $B$, uniqueness makes $\lambda_\alpha$ constant.
The unique zero of the entropy is the fixed point
$\zeta'(0)=a/(1-b)$, proving the last assertion.
\end{proof}

For the constant urn in Figure~\ref{fig:constant_urn_paths}, $B\equiv0$;
the optimal means are therefore constant for every terminal state.
The four nonconstant affine examples below all satisfy the finite-curvature
assumption of the lemma, by Theorem~\ref{thm:affine_second_derivatives}.

\subsection{Binary increments}
\label{subsec:binary_increments}

\subsubsection{General urn functions}
\label{subsec:binary_general}

Let $K=1$ and let $p:[0,1]\to(0,1)$ be continuous. Set
$\pi_1=p$ and $\pi_0=1-p$. Thus the next increment equals one
with probability $p(M_n/n)$ and zero otherwise. Theorem~\ref{thm:intro_path_ldp}
reduces to the Bernoulli relative entropy:
\begin{equation}
L_p(\alpha,v)=v\log\frac{v}{p(\alpha)}
+(1-v)\log\frac{1-v}{1-p(\alpha)},\qquad0\leq\alpha,v\leq1,
\label{eq:binary_general_action}
\end{equation}
\begin{equation}
I(\varphi)=\int_0^1L_p\!\left(\frac{\varphi(\tau)}{\tau},\partial_\tau\varphi(\tau)\right)d\tau.
\label{eq:style_s9_1_general_action}
\end{equation}
Here $\varphi\in Q$ and $0\log0=0$. For $0<v<1$, the maximizing parameter in
\eqref{eq:intro_path_rate} is
\begin{equation}
\lambda=\log\frac{v[1-p(\alpha)]}{(1-v)p(\alpha)};
\label{eq:style_s9_1_general_tilt}
\end{equation}
the endpoints follow by continuity. Unlike the case $K>1$,
the cumulative path determines all counts:
$\gamma_1(\tau)=\varphi(\tau)$ and
$\gamma_0(\tau)=\tau-\varphi(\tau)$, so no minimization over
different count vectors remains.

By Theorem~\ref{thm:intro_terminal_entropy},
\begin{equation}
\phi(\alpha)=-\min_{\varphi\in Q(\alpha)}I(\varphi),\qquad0\leq\alpha\leq1,
\label{eq:binary_general_entropy}
\end{equation}
\begin{equation}
\phi(0)=\log(1-p(0)),
\label{eq:style_s9_1_general_entropy_zero}
\end{equation}
\begin{equation}
\phi(1)=\log p(1),
\label{eq:style_s9_1_general_entropy_one}
\end{equation}
\begin{equation}
\lim_{N\to\infty}\max_{r\in\{0,\ldots,N\}}
\left|\frac1N\log\mathbb P(M_N=r)-\phi(r/N)\right|=0.
\label{eq:style_s9_1_general_point_probabilities}
\end{equation}
The entropy is continuous and nonpositive. Write
$C=\{\alpha\in[0,1]:p(\alpha)=\alpha\}$,
$c_- =\min C$ and $c_+=\max C$.
Theorem~\ref{thm:intro_zero_cost} gives
\begin{equation}
\phi(\alpha)=0\quad\Longleftrightarrow\quad c_-\leq\alpha\leq c_+,
\label{eq:binary_general_zero_cost}
\end{equation}
\begin{equation}
I(\varphi)=0\quad\Longleftrightarrow\quad
\partial_\tau\varphi(\tau)=p\!\left(\frac{\varphi(\tau)}{\tau}\right)
\qquad\text{for almost every }\tau\in(0,1).
\label{eq:style_s9_1_general_zero_cost_paths}
\end{equation}
Equivalently, $\psi=\varphi/\tau$ solves
$\tau\partial_\tau\psi=p(\psi)-\psi$.
These zero-cost means are monotone, remain in $[c_-,c_+]$ and
approach a fixed point as $\tau\downarrow0$.
At a fixed point $c\in C$, the complete criterion for uniqueness
of the straight optimal path $\varphi(\tau)=c\tau$ is given by
\eqref{eq:intro_fixed_point_arrival_left}--\eqref{eq:style_s2_intro_fixed_point_arrival_right_2},
with $\bar\pi=p$; no finiteness assumption on $C$ is needed.

Theorem~\ref{thm:intro_cumulant} applies without monotonicity
of $p$. In particular, the limiting cumulant generating function
is strictly convex, nondecreasing and $1$-Lipschitz, with $\zeta(0)=0$.
The convergence defining $\zeta$ is uniform on compact parameter intervals.
For every $\lambda\in\mathbb R$, it satisfies
\begin{equation}
\zeta(\lambda)=\max_{0\leq\alpha\leq1}\{\lambda\alpha+\phi(\alpha)\},
\label{eq:binary_general_cumulant}
\end{equation}
\begin{equation}
e^{\zeta(\lambda)}=1+(e^\lambda-1)p(\zeta'_-(\lambda)),
\label{eq:style_s9_1_cumulant_left_equation}
\end{equation}
\begin{equation}
e^{\zeta(\lambda)}=1+(e^\lambda-1)p(\zeta'_+(\lambda)),
\label{eq:style_s9_1_cumulant_right_equation}
\end{equation}
\begin{equation}
\zeta'_-(0)=c_-,
\label{eq:style_s9_1_cumulant_left_slope}
\end{equation}
\begin{equation}
\zeta'_+(0)=c_+,
\label{eq:style_s9_1_cumulant_right_slope}
\end{equation}
\begin{equation}
\lim_{\lambda\to-\infty}\zeta(\lambda)=\log(1-p(0)),
\label{eq:style_s9_1_cumulant_negative_limit}
\end{equation}
\begin{equation}
\lim_{\lambda\to+\infty}\bigl(\zeta(\lambda)-\lambda\bigr)=\log p(1).
\label{eq:style_s9_1_cumulant_positive_limit}
\end{equation}
For every $\lambda\neq0$, injectivity of the exponential sum in
Theorem~\ref{thm:intro_cumulant_duality} is exactly injectivity
of $p$ on the interval in question. If $p$ is injective on
$[0,c_-]$, $\zeta$ is the unique $C^1$ solution on $(-\infty,0)$ of
\begin{equation}
e^{g(\lambda)}=1+(e^\lambda-1)p(g'(\lambda)),
\qquad 0\leq g'(\lambda)\leq c_-.
\label{eq:binary_general_negative_problem}
\end{equation}
Its derivative ranges over $(0,c_-)$, and
\begin{equation}
\phi(\alpha)=\inf_{\lambda\leq0}\{\zeta(\lambda)-\lambda\alpha\},\qquad\alpha\in[0,c_-].
\label{eq:style_s9_1_negative_duality}
\end{equation}
Independently, injectivity of $p$ on $[c_+,1]$ gives the
corresponding unique solution on $(0,\infty)$, with the same
equation, slope constraint $c_+\leq g'\leq1$, derivative range
$(c_+,1)$, and dual formula with $\lambda\geq0$ on $[c_+,1]$.
On each half-line where its injectivity assumption holds, the equation
also gives the explicit derivative
\begin{equation}
\zeta'(\lambda)=
\begin{cases}
\displaystyle
\left(p\big|_{[0,c_-]}\right)^{-1}
\!\left(\frac{e^{\zeta(\lambda)}-1}{e^\lambda-1}\right),&\lambda<0,\\[3mm]
\displaystyle
\left(p\big|_{[c_+,1]}\right)^{-1}
\!\left(\frac{e^{\zeta(\lambda)}-1}{e^\lambda-1}\right),&\lambda>0.
\end{cases}
\label{eq:binary_general_cumulant_ode}
\end{equation}
The inverse in each row is that of the indicated restriction of $p$.
The slope constraints are imposed on the entire half-line;
the values at zero and infinity in
\eqref{eq:style_s9_1_cumulant_left_slope}--\eqref{eq:style_s9_1_cumulant_positive_limit}
follow from the equation, without assuming convexity of a candidate.
If both assumptions hold, then
\begin{equation}
\phi(\alpha)=\inf_{\lambda\in\mathbb R}\{\zeta(\lambda)-\lambda\alpha\},
\label{eq:binary_general_duality}
\end{equation}
\begin{equation}
\phi'(\zeta'(\lambda))=-\lambda,\qquad\lambda\neq0.
\label{eq:style_s9_1_general_dual_derivative}
\end{equation}
The entropy is then concave on $[0,1]$, absolutely continuous,
$C^1$ on $(0,1)$ and strictly concave on the two exterior intervals.
It is strictly concave on all of $[0,1]$ exactly when $c_-=c_+$.
These are sufficient injectivity assumptions, not restrictions
needed for \eqref{eq:binary_general_action}--\eqref{eq:style_s9_1_cumulant_positive_limit}.

Finally, for a prescribed continuous concave $f$ on $[0,1]$,
Theorem~\ref{thm:intro_inverse_entropy} gives a realization
exactly when $f\in C^1((0,1))$, $\max f=0$,
$f'(0+)=+\infty$ and $f'(1-)=-\infty$. One such realization is
\begin{equation}
p_f(\alpha)=\frac{e^{f(\alpha)-\alpha f'(\alpha)}-1}{e^{-f'(\alpha)}-1},
\qquad0<\alpha<1,\qquad f'(\alpha)\neq0,
\label{eq:binary_general_inverse}
\end{equation}
\begin{equation}
p_f(\alpha)=
\begin{cases}
\alpha,&0<\alpha<1,\qquad f'(\alpha)=0,\\[1mm]
1-e^{f(0)},&\alpha=0,\\[1mm]
e^{f(1)},&\alpha=1.
\end{cases}
\label{eq:style_s9_1_inverse_extensions}
\end{equation}
This recovers \eqref{eq:inverse_binary_formula} with its full
domain and boundary extensions. The inverse equation determines
$p$ uniquely where $f'\neq0$; on a nondegenerate interval where
$f=f'=0$, it reduces to normalization. Thus the value
$p_f(\alpha)=\alpha$ on that interval is a choice of realization,
not a general uniqueness conclusion. Continuity still forces
$p(c)=c$ at any zero of $f$ approached through states where $f'\neq0$.

\subsubsection{Affine urns}
\label{subsec:binary_affine_examples}

For $K=1$, we consider
\begin{equation}
\pi_1(\alpha)=a+b\alpha,\qquad0\leq\alpha\leq1,
\label{eq:binary_urn_functions}
\end{equation}
\begin{equation}
\pi_0(\alpha)=1-a-b\alpha,\qquad0\leq\alpha\leq1.
\label{eq:style_s9_1_affine_probability_zero}
\end{equation}
The strict positivity
conditions are
\begin{equation}
0<a<1,
\label{eq:binary_admissibility}
\end{equation}
\begin{equation}
0<a+b<1.
\label{eq:style_s9_1_affine_endpoint_admissibility}
\end{equation}
In particular, $-1<b<1$.
The sums $A,B$ are
those in
\eqref{eq:affine_binary_sums}--\eqref{eq:style_s8_binary_slope_sum},
and the unique fixed point
is $a/(1-b)$.

If $b=0$, we have
\begin{equation}
\zeta(\lambda)=\log(1-a+ae^\lambda).
\label{eq:binary_constant_cumulant}
\end{equation}
For $0<\alpha<1$,
the equation
$\zeta'(\lambda)=\alpha$
is solved as
\begin{equation}
\lambda
=\log\frac{\alpha(1-a)}{a(1-\alpha)}.
\label{eq:binary_constant_dual_parameter}
\end{equation}
Substitution into the dual
formula gives
\begin{equation}
\phi(\alpha)
=-\alpha\log\frac{\alpha}{a}
-(1-\alpha)\log\frac{1-\alpha}{1-a},
\qquad 0\leq\alpha\leq1.
\label{eq:binary_constant_entropy}
\end{equation}
The terms at the endpoints
are interpreted by
continuity, with the
convention $0\log0=0$.

In the case $b\neq0$,
the integral
representations of the preceding
section take
the following form.

\begin{theorem}[Integral formulas in the binary case]
\label{thm:binary_integral_formulas}
Under
\eqref{eq:binary_admissibility}--\eqref{eq:style_s9_1_affine_endpoint_admissibility},
assume $b\neq0$.
For $\lambda>0$ and $b>0$,
\begin{equation}
e^{-\zeta(\lambda)}
=1-\frac{a}{b}
e^{a\lambda/b}
(1-e^{-\lambda})^{1/b}\,
\int_{1-e^{-\lambda}}^1
(1-t)^{a/b-1}t^{-1/b}\,dt.
\label{eq:binary_positive_lambda_positive_b}
\end{equation}
For $\lambda>0$ and $b<0$,
\begin{equation}
e^{-\zeta(\lambda)}
=1+\frac{a}{b}
e^{a\lambda/b}
(1-e^{-\lambda})^{1/b}\,
\int_0^{1-e^{-\lambda}}
(1-t)^{a/b-1}t^{-1/b}\,dt.
\label{eq:binary_positive_lambda_negative_b}
\end{equation}
For $\lambda<0$ and $b>0$,
\begin{equation}
e^{-\zeta(\lambda)}
=1+\frac{a}{b}
e^{-(1-a)\lambda/b}
(1-e^\lambda)^{1/b}\,
\int_{1-e^\lambda}^1
(1-t)^{(1-a)/b-1}t^{-1/b}\,dt.
\label{eq:binary_negative_lambda_positive_b}
\end{equation}
For $\lambda<0$ and $b<0$,
\begin{equation}
e^{-\zeta(\lambda)}
=1-\frac{a}{b}
e^{-(1-a)\lambda/b}
(1-e^\lambda)^{1/b}\,
\int_0^{1-e^\lambda}
(1-t)^{(1-a)/b-1}t^{-1/b}\,dt.
\label{eq:binary_negative_lambda_negative_b}
\end{equation}
All the integrals
converge, and the expressions
for $e^{-\zeta}$ are
strictly positive.
The two branches join
with $C^1$ regularity and
\begin{equation}
\zeta(0)=0,
\label{eq:binary_origin_matching}
\end{equation}
\begin{equation}
\zeta'(0)=\frac{a}{1-b}.
\label{eq:style_s9_1_affine_origin_slope}
\end{equation}
At the endpoints of the real line,
we have
\begin{equation}
\lim_{\lambda\to-\infty}\zeta(\lambda)=\log(1-a),
\label{eq:binary_cumulant_boundary_limits}
\end{equation}
\begin{equation}
\lim_{\lambda\to+\infty}\bigl(\zeta(\lambda)-\lambda\bigr)=\log(a+b).
\label{eq:style_s9_1_affine_positive_limit}
\end{equation}
\end{theorem}

\begin{proof}
By
\eqref{eq:style_s8_binary_inhomogeneous_coefficient},
the equation for $y$ is
\begin{equation}
y'(\lambda)
-\frac{A(\lambda)}{B(\lambda)}y(\lambda)
=\frac{a}{b}.
\label{eq:binary_linear_equation}
\end{equation}
With $y=e^{-\zeta}-1$, we use the integrating factor
\eqref{eq:affine_binary_integrating_factor} and set
$z=1-e^{-\lambda}$ on the positive half-line, $z=1-e^\lambda$
on the negative one; in both cases, $0<z<1$.
On the negative branch, the factor can be rewritten as
\begin{equation}
e^{-a\lambda/b}
|1-e^{-\lambda}|^{-1/b}
=
e^{(1-a)\lambda/b}
(1-e^\lambda)^{-1/b},
\qquad \lambda<0.
\label{eq:binary_negative_integrating_factor}
\end{equation}
The chain rule gives the two identities
\begin{equation}
\frac d{dz}\left\{(1-z)^{a/b}z^{-1/b}
y\bigl(-\log(1-z)\bigr)\right\}
=\frac ab(1-z)^{a/b-1}z^{-1/b},
\label{eq:binary_positive_change_of_variable}
\end{equation}
\begin{equation}
\frac d{dz}\left\{(1-z)^{(1-a)/b}z^{-1/b}
y(\log(1-z))\right\}
=-\frac ab(1-z)^{(1-a)/b-1}z^{-1/b}.
\label{eq:binary_negative_change_of_variable}
\end{equation}
For $b>0$, Theorem~\ref{thm:affine_integral_representation}
selects $+\infty$ and $-\infty$, respectively, that is, $z=1$.
The two factors tend to zero because $a,1-a>0$, while $y$ remains
bounded by the data \eqref{eq:affine_cumulant_boundaries}--\eqref{eq:style_s8_cumulant_right_boundary}.
Their product with $y$ thus has limit zero; integrating from $z$
to $1$ gives \eqref{eq:binary_positive_lambda_positive_b}
and \eqref{eq:binary_negative_lambda_positive_b}, with the opposite signs
imposed by the two identities.
For $b<0$, the selected boundaries are $0^+$ and $0^-$, both $z=0$.
Here $z^{-1/b}\to0$, the other factors tend to one, and $y$ is bounded:
the product again has limit zero. Integrating from $0$ to $z$
gives \eqref{eq:binary_positive_lambda_negative_b}
and \eqref{eq:binary_negative_lambda_negative_b}.

The integrals converge: for $b>0$, the only improper endpoint is
$t=1$, with exponents $a/b-1>-1$ and $(1-a)/b-1>-1$; for $b<0$,
it is $t=0$, with exponent $-1/b>0$. The other endpoint lies inside
$(0,1)$ for every finite $\lambda\neq0$.
Identification with
Theorems~\ref{thm:affine_cumulant_equation}
and~\ref{thm:affine_integral_representation}, applicable under
\eqref{eq:binary_admissibility}--\eqref{eq:style_s9_1_affine_endpoint_admissibility}, gives positivity, $C^1$ matching,
and all the data
\eqref{eq:binary_origin_matching}--\eqref{eq:style_s9_1_affine_positive_limit}.
It remains to check the negative limit directly, which distinguishes
the prefactor of the corresponding branch.

For the limit at
$-\infty$, set
$z=1-e^\lambda$ in the
formulas for the negative branch.
If $b>0$, then as $z\uparrow1$,
\begin{equation}
\int_z^1(1-t)^{(1-a)/b-1}t^{-1/b}\,dt
\sim\frac{b}{1-a}(1-z)^{(1-a)/b}.
\label{eq:style_s9_1_positive_slope_integral_asymptotic}
\end{equation}
If $b<0$, then as $z\uparrow1$,
\begin{equation}
\int_0^z(1-t)^{(1-a)/b-1}t^{-1/b}\,dt
\sim-\frac{b}{1-a}(1-z)^{(1-a)/b}.
\label{eq:style_s9_1_negative_slope_integral_asymptotic}
\end{equation}
Multiplying by the
respective prefactor
gives
\begin{equation}
\lim_{\lambda\to-\infty}y(\lambda)=\frac{a}{1-a},
\label{eq:binary_negative_infinity_asymptotic}
\end{equation}
\begin{equation}
\lim_{\lambda\to-\infty}e^{-\zeta(\lambda)}=\frac1{1-a}.
\label{eq:style_s9_1_affine_exponential_negative_limit}
\end{equation}
This recovers the limit in
\eqref{eq:binary_cumulant_boundary_limits}.
\end{proof}

For $\lambda>0$, the
formulas
\eqref{eq:binary_positive_lambda_positive_b}
and
\eqref{eq:binary_positive_lambda_negative_b}
coincide with those
of Ref.~\cite{Franchini_URNS},
Corollary~12,
Eq.~(2.44)--(2.45),
in the form corrected
by Eq.~(1.1)--(1.2)
of the Corrigendum.
The first also recovers
Ref.~\cite{FB},
Section~VIII, Eq.~(93),
in version
arXiv:2210.12585v5,
with $a=1-p$,
$b=2p-1$ and $1/2<p<1$.

For $\lambda<0$, the
exponential prefactor
is $e^{-(1-a)\lambda/b}$.
In Eq.~(2.47)--(2.48)
of version
arXiv:1412.5762v7 of
Ref.~\cite{Franchini_URNS},
$e^{-(1-a+b)\lambda/b}$
appears instead.
The form obtained here
follows from
\eqref{eq:binary_negative_integrating_factor}
and satisfies
\eqref{eq:binary_negative_infinity_asymptotic}--\eqref{eq:style_s9_1_affine_exponential_negative_limit}.
For $b>0$, the additional
factor
$e^{-\lambda}$ in the
cited formula would make
the term that must tend
to $a/(1-a)$
diverge.

The entropy is computed
from the four formulas
using
\begin{equation}
\alpha=\frac1b\left\{\frac{e^{\zeta(\lambda)}-1}{e^\lambda-1}-a\right\},
\qquad\lambda\neq0,
\label{eq:binary_entropy_parametrization}
\end{equation}
\begin{equation}
\phi(\alpha)=\zeta(\lambda)-\lambda\alpha,\qquad\lambda\neq0.
\label{eq:style_s9_1_affine_entropy_parametrization}
\end{equation}
The parameter traverses
each branch of the entropy
exactly once.
The central value and
the endpoint values are
\begin{equation}
\phi\!\left(\frac{a}{1-b}\right)=0,
\label{eq:binary_entropy_values}
\end{equation}
\begin{equation}
\phi(0)=\log(1-a),
\label{eq:style_s9_1_affine_entropy_zero}
\end{equation}
\begin{equation}
\phi(1)=\log(a+b).
\label{eq:style_s9_1_affine_entropy_one}
\end{equation}

For $b<1/2$, the formula
\eqref{eq:affine_origin_finite_curvature}
reduces to
\begin{equation}
\zeta''(0)=\frac{a(1-a-b)}{(1-b)^2(1-2b)},
\label{eq:binary_origin_curvatures}
\end{equation}
\begin{equation}
\phi''\!\left(\frac{a}{1-b}\right)=-\frac{(1-b)^2(1-2b)}{a(1-a-b)}.
\label{eq:style_s9_1_affine_entropy_curvature}
\end{equation}
For $b=1/2$,
\begin{equation}
\zeta''(\lambda)
\sim4a(1-2a)
\log\frac{1}{|\lambda|},
\qquad \lambda\to0.
\label{eq:binary_critical_curvature}
\end{equation}
For $b>1/2$, the
divergence is that
of
\eqref{eq:affine_root_power_curvature},
with exponent $1/b-2$.
In both of these
latter cases,
$\phi''(a/(1-b))=0$,
as proved in
Theorem~\ref{thm:affine_second_derivatives}.

Regularity can
be specified
further when
$b<0$.

\begin{lemma}[Analyticity in the binary case with negative slope]
\label{lem:binary_negative_slope_analyticity}
Under
\eqref{eq:binary_admissibility}--\eqref{eq:style_s9_1_affine_endpoint_admissibility},
if $b<0$, the cumulant generating function
is real analytic on
$\mathbb{R}$ and the entropy
is real analytic on
$(0,1)$.
\end{lemma}

\begin{proof}
Theorem~\ref{thm:affine_second_derivatives} gives
analyticity away from the origin. To extend across zero,
we use $z=1-e^{-\lambda}$ and
$y=e^{-\zeta}-1$ on both sides. Equation \eqref{eq:binary_linear_equation} becomes
\begin{equation}
bz(1-z)\frac{d}{dz}y\bigl(-\log(1-z)\bigr)
=\{1+(a-1)z\}y\bigl(-\log(1-z)\bigr)+az.
\label{eq:binary_analytic_equation}
\end{equation}
The integrating factor is $(1-z)^{a/b}|z|^{-1/b}$.
For $b<0$, its product with $y$ tends to zero as
$z\to0^\pm$, since the factor tends to zero and $y$ is
bounded by the affine theorem. Integrating from the origin
and setting $t=zu$ gives a single formula:
\begin{equation}
y\bigl(-\log(1-z)\bigr)=
\frac ab\,z(1-z)^{-a/b}
\int_0^1(1-zu)^{a/b-1}u^{-1/b}\,du,
\qquad |z|<1.
\label{eq:binary_negative_slope_common_integral}
\end{equation}
For $z<0$, the factor $z$ preserves the orientation
of the integral; at zero, the formula holds by continuity.
For $|z|\leq1/2$, the binomial expansion converges
absolutely and uniformly in $u\in[0,1]$, while
$u^{-1/b}$ is integrable. Term-by-term integration
proves analyticity at zero of $y$ and, since $1+y(0)=1$,
of $\zeta=-\log(1+y)$.

Finally, $\zeta''>0$ away from the origin by
\eqref{eq:affine_strict_second_derivative}, and at the origin
by \eqref{eq:binary_origin_curvatures}. The analytic inverse
of $\zeta'$, whose image is $(0,1)$, and the duality
\eqref{eq:affine_entropy_parametric_relations} give
analyticity of $\phi$ on $(0,1)$.
\end{proof}

Finally, we consider
a case with positive
slope in which the
parametric inversion
is also elementary:
\begin{equation}
a=b=\frac13,
\label{eq:binary_third_parameters}
\end{equation}
\begin{equation}
\pi_1(\alpha)=\frac{1+\alpha}{3},
\label{eq:style_s9_1_third_probability_one}
\end{equation}
\begin{equation}
\pi_0(\alpha)=\frac{2-\alpha}{3}.
\label{eq:style_s9_1_third_probability_zero}
\end{equation}
The two probabilities are
strictly positive
on $[0,1]$.
On the two branches with $b>0$, the integrals reduce,
respectively, to
\begin{equation}
\int_z^1t^{-3}\,dt=\frac{1-z^2}{2z^2},
\label{eq:display_031}
\end{equation}
\begin{equation}
\int_z^1(1-t)t^{-3}\,dt=\frac{(1-z)^2}{2z^2}.
\label{eq:style_s9_1_third_negative_integral}
\end{equation}
With $z=1-e^{-\lambda}$ in the first and $z=1-e^\lambda$
in the second, the prefactors of the quadratures give
\begin{equation}
e^{-\zeta(\lambda)}
=
\begin{cases}
\displaystyle\frac{3-e^\lambda}{2},
&\lambda\leq0,\\[3mm]
\displaystyle\frac{3e^\lambda-1}{2e^{2\lambda}},
&\lambda\geq0.
\end{cases}
\label{eq:binary_third_cumulant_exponential}
\end{equation}
The expressions coincide at zero. Differentiating,
\begin{equation}
\zeta'(\lambda)
=
\begin{cases}
\displaystyle\frac{e^\lambda}{3-e^\lambda},
&\lambda\leq0,\\[3mm]
\displaystyle\frac{3e^\lambda-2}{3e^\lambda-1},
&\lambda\geq0.
\end{cases}
\label{eq:binary_third_cumulant_derivative}
\end{equation}
The equation
$\zeta'(\lambda)=\alpha$
is solved explicitly:
\begin{equation}
\lambda=
\begin{cases}
\displaystyle\log\frac{3\alpha}{1+\alpha},
&0<\alpha\leq\frac12,\\[3mm]
\displaystyle\log\frac{2-\alpha}{3(1-\alpha)},
&\frac12\leq\alpha<1.
\end{cases}
\label{eq:binary_third_dual_parameter}
\end{equation}
Substituting into the
dual formula, we obtain
for $0\leq\alpha\leq1/2$
\begin{equation}
\phi(\alpha)=\log\frac23+(1+\alpha)\log(1+\alpha)-\alpha\log(3\alpha),
\label{eq:binary_third_entropy}
\end{equation}
and for $1/2\leq\alpha\leq1$
\begin{equation}
\phi(\alpha)=\log\frac23+(2-\alpha)\log(2-\alpha)-(1-\alpha)\log\bigl(3(1-\alpha)\bigr).
\label{eq:style_s9_1_third_entropy_upper}
\end{equation}
The terms at the endpoints
are again interpreted
by continuity.
In particular,
\begin{equation}
\phi(0)=\log\frac23,
\label{eq:binary_third_entropy_values}
\end{equation}
\begin{equation}
\phi(1)=\log\frac23,
\label{eq:style_s9_1_third_entropy_one}
\end{equation}
\begin{equation}
\phi\!\left(\frac12\right)=0,
\label{eq:style_s9_1_third_entropy_center}
\end{equation}
\begin{equation}
\phi(\alpha)=\phi(1-\alpha).
\label{eq:style_s9_1_third_entropy_symmetry}
\end{equation}
Theorem~\ref{thm:affine_second_derivatives} gives $C^2$
matching of the entropy, with $\phi''(1/2)=-4/3$ by
\eqref{eq:style_s9_1_affine_entropy_curvature}.

The associated urn and its optimal means are shown in
Figure~\ref{fig:binary_affine_paths}. For $0<\alpha<1/2$, the mean
$u=\psi_\alpha^*(\tau)$ is determined by the monotone relation
\begin{equation}
\tau=\left(\frac{u}{\alpha}\right)^2
\left(\frac{1-2\alpha}{1-2u}\right)^3,
\qquad0<u\leq\alpha.
\label{eq:binary_affine_optimal_mean}
\end{equation}
This follows from Lemma~\ref{lem:affine_plot_paths} and the explicit
cumulant above. For $1/2<\alpha<1$, reflection gives
$\psi_\alpha^*=1-\psi_{1-\alpha}^*$; the central mean is $1/2$.

The example also
allows us to check
higher-order
regularity. The explicit
formulas yield
\begin{equation}
\zeta'(0)=\frac12,
\label{eq:binary_third_derivative_matching}
\end{equation}
\begin{equation}
\zeta''(0)=\frac34,
\label{eq:style_s9_1_third_second_derivative}
\end{equation}
\begin{equation}
\lim_{\lambda\uparrow0}\zeta'''(\lambda)=\frac32,
\label{eq:style_s9_1_third_left_third_derivative}
\end{equation}
\begin{equation}
\lim_{\lambda\downarrow0}\zeta'''(\lambda)=-\frac32.
\label{eq:style_s9_1_third_right_third_derivative}
\end{equation}
Thus $\zeta$ is
$C^2$ but not $C^3$
at the origin, and
\begin{equation}
\zeta(\lambda)
=\frac{\lambda}{2}
+\frac{3\lambda^2}{8}
-\frac{|\lambda|^3}{4}
+O(\lambda^4).
\label{eq:binary_third_origin_expansion}
\end{equation}
For the entropy,
we also have
\begin{equation}
\lim_{\alpha\uparrow1/2}\phi'''(\alpha)=\frac{32}{9},
\label{eq:binary_third_entropy_third_derivatives}
\end{equation}
\begin{equation}
\lim_{\alpha\downarrow1/2}\phi'''(\alpha)=-\frac{32}{9}.
\label{eq:style_s9_1_third_entropy_right_third_derivative}
\end{equation}
The matching is therefore
$C^2$ but not $C^3$
for $\phi$ as well.

In this case, $1/b=3$
is an integer, but no
logarithmic term appears
in the expansion of the
cumulant generating function.
Thus integrality of $1/b$
is not sufficient to
guarantee a nonzero
logarithmic coefficient,
as suggested by the
discussion following
Corollary~12 of
Ref.~\cite{Franchini_URNS}
and at the end of
Section~VIII of
Ref.~\cite{FB}, after
Eq.~(93).
The lack of third-order
regularity here results
from the different
one-sided limits
of the derivatives.

\subsection{Three possible increments}
\label{subsec:three_increments}

\subsubsection{General urn functions}
\label{subsec:ternary_general}

Let $\pi_0,\pi_1,\pi_2:[0,2]\to(0,1)$ be continuous,
with $\pi_0+\pi_1+\pi_2=1$. No affine or monotonicity
assumption is imposed here. Write
\begin{equation}
P_\alpha(z)=\pi_0(\alpha)+\pi_1(\alpha)z+\pi_2(\alpha)z^2,
\label{eq:ternary_general_polynomial_mean}
\end{equation}
\begin{equation}
\bar\pi(\alpha)=\pi_1(\alpha)+2\pi_2(\alpha).
\label{eq:style_s9_2_001}
\end{equation}
For $0<v<2$, the equation $zP_\alpha'(z)/P_\alpha(z)=v$
has the unique positive solution
\begin{equation}
z(\alpha,v)=
\frac{(v-1)\pi_1(\alpha)+
\sqrt{(1-v)^2\pi_1(\alpha)^2+
4v(2-v)\pi_0(\alpha)\pi_2(\alpha)}}
{2(2-v)\pi_2(\alpha)}.
\label{eq:ternary_general_tilt}
\end{equation}
Multiplication by $P_\alpha(z)$ gives a quadratic equation
whose leading and constant coefficients have opposite signs;
the positive root is \eqref{eq:ternary_general_tilt}.
Theorem~\ref{thm:intro_path_ldp} gives the cost
\begin{equation}
L_2(\alpha,v)=
\begin{cases}
v\log z(\alpha,v)-\log P_\alpha\bigl(z(\alpha,v)\bigr),&0<v<2,\\
-\log\pi_0(\alpha),& v=0,\\
-\log\pi_2(\alpha),& v=2.
\end{cases}
\label{eq:ternary_general_local_cost}
\end{equation}
The minimizing local frequencies are unique and, for $0<v<2$, are
\begin{equation}
\bigl(q_0^*,q_1^*,q_2^*\bigr)(\alpha,v)
=\frac{\bigl(\pi_0(\alpha),\,
\pi_1(\alpha)z(\alpha,v),\,
\pi_2(\alpha)z(\alpha,v)^2\bigr)}
{P_\alpha\bigl(z(\alpha,v)\bigr)}.
\label{eq:ternary_general_optimal_frequencies}
\end{equation}
At $v=0$ and $v=2$ they are $(1,0,0)$ and $(0,0,1)$,
respectively. Indeed, all frequencies of mean $v$ have the form
$(1-v+w,v-2w,w)$ with
$\max\{0,v-1\}\leq w\leq v/2$;
minimizing their relative entropy against
$(\pi_0(\alpha),\pi_1(\alpha),\pi_2(\alpha))$ gives
\eqref{eq:ternary_general_optimal_frequencies} by
Lemma~\ref{lem:local_frequency_minimum}.
In particular, the count and path rates are
\begin{equation}
I_{\mathrm c}(\gamma)=\int_0^1\sum_{j=0}^2\partial_\tau\gamma_j(\tau)
\log\frac{\partial_\tau\gamma_j(\tau)}{\pi_j\bigl((\gamma_1(\tau)+2\gamma_2(\tau))/\tau\bigr)}\,d\tau,
\label{eq:ternary_general_rates}
\end{equation}
\begin{equation}
I(\varphi)=\int_0^1 L_2\!\left(\frac{\varphi(\tau)}\tau,\partial_\tau\varphi(\tau)\right)d\tau
=\min_{\gamma_1+2\gamma_2=\varphi}I_{\mathrm c}(\gamma),
\label{eq:style_s9_2_002}
\end{equation}
where the last minimum is over admissible counts and is attained
by integrating the frequencies in
\eqref{eq:ternary_general_optimal_frequencies}, with their endpoint
values. The convention $0\log0=0$ applies.

Theorem~\ref{thm:intro_terminal_entropy} now gives
\begin{equation}
\phi(\alpha)=-\min_{\varphi\in Q(\alpha)}\int_0^1L_2\!\left(\frac{\varphi(\tau)}\tau,\partial_\tau\varphi(\tau)\right)d\tau,
\qquad0\leq\alpha\leq2,
\label{eq:ternary_general_entropy}
\end{equation}
\begin{equation}
\phi(0)=\log\pi_0(0),
\label{eq:style_s9_2_003}
\end{equation}
\begin{equation}
\phi(2)=\log\pi_2(2),
\label{eq:style_s9_2_004}
\end{equation}
\begin{equation}
\lim_{N\to\infty}\max_{0\leq r\leq2N,\ r\in\mathbb Z}
\left|N^{-1}\log\mathbb P(M_N=r)-\phi(r/N)\right|=0.
\label{eq:style_s9_2_005}
\end{equation}
Put $C=\{\alpha\in[0,2]:\pi_1(\alpha)+2\pi_2(\alpha)=\alpha\}$,
$c_-=\inf C$ and $c_+=\sup C$. By
Theorem~\ref{thm:intro_zero_cost},
$\phi(\alpha)=0$ exactly on $[c_-,c_+]$.
For endpoints in this interval, the optimal paths are precisely those
satisfying the following differential equation, with their unique
minimizing counts given by the subsequent formula:
\begin{equation}
\partial_\tau\varphi(\tau)=\pi_1\bigl(\varphi(\tau)/\tau\bigr)+2\pi_2\bigl(\varphi(\tau)/\tau\bigr),
\label{eq:ternary_general_zero_cost}
\end{equation}
\begin{equation}
\gamma_j^*(\tau)=\int_0^\tau\pi_j\bigl(\varphi(s)/s\bigr)\,ds,\qquad j=0,1,2.
\label{eq:style_s9_2_006}
\end{equation}
The differential equation holds almost everywhere.
Equivalently, $\tau\partial_\tau\psi=\pi_1(\psi)+2\pi_2(\psi)-\psi$;
the separation and fixed-point arrival criteria in
Theorem~\ref{thm:intro_zero_cost} therefore depend only on this mean.

For arbitrary functions in the stated class,
Theorem~\ref{thm:intro_cumulant} gives a strictly convex, nondecreasing,
$2$-Lipschitz limiting cumulant generating function, with $\zeta(0)=0$
and uniform convergence on compact parameter intervals. For every
$\lambda\in\mathbb R$,
\begin{equation}
\zeta(\lambda)=\max_{0\leq\alpha\leq2}\{\lambda\alpha+\phi(\alpha)\},
\label{eq:ternary_general_cumulant}
\end{equation}
\begin{equation}
e^{\zeta(\lambda)}=P_{\zeta'_-(\lambda)}(e^\lambda)=P_{\zeta'_+(\lambda)}(e^\lambda),
\label{eq:style_s9_2_007}
\end{equation}
\begin{equation}
\zeta'_-(0)=c_-,
\label{eq:style_s9_2_008}
\end{equation}
\begin{equation}
\zeta'_+(0)=c_+,
\label{eq:style_s9_2_009}
\end{equation}
\begin{equation}
\lim_{\lambda\to-\infty}\zeta(\lambda)=\log\pi_0(0),
\label{eq:style_s9_2_010}
\end{equation}
\begin{equation}
\lim_{\lambda\to+\infty}\bigl\{\zeta(\lambda)-2\lambda\bigr\}=\log\pi_2(2).
\label{eq:style_s9_2_011}
\end{equation}
To specify exactly when Theorem~\ref{thm:intro_cumulant_duality}
applies, observe that
\begin{equation}
P_\alpha(t)=t+(1-t)\bigl(\pi_0(\alpha)-t\pi_2(\alpha)\bigr).
\label{eq:ternary_general_injectivity_identity}
\end{equation}
Set $h_t(\alpha)=\pi_0(\alpha)-t\pi_2(\alpha)$.
The two injectivity assumptions are consequently equivalent to
requiring $h_t$ to be
injective on $[0,c_-]$ for every $0<t<1$, and on $[c_+,2]$
for every $t>1$. Each condition can be imposed separately.
A sufficient condition on either interval is that $\pi_0$ is
nonincreasing and $\pi_2$ nondecreasing, with at least one strict
inequality for each pair of distinct ordered states; the reversed
inequalities also suffice (Lemma~\ref{lem:cumulant_probability_sums}). Under both injectivity conditions,
$\zeta\in C^1(\mathbb R\setminus\{0\})$ is uniquely determined
by the equation \eqref{eq:style_s9_2_007} and the
slope constraints $0\leq\zeta'\leq c_-$ on $(-\infty,0)$ and
$c_+\leq\zeta'\leq2$ on $(0,\infty)$, in the continuous class
specified in Theorem~\ref{thm:intro_cumulant_duality}.
With $J_\lambda=[0,c_-]$ for $\lambda<0$ and
$J_\lambda=[c_+,2]$ for $\lambda>0$, this equation can be written as
\begin{equation}
\zeta'(\lambda)=
\bigl(h_{e^\lambda}|_{J_\lambda}\bigr)^{-1}
\!\left(\frac{e^{\zeta(\lambda)}-e^\lambda}{1-e^\lambda}\right),
\qquad \lambda\neq0.
\label{eq:ternary_general_inverse_cumulant_equation}
\end{equation}
Here the inverse is that of the indicated restriction,
whose range contains the displayed argument by
\eqref{eq:style_s9_2_007}. Moreover,
\begin{equation}
\phi(\alpha)=\inf_{\lambda\in\mathbb R}\{\zeta(\lambda)-\lambda\alpha\},
\label{eq:ternary_general_duality}
\end{equation}
\begin{equation}
\phi'\bigl(\zeta'(\lambda)\bigr)=-\lambda,\qquad\lambda\neq0.
\label{eq:style_s9_2_012}
\end{equation}
The entropy is concave and absolutely continuous on $[0,2]$,
$C^1$ on $(0,2)$, and strictly concave on $[0,c_-]$ and $[c_+,2]$;
it is strictly concave on all of $[0,2]$ exactly when $c_-=c_+$.
If only the condition for $0<t<1$ holds, uniqueness and the inverse
equation apply on $(-\infty,0)$, with duality on $[0,c_-]$;
the condition for $t>1$ independently gives the corresponding
conclusions on $(0,\infty)$ and $[c_+,2]$.
Without the injectivity assumptions, the variational formula
\eqref{eq:ternary_general_cumulant} still holds, but the reverse
formula \eqref{eq:ternary_general_duality} is not asserted.

Finally, let $f:[0,2]\to\mathbb R$ be continuous and concave,
and satisfy the realizability conditions \eqref{eq:intro_inverse_conditions}--\eqref{eq:style_s2_intro_inverse_conditions_2} with $K=2$.
At any $0<\alpha<2$ for which $f'(\alpha)\neq0$, set
$t=e^{-f'(\alpha)}$ and $E=e^{f(\alpha)-\alpha f'(\alpha)}$.
All positive probability vectors solving the inverse equation
are obtained by choosing $q=\pi_1(\alpha)>0$ and setting
\begin{equation}
\pi_2(\alpha)=\frac{E-1-q(t-1)}{t^2-1},
\label{eq:ternary_general_inverse_free}
\end{equation}
\begin{equation}
\pi_0(\alpha)=\frac{t^2-E-qt(t-1)}{t^2-1},
\label{eq:style_s9_2_013}
\end{equation}
subject to the equivalent strict positivity conditions
\begin{equation}
\displaystyle
0<q<\min\left\{\frac{E-1}{t-1},\,
\frac{t^2-E}{t(t-1)}\right\}.
\label{eq:ternary_general_inverse_positive}
\end{equation}
This is Lemma~\ref{lem:inverse_probability_freedom} for $K=2$.
Where $f'=0$, we have $f=0$ and the inverse equation reduces
to normalization, leaving any positive normalized vector as a
pointwise choice. These choices must extend to continuous,
strictly positive functions on $[0,2]$; if they do,
Theorem~\ref{thm:inverse_verification} gives $\phi=f$.
Continuity forces $\bar\pi(\alpha)=\alpha$ at any point of
$\{f=0\}$ approached by states where $f'\neq0$.

A global realization is always supplied by
Theorem~\ref{thm:intro_inverse_entropy}:
\begin{equation}
\bigl(\pi_0,\pi_1,\pi_2\bigr)=\bigl((1-p)^2,2p(1-p),p^2\bigr),
\label{eq:ternary_general_inverse_binomial}
\end{equation}
\begin{equation}
p(\alpha)=\frac{e^{[f(\alpha)-\alpha f'(\alpha)]/2}-1}{e^{-f'(\alpha)}-1},\qquad0<\alpha<2,\qquad f'(\alpha)\neq0,
\label{eq:style_s9_2_014}
\end{equation}
\begin{equation}
p(\alpha)=\alpha/2,\qquad0<\alpha<2,\qquad f'(\alpha)=0,
\label{eq:style_s9_2_015}
\end{equation}
\begin{equation}
p(0)=1-e^{f(0)/2},
\label{eq:style_s9_2_016}
\end{equation}
\begin{equation}
p(2)=e^{f(2)/2}.
\label{eq:style_s9_2_017}
\end{equation}
The function $p$ is continuous and belongs to $(0,1)$.
The affine urns studied next form a subclass of this general
three-increment model, to which Theorem~\ref{thm:intro_affine}
adds explicit quadratures and curvature formulas.

\subsubsection{Affine urns}
\label{subsec:ternary_affine_examples}

For $K=2$, the coefficient $B$ may have a zero other
than the origin, or a double zero at the origin.
The slope $b$ of the mean function no longer determines
all the coefficients $b_k$. In this subsection,
we compute examples of both types and a family with $b>0$
for which the cumulant generating function is analytic on the whole real line.
We begin with a restriction specific to the case
with three increments.

\begin{lemma}[Curvature near nonzero roots]
\label{lem:three_increment_nonzero_roots}
Let $K=2$ and assume the conditions
\eqref{eq:affine_admissibility}--\eqref{eq:style_s8_right_probabilities_positive}.
If $B\not\equiv0$, each of its real zeros
$\lambda_0\neq0$ is simple and satisfies
\begin{equation}
A(\lambda_0)-2B'(\lambda_0)>0.
\label{eq:three_increment_root_curvature_bound}
\end{equation}
For every affine urn in the class considered, we have
\begin{equation}
\zeta\in C^2(\mathbb R\setminus\{0\}),
\label{eq:three_increment_global_curvature}
\end{equation}
\begin{equation}
\zeta''(\lambda)>0,\qquad\lambda\neq0,
\label{eq:style_s9_2_018}
\end{equation}
\begin{equation}
\zeta\in C^2(\mathbb R)\quad\Longleftrightarrow\quad b<\frac12.
\label{eq:style_s9_2_019}
\end{equation}
If $b<1/2$, the second derivative of the entropy is
strictly negative on $(0,2)$. If $b\geq1/2$,
it vanishes only at $a/(1-b)$.
\end{lemma}

\begin{proof}
Assume $B\not\equiv0$ and fix a zero
$\lambda_0\neq0$. Set locally $x=e^{\lambda_0}$,
so that $x>0$ and $x\neq1$.
The polynomial $b_0+b_1u+b_2u^2$ has the two distinct
zeros $1$ and $x$. Thus $b_2\neq0$, the two zeros
are simple, and
\begin{equation}
B(\lambda)=b_2(e^\lambda-1)(e^\lambda-x),
\label{eq:three_increment_root_factorization}
\end{equation}
\begin{equation}
b_0=b_2x,
\label{eq:style_s9_2_020}
\end{equation}
\begin{equation}
b_1=-b_2(1+x),
\label{eq:style_s9_2_021}
\end{equation}
\begin{equation}
B'(\lambda_0)=b_2x(x-1).
\label{eq:style_s9_2_022}
\end{equation}
The following identities use only positivity
of the probabilities at the two boundaries. If $b_2>0$, then
\begin{equation}
A(\lambda_0)-2B'(\lambda_0)=a_0+(a_1+2b_1)x+a_2x^2+4b_2x>0,
\label{eq:three_increment_positive_root_identities}
\end{equation}
whereas for $b_2<0$,
\begin{equation}
A(\lambda_0)-2B'(\lambda_0)=(a_0+2b_0)+a_1x+(a_2+2b_2)x^2-4b_2x^2>0.
\label{eq:style_s9_2_023}
\end{equation}
Each summand is positive in the relevant case; this
proves \eqref{eq:three_increment_root_curvature_bound}.

Theorem~\ref{thm:affine_second_derivatives}
therefore gives a continuous, finite, and positive
second derivative also near the nonzero zeros of $B$.
At the origin, $A(0)=1$ and $B'(0)=b$:
the same theorem gives finite curvature
if and only if $b<1/2$.
If $B\equiv0$, we have $b=0$ and $\zeta=\log A$,
which is analytic with positive second derivative.
The conclusions for the entropy follow from
\eqref{eq:affine_entropy_curvature_reciprocal} and
\eqref{eq:affine_entropy_root_curvature}--\eqref{eq:style_s8_entropy_root_zero_curvature}.
\end{proof}

We return to the example \eqref{eq:affine_extra_root_example}--\eqref{eq:style_s8_extra_root_slopes}.
Its probabilities are
\begin{equation}
\pi_0(\alpha)=\frac13+\frac\alpha{12},\qquad0\leq\alpha\leq2,
\label{eq:three_increment_extra_root_probabilities}
\end{equation}
\begin{equation}
\pi_1(\alpha)=\frac13-\frac\alpha8,\qquad0\leq\alpha\leq2,
\label{eq:style_s9_2_024}
\end{equation}
\begin{equation}
\pi_2(\alpha)=\frac13+\frac\alpha{24},\qquad0\leq\alpha\leq2.
\label{eq:style_s9_2_025}
\end{equation}
At the endpoint $\alpha=2$, their values are respectively
$1/2$, $1/12$, and $5/12$; they are thus strictly
positive on the entire interval.
In this case,
\begin{equation}
a=1,
\label{eq:three_increment_extra_root_coefficients}
\end{equation}
\begin{equation}
b=-\frac1{24},
\label{eq:style_s9_2_026}
\end{equation}
\begin{equation}
\bar\pi(\alpha)=1-\frac\alpha{24},
\label{eq:style_s9_2_027}
\end{equation}
\begin{equation}
C=\left\{\frac{24}{25}\right\},
\label{eq:style_s9_2_028}
\end{equation}
\begin{equation}
A(\lambda)=\frac{1+e^\lambda+e^{2\lambda}}3,
\label{eq:style_s9_2_029}
\end{equation}
\begin{equation}
B(\lambda)=\frac{(e^\lambda-1)(e^\lambda-2)}{24}.
\label{eq:style_s9_2_030}
\end{equation}
The cumulant generating function is therefore $C^2$ on the whole real line.
We compute its integral representation
and the matching at the nonzero root.
Figure~\ref{fig:extra_root_affine_paths} displays the urn and its
optimal means, including the constant mean corresponding to the
nonzero root. Its positive cost is the negative of
\eqref{eq:style_s9_2_037}.

Setting $x=e^\lambda$, we find
\begin{equation}
\frac{A(\lambda)}{B(\lambda)}\,d\lambda
=\left\{\frac4x-\frac{24}{x-1}+\frac{28}{x-2}\right\}\,dx.
\label{eq:three_increment_extra_root_partial_fractions}
\end{equation}
An integrating factor for the linear equation
for $e^{-\zeta}$ is therefore
\begin{equation}
\frac{(x-1)^{24}}{x^4(x-2)^{28}},
\qquad x>0,\qquad x\neq1,2.
\label{eq:three_increment_extra_root_integrating_factor}
\end{equation}
The sign of $B$ is positive for $0<x<1$ and $x>2$,
negative for $1<x<2$.
Theorem~\ref{thm:affine_integral_representation}
therefore selects the endpoint $x=1$ in the first
two components and $x=+\infty$ in the third. For $0<x<2$ with $x\neq1$, this gives
\begin{equation}
e^{-\zeta(\log x)}=24\frac{x^4(x-2)^{28}}{(x-1)^{24}}\int_x^1\frac{(t-1)^{23}}{t^5(t-2)^{29}}\,dt,
\label{eq:three_increment_extra_root_lower_integral}
\end{equation}
where the integral is oriented, and
\begin{equation}
e^{-\zeta(\log x)}
=24\frac{x^4(x-2)^{28}}{(x-1)^{24}}
\int_x^\infty\frac{(t-1)^{23}}{t^5(t-2)^{29}}\,dt,
\qquad x>2.
\label{eq:three_increment_extra_root_upper_integral}
\end{equation}
For $1<x<2$, the integrand is negative and
the orientation of the integral makes the result positive.
Near $t=1$, the integrand vanishes to order $23$,
while as $t\to+\infty$ it is asymptotically equivalent to $t^{-11}$;
the integrals therefore converge.
The points $x=1$ and $x=2$ are filled in using
the matching data.

In particular, \eqref{eq:affine_root_matching}--\eqref{eq:style_s8_root_matching_derivative} and
\eqref{eq:affine_root_finite_curvature} give
\begin{equation}
\zeta(0)=0,
\label{eq:three_increment_extra_root_matching}
\end{equation}
\begin{equation}
\zeta'(0)=\frac{24}{25},
\label{eq:style_s9_2_031}
\end{equation}
\begin{equation}
\zeta''(0)=\frac{5888}{8125},
\label{eq:style_s9_2_032}
\end{equation}
\begin{equation}
\zeta(\log2)=\log\frac73,
\label{eq:style_s9_2_033}
\end{equation}
\begin{equation}
\zeta'(\log2)=\frac{40}{27},
\label{eq:style_s9_2_034}
\end{equation}
\begin{equation}
\zeta''(\log2)=\frac{6544}{9477}.
\label{eq:style_s9_2_035}
\end{equation}
Duality and the reciprocal of the curvatures give
\begin{equation}
\phi\!\left(\frac{24}{25}\right)=0,
\label{eq:three_increment_extra_root_entropy_data}
\end{equation}
\begin{equation}
\phi''\!\left(\frac{24}{25}\right)=-\frac{8125}{5888},
\label{eq:style_s9_2_036}
\end{equation}
\begin{equation}
\phi\!\left(\frac{40}{27}\right)=\log\frac73-\frac{40}{27}\log2,
\label{eq:style_s9_2_037}
\end{equation}
\begin{equation}
\phi'\!\left(\frac{40}{27}\right)=-\log2,
\label{eq:style_s9_2_038}
\end{equation}
\begin{equation}
\phi''\!\left(\frac{40}{27}\right)=-\frac{9477}{6544}.
\label{eq:style_s9_2_039}
\end{equation}
The other values of $\phi$ are uniquely determined
by the parametrization
\eqref{eq:affine_entropy_parametric_relations}
and the two integral formulas.
The boundary limits are
\begin{equation}
\lim_{\lambda\to-\infty}\zeta(\lambda)=\phi(0)=\log\frac13,
\label{eq:three_increment_extra_root_boundaries}
\end{equation}
\begin{equation}
\lim_{\lambda\to+\infty}\bigl\{\zeta(\lambda)-2\lambda\bigr\}=\phi(2)=\log\frac5{12}.
\label{eq:style_s9_2_040}
\end{equation}
The zero $\log2$ is not the selected endpoint
in either of the two adjacent components.
Matching at that point follows from the general
theorems, without imposing a new constant of integration.

Finite curvature does not imply analyticity at the
nonzero root. In this example, it is possible
to specify the order of regularity.
The kernel 
\begin{equation}
\frac{24(t-1)^{23}}{t^5(t-2)^{29}}
\label{eq:style_s9_2_065}
\end{equation}
 has
only the poles $0,2$ and is $O(t^{-11})$ at infinity.
In the partial fractions, the coefficients of $t^{-1}$ and
$(t-2)^{-1}$ are therefore opposite. The first is computed
using a fourth-order derivative:
\begin{equation}
\left.\frac1{4!}\frac{d^4}{dt^4}
\left\{\frac{24(t-1)^{23}}{(t-2)^{29}}\right\}
\right|_{t=0}=\frac{327}{2^{29}}.
\label{eq:three_increment_extra_root_log_coefficient}
\end{equation}
Indeed, the coefficient of $t^4$ in
$(1-t)^{23}(1-t/2)^{-29}$ is $109/8$, to be multiplied
by $24\,2^{-29}$. The coefficient of $(t-2)^{-1}$
is therefore $-327/2^{29}$.
Integrating with respect to the lower endpoint $x$,
the logarithmic term of $e^{-\zeta(\log x)}$ is therefore
\begin{equation}
\frac{327}{2^{29}}
\frac{x^4(x-2)^{28}}{(x-1)^{24}}\log|x-2|.
\label{eq:three_increment_extra_root_exponential_log_term}
\end{equation}
The other terms of the antiderivative, after
multiplication by the prefactor, are analytic
in a neighborhood of $2$. The constants fixed in the two
intervals add only multiples,
possibly different, of

\begin{equation}
\frac{x^4(x-2)^{28}}{(x-1)^{24}}.
\label{eq:style_s9_2_066}
\end{equation}
Thus the two branches have the same derivatives
up to order $27$, while the term
\eqref{eq:three_increment_extra_root_exponential_log_term}
prevents $C^{28}$ matching.
Since $e^{-\zeta(\log2)}=3/7>0$, the logarithm
preserves this order of regularity. The same holds
for the analytic change of variable $x=e^\lambda$, with nonzero
derivative at $\log2$.
The cumulant generating function is therefore $C^{27}$ but not $C^{28}$
at that point, despite having the finite curvature
computed in \eqref{eq:style_s9_2_035}.

We now consider a double zero at the origin. Let
\begin{equation}
\pi_0(\alpha)=\frac13+c\alpha,\qquad0<c<\frac1{12},\qquad 0\leq\alpha\leq2,
\label{eq:three_increment_double_root_probabilities}
\end{equation}
\begin{equation}
\pi_1(\alpha)=\frac13-2c\alpha,\qquad0<c<\frac1{12},\qquad 0\leq\alpha\leq2,
\label{eq:style_s9_2_041}
\end{equation}
\begin{equation}
\pi_2(\alpha)=\frac13+c\alpha,\qquad0<c<\frac1{12},\qquad 0\leq\alpha\leq2.
\label{eq:style_s9_2_042}
\end{equation}
The constraint on $c$ ensures strict positivity
also at $\alpha=2$.
The mean function is constant, but the probabilities
of the individual increments depend on the state:
\begin{equation}
a=1,
\label{eq:three_increment_double_root_mean}
\end{equation}
\begin{equation}
b=0,
\label{eq:style_s9_2_043}
\end{equation}
\begin{equation}
\bar\pi(\alpha)=1,
\label{eq:style_s9_2_044}
\end{equation}
\begin{equation}
C=\{1\}.
\label{eq:style_s9_2_045}
\end{equation}
In particular, $b=0$ does not reduce the urn
to the case of constant probabilities.
The functions in the equation for the cumulant generating function are
\begin{equation}
A(\lambda)=\frac{1+e^\lambda+e^{2\lambda}}3,
\label{eq:three_increment_double_root_coefficients}
\end{equation}
\begin{equation}
B(\lambda)=c(e^\lambda-1)^2.
\label{eq:style_s9_2_046}
\end{equation}
For $c>0$, $B$ vanishes only at the origin,
where the zero has multiplicity two.

Setting $x=e^\lambda$ again,
\begin{equation}
\frac{A(\lambda)}{B(\lambda)}\,d\lambda
=\frac1{3c}
\left\{\frac1x+\frac3{(x-1)^2}\right\}\,dx.
\label{eq:three_increment_double_root_partial_fractions}
\end{equation}
An integrating factor is therefore
\begin{equation}
x^{-1/(3c)}
\exp\!\left\{\frac{1}{c(x-1)}\right\}.
\label{eq:three_increment_double_root_integrating_factor}
\end{equation}
Since $B>0$ away from the origin, the integral
representation uses the right endpoint of each component.
This gives
\begin{multline}
e^{-\zeta(\log x)}=\frac{x^{1/(3c)}}c\exp\!\left\{-\frac1{c(x-1)}\right\}
\\
{}\times\int_x^{1}\frac{t^{-1/(3c)-1}}{(t-1)^2}\exp\!\left\{\frac1{c(t-1)}\right\}\,dt,
\qquad0<x<1,
\label{eq:three_increment_double_root_lower_integral}
\end{multline}
and
\begin{multline}
e^{-\zeta(\log x)}=\frac{x^{1/(3c)}}c\exp\!\left\{-\frac1{c(x-1)}\right\}
\\
{}\times\int_x^{\infty}\frac{t^{-1/(3c)-1}}{(t-1)^2}\exp\!\left\{\frac1{c(t-1)}\right\}\,dt,
\qquad x>1.
\label{eq:three_increment_double_root_upper_integral}
\end{multline}
The first integral converges because of the exponential
decay as $t\uparrow1$.
In the second, the integrand is asymptotically equivalent to
$t^{-1/(3c)-3}$ at infinity, and is therefore integrable.
Neither integral crosses $t=1$.

Theorem~\ref{thm:affine_second_derivatives}
ensures $C^2$ matching.
The formulas at the origin are
\begin{equation}
\zeta(0)=0,
\label{eq:three_increment_double_root_curvature}
\end{equation}
\begin{equation}
\zeta'(0)=1,
\label{eq:style_s9_2_047}
\end{equation}
\begin{equation}
\zeta''(0)=\frac23+2c,
\label{eq:style_s9_2_048}
\end{equation}
\begin{equation}
\phi''(1)=-\frac1{2/3+2c}.
\label{eq:style_s9_2_049}
\end{equation}
Indeed, $B'(0)=0$ and the variance at the fixed state is
$\pi_0(1)+\pi_2(1)=2/3+2c$.
The cumulant generating function is analytic on each open half-line.
The preceding formulas do not require assuming
analyticity across the double zero.
The boundary data and the maximum of the entropy are
\begin{equation}
\lim_{\lambda\to-\infty}\zeta(\lambda)=\phi(0)=\log\frac13,
\label{eq:three_increment_double_root_boundaries}
\end{equation}
\begin{equation}
\lim_{\lambda\to+\infty}\bigl\{\zeta(\lambda)-2\lambda\bigr\}=\phi(2)=\log\!\left(\frac13+2c\right),
\label{eq:style_s9_2_050}
\end{equation}
\begin{equation}
\phi(1)=0.
\label{eq:style_s9_2_051}
\end{equation}
The quadratures and \eqref{eq:affine_entropy_parametric_relations}
determine the entire entropy. As $c$ varies, the mean
remains $1$, while both $\phi(2)$ and $\phi''(1)$ change:
the same mean therefore does not determine the entropy.
Figure~\ref{fig:double_root_affine_paths} shows the urn and optimal
means for $c=1/24$: the constant mean increment does not force the
optimal paths to be straight.

The presence of a slope $b>0$ does not impose,
for $K=2$, a loss of analyticity of the cumulant generating function.
The following result provides a family
for which the variational principle also
admits an elementary solution. We consider
\begin{equation}
\pi_0(\alpha)=\frac14-\frac c2\alpha,\qquad0<c<\frac14,\qquad 0\leq\alpha\leq2,
\label{eq:three_increment_analytic_family}
\end{equation}
\begin{equation}
\pi_1(\alpha)=\frac12+c,\qquad0<c<\frac14,\qquad 0\leq\alpha\leq2,
\label{eq:style_s9_2_052}
\end{equation}
\begin{equation}
\pi_2(\alpha)=\frac14-c+\frac c2\alpha,\qquad0<c<\frac14,\qquad 0\leq\alpha\leq2.
\label{eq:style_s9_2_053}
\end{equation}
The probabilities are normalized and strictly
positive on $[0,2]$. The coefficients of the
mean function and the fixed point are
\begin{equation}
a=1-c,
\label{eq:three_increment_analytic_mean}
\end{equation}
\begin{equation}
b=c,
\label{eq:style_s9_2_054}
\end{equation}
\begin{equation}
\bar\pi(\alpha)=1-c+c\alpha,
\label{eq:style_s9_2_055}
\end{equation}
\begin{equation}
C=\{1\}.
\label{eq:style_s9_2_056}
\end{equation}

\begin{theorem}[A family with explicit cumulant generating function and entropy]
\label{thm:three_increment_analytic_family}
For the urn functions
\eqref{eq:three_increment_analytic_family}--\eqref{eq:style_s9_2_053},
the cumulant generating function and the terminal entropy are
\begin{equation}
\zeta(\lambda)=2\log\!\left(\frac{1+e^\lambda}2\right),\qquad\lambda\in\mathbb R,
\label{eq:three_increment_analytic_solution}
\end{equation}
\begin{equation}
\phi(\alpha)=-\alpha\log\alpha-(2-\alpha)\log(2-\alpha),\qquad0\leq\alpha\leq2,
\label{eq:style_s9_2_057}
\end{equation}
with the convention $0\log0=0$.
Both are independent of $c$.
The cumulant generating function is analytic and strictly convex
on $\mathbb{R}$; the entropy is analytic and
strictly concave on $(0,2)$, with a unique zero
at $\alpha=1$.
The dual parameter and the derivatives are
\begin{equation}
\lambda=\log\frac\alpha{2-\alpha},\qquad0<\alpha<2,
\label{eq:three_increment_analytic_derivatives}
\end{equation}
\begin{equation}
\zeta'(\lambda)=\frac{2e^\lambda}{1+e^\lambda},
\label{eq:style_s9_2_058}
\end{equation}
\begin{equation}
\zeta''(\lambda)=\frac{2e^\lambda}{(1+e^\lambda)^2},
\label{eq:style_s9_2_059}
\end{equation}
\begin{equation}
\phi'(\alpha)=\log\frac{2-\alpha}\alpha,
\label{eq:style_s9_2_060}
\end{equation}
\begin{equation}
\phi''(\alpha)=-\frac2{\alpha(2-\alpha)}.
\label{eq:style_s9_2_061}
\end{equation}
The parameter is the unique minimizer in
\eqref{eq:affine_entropy_duality}. In particular,
$\zeta'(0)=1$, $\zeta''(0)=1/2$, $\phi(1)=0$,
$\phi''(1)=-2$ and $\phi(0)=\phi(2)=-2\log2$.
\end{theorem}

\begin{proof}
The functions $A$ and $B$ are
\begin{equation}
A(\lambda)=\frac14+\left(\frac12+c\right)e^\lambda+\left(\frac14-c\right)e^{2\lambda},
\label{eq:three_increment_analytic_AB}
\end{equation}
\begin{equation}
A(\lambda)=\left(\frac{1+e^\lambda}2\right)^2-ce^\lambda(e^\lambda-1),
\label{eq:style_s9_2_067}
\end{equation}
\begin{equation}
B(\lambda)=\frac c2(e^{2\lambda}-1).
\label{eq:style_s9_2_062}
\end{equation}
The proposed function for $\zeta$ satisfies
\begin{equation}
B(\lambda)\zeta'(\lambda)=ce^\lambda(e^\lambda-1),
\label{eq:three_increment_analytic_equation_check}
\end{equation}
\begin{equation}
A(\lambda)+B(\lambda)\zeta'(\lambda)=\left(\frac{1+e^\lambda}2\right)^2=e^{\zeta(\lambda)}.
\label{eq:style_s9_2_063}
\end{equation}
Moreover,
\begin{equation}
\lim_{\lambda\to-\infty}\zeta(\lambda)=\log\frac14,
\label{eq:three_increment_analytic_boundaries}
\end{equation}
\begin{equation}
\lim_{\lambda\to+\infty}\bigl\{\zeta(\lambda)-2\lambda\bigr\}=\log\frac14.
\label{eq:style_s9_2_064}
\end{equation}
Since $a_0=a_2+2b_2=1/4$, these are
the data \eqref{eq:affine_cumulant_boundaries}--\eqref{eq:style_s8_cumulant_right_boundary}.
Uniqueness in
Theorem~\ref{thm:affine_cumulant_equation}
therefore identifies the cumulant generating function on the whole real line.

The derivative $\zeta'$ is strictly increasing
and takes every value in $(0,2)$.
Solving $\zeta'(\lambda)=\alpha$ gives
the parameter in \eqref{eq:three_increment_analytic_derivatives}.
Substituting into the duality
\eqref{eq:affine_entropy_duality} directly gives
the entropy formula \eqref{eq:style_s9_2_057}.
Continuity extends the formula to both endpoints.
The derivative expressions follow by direct
differentiation. Strict concavity and $\phi(1)=0$
show that $1$ is the unique zero.
\end{proof}

For $c=0$, the probabilities become constant
and equal to $(1/4,1/2,1/4)$; the same formulas
follow from \eqref{eq:affine_constant_probabilities}.
For every admissible $c>0$, however, both $\pi_0$
and $\pi_2$ do depend on the state,
and $b=c>0$.
The variance at the fixed state is
$\pi_0(1)+\pi_2(1)=1/2-c$; therefore
\begin{equation}
\frac{\sum_{k=0}^2(k-1)^2\pi_k(1)}{1-2b}
=\frac{1/2-c}{1-2c}=\frac12,
\label{eq:three_increment_analytic_variance_compensation}
\end{equation}
in agreement with \eqref{eq:affine_origin_finite_curvature}.

Equality of the terminal entropies does not imply
equality of the functionals on paths.
Indeed, consider the path
\begin{equation}
\varphi(\tau)=
\begin{cases}
0,&0\leq\tau\leq\frac12,\\
2\tau-1,&\frac12\leq\tau\leq1.
\end{cases}
\label{eq:three_increment_comparison_path}
\end{equation}
It belongs to $Q(1)$.
In the first half of the interval,
$\partial_\tau\varphi=0$ and $\varphi(\tau)/\tau=0$;
in the second, $\partial_\tau\varphi=2$ and
$\varphi(\tau)/\tau=2-1/\tau$.
The local cost at velocity $0$ is
$-\log\pi_0$, while at velocity $2$ it is
$-\log\pi_2$.
Definition \eqref{eq:path_rate} therefore gives
\begin{equation}
I(\varphi)
=\log2-\int_{1/2}^{1}
\log\left(\frac14-\frac{c}{2\tau}\right)\,d\tau.
\label{eq:three_increment_comparison_path_cost}
\end{equation}
The term in the logarithm is strictly
positive for $0<c<1/4$ and decreases
strictly as $c$ increases.
The cost therefore increases strictly with $c$.
The family \eqref{eq:three_increment_analytic_family}--\eqref{eq:style_s9_2_053}
has the same terminal entropy for all
admissible parameters, but distinct rate functionals
on paths.

Lemma~\ref{lem:affine_plot_paths} gives a closed formula for the optimal
means of this family. For $0<\alpha<2$, write $d=\alpha-1$. Then
\begin{equation}
\psi_\alpha^*(\tau)=1+\frac{d}{\sqrt{d^2+(1-d^2)\tau^{2c}}},
\qquad0<\tau\leq1.
\label{eq:analytic_affine_optimal_mean}
\end{equation}
Indeed, the characteristic equation reduces to
\begin{equation}
\tau\partial_\tau\psi_\alpha^*
=c\psi_\alpha^*(\psi_\alpha^*-1)(\psi_\alpha^*-2).
\label{eq:analytic_affine_optimal_flow}
\end{equation}
Formula~\eqref{eq:analytic_affine_optimal_mean} satisfies this equation
and the terminal condition, including $d=0$.
Figure~\ref{fig:analytic_affine_paths} pairs these means with the urn
functions. Their dependence on $c$ contrasts with the common terminal
entropy, while every corresponding cumulative path has the minimum
cost $-\phi(\alpha)$.

\subsection{A strictly positive cubic example}
\label{subsec:cubic_positive_example}

We again consider $K=2$, with the following urn probabilities for $0\leq\alpha\leq2$:
\begin{equation}
\pi_0(\alpha)=\frac{9-6\alpha^2+2\alpha^3}{11},
\label{eq:cubic_urn_probabilities}
\end{equation}
\begin{equation}
\pi_1(\alpha)=\frac1{11},
\label{eq:style_s9_cubic_urn_probabilities_line2}
\end{equation}
\begin{equation}
\pi_2(\alpha)=\frac{1+6\alpha^2-2\alpha^3}{11}.
\label{eq:style_s9_cubic_urn_probabilities_line3}
\end{equation}
The probabilities sum to one and satisfy
\begin{equation}
\pi_0'(\alpha)=-\frac{6\alpha(2-\alpha)}{11},
\qquad
\pi_2'(\alpha)=\frac{6\alpha(2-\alpha)}{11}.
\label{eq:cubic_probability_derivatives}
\end{equation}
The functions $\pi_0,\pi_1,\pi_2$ are respectively
strictly decreasing, constant, and strictly increasing,
with boundary data
\begin{equation}
\bigl(\pi_0(0),\pi_1(0),\pi_2(0)\bigr)=\frac1{11}\bigl(9,1,1\bigr),
\label{eq:cubic_strict_positivity}
\end{equation}
\begin{equation}
\bigl(\pi_0(2),\pi_1(2),\pi_2(2)\bigr)=\frac1{11}\bigl(1,1,9\bigr).
\label{eq:style_s9_cubic_strict_positivity_line2}
\end{equation}
\begin{equation}
\pi_k(\alpha)\geq\frac1{11},\qquad k=0,1,2,\qquad 0\leq\alpha\leq2.
\label{eq:style_s9_cubic_strict_positivity_2}
\end{equation}
The mean is strictly increasing:
\begin{equation}
\bar\pi(\alpha)=\frac{3+12\alpha^2-4\alpha^3}{11},
\qquad \bar\pi'(\alpha)=\frac{12\alpha(2-\alpha)}{11}.
\label{eq:cubic_mean_function}
\end{equation}
The factorization
\begin{equation}
\bar\pi(\alpha)-\alpha=-\frac4{11}
\Bigl(\alpha-\frac12\Bigr)\Bigl(\alpha-1\Bigr)
\Bigl(\alpha-\frac32\Bigr)
\label{eq:cubic_mean_drift}
\end{equation}
gives
\begin{equation}
C=\left\{\frac12,1,\frac32\right\}.
\label{eq:cubic_fixed_points}
\end{equation}
We can therefore make the paths between these points explicit
and apply the nonaffine characterization of the cumulant generating function.

\begin{theorem}[Zero entropy and optimal paths in the cubic example]
\label{thm:cubic_zero_cost_paths}
For the probabilities \eqref{eq:cubic_urn_probabilities}--\eqref{eq:style_s9_cubic_urn_probabilities_line3},
\begin{equation}
\phi(\alpha)=0
\qquad\Longleftrightarrow\qquad
\frac12\leq\alpha\leq\frac32,
\qquad
\phi(0)=\phi(2)=\log\frac9{11}.
\label{eq:cubic_entropy_zero_set}
\end{equation}
For every $\alpha\in[1/2,3/2]$, there exists
a unique optimal path in $Q(\alpha)$.
Its mean is
\begin{equation}
\psi^*(\tau)
=
1+
\frac{(\alpha-1)\tau^{1/11}}
{\sqrt{1-4(\alpha-1)^2+
4(\alpha-1)^2\tau^{2/11}}},
\qquad 0<\tau\leq1,
\label{eq:cubic_optimal_average}
\end{equation}
and the cumulative path is
\begin{equation}
\varphi^*(\tau)=\tau\psi^*(\tau)
\qquad(0<\tau\leq1),
\qquad \varphi^*(0)=0.
\label{eq:cubic_optimal_path}
\end{equation}
The formula also holds at the three fixed points,
where it gives $\varphi^*(\tau)=\alpha\tau$.
For $0\leq\tau\leq1$, the optimal increment counts are
\begin{equation}
\gamma_0^*(\tau)=\frac{21\tau}{22}-\frac{\varphi^*(\tau)}2,\qquad
\gamma_1^*(\tau)=\frac{\tau}{11},\qquad
\gamma_2^*(\tau)=\frac{\varphi^*(\tau)}2-\frac{\tau}{22}.
\label{eq:cubic_optimal_counts}
\end{equation}
\end{theorem}

\begin{proof}
The first equivalence in
\eqref{eq:cubic_entropy_zero_set}
follows from Theorem~\ref{thm:zero_entropy_interval}
and \eqref{eq:cubic_fixed_points}.
The endpoint values follow from
Theorem~\ref{thm:terminal_endpoints}
and \eqref{eq:cubic_strict_positivity}--\eqref{eq:style_s9_cubic_strict_positivity_line2}.

To verify the formula for the mean, set locally
\begin{equation}
d=\alpha-1,\qquad
D(\tau)=1-4d^2+4d^2\tau^{2/11},\qquad
r(\tau)=\frac{d\tau^{1/11}}{\sqrt{D(\tau)}}.
\label{eq:display_032}
\end{equation}
For every $\tau>0$ with $D(\tau)>0$, algebraic differentiation,
valid also for $d=0$, gives the first identity below;
the other two follow from \eqref{eq:cubic_mean_drift}
where $\psi=1+r\in[0,2]$:
\begin{equation}
\tau r'=\frac r{11}-\frac{4d^2\tau^{2/11}r}{11D}
=\frac{r(1-4d^2)}{11D}=\frac{r(1-4r^2)}{11},
\label{eq:cubic_average_differential_identity}
\end{equation}
\begin{equation}
\tau\psi'=\bar\pi(\psi)-\psi,\qquad (\tau\psi)'=\bar\pi(\psi).
\label{eq:style_s9_cubic_average_differential_identity_2}
\end{equation}
For $\alpha\in[1/2,3/2]$, we have $|d|\leq1/2$.
Thus $D(\tau)>0$ for $\tau>0$, even when $|d|=1/2$,
and, for $0<\tau\leq1$,
\begin{equation}
r(\tau)^2=
\frac{d^2\tau^{2/11}}{1-4d^2+4d^2\tau^{2/11}}\leq d^2,
\label{eq:display_033}
\end{equation}
since $(1-4d^2)\tau^{2/11}\leq1-4d^2$.
The sign of $r$ is that of $d$: the mean remains between $1$ and
$\alpha$, with $\psi(1)=\alpha$, and is constant at the three fixed points.
Boundedness of the mean gives $\tau\psi(\tau)\to0$ at zero;
\eqref{eq:style_s9_cubic_average_differential_identity_2} gives the slope
$\bar\pi(\psi)\in(0,2)$. Integrating between positive times and passing
to the limit at zero gives the constraints of $Q(\alpha)$.
Lemma~\ref{lem:zero_cost_paths} therefore proves zero cost.

The quadrature remains explicit: for $\alpha,u$ in the same
component $(1/2,1)$ or $(1,3/2)$,
\begin{equation}
F_\pi(\alpha,u)
=
\frac{11}{2}
\log\left[
\frac{(\alpha-1)^2\{1-4(u-1)^2\}}
{(u-1)^2\{1-4(\alpha-1)^2\}}
\right].
\label{eq:cubic_separation_integral}
\end{equation}
Indeed, an antiderivative of $1/(\bar\pi(u)-u)$ on each component
is $11\log|u-1|-(11/2)\log\{1-4(u-1)^2\}$.

For every $\varepsilon>0$, the equation
$\partial_\tau\psi=(\bar\pi(\psi)-\psi)/\tau$ has a right-hand
side that is Lipschitz in the mean on
$[\varepsilon,1]\times[0,2]$, since $\bar\pi$ is polynomial.
Uniqueness with the terminal datum $\psi(1)=\alpha$ therefore holds
on every $[\varepsilon,1]$, and hence on $(0,1]$.
This proves uniqueness of the zero-cost path for all
$\alpha\in[1/2,3/2]$, including the fixed points.

Finally, Lemma~\ref{lem:zero_cost_paths}
gives $\partial_\tau\gamma_k^*=\pi_k(\psi^*)$
and $\partial_\tau\varphi^*=\bar\pi(\psi^*)$.
The identities
\begin{equation}
\pi_0(u)=\frac{21}{22}-\frac{\bar\pi(u)}2,
\qquad
\pi_1(u)=\frac1{11},
\qquad
\pi_2(u)=\frac{\bar\pi(u)}2-\frac1{22}
\label{eq:cubic_probabilities_from_mean}
\end{equation}
can be integrated directly, with zero increment counts
at zero, and give \eqref{eq:cubic_optimal_counts}.
\end{proof}

The initial behavior of the optimal means is
\begin{equation}
\lim_{\tau\downarrow0}\psi^*(\tau)=\partial_\tau\varphi^*(0)=1,
\qquad \frac12<\alpha<\frac32,
\label{eq:cubic_initial_optimal_mean}
\end{equation}
\begin{equation}
\lim_{\tau\downarrow0}\psi^*(\tau)=\partial_\tau\varphi^*(0)=\alpha,
\qquad \alpha\in\Bigl\{\frac12,\frac32\Bigr\}.
\label{eq:style_s9_cubic_initial_optimal_mean_2}
\end{equation}
Thus states in the interior of the zero interval
of the entropy are reached by paths
whose initial mean is $1$.
For $\alpha\notin C$, these paths
are not straight lines. The straight line
$\varphi(\tau)=\alpha\tau$, however,
has strictly positive cost, since
$\alpha\neq\bar\pi(\alpha)$.
In particular, the three fixed points
and the entire interval of zero entropy
have distinct roles in the description
of the paths.

Theorem~\ref{thm:point_probabilities}
makes this distinction explicit also
for the point probabilities:
\begin{equation}
\lim_{N\to\infty}\frac1N
\log\mathbb{P}\bigl(M_N=\lfloor N\alpha\rfloor\bigr)
=0,
\qquad \frac12\leq\alpha\leq\frac32.
\label{eq:cubic_zero_point_probability_rate}
\end{equation}
The statement concerns the exponential scale $N$
even when the final state is not a fixed point.

We turn to the cumulant generating function and the two intervals
on which the entropy is negative.

\begin{theorem}[Cumulant generating function and entropy of the cubic example]
\label{thm:cubic_cumulant_entropy}
For the probabilities \eqref{eq:cubic_urn_probabilities}--\eqref{eq:style_s9_cubic_urn_probabilities_line3},
the cumulant generating function is real analytic on each
of the two open half-lines, with
$\zeta''(\lambda)>0$ for $\lambda\neq0$.
It satisfies
\begin{equation}
11e^{\zeta(\lambda)}=9+e^\lambda+e^{2\lambda}
+(e^{2\lambda}-1)\bigl\{6\zeta'(\lambda)^2-2\zeta'(\lambda)^3\bigr\},
\qquad \lambda\neq0,
\label{eq:cubic_cumulant_equation}
\end{equation}
and, as consequences of the equation and the global constraints
on the slopes, we have
\begin{equation}
\lim_{\lambda\to-\infty}\zeta(\lambda)=\log\frac9{11},
\label{eq:cubic_cumulant_boundary_data}
\end{equation}
\begin{equation}
\lim_{\lambda\to+\infty}\{\zeta(\lambda)-2\lambda\}=\log\frac9{11},
\label{eq:style_s9_cubic_cumulant_boundary_data_2}
\end{equation}
\begin{equation}
\zeta(0)=0,\qquad \zeta'_-(0)=\frac12,\qquad \zeta'_+(0)=\frac32.
\label{eq:style_s9_cubic_cumulant_boundary_data_3}
\end{equation}
The derivative is strictly increasing
on each half-line, with images
\begin{equation}
\{\zeta'(\lambda):\lambda<0\}=(0,1/2),
\qquad
\{\zeta'(\lambda):\lambda>0\}=(3/2,2).
\label{eq:cubic_cumulant_derivative_images}
\end{equation}
On each half-line, the equation
\eqref{eq:cubic_cumulant_equation}
with $\zeta$ replaced by $f$ and the constraints
\begin{equation}
0\leq f'(\lambda)\leq\frac12,\qquad \lambda<0,
\label{eq:cubic_cumulant_candidate_slopes}
\end{equation}
\begin{equation}
\frac32\leq f'(\lambda)\leq2,\qquad \lambda>0
\label{eq:style_s9_cubic_cumulant_candidate_slopes_2}
\end{equation}
uniquely characterize the restriction of $\zeta$ among functions
$f$ of class $C^1$ on the entire half-line considered. The value zero
at the origin and the corresponding asymptotes in
\eqref{eq:cubic_cumulant_boundary_data}--\eqref{eq:style_s9_cubic_cumulant_boundary_data_3} follow from the equation and the
global constraints, without assuming convexity of the candidates.

The entropy is concave and absolutely continuous
on $[0,2]$, of class $C^1$ on $(0,2)$,
and satisfies
\begin{equation}
\phi(\alpha)
=\inf_{\lambda\in\mathbb{R}}
\{\zeta(\lambda)-\lambda\alpha\},
\qquad 0\leq\alpha\leq2.
\label{eq:cubic_entropy_duality}
\end{equation}
It is strictly increasing and strictly
concave on $[0,1/2]$, identically zero
on $[1/2,3/2]$, strictly decreasing
and strictly concave on $[3/2,2]$.
On the two outer open intervals,
it is real analytic and has strictly
negative second derivative.
The derivatives at the boundaries of the three portions are
\begin{equation}
\lim_{\alpha\downarrow0}\phi'(\alpha)=+\infty,\qquad
\lim_{\alpha\uparrow2}\phi'(\alpha)=-\infty,
\label{eq:cubic_entropy_derivative_boundaries}
\end{equation}
\begin{equation}
\phi'(1/2)=\phi'(3/2)=0.
\label{eq:style_s9_cubic_entropy_derivative_boundaries_2}
\end{equation}
The following symmetries also hold:
\begin{equation}
\phi(2-\alpha)=\phi(\alpha),
\qquad
\zeta(-\lambda)=\zeta(\lambda)-2\lambda.
\label{eq:cubic_entropy_cumulant_symmetry}
\end{equation}
\end{theorem}

\begin{proof}
The exponential sum of the probabilities
can be rewritten as
\begin{equation}
\sum_{k=0}^2\pi_k(u)e^{k\lambda}
=
1-\frac{(e^\lambda-1)^2}{22}
+\frac{\bar\pi(u)}2(e^{2\lambda}-1).
\label{eq:cubic_exponential_sum}
\end{equation}
Since $\bar\pi$ is strictly increasing,
this sum is strictly increasing
in $u$ for $\lambda>0$ and strictly
decreasing for $\lambda<0$.
Thus both hypotheses of
Theorem~\ref{thm:cumulant_boundary_uniqueness}
are satisfied.
The theorem gives the equation, uniqueness
under the stated constraints, and
\eqref{eq:cubic_cumulant_derivative_images}.
The boundary conditions follow from
Lemma~\ref{lem:cumulant_boundary_values}
and \eqref{eq:cubic_fixed_points}.
Expanding the sum with
$u=\zeta'(\lambda)$ gives
\eqref{eq:cubic_cumulant_equation}.

To prove analyticity,
observe that for $\lambda\neq0$
\begin{equation}
\sum_{k=0}^2
\pi_k'(\zeta'(\lambda))e^{k\lambda}
=
\frac{6\zeta'(\lambda)\{2-\zeta'(\lambda)\}}{11}
(e^{2\lambda}-1)\neq0.
\label{eq:cubic_cumulant_implicit_derivative}
\end{equation}
The implicit function theorem
locally determines $\zeta'$ as
an analytic function of $\lambda$ and $\zeta$.
Thus the cumulant generating function solves an analytic
differential equation in a neighborhood
of every nonzero point, and is analytic
on both half-lines.

Convexity gives $\zeta''\geq0$.
If $\zeta''(\lambda_0)=0$ for
$\lambda_0\neq0$, then
$\zeta'''(\lambda_0)=0$, since
$\zeta''$ would have a local minimum.
Setting $\alpha=\zeta'(\lambda_0)$
and differentiating the equation
for the cumulant generating function twice,
we instead obtain
\begin{equation}
\left\{\sum_{k=0}^2
\pi_k'(\alpha)e^{k\lambda_0}\right\}
\zeta'''(\lambda_0)
=
-\sum_{k=0}^2
(k-\alpha)^2\pi_k(\alpha)e^{k\lambda_0}<0.
\label{eq:cubic_positive_curvature_argument}
\end{equation}
The identity also uses the
first derivative of the equation
at the point where $\zeta''=0$.
The strict inequality follows
from positivity of all the probabilities.
The contradiction proves
$\zeta''(\lambda)>0$ for $\lambda\neq0$.

Theorem~\ref{thm:entropy_shape}
gives duality, concavity,
absolute continuity, and the limits
of the derivatives of the entropy.
On the two outer intervals,
the inverse of $\zeta'$ is analytic,
and the relations
\begin{equation}
\phi'\bigl(\zeta'(\lambda)\bigr)=-\lambda,
\qquad
\phi''\bigl(\zeta'(\lambda)\bigr)
=-\frac1{\zeta''(\lambda)}<0,
\qquad \lambda\neq0,
\label{eq:cubic_entropy_parametric_derivatives}
\end{equation}
give the further properties asserted.

Finally, the probabilities satisfy
$\pi_k(2-\alpha)=\pi_{2-k}(\alpha)$.
Reflection of the history,
which replaces each $\sigma_n$
by $2-\sigma_n$, preserves
its probability, since the
uniform initialization is also symmetric.
For every $N$, we therefore have
\begin{equation}
\mathbb{P}(M_N=r)
=
\mathbb{P}(M_N=2N-r),
\qquad r=0,\ldots,2N.
\label{eq:cubic_finite_reflection}
\end{equation}
Passing to the limits defining
$\phi$ and $\zeta$ gives
\eqref{eq:cubic_entropy_cumulant_symmetry}.
\end{proof}

The cumulant generating function therefore has a corner
at the origin:
\begin{equation}
\zeta(\lambda)
=
\lambda+\frac12|\lambda|+o(|\lambda|),
\qquad \lambda\to0.
\label{eq:cubic_cumulant_corner}
\end{equation}
The two slopes correspond
to the endpoints of the interval of zero
entropy. The entropy, in turn, joins
that interval with zero derivative,
as specified in
\eqref{eq:style_s9_cubic_entropy_derivative_boundaries_2}.

The outer branches of the entropy are
uniquely determined by the boundary
problem and the parametrization
\begin{equation}
\alpha=\zeta'(\lambda),
\qquad
\phi(\alpha)=\zeta(\lambda)-\lambda\alpha,
\qquad \lambda\neq0.
\label{eq:cubic_entropy_parametrization}
\end{equation}
We conclude with an explicit estimate
of the cost on the outer intervals,
obtained by constructing an admissible
path. Fix $0\leq\alpha<1/2$
and set locally
\begin{equation}
\tau_0:=\exp\left\{-\int_0^\alpha\frac{du}{\bar\pi(u)-u}\right\}
=\left[\frac{4(1/2-\alpha)(3/2-\alpha)}{3(1-\alpha)^2}\right]^{11/2}.
\label{eq:cubic_external_departure_time}
\end{equation}
For $0<\alpha<1/2$, we have $0<\tau_0<1$;
for $\alpha=0$, we have $\tau_0=1$.
Evaluation of the integral uses
the antiderivative
$11\log|u-1|-(11/2)\log\{4(u-1)^2-1\}$
on the interval $[0,1/2)$.

For $\tau\in[\tau_0,1]$, let $\psi_\alpha(\tau)$ denote
the algebraic expression on the right-hand side of
\eqref{eq:cubic_optimal_average}, now with $0\leq\alpha<1/2$.
With $d=\alpha-1\in[-1,-1/2)$ and $D,r$ as in the preceding
verification, the choice of $\tau_0$ gives
\begin{equation}
\tau_0^{2/11}=\frac{4d^2-1}{3d^2},\qquad D(\tau_0)=\frac{4d^2-1}{3}>0,
\label{eq:display_034}
\end{equation}
\begin{equation}
\psi_\alpha(\tau_0)=0,\qquad \psi_\alpha(1)=\alpha.
\label{eq:style_s9_display_034_2}
\end{equation}
Since $D$ increases, the radicand is positive throughout $[\tau_0,1]$.
Moreover, $r<0$ and $1-4d^2<0$: the first identity in
\eqref{eq:cubic_average_differential_identity} shows that the mean
increases from $0$ to $\alpha$. It therefore remains in $[0,\alpha]$;
\eqref{eq:style_s9_cubic_average_differential_identity_2} gives the equation for the mean and the derivative
of the cumulative path on this new domain.
Set
\begin{equation}
\varphi(\tau)=
\begin{cases}
0,&0\leq\tau\leq\tau_0,\\
\tau\psi_\alpha(\tau),&\tau_0<\tau\leq1.
\end{cases}
\label{eq:cubic_external_trial_path}
\end{equation}
The matching is continuous, and on the second portion
$\partial_\tau\varphi=\bar\pi(\psi_\alpha)\in(0,2)$.
The path is therefore absolutely continuous and belongs
to $Q(\alpha)$; for $\alpha=0$, the second portion is empty.
This verification reuses only the algebraic formula,
not Theorem~\ref{thm:cubic_zero_cost_paths} outside
its interval of final states.

The cost of the portion after
$\tau_0$ is zero.
On the initial portion, both the state and the
velocity are zero,
so the local cost is
$-\log\pi_0(0)=\log(11/9)$.
Therefore
\begin{equation}
I(\varphi)=\tau_0\log\frac{11}{9}.
\label{eq:cubic_external_trial_cost}
\end{equation}
The variational formula
\eqref{eq:terminal_entropy_variational}
gives the explicit upper bound
\begin{equation}
0\leq-\phi(\alpha)\leq\log\frac{11}{9}
\left[\frac{4(1/2-\alpha)(3/2-\alpha)}{3(1-\alpha)^2}\right]^{11/2},
\qquad 0\leq\alpha\leq\frac12.
\label{eq:cubic_external_entropy_bound}
\end{equation}
The value at $1/2$ follows by continuity.
At $\alpha=0$, the estimate is an equality.
Symmetry provides the corresponding
estimate on $[3/2,2]$; in particular,
\begin{equation}
-\phi\!\left(\frac12-\varepsilon\right)
=
-\phi\!\left(\frac32+\varepsilon\right)
=
O(\varepsilon^{11/2}),
\qquad \varepsilon\downarrow0.
\label{eq:cubic_entropy_flat_boundary_bound}
\end{equation}
These are estimates of the minimum cost
obtained from the constructed path.
The two branches of the entropy
are determined by
\eqref{eq:cubic_cumulant_equation}--%
\eqref{eq:style_s9_cubic_cumulant_candidate_slopes_2}
and \eqref{eq:cubic_entropy_parametrization}.

\subsection{Finite-time arrival and increment counts}
\label{subsec:nonlipschitz_example}

The following fixed urn with $K=2$ makes the waiting and arrival
mechanism of Theorem~\ref{thm:fixed_point_component_arrival} explicit.
Its mean is obtained by rescaling the binary urn function in
Figure~2.2 of Ref.~\cite{Franchini_URNS}:
\begin{equation}
\bar\pi(\alpha)=
\begin{cases}
\alpha+\sqrt{1-2\alpha},&0\leq\alpha\leq1/2,\\
\alpha-\sqrt{(\alpha-1/2)(1-\alpha)},&1/2<\alpha\leq1,\\
\alpha+\sqrt{(\alpha-1)(3/2-\alpha)},&1<\alpha\leq3/2,\\
\alpha-\sqrt{2\alpha-3},&3/2<\alpha\leq2.
\end{cases}
\label{eq:nonlipschitz_mean}
\end{equation}
We choose the binomial increment probabilities
\begin{equation}
\pi_0(\alpha)=\Bigl(1-\frac{\bar\pi(\alpha)}2\Bigr)^2,
\label{eq:nonlipschitz_probabilities}
\end{equation}
\begin{equation}
\pi_1(\alpha)=\bar\pi(\alpha)\Bigl(1-\frac{\bar\pi(\alpha)}2\Bigr),
\label{eq:style_s9_nonlipschitz_probabilities_line2}
\end{equation}
\begin{equation}
\pi_2(\alpha)=\frac{\bar\pi(\alpha)^2}{4}.
\label{eq:style_s9_nonlipschitz_probabilities_line3}
\end{equation}
They sum to one and have mean $\pi_1+2\pi_2=\bar\pi$.
The formulas in \eqref{eq:nonlipschitz_mean} agree at their joining
points. Direct minimization and the symmetry
$\bar\pi(2-\alpha)=2-\bar\pi(\alpha)$ give
\begin{equation}
\frac{3-\sqrt2}{4}\leq\bar\pi(\alpha)\leq2-\frac{3-\sqrt2}{4},\qquad0\leq\alpha\leq2,
\label{eq:nonlipschitz_positivity}
\end{equation}
\begin{equation}
\min_{\substack{0\leq\alpha\leq2\\0\leq k\leq2}}\pi_k(\alpha)
=\Bigl(\frac{3-\sqrt2}{8}\Bigr)^2>0,
\label{eq:style_s9_nonlipschitz_positivity_2}
\end{equation}
\begin{equation}
C=\Bigl\{\frac12,1,\frac32\Bigr\}.
\label{eq:style_s9_nonlipschitz_positivity_3}
\end{equation}
Thus the general assumptions hold, and
Theorem~\ref{thm:zero_entropy_interval} gives zero entropy exactly
on $[1/2,3/2]$.

For every arrival time $0<t\leq1$, set $s=e^{-\pi}t$ and define
\begin{equation}
\psi_t^-(\tau)=1,\qquad0\leq\tau\leq s,
\label{eq:nonlipschitz_paths}
\end{equation}
\begin{equation}
\psi_t^-(\tau)=\frac12+\frac12\sin^2\!\Bigl(\frac12\log\frac t\tau\Bigr),\qquad s<\tau<t,
\label{eq:style_s9_nonlipschitz_paths_2}
\end{equation}
\begin{equation}
\psi_t^-(\tau)=\frac12,\qquad t\leq\tau\leq1,
\label{eq:style_s9_nonlipschitz_paths_3}
\end{equation}
\begin{equation}
\psi_t^+(\tau)=2-\psi_t^-(\tau),\qquad\varphi_t^\pm(\tau)=\tau\psi_t^\pm(\tau).
\label{eq:style_s9_nonlipschitz_paths_4}
\end{equation}
The means wait at $1$ until $s$, move to $1/2$ or $3/2$, and remain
there after time $t$. The finite transit integral
\begin{equation}
\int_{1/2}^1\frac{du}{\sqrt{(u-1/2)(1-u)}}=\pi
\label{eq:nonlipschitz_transit}
\end{equation}
gives $s/t=e^{-\pi}$ in
\eqref{eq:style_s6_full_transit_time}.
For a direct verification, put $z=\log(t/\tau)/2$ on $(s,t)$.
Then
\begin{equation}
\tau\partial_\tau\psi_t^-=-\frac12\sin z\cos z
=-\sqrt{(\psi_t^--1/2)(1-\psi_t^-)}=\bar\pi(\psi_t^-)-\psi_t^-.
\label{eq:nonlipschitz_path_verification}
\end{equation}
Both sides vanish at the joining times. On the constant portions
the values belong to $C$, and the right branch follows by symmetry.
Consequently $\partial_\tau\varphi_t^\pm=\bar\pi(\psi_t^\pm)$
on $[0,1]$, with one-sided derivatives at the endpoints.
The paths satisfy the increment constraints and have zero cost by
Lemma~\ref{lem:zero_cost_paths}. They are optimal for their terminal
states $1/2$ and $3/2$, respectively. The straight paths
$\varphi(\tau)=\tau/2$ and $\varphi(\tau)=3\tau/2$ correspond
separately to $t=0$. Figure~\ref{fig:finite_arrival_paths} shows the
mean function and the two nonconstant families.

The increment counts retain information that the common terminal
portion of the scalar paths no longer records.
For the left branch, let
\begin{equation}
J_t(\tau)=\int_0^\tau\bar\pi\bigl(\psi_t^-(v)\bigr)^2\,dv.
\label{eq:style_s9_integrated_square_definition}
\end{equation}
Writing $u=\log(t/\tau)$ in the middle portion gives the explicit formula
\begin{equation}
J_t(\tau)=\tau,\qquad0\leq\tau\leq s,
\label{eq:nonlipschitz_integrated_square}
\end{equation}
\begin{equation}
J_t(\tau)=\tau\Bigl[\frac58-\frac38\cos u+\frac{2\cos(2u)+\sin(2u)}{80}\Bigr]-\frac{s}{40},
\qquad s<\tau<t,
\label{eq:style_s9_nonlipschitz_integrated_square_2}
\end{equation}
\begin{equation}
J_t(\tau)=\frac\tau4+\frac{t-s}{40},\qquad t\leq\tau\leq1.
\label{eq:style_s9_nonlipschitz_integrated_square_3}
\end{equation}
Indeed, on $(s,t)$ one has
$\bar\pi(\psi_t^-)=3/4-(\cos u+\sin u)/4$;
differentiation of the middle expression yields its square.
The values match at $s$ and $t$, and $J_t(0)=0$.
Integration of \eqref{eq:nonlipschitz_probabilities}--\eqref{eq:style_s9_nonlipschitz_probabilities_line3} therefore gives
the minimizing counts
\begin{equation}
\gamma_{0,t}^*(\tau)=\tau-\varphi_t^-(\tau)+\frac{J_t(\tau)}4,
\label{eq:nonlipschitz_counts}
\end{equation}
\begin{equation}
\gamma_{1,t}^*(\tau)=\varphi_t^-(\tau)-\frac{J_t(\tau)}2,
\label{eq:style_s9_nonlipschitz_counts_line2}
\end{equation}
\begin{equation}
\gamma_{2,t}^*(\tau)=\frac{J_t(\tau)}4.
\label{eq:style_s9_nonlipschitz_counts_line3}
\end{equation}
They start at zero and satisfy
$\partial_\tau\gamma_{k,t}^*=\pi_k(\psi_t^-)$,
$\sum_k\gamma_{k,t}^*(\tau)=\tau$ and
$\sum_k k\gamma_{k,t}^*(\tau)=\varphi_t^-(\tau)$.
For the right branch, interchange the counts for $k=0$ and $k=2$.
After arrival, the left scalar path equals $\tau/2$, whereas
\begin{equation}
\Bigl(\gamma_{0,t}^*,\gamma_{1,t}^*,\gamma_{2,t}^*\Bigr)(\tau)
=\tau\Bigl(\frac9{16},\frac38,\frac1{16}\Bigr)
+\frac{t-s}{160}\Bigl(1,-2,1\Bigr),\qquad t\leq\tau\leq1.
\label{eq:nonlipschitz_count_memory}
\end{equation}
Different arrival times thus give distinct constant offsets in the
counts, while their sums and weighted sums agree. This realizes the
dependence on the earlier portion discussed after
\eqref{eq:style_global_component_optimal_counts_2}.
Figure~\ref{fig:component_memory} displays the normalized counts
$\gamma_{k,t}^*(\tau)/\tau$ and their common scalar mean after arrival;
the normalized counts still vary with time and with $t$.

\section*{Declaration of Generative AI}
Research, writing, and proofreading of the present manuscript were performed with AI assistance (GPT-6 Astra) under strict human guidance and scrutiny. We independently checked all mathematical arguments and take full responsibility for the content of the publication.


\begin{thebibliography}{54}
\bibitem{HLS}
B.~M.~Hill, D.~Lane, W.~Sudderth,
\textit{A strong law for some generalized urn processes},
Ann. Probab. 8(2), 214--226 (1980).
\url{https://doi.org/10.1214/aop/1176994772}.

\bibitem{HLS_orig_02}
J.~R.~Blum, M.~Brennan,
\textit{On the strong law of large numbers for dependent random variables},
Israel Journal of Mathematics 37(3), 241--245 (1980).
\url{https://doi.org/10.1007/BF02760965}.

\bibitem{HLS_orig_03}
W.~B.~Arthur, Yu.~M.~Ermoliev, Yu.~M.~Kaniovski,
\textit{A generalized urn problem and its applications},
Kibernetika, no. 1, 49--56 (1983).

\bibitem{HLS_orig_04}
R.~Pemantle,
\textit{When are Touchpoints Limits for Generalized Pólya Urns?},
Proceedings of the American Mathematical Society 113(1), 235--243 (1991).

\bibitem{HLS_orig_05}
R.~Gouet,
\textit{Martingale Functional Central Limit Theorems for a Generalized Pólya Urn},
Annals of Probability 21, 1624--1639 (1993).

\bibitem{HLS_orig_06}
K.~Nakayama, S.~Mori,
\textit{Universal function of the nonequilibrium phase transition of a nonlinear Pólya urn},
Physical Review E 104, 014109 (2021).

\bibitem{Dembo_Zeitouni}
A.~Dembo, O.~Zeitouni,
\textit{Large Deviations Techniques and Applications},
second edition, corrected printing, Springer, Berlin (2010).
\url{https://doi.org/10.1007/978-3-642-03311-7}.

\bibitem{Franchini_URNS}
S.~Franchini,
\textit{Large deviations for generalized Polya urns with arbitrary urn function},
Stochastic Process. Appl. 127(10), 3372--3411 (2017).
\url{https://doi.org/10.1016/j.spa.2017.02.010}.
arXiv version consulted: \url{https://arxiv.org/abs/1412.5762v7}.
\textit{Corrigendum}, Stochastic Process. Appl. 189, 104745 (2025).
\url{https://doi.org/10.1016/j.spa.2025.104745}.

\bibitem{FB}
S.~Franchini, R.~Balzan,
\textit{Large-deviation theory of increasing returns},
Phys. Rev. E 107, 064142 (2023).
\url{https://doi.org/10.1103/PhysRevE.107.064142}.
Corrected version: \url{https://arxiv.org/abs/2210.12585v5} (2025),
with a correction to equation~(93).

\bibitem{HLS_orig_10}
P.~Flajolet, J.~Gabarró, H.~Pekari,
\textit{Analytic urns},
Annals of Probability 33(3), 1200--1233 (2005).

\bibitem{HLS_orig_11}
P.~Flajolet, V.~Puyhaubert,
\textit{Analytic combinatorics at OK Corral},
Unpublished technical report (2005).

\bibitem{HLS_orig_12}
P.~Flajolet, P.~Dumas, V.~Puyhaubert,
\textit{Some exactly solvable models of urn process theory},
Fourth Colloquium on Mathematics and Computer Science, DMTCS Proceedings AG, 59--118 (2006).

\bibitem{HLS_orig_13}
W.~Bryc, D.~Minda, S.~Sethuraman,
\textit{Large deviations for the leaves in some random trees},
Advances in Applied Probability 41(3), 845--873 (2009).
\url{https://doi.org/10.1239/aap/1253281066}.

\bibitem{HLS_orig_14}
B.~Morcrette, H.~Mahmoud,
\textit{Exactly solvable balanced tenable urns with random entries via the analytic methodology},
Discrete Mathematics and Theoretical Computer Science Proceedings AQ, 219--232 (2012).

\bibitem{HLS_orig_15}
S.~Franchini,
\textit{Large deviations for generalized Polya urns with general urn functions},
PhD thesis, Università degli Studi Roma Tre (2015).
\url{https://arcadia.sba.uniroma3.it/handle/2307/5212}.

\bibitem{HLS_orig_16}
S.~Franchini,
\textit{Catene ideali con numero fissato di auto-intersezioni},
Master's thesis, Sapienza Università di Roma (2011).
\url{https://arxiv.org/abs/2412.10485}.
Deposited on arXiv in 2024.

\bibitem{HLS_orig_17}
S.~Franchini,
\textit{The Urn of Hill, Lane and Sudderth},
arXiv:2506.20826 (2025).
\url{https://arxiv.org/abs/2506.20826}.

\bibitem{HLS_orig_18}
R.~Pemantle,
\textit{A survey of random processes with reinforcement},
Probability Surveys 4 (2007).

\bibitem{HLS_orig_19}
H.~M.~Mahmoud,
\textit{Pólya Urn Models},
Chapman \& Hall/CRC, Taylor \& Francis, first edition (2008).

\bibitem{HLS_orig_20}
A.~Bagchi, A.~K.~Pal,
\textit{Asymptotic Normality in the Generalized Polya--Eggenberger Urn Model, with an Application to Computer Data Structures},
SIAM Journal on Algebraic and Discrete Methods 6(3), 394--405 (1985).
\url{https://doi.org/10.1137/0606041}.

\bibitem{ERW=000020shcutz=000020trimper}
G.~M.~Schütz, S.~Trimper,
\textit{Elephants can always remember: Exact long-range memory effects
in a non-Markovian random walk},
Phys. Rev. E 70, 045101(R) (2004).
\url{https://doi.org/10.1103/PhysRevE.70.045101}.

\bibitem{ERW=000020UM=000020Baur=000020Berton}
E.~Baur, J.~Bertoin,
\textit{Elephant random walks and their connection to Pólya-type urns},
Phys. Rev. E 94, 052134 (2016).
\url{https://doi.org/10.1103/PhysRevE.94.052134}.

\bibitem{HLS_orig_23}
A.~Gut, U.~Stadtmüller,
\textit{Variations of the elephant random walk},
Journal of Applied Probability 58(3), 805--829 (2021).

\bibitem{HLS_orig_24}
B.~Bercu, L.~Laulin,
\textit{On the center of mass of the elephant random walk},
Stochastic Processes and their Applications 133, 111--128 (2021).

\bibitem{HLS_orig_25}
K.~Maulik, P.~Roy, T.~Sadhukhan,
\textit{Asymptotic Properties of Generalized Elephant Random Walks},
arXiv:2406.19383 (2024).
\url{https://arxiv.org/abs/2406.19383v1}.

\bibitem{HLS_orig_26}
M.~Podder, A.~Roy,
\textit{Elephant random walks with multiple extractions and general reinforcement functions},
arXiv:2507.14626 (2025).
\url{https://arxiv.org/abs/2507.14626v1}.

\bibitem{HLS_orig_27}
R.~L.~Jack, R.~J.~Harris,
\textit{Giant leaps and long excursions: Fluctuation mechanisms in systems with long-range memory},
Physical Review E 102, 012154 (2020).
\url{https://doi.org/10.1103/PhysRevE.102.012154}.

\bibitem{HLS_orig_28}
{\raggedright
S.~Franchini,
\textit{Elephant Random Walk with multiple extractions},
arXiv:2210.12585v1 (2022), first version of the work
subsequently developed in ref.~\cite{FB}.
\url{https://arxiv.org/abs/2210.12585v1}.
\par}

\bibitem{HLS_orig_29}
W.~B.~Arthur,
\textit{Foundations of complexity economics},
Nature Reviews Physics 3 (2021).
\url{https://doi.org/10.1038/s42254-020-00273-3}.

\bibitem{HLS_orig_30}
T.~Gottfried, S.~Grosskinsky,
\textit{Asymptotics of generalized Pólya urns with non-linear feedback},
Electronic Journal of Probability 29, 1--56 (2024).

\bibitem{HLS_orig_31}
W.~B.~Arthur, Yu.~M.~Ermoliev, Yu.~M.~Kaniovski,
\textit{Path dependent processes and the emergence of macro-structure},
European Journal of Operational Research 30, 294--303 (1987).

\bibitem{HLS_orig_32}
W.~B.~Arthur,
\textit{Competing Technologies, Increasing Returns, and Lock-In by Historical Events},
The Economic Journal 99, 116--131 (1989).

\bibitem{HLS_orig_33}
W.~B.~Arthur, Yu.~M.~Ermoliev, Yu.~M.~Kaniovski,
\textit{Strong laws for a class of path-dependent stochastic processes with applications},
Stochastic Optimization, Lecture Notes in Control and Information Sciences 81, 287--300 (1986).

\bibitem{HLS_orig_34}
W.~B.~Arthur, Yu.~M.~Ermoliev, Yu.~M.~Kaniovski,
\textit{Limit Theorems for Proportions of Balls in a Generalized Urn Scheme},
IIASA Working Paper WP-87-111 (1987).

\bibitem{HLS_orig_35}
W.~B.~Arthur, Yu.~M.~Ermoliev, Yu.~M.~Kaniovski,
\textit{Non-Linear Urn Processes: Asymptotic Behavior and Applications},
IIASA Working Paper WP-87-085 (1987).

\bibitem{HLS_orig_36}
W.~B.~Arthur,
\textit{Increasing Returns and Path Dependence in the Economy},
University of Michigan Press, Ann Arbor (1994).

\bibitem{HLS_orig_37}
G.~Dosi, Yu.~M.~Ermoliev, Yu.~M.~Kaniovski,
\textit{Generalized urn schemes and technological dynamics},
Journal of Mathematical Economics 23, 1--19 (1994).

\bibitem{HLS_orig_38}
A.~M.~Espinosa, L.~Horna,
\textit{The statistical properties of the threshold model and the feedback leadership condition},
Journal of Applied Statistics 47(5), 844--864 (2020).
\url{https://doi.org/10.1080/02664763.2019.1658728}.

\bibitem{HLS_orig_39}
T.~Iyer,
\textit{Fixation of leadership in non-Markovian growth processes},
arXiv:2408.11516 (2024).
\url{https://arxiv.org/abs/2408.11516}.

\bibitem{HLS_orig_40}
T.~Gottfried, S.~Grosskinsky,
\textit{Wages and capital returns in a generalized Pólya urn},
J. Econ. Interact. Coord., online publication (2024).
\url{https://doi.org/10.1007/s11403-024-00437-9}.

\bibitem{HLS_orig_41}
T.~Gottfried,
\textit{Theory and application of a Pólya urn with non-linear feedback},
PhD thesis, University of Augsburg (2024).
\url{https://opus.bibliothek.uni-augsburg.de/opus4/frontdoor/deliver/index/docId/116236/file/Gottfried_Diss.pdf}.

\bibitem{HLS_orig_42}
A.~van de Rijt,
\textit{Self-Correcting Dynamics in Social Influence Processes},
American Journal of Sociology 124 (2019).
\url{https://doi.org/10.1086/702899}.

\bibitem{HLS_orig_43}
A.~Gelastopoulos, P.~P.~Analytis, G.~Le Mens, A.~van de Rijt,
\textit{The marginal majority effect: when social influence produces lock-in},
arXiv:2408.03952 (2024).
\url{https://arxiv.org/abs/2408.03952}.

\bibitem{HLS_orig_44}
R.~L.~Jack,
\textit{Large deviations in models of growing clusters with symmetry-breaking transitions},
Physical Review E 100, 012140 (2019).

\bibitem{HLS_orig_45}
R.~L.~Jack,
\textit{Ergodicity and large deviations in physical systems with stochastic dynamics},
European Physical Journal B 93, 74 (2020).

\bibitem{HLS_orig_46}
K.~Klymko, J.~P.~Garrahan, S.~Whitelam,
\textit{Similarity of ensembles of trajectories of reversible and irreversible growth processes},
Physical Review E 96, 042126 (2017).

\bibitem{HLS_orig_47}
K.~Klymko, P.~L.~Geissler, J.~P.~Garrahan, S.~Whitelam,
\textit{Rare behavior of growth processes via umbrella sampling of trajectories},
Physical Review E 97, 032123 (2018).

\bibitem{HLS_orig_48}
K.~Khanin, R.~Khanin,
\textit{A probabilistic model for establishment of neuron polarity},
Journal of Mathematical Biology 42, 26--40 (2001).

\bibitem{HLS_orig_49}
S.~Franchini, R.~Balzan,
\textit{Random Polymers and Generalized Urn Processes},
Physical Review E 98, 042502 (2018).

\bibitem{HLS_orig_50}
S.~Franchini,
\textit{Ideal chains with fixed self-intersection rate},
Physical Review E 84, 051104 (2011).

\bibitem{HLS_orig_51}
B.~D.~Hughes,
\textit{Random Walks and Random Environments, Volume 1},
Clarendon Press, Oxford (1995).

\bibitem{HLS_orig_52}
S.~Franchini, R.~Balzan,
\textit{Energy of the interacting self-avoiding walk at the $\theta$ point},
Physical Review E 102, 032143 (2020).

\bibitem{HLS_orig_53}
M.~van den Berg, E.~Bolthausen, F.~den Hollander,
\textit{Moderate deviations for the volume of the Wiener sausage},
Annals of Mathematics 153, 355--406 (2001).

\bibitem{HLS_orig_54}
G.~Dosi, A.~Moneta, E.~Stepanova,
\textit{Dynamic increasing returns and innovation diffusion: bringing Polya Urn processes to the empirical data},
Industry and Innovation 26(4), 461--478 (2019).
\url{https://doi.org/10.1080/13662716.2018.1444978}.
Published online in 2018.

\bibitem{Launay_simultaneous}
M.~Launay,
\textit{Urns with simultaneous drawing},
arXiv:1201.3495 (2012).
\url{https://arxiv.org/abs/1201.3495}.

\bibitem{Kuba_Mahmoud_multidrawing}
M.~Kuba, H.~M.~Mahmoud,
\textit{Two-color balanced affine urn models with multiple drawings},
Advances in Applied Mathematics 90, 1--26 (2017).
\url{https://doi.org/10.1016/j.aam.2017.04.004}.

\bibitem{Lasmar_Mailler_Selmi}
N.~Lasmar, C.~Mailler, O.~Selmi,
\textit{Multiple drawing multi-colour urns by stochastic approximation},
Journal of Applied Probability 55(1), 254--281 (2018).
\url{https://doi.org/10.1017/jpr.2018.16}.
Corrected version: \url{https://arxiv.org/abs/1611.09090v6}.

\bibitem{Steffens_multiple_units}
P.~R.~Steffens,
\textit{A model of multiple-unit ownership as a diffusion process},
Technological Forecasting and Social Change 70(9), 901--917 (2003).
\url{https://doi.org/10.1016/S0040-1625(02)00398-0}.

\bibitem{Budhiraja_Waterbury}
A.~Budhiraja, A.~Waterbury,
\textit{Empirical measure large deviations for reinforced chains on finite spaces},
Systems \& Control Letters 169, 105379 (2022).
\url{https://doi.org/10.1016/j.sysconle.2022.105379}.

\bibitem{Budhiraja_Waterbury_Zoubouloglou}
A.~Budhiraja, A.~Waterbury, P.~Zoubouloglou,
\textit{Large deviations for empirical measures of self-interacting Markov chains},
Stochastic Processes and their Applications 186, 104640 (2025).
\url{https://doi.org/10.1016/j.spa.2025.104640}.

\bibitem{Budhiraja_Coghi}
A.~Budhiraja, F.~Coghi,
\textit{Jump processes with self-interactions: large deviation asymptotics},
arXiv:2510.14659 (2025), revised version~2 (2026).
\url{https://arxiv.org/abs/2510.14659v2}.

\bibitem{Coghi_Garrahan_SIJP}
F.~Coghi, J.~P.~Garrahan,
\textit{Level-2.5 large deviations and uncertainty relations for
self-interacting jump processes: tilting constructions and the
emergence of timescale separation},
Phys. Rev. E 114, 054112 (2026).
\url{https://doi.org/10.1103/vhv8-tcrs}.

\bibitem{Coghi_Garrahan_spectral}
F.~Coghi, J.~P.~Garrahan,
\textit{Spectral theory for dynamical large deviations in non-Markov
self-interacting processes},
arXiv:2609.40329 (2026), revised version~2.
\url{https://arxiv.org/abs/2609.40329v2}.
\end{thebibliography}
\end{document}